
\ifx\shlhetal\undefinedcontrolsequence\let\shlhetal\relax\fi
\def\fmtname{AmS-TeX}

\def\fmtversion{2.1}
\catcode`\@=11
\ifx\amstexloaded@\relax\catcode`\@=\active
  \endinput\else\let\amstexloaded@\relax\fi
\newlinechar=`\^^J
\def\W@{\immediate\write\sixt@@n}
\def\CR@{\W@{^^J\fmtname - Version \fmtversion^^J%
COPYRIGHT 1985, 1990, 1991 - AMERICAN MATHEMATICAL SOCIETY^^J%
Use of this macro package is not restricted provided^^J%
each use is acknowledged upon publication.^^J}}
\CR@ \everyjob{\CR@}
\message{Loading definitions for}
\message{misc utility macros,}
\toksdef\toks@@=2
\long\def\rightappend@#1\to#2{\toks@{\\{#1}}\toks@@
 =\expandafter{#2}\xdef#2{\the\toks@@\the\toks@}\toks@{}\toks@@{}}
\def\alloclist@{}
\newif\ifalloc@
\def\showallocations{{\def\\{\immediate\write\m@ne}\alloclist@}\alloc@true}
\def\alloc@#1#2#3#4#5{\global\advance\count1#1by\@ne
 \ch@ck#1#4#2\allocationnumber=\count1#1
 \global#3#5=\allocationnumber
 \edef\next@{\string#5=\string#2\the\allocationnumber}%
 \expandafter\rightappend@\next@\to\alloclist@}
\newcount\count@@
\newcount\count@@@
\def\FN@{\futurelet\next}
\def\DN@{\def\next@}
\def\DNii@{\def\nextii@}
\def\RIfM@{\relax\ifmmode}
\def\RIfMIfI@{\relax\ifmmode\ifinner}
\def\setboxz@h{\setbox\z@\hbox}
\def\wdz@{\wd\z@}
\def\boxz@{\box\z@}
\def\setbox@ne{\setbox\@ne}
\def\wd@ne{\wd\@ne}
\def\iterate{\body\expandafter\iterate\else\fi}
\def\err@#1{\errmessage{AmS-TeX error: #1}}
\newhelp\defaulthelp@{Sorry, I already gave what help I could...^^J
Maybe you should try asking a human?^^J
An error might have occurred before I noticed any problems.^^J
``If all else fails, read the instructions.''}
\def\Err@{\errhelp\defaulthelp@\err@}
\def\eat@#1{}
\def\in@#1#2{\def\in@@##1#1##2##3\in@@{\ifx\in@##2\in@false\else\in@true\fi}%
 \in@@#2#1\in@\in@@}
\newif\ifin@
\def\space@.{\futurelet\space@\relax}
\space@. %
\newhelp\athelp@
{Only certain combinations beginning with @ make sense to me.^^J
Perhaps you wanted \string\@\space for a printed @?^^J
I've ignored the character or group after @.}
{\catcode`\~=\active 
 \lccode`\~=`\@ \lowercase{\gdef~{\FN@\at@}}}
\def\at@{\let\next@\at@@
 \ifcat\noexpand\next a\else\ifcat\noexpand\next0\else
 \ifcat\noexpand\next\relax\else
   \let\next\at@@@\fi\fi\fi
 \next@}
\def\at@@#1{\expandafter
 \ifx\csname\space @\string#1\endcsname\relax
  \expandafter\at@@@ \else
  \csname\space @\string#1\expandafter\endcsname\fi}
\def\at@@@#1{\errhelp\athelp@ \err@{\Invalid@@ @}}
\def\atdef@#1{\expandafter\def\csname\space @\string#1\endcsname}
\newhelp\defahelp@{If you typed \string\define\space cs instead of
\string\define\string\cs\space^^J
I've substituted an inaccessible control sequence so that your^^J
definition will be completed without mixing me up too badly.^^J
If you typed \string\define{\string\cs} the inaccessible control sequence^^J
was defined to be \string\cs, and the rest of your^^J
definition appears as input.}
\newhelp\defbhelp@{I've ignored your definition, because it might^^J
conflict with other uses that are important to me.}
\def\define{\FN@\define@}
\def\define@{\ifcat\noexpand\next\relax
 \expandafter\define@@\else\errhelp\defahelp@                               
 \err@{\string\define\space must be followed by a control
 sequence}\expandafter\def\expandafter\nextii@\fi}                          
\def\undefined@@@@@@@@@@{}
\def\preloaded@@@@@@@@@@{}
\def\next@@@@@@@@@@{}
\def\define@@#1{\ifx#1\relax\errhelp\defbhelp@                              
 \err@{\string#1\space is already defined}\DN@{\DNii@}\else
 \expandafter\ifx\csname\expandafter\eat@\string                            
 #1@@@@@@@@@@\endcsname\undefined@@@@@@@@@@\errhelp\defbhelp@
 \err@{\string#1\space can't be defined}\DN@{\DNii@}\else
 \expandafter\ifx\csname\expandafter\eat@\string#1\endcsname\relax          
 \global\let#1\undefined\DN@{\def#1}\else\errhelp\defbhelp@
 \err@{\string#1\space is already defined}\DN@{\DNii@}\fi
 \fi\fi\next@}

\def\predefine#1#2{\let#1#2}
\def\undefine#1{\let#1\undefined}
\message{page layout,}
\newdimen\captionwidth@
\captionwidth@\hsize
\advance\captionwidth@-1.5in
\def\pagewidth#1{\hsize#1\relax
 \captionwidth@\hsize\advance\captionwidth@-1.5in}
\def\pageheight#1{\vsize#1\relax}
\def\hcorrection#1{\advance\hoffset#1\relax}
\def\vcorrection#1{\advance\voffset#1\relax}
\message{accents/punctuation,}

\let\graveaccent\`
\let\acuteaccent\'
\let\tildeaccent\~
\let\hataccent\^
\let\underscore\_
\let\B\=
\let\D\.
\let\ic@\/
\def\/{\unskip\ic@}
\def\textfonti{\the\textfont\@ne}
\def\t#1#2{{\edef\next@{\the\font}\textfonti\accent"7F \next@#1#2}}
\def~{\unskip\nobreak\ \ignorespaces}
\def\.{.\spacefactor\@m}
\atdef@;{\leavevmode\null;}
\atdef@:{\leavevmode\null:}
\atdef@?{\leavevmode\null?}
\edef\@{\string @}
\def\relaxnext@{\let\next\relax}
\atdef@-{\relaxnext@\leavevmode
 \DN@{\ifx\next-\DN@-{\FN@\nextii@}\else
  \DN@{\leavevmode\hbox{-}}\fi\next@}%
 \DNii@{\ifx\next-\DN@-{\leavevmode\hbox{---}}\else
  \DN@{\leavevmode\hbox{--}}\fi\next@}%
 \FN@\next@}
\def\srdr@{\kern.16667em}
\def\drsr@{\kern.02778em}
\def\sldl@{\drsr@}
\def\dlsl@{\srdr@}
\atdef@"{\unskip\relaxnext@
 \DN@{\ifx\next\space@\DN@. {\FN@\nextii@}\else
  \DN@.{\FN@\nextii@}\fi\next@.}%
 \DNii@{\ifx\next`\DN@`{\FN@\nextiii@}\else
  \ifx\next\lq\DN@\lq{\FN@\nextiii@}\else
  \DN@####1{\FN@\nextiv@}\fi\fi\next@}%
 \def\nextiii@{\ifx\next`\DN@`{\sldl@``}\else\ifx\next\lq
  \DN@\lq{\sldl@``}\else\DN@{\dlsl@`}\fi\fi\next@}%
 \def\nextiv@{\ifx\next'\DN@'{\srdr@''}\else
  \ifx\next\rq\DN@\rq{\srdr@''}\else\DN@{\drsr@'}\fi\fi\next@}%
 \FN@\next@}

\def\textfontii{\the\textfont\tw@}
\def\lbrace@{\delimiter"4266308 }
\def\rbrace@{\delimiter"5267309 }
\def\{{\RIfM@\lbrace@\else{\textfontii f}\spacefactor\@m\fi}
\def\}{\RIfM@\rbrace@\else
 \let\@sf\empty\ifhmode\edef\@sf{\spacefactor\the\spacefactor}\fi
 {\textfontii g}\@sf\relax\fi}
\let\lbrace\{
\let\rbrace\}
\def\AmSTeX{{\textfontii A\kern-.1667em%
  \lower.5ex\hbox{M}\kern-.125emS}-\TeX}
\message{line and page breaks,}
\def\vmodeerr@#1{\Err@{\string#1\space not allowed between paragraphs}}
\def\mathmodeerr@#1{\Err@{\string#1\space not allowed in math mode}}
\def\linebreak{\RIfM@\mathmodeerr@\linebreak\else
 \ifhmode\unskip\unkern\break\else\vmodeerr@\linebreak\fi\fi}

\newskip\saveskip@
\def\allowlinebreak{\RIfM@\mathmodeerr@\allowlinebreak\else
 \ifhmode\saveskip@\lastskip\unskip
 \allowbreak\ifdim\saveskip@>\z@\hskip\saveskip@\fi
 \else\vmodeerr@\allowlinebreak\fi\fi}
\def\nolinebreak{\RIfM@\mathmodeerr@\nolinebreak\else
 \ifhmode\saveskip@\lastskip\unskip
 \nobreak\ifdim\saveskip@>\z@\hskip\saveskip@\fi
 \else\vmodeerr@\nolinebreak\fi\fi}
\def\newline{\relaxnext@
 \DN@{\RIfM@\expandafter\mathmodeerr@\expandafter\newline\else
  \ifhmode\ifx\next\par\else
  \expandafter\unskip\expandafter\null\expandafter\hfill\expandafter\break\fi
  \else
  \expandafter\vmodeerr@\expandafter\newline\fi\fi}%
 \FN@\next@}
\def\dmatherr@#1{\Err@{\string#1\space not allowed in display math mode}}
\def\nondmatherr@#1{\Err@{\string#1\space not allowed in non-display math
 mode}}
\def\onlydmatherr@#1{\Err@{\string#1\space allowed only in display math mode}}
\def\nonmatherr@#1{\Err@{\string#1\space allowed only in math mode}}
\def\mathbreak{\RIfMIfI@\break\else
 \dmatherr@\mathbreak\fi\else\nonmatherr@\mathbreak\fi}
\def\nomathbreak{\RIfMIfI@\nobreak\else
 \dmatherr@\nomathbreak\fi\else\nonmatherr@\nomathbreak\fi}
\def\allowmathbreak{\RIfMIfI@\allowbreak\else
 \dmatherr@\allowmathbreak\fi\else\nonmatherr@\allowmathbreak\fi}
\def\pagebreak{\RIfM@
 \ifinner\nondmatherr@\pagebreak\else\postdisplaypenalty-\@M\fi
 \else\ifvmode\removelastskip\break\else\vadjust{\break}\fi\fi}
\def\nopagebreak{\RIfM@
 \ifinner\nondmatherr@\nopagebreak\else\postdisplaypenalty\@M\fi
 \else\ifvmode\nobreak\else\vadjust{\nobreak}\fi\fi}
\def\nonvmodeerr@#1{\Err@{\string#1\space not allowed within a paragraph
 or in math}}
\def\vnonvmode@#1#2{\relaxnext@\DNii@{\ifx\next\par\DN@{#1}\else
 \DN@{#2}\fi\next@}%
 \ifvmode\DN@{#1}\else
 \DN@{\FN@\nextii@}\fi\next@}
\def\newpage{\vnonvmode@{\vfill\break}{\nonvmodeerr@\newpage}}
\def\smallpagebreak{\vnonvmode@\smallbreak{\nonvmodeerr@\smallpagebreak}}
\def\medpagebreak{\vnonvmode@\medbreak{\nonvmodeerr@\medpagebreak}}
\def\bigpagebreak{\vnonvmode@\bigbreak{\nonvmodeerr@\bigpagebreak}}
\def\NoBlackBoxes{\global\overfullrule\z@}
\def\BlackBoxes{\global\overfullrule5\p@}
\def\Invalid@#1{\def#1{\Err@{\Invalid@@\string#1}}}
\def\Invalid@@{Invalid use of }
\message{figures,}
\Invalid@\caption
\Invalid@\captionwidth
\newdimen\smallcaptionwidth@
\def\topspace{\mid@false\ins@}
\def\midspace{\mid@true\ins@}
\newif\ifmid@
\def\captionfont@{}
\def\ins@#1{\relaxnext@\allowbreak
 \smallcaptionwidth@\captionwidth@\gdef\thespace@{#1}%
 \DN@{\ifx\next\space@\DN@. {\FN@\nextii@}\else
  \DN@.{\FN@\nextii@}\fi\next@.}%
 \DNii@{\ifx\next\caption\DN@\caption{\FN@\nextiii@}%
  \else\let\next@\nextiv@\fi\next@}%
 \def\nextiv@{\vnonvmode@
  {\ifmid@\expandafter\midinsert\else\expandafter\topinsert\fi
   \vbox to\thespace@{}\endinsert}
  {\ifmid@\nonvmodeerr@\midspace\else\nonvmodeerr@\topspace\fi}}%
 \def\nextiii@{\ifx\next\captionwidth\expandafter\nextv@
  \else\expandafter\nextvi@\fi}%
 \def\nextv@\captionwidth##1##2{\smallcaptionwidth@##1\relax\nextvi@{##2}}%
 \def\nextvi@##1{\def\thecaption@{\captionfont@##1}%
  \DN@{\ifx\next\space@\DN@. {\FN@\nextvii@}\else
   \DN@.{\FN@\nextvii@}\fi\next@.}%
  \FN@\next@}%
 \def\nextvii@{\vnonvmode@
  {\ifmid@\expandafter\midinsert\else
  \expandafter\topinsert\fi\vbox to\thespace@{}\nobreak\smallskip
  \setboxz@h{\noindent\ignorespaces\thecaption@\unskip}%
  \ifdim\wdz@>\smallcaptionwidth@\centerline{\vbox{\hsize\smallcaptionwidth@
   \noindent\ignorespaces\thecaption@\unskip}}%
  \else\centerline{\boxz@}\fi\endinsert}
  {\ifmid@\nonvmodeerr@\midspace
  \else\nonvmodeerr@\topspace\fi}}%
 \FN@\next@}
\message{comments,}
\def\newcodes@{\catcode`\\12\catcode`\{12\catcode`\}12\catcode`\#12%
 \catcode`\%12\relax}
\def\oldcodes@{\catcode`\\0\catcode`\{1\catcode`\}2\catcode`\#6%
 \catcode`\%14\relax}
\def\comment{\newcodes@\endlinechar=10 \comment@}
{\lccode`\0=`\\
\lowercase{\gdef\comment@#1^^J{\comment@@#10endcomment\comment@@@}%
\gdef\comment@@#10endcomment{\FN@\comment@@@}%
\gdef\comment@@@#1\comment@@@{\ifx\next\comment@@@\let\next\comment@
 \else\def\next{\oldcodes@\endlinechar=`\^^M\relax}%
 \fi\next}}}
\def\pr@m@s{\ifx'\next\DN@##1{\prim@s}\else\let\next@\egroup\fi\next@}
\def\prime{{\null\prime@\null}}
\mathchardef\prime@="0230
\let\dsize\displaystyle

\let\ssize\scriptstyle

\message{math spacing,}
\def\,{\RIfM@\mskip\thinmuskip\relax\else\kern.16667em\fi}
\def\!{\RIfM@\mskip-\thinmuskip\relax\else\kern-.16667em\fi}
\let\thinspace\,
\let\negthinspace\!
\def\medspace{\RIfM@\mskip\medmuskip\relax\else\kern.222222em\fi}
\def\negmedspace{\RIfM@\mskip-\medmuskip\relax\else\kern-.222222em\fi}
\def\thickspace{\RIfM@\mskip\thickmuskip\relax\else\kern.27777em\fi}
\let\;\thickspace
\def\negthickspace{\RIfM@\mskip-\thickmuskip\relax\else
 \kern-.27777em\fi}
\atdef@,{\RIfM@\mskip.1\thinmuskip\else\leavevmode\null,\fi}
\atdef@!{\RIfM@\mskip-.1\thinmuskip\else\leavevmode\null!\fi}
\atdef@.{\RIfM@&&\else\leavevmode.\spacefactor3000 \fi}
\def\and{\DOTSB\;\mathchar"3026 \;}

\message{fractions,}
\def\frac#1#2{{#1\over#2}}

\newdimen\ex@
\ex@.2326ex
\Invalid@\thickness
\def\thickfrac{\relaxnext@
 \DN@{\ifx\next\thickness\let\next@\nextii@\else
 \DN@{\nextii@\thickness1}\fi\next@}%
 \DNii@\thickness##1##2##3{{##2\above##1\ex@##3}}%
 \FN@\next@}

\def\thickfracwithdelims#1#2{\relaxnext@\def\ldelim@{#1}\def\rdelim@{#2}%
 \DN@{\ifx\next\thickness\let\next@\nextii@\else
 \DN@{\nextii@\thickness1}\fi\next@}%
 \DNii@\thickness##1##2##3{{##2\abovewithdelims
 \ldelim@\rdelim@##1\ex@##3}}%
 \FN@\next@}

\def\:{\nobreak\hskip.1111em\mathpunct{}\nonscript\mkern-\thinmuskip{:}\hskip
 .3333emplus.0555em\relax}
\def\snug{\unskip\kern-\mathsurround}
\message{smash commands,}
\def\topsmash{\top@true\bot@false\smash@}
\def\botsmash{\top@false\bot@true\smash@}
\newif\iftop@
\newif\ifbot@
\def\smash{\top@true\bot@true\smash@}
\def\smash@{\RIfM@\expandafter\mathpalette\expandafter\mathsm@sh\else
 \expandafter\makesm@sh\fi}
\def\finsm@sh{\iftop@\ht\z@\z@\fi\ifbot@\dp\z@\z@\fi\leavevmode\boxz@}
\message{large operator symbols,}
\def\LimitsOnSums{\global\let\slimits@\displaylimits}
\def\NoLimitsOnSums{\global\let\slimits@\nolimits}
\LimitsOnSums
\mathchardef\coprod@="1360       \def\coprod{\DOTSB\coprod@\slimits@}
\mathchardef\bigvee@="1357       \def\bigvee{\DOTSB\bigvee@\slimits@}
\mathchardef\bigwedge@="1356     \def\bigwedge{\DOTSB\bigwedge@\slimits@}
\mathchardef\biguplus@="1355     \def\biguplus{\DOTSB\biguplus@\slimits@}
\mathchardef\bigcap@="1354       \def\bigcap{\DOTSB\bigcap@\slimits@}
\mathchardef\bigcup@="1353       \def\bigcup{\DOTSB\bigcup@\slimits@}
\mathchardef\prod@="1351         \def\prod{\DOTSB\prod@\slimits@}
\mathchardef\sum@="1350          \def\sum{\DOTSB\sum@\slimits@}
\mathchardef\bigotimes@="134E    \def\bigotimes{\DOTSB\bigotimes@\slimits@}
\mathchardef\bigoplus@="134C     \def\bigoplus{\DOTSB\bigoplus@\slimits@}
\mathchardef\bigodot@="134A      \def\bigodot{\DOTSB\bigodot@\slimits@}
\mathchardef\bigsqcup@="1346     \def\bigsqcup{\DOTSB\bigsqcup@\slimits@}
\message{integrals,}
\def\LimitsOnInts{\global\let\ilimits@\displaylimits}
\def\NoLimitsOnInts{\global\let\ilimits@\nolimits}
\NoLimitsOnInts
\def\int{\DOTSI\intop\ilimits@}
\def\oint{\DOTSI\ointop\ilimits@}
\def\intic@{\mathchoice{\hskip.5em}{\hskip.4em}{\hskip.4em}{\hskip.4em}}
\def\negintic@{\mathchoice
 {\hskip-.5em}{\hskip-.4em}{\hskip-.4em}{\hskip-.4em}}
\def\intkern@{\mathchoice{\!\!\!}{\!\!}{\!\!}{\!\!}}
\def\intdots@{\mathchoice{\plaincdots@}
 {{\cdotp}\mkern1.5mu{\cdotp}\mkern1.5mu{\cdotp}}
 {{\cdotp}\mkern1mu{\cdotp}\mkern1mu{\cdotp}}
 {{\cdotp}\mkern1mu{\cdotp}\mkern1mu{\cdotp}}}
\newcount\intno@
\def\iint{\DOTSI\intno@\tw@\FN@\ints@}
\def\iiint{\DOTSI\intno@\thr@@\FN@\ints@}
\def\iiiint{\DOTSI\intno@4 \FN@\ints@}
\def\idotsint{\DOTSI\intno@\z@\FN@\ints@}
\def\ints@{\findlimits@\ints@@}
\newif\iflimtoken@
\newif\iflimits@
\def\findlimits@{\limtoken@true\ifx\next\limits\limits@true
 \else\ifx\next\nolimits\limits@false\else
 \limtoken@false\ifx\ilimits@\nolimits\limits@false\else
 \ifinner\limits@false\else\limits@true\fi\fi\fi\fi}
\def\multint@{\int\ifnum\intno@=\z@\intdots@                                
 \else\intkern@\fi                                                          
 \ifnum\intno@>\tw@\int\intkern@\fi                                         
 \ifnum\intno@>\thr@@\int\intkern@\fi                                       
 \int}                                                                      
\def\multintlimits@{\intop\ifnum\intno@=\z@\intdots@\else\intkern@\fi
 \ifnum\intno@>\tw@\intop\intkern@\fi
 \ifnum\intno@>\thr@@\intop\intkern@\fi\intop}
\def\ints@@{\iflimtoken@                                                    
 \def\ints@@@{\iflimits@\negintic@\mathop{\intic@\multintlimits@}\limits    
  \else\multint@\nolimits\fi                                                
  \eat@}                                                                    
 \else                                                                      
 \def\ints@@@{\iflimits@\negintic@
  \mathop{\intic@\multintlimits@}\limits\else
  \multint@\nolimits\fi}\fi\ints@@@}
\def\LimitsOnNames{\global\let\nlimits@\displaylimits}
\def\NoLimitsOnNames{\global\let\nlimits@\nolimits@}
\LimitsOnNames
\def\nolimits@{\relaxnext@
 \DN@{\ifx\next\limits\DN@\limits{\nolimits}\else
  \let\next@\nolimits\fi\next@}%
 \FN@\next@}
\message{operator names,}
\def\newmcodes@{\mathcode`\'"27\mathcode`\*"2A\mathcode`\."613A%
 \mathcode`\-"2D\mathcode`\/"2F\mathcode`\:"603A }
\def\operatorname#1{\mathop{\newmcodes@\kern\z@\fam\z@#1}\nolimits@}
\def\operatornamewithlimits#1{\mathop{\newmcodes@\kern\z@\fam\z@#1}\nlimits@}
\def\qopname@#1{\mathop{\fam\z@#1}\nolimits@}
\def\qopnamewl@#1{\mathop{\fam\z@#1}\nlimits@}
\def\arccos{\qopname@{arccos}}
\def\arcsin{\qopname@{arcsin}}
\def\arctan{\qopname@{arctan}}
\def\arg{\qopname@{arg}}
\def\cos{\qopname@{cos}}
\def\cosh{\qopname@{cosh}}
\def\cot{\qopname@{cot}}
\def\coth{\qopname@{coth}}
\def\csc{\qopname@{csc}}
\def\deg{\qopname@{deg}}
\def\det{\qopnamewl@{det}}
\def\dim{\qopname@{dim}}
\def\exp{\qopname@{exp}}
\def\gcd{\qopnamewl@{gcd}}
\def\hom{\qopname@{hom}}
\def\inf{\qopnamewl@{inf}}
\def\injlim{\qopnamewl@{inj\,lim}}
\def\ker{\qopname@{ker}}
\def\lg{\qopname@{lg}}
\def\lim{\qopnamewl@{lim}}
\def\liminf{\qopnamewl@{lim\,inf}}
\def\limsup{\qopnamewl@{lim\,sup}}
\def\ln{\qopname@{ln}}
\def\log{\qopname@{log}}
\def\max{\qopnamewl@{max}}
\def\min{\qopnamewl@{min}}
\def\Pr{\qopnamewl@{Pr}}
\def\projlim{\qopnamewl@{proj\,lim}}
\def\sec{\qopname@{sec}}
\def\sin{\qopname@{sin}}
\def\sinh{\qopname@{sinh}}
\def\sup{\qopnamewl@{sup}}
\def\tan{\qopname@{tan}}
\def\tanh{\qopname@{tanh}}
\def\varinjlim{\mathop{\vtop{\ialign{##\crcr
 \hfil\rm lim\hfil\crcr\noalign{\nointerlineskip}\rightarrowfill\crcr
 \noalign{\nointerlineskip\kern-\ex@}\crcr}}}}
\def\varprojlim{\mathop{\vtop{\ialign{##\crcr
 \hfil\rm lim\hfil\crcr\noalign{\nointerlineskip}\leftarrowfill\crcr
 \noalign{\nointerlineskip\kern-\ex@}\crcr}}}}
\def\varliminf{\mathop{\underline{\vrule height\z@ depth.2exwidth\z@
 \hbox{\rm lim}}}}

\newdimen\buffer@
\buffer@\fontdimen13 \tenex
\newdimen\buffer
\buffer\buffer@

\def\ResetBuffer{\fontdimen13 \tenex\buffer@\global\buffer\buffer@}
\def\shave#1{\mathop{\hbox{$\m@th\fontdimen13 \tenex\z@                     
 \displaystyle{#1}$}}\fontdimen13 \tenex\buffer}

\message{multilevel sub/superscripts,}
\Invalid@\\
\def\Let@{\relax\iffalse{\fi\let\\=\cr\iffalse}\fi}
\Invalid@\vspace
\def\vspace@{\def\vspace##1{\crcr\noalign{\vskip##1\relax}}}
\def\multilimits@{\bgroup\vspace@\Let@
 \baselineskip\fontdimen10 \scriptfont\tw@
 \advance\baselineskip\fontdimen12 \scriptfont\tw@
 \lineskip\thr@@\fontdimen8 \scriptfont\thr@@
 \lineskiplimit\lineskip
 \vbox\bgroup\ialign\bgroup\hfil$\m@th\scriptstyle{##}$\hfil\crcr}
\def\Sb{_\multilimits@}
\def\endSb{\crcr\egroup\egroup\egroup}
\def\Sp{^\multilimits@}

\def\spreadlines#1{\RIfMIfI@\onlydmatherr@\spreadlines\else
 \openup#1\relax\fi\else\onlydmatherr@\spreadlines\fi}
\def\Mathstrut@{\copy\Mathstrutbox@}
\newbox\Mathstrutbox@
\setbox\Mathstrutbox@\null
\setboxz@h{$\m@th($}
\ht\Mathstrutbox@\ht\z@
\dp\Mathstrutbox@\dp\z@
\message{matrices,}
\newdimen\spreadmlines@
\def\spreadmatrixlines#1{\RIfMIfI@
 \onlydmatherr@\spreadmatrixlines\else
 \spreadmlines@#1\relax\fi\else\onlydmatherr@\spreadmatrixlines\fi}
\def\matrix{\null\,\vcenter\bgroup\Let@\vspace@
 \normalbaselines\openup\spreadmlines@\ialign
 \bgroup\hfil$\m@th##$\hfil&&\quad\hfil$\m@th##$\hfil\crcr
 \Mathstrut@\crcr\noalign{\kern-\baselineskip}}
\def\endmatrix{\crcr\Mathstrut@\crcr\noalign{\kern-\baselineskip}\egroup
 \egroup\,}
\def\format{\crcr\egroup\iffalse{\fi\ifnum`}=0 \fi\format@}
\newtoks\hashtoks@
\hashtoks@{#}
\def\format@#1\\{\def\preamble@{#1}%
 \def\l{$\m@th\the\hashtoks@$\hfil}%
 \def\c{\hfil$\m@th\the\hashtoks@$\hfil}%
 \def\r{\hfil$\m@th\the\hashtoks@$}%
 \edef\preamble@@{\preamble@}\ifnum`{=0 \fi\iffalse}\fi
 \ialign\bgroup\span\preamble@@\crcr}
\def\smallmatrix{\null\,\vcenter\bgroup\vspace@\Let@
 \baselineskip9\ex@\lineskip\ex@
 \ialign\bgroup\hfil$\m@th\scriptstyle{##}$\hfil&&\thickspace\hfil
 $\m@th\scriptstyle{##}$\hfil\crcr}
\def\endsmallmatrix{\crcr\egroup\egroup\,}

\newmuskip\dotsspace@
\dotsspace@1.5mu
\def\strip@#1 {#1}
\def\spacehdots#1\for#2{\multispan{#2}\xleaders
 \hbox{$\m@th\mkern\strip@#1 \dotsspace@.\mkern\strip@#1 \dotsspace@$}\hfill}
\def\hdotsfor#1{\spacehdots\@ne\for{#1}}
\def\multispan@#1{\omit\mscount#1\unskip\loop\ifnum\mscount>\@ne\sp@n\repeat}
\def\spaceinnerhdots#1\for#2\after#3{\multispan@{\strip@#2 }#3\xleaders
 \hbox{$\m@th\mkern\strip@#1 \dotsspace@.\mkern\strip@#1 \dotsspace@$}\hfill}
\def\innerhdotsfor#1\after#2{\spaceinnerhdots\@ne\for#1\after{#2}}
\def\cases{\bgroup\spreadmlines@\jot\left\{\,\matrix\format\l&\quad\l\\}
\def\endcases{\endmatrix\right.\egroup}
\message{multiline displays,}
\newif\ifinany@
\newif\ifinalign@
\newif\ifingather@
\def\strut@{\copy\strutbox@}
\newbox\strutbox@
\setbox\strutbox@\hbox{\vrule height8\p@ depth3\p@ width\z@}
\def\topaligned{\null\,\vtop\aligned@}
\def\botaligned{\null\,\vbox\aligned@}
\def\aligned{\null\,\vcenter\aligned@}
\def\aligned@{\bgroup\vspace@\Let@
 \ifinany@\else\openup\jot\fi\ialign
 \bgroup\hfil\strut@$\m@th\displaystyle{##}$&
 $\m@th\displaystyle{{}##}$\hfil\crcr}
\def\endaligned{\crcr\egroup\egroup}

\def\alignedat#1{\null\,\vcenter\bgroup\doat@{#1}\vspace@\Let@
 \ifinany@\else\openup\jot\fi\ialign\bgroup\span\preamble@@\crcr}
\newcount\atcount@
\def\doat@#1{\toks@{\hfil\strut@$\m@th
 \displaystyle{\the\hashtoks@}$&$\m@th\displaystyle
 {{}\the\hashtoks@}$\hfil}
 \atcount@#1\relax\advance\atcount@\m@ne                                    
 \loop\ifnum\atcount@>\z@\toks@=\expandafter{\the\toks@&\hfil$\m@th
 \displaystyle{\the\hashtoks@}$&$\m@th
 \displaystyle{{}\the\hashtoks@}$\hfil}\advance
  \atcount@\m@ne\repeat                                                     
 \xdef\preamble@{\the\toks@}\xdef\preamble@@{\preamble@}}

\def\gathered{\null\,\vcenter\bgroup\vspace@\Let@
 \ifinany@\else\openup\jot\fi\ialign
 \bgroup\hfil\strut@$\m@th\displaystyle{##}$\hfil\crcr}
\def\endgathered{\crcr\egroup\egroup}
\newif\iftagsleft@
\def\TagsOnLeft{\global\tagsleft@true}
\def\TagsOnRight{\global\tagsleft@false}
\TagsOnLeft
\newif\ifmathtags@
\def\TagsAsMath{\global\mathtags@true}
\def\TagsAsText{\global\mathtags@false}
\TagsAsText
\def\tagform@#1{\hbox{\rm(\ignorespaces#1\unskip)}}
\def\thetag{\leavevmode\tagform@}
\def\tag#1$${\iftagsleft@\leqno\else\eqno\fi                                
 \maketag@#1\maketag@                                                       
 $$}                                                                        
\def\maketag@{\FN@\maketag@@}
\def\maketag@@{\ifx\next"\expandafter\maketag@@@\else\expandafter\maketag@@@@
 \fi}
\def\maketag@@@"#1"#2\maketag@{\hbox{\rm#1}}                                
\def\maketag@@@@#1\maketag@{\ifmathtags@\tagform@{$\m@th#1$}\else
 \tagform@{#1}\fi}
\interdisplaylinepenalty\@M
\def\allowdisplaybreaks{\RIfMIfI@
 \onlydmatherr@\allowdisplaybreaks\else
 \interdisplaylinepenalty\z@\fi\else\onlydmatherr@\allowdisplaybreaks\fi}
\Invalid@\allowdisplaybreak
\Invalid@\displaybreak
\Invalid@\intertext
\def\allowdisplaybreak@{\def\allowdisplaybreak{\crcr\noalign{\allowbreak}}}
\def\displaybreak@{\def\displaybreak{\crcr\noalign{\break}}}
\def\intertext@{\def\intertext##1{\crcr\noalign{%
 \penalty\postdisplaypenalty \vskip\belowdisplayskip
 \vbox{\normalbaselines\noindent##1}%
 \penalty\predisplaypenalty \vskip\abovedisplayskip}}}
\newskip\centering@
\centering@\z@ plus\@m\p@
\def\align{\relax\ifingather@\DN@{\csname align (in
  \string\gather)\endcsname}\else
 \ifmmode\ifinner\DN@{\onlydmatherr@\align}\else
  \let\next@\align@\fi
 \else\DN@{\onlydmatherr@\align}\fi\fi\next@}
\newhelp\andhelp@
{An extra & here is so disastrous that you should probably exit^^J
and fix things up.}
\newif\iftag@
\newcount\and@
\def\align@{\inalign@true\inany@true
 \vspace@\allowdisplaybreak@\displaybreak@\intertext@
 \def\tag{\global\tag@true\ifnum\and@=\z@\DN@{&&}\else
          \DN@{&}\fi\next@}%
 \iftagsleft@\DN@{\csname align \endcsname}\else
  \DN@{\csname align \space\endcsname}\fi\next@}
\def\Tag@{\iftag@\else\errhelp\andhelp@\err@{Extra & on this line}\fi}
\newdimen\lwidth@
\newdimen\rwidth@
\newdimen\maxlwidth@
\newdimen\maxrwidth@
\newdimen\totwidth@
\def\measure@#1\endalign{\lwidth@\z@\rwidth@\z@\maxlwidth@\z@\maxrwidth@\z@
 \global\and@\z@                                                            
 \setbox@ne\vbox                                                            
  {\everycr{\noalign{\global\tag@false\global\and@\z@}}\Let@                
  \halign{\setboxz@h{$\m@th\displaystyle{\@lign##}$}
   \global\lwidth@\wdz@                                                     
   \ifdim\lwidth@>\maxlwidth@\global\maxlwidth@\lwidth@\fi                  
   \global\advance\and@\@ne                                                 
   &\setboxz@h{$\m@th\displaystyle{{}\@lign##}$}\global\rwidth@\wdz@        
   \ifdim\rwidth@>\maxrwidth@\global\maxrwidth@\rwidth@\fi                  
   \global\advance\and@\@ne                                                
   &\Tag@
   \eat@{##}\crcr#1\crcr}}
 \totwidth@\maxlwidth@\advance\totwidth@\maxrwidth@}                       
\def\displ@y@{\global\dt@ptrue\openup\jot
 \everycr{\noalign{\global\tag@false\global\and@\z@\ifdt@p\global\dt@pfalse
 \vskip-\lineskiplimit\vskip\normallineskiplimit\else
 \penalty\interdisplaylinepenalty\fi}}}
\def\black@#1{\noalign{\ifdim#1>\displaywidth
 \dimen@\prevdepth\nointerlineskip                                          
 \vskip-\ht\strutbox@\vskip-\dp\strutbox@                                   
 \vbox{\noindent\hbox to#1{\strut@\hfill}}
 \prevdepth\dimen@                                                          
 \fi}}
\expandafter\def\csname align \space\endcsname#1\endalign
 {\measure@#1\endalign\global\and@\z@                                       
 \ifingather@\everycr{\noalign{\global\and@\z@}}\else\displ@y@\fi           
 \Let@\tabskip\centering@                                                   
 \halign to\displaywidth
  {\hfil\strut@\setboxz@h{$\m@th\displaystyle{\@lign##}$}
  \global\lwidth@\wdz@\boxz@\global\advance\and@\@ne                        
  \tabskip\z@skip                                                           
  &\setboxz@h{$\m@th\displaystyle{{}\@lign##}$}
  \global\rwidth@\wdz@\boxz@\hfill\global\advance\and@\@ne                  
  \tabskip\centering@                                                       
  &\setboxz@h{\@lign\strut@\maketag@##\maketag@}
  \dimen@\displaywidth\advance\dimen@-\totwidth@
  \divide\dimen@\tw@\advance\dimen@\maxrwidth@\advance\dimen@-\rwidth@     
  \ifdim\dimen@<\tw@\wdz@\llap{\vtop{\normalbaselines\null\boxz@}}
  \else\llap{\boxz@}\fi                                                    
  \tabskip\z@skip                                                          
  \crcr#1\crcr                                                             
  \black@\totwidth@}}                                                      
\newdimen\lineht@
\expandafter\def\csname align \endcsname#1\endalign{\measure@#1\endalign
 \global\and@\z@
 \ifdim\totwidth@>\displaywidth\let\displaywidth@\totwidth@\else
  \let\displaywidth@\displaywidth\fi                                        
 \ifingather@\everycr{\noalign{\global\and@\z@}}\else\displ@y@\fi
 \Let@\tabskip\centering@\halign to\displaywidth
  {\hfil\strut@\setboxz@h{$\m@th\displaystyle{\@lign##}$}%
  \global\lwidth@\wdz@\global\lineht@\ht\z@                                 
  \boxz@\global\advance\and@\@ne
  \tabskip\z@skip&\setboxz@h{$\m@th\displaystyle{{}\@lign##}$}%
  \global\rwidth@\wdz@\ifdim\ht\z@>\lineht@\global\lineht@\ht\z@\fi         
  \boxz@\hfil\global\advance\and@\@ne
  \tabskip\centering@&\kern-\displaywidth@                                  
  \setboxz@h{\@lign\strut@\maketag@##\maketag@}%
  \dimen@\displaywidth\advance\dimen@-\totwidth@
  \divide\dimen@\tw@\advance\dimen@\maxlwidth@\advance\dimen@-\lwidth@
  \ifdim\dimen@<\tw@\wdz@
   \rlap{\vbox{\normalbaselines\boxz@\vbox to\lineht@{}}}\else
   \rlap{\boxz@}\fi
  \tabskip\displaywidth@\crcr#1\crcr\black@\totwidth@}}
\expandafter\def\csname align (in \string\gather)\endcsname
  #1\endalign{\vcenter{\align@#1\endalign}}
\Invalid@\endalign
\newif\ifxat@
\def\alignat{\RIfMIfI@\DN@{\onlydmatherr@\alignat}\else
 \DN@{\csname alignat \endcsname}\fi\else
 \DN@{\onlydmatherr@\alignat}\fi\next@}
\newif\ifmeasuring@
\newbox\savealignat@
\expandafter\def\csname alignat \endcsname#1#2\endalignat                   
 {\inany@true\xat@false
 \def\tag{\global\tag@true\count@#1\relax\multiply\count@\tw@
  \xdef\tag@{}\loop\ifnum\count@>\and@\xdef\tag@{&\tag@}\advance\count@\m@ne
  \repeat\tag@}%
 \vspace@\allowdisplaybreak@\displaybreak@\intertext@
 \displ@y@\measuring@true                                                   
 \setbox\savealignat@\hbox{$\m@th\displaystyle\Let@
  \attag@{#1}
  \vbox{\halign{\span\preamble@@\crcr#2\crcr}}$}%
 \measuring@false                                                           
 \Let@\attag@{#1}
 \tabskip\centering@\halign to\displaywidth
  {\span\preamble@@\crcr#2\crcr                                             
  \black@{\wd\savealignat@}}}                                               
\Invalid@\endalignat
\def\xalignat{\RIfMIfI@
 \DN@{\onlydmatherr@\xalignat}\else
 \DN@{\csname xalignat \endcsname}\fi\else
 \DN@{\onlydmatherr@\xalignat}\fi\next@}
\expandafter\def\csname xalignat \endcsname#1#2\endxalignat
 {\inany@true\xat@true
 \def\tag{\global\tag@true\def\tag@{}\count@#1\relax\multiply\count@\tw@
  \loop\ifnum\count@>\and@\xdef\tag@{&\tag@}\advance\count@\m@ne\repeat\tag@}%
 \vspace@\allowdisplaybreak@\displaybreak@\intertext@
 \displ@y@\measuring@true\setbox\savealignat@\hbox{$\m@th\displaystyle\Let@
 \attag@{#1}\vbox{\halign{\span\preamble@@\crcr#2\crcr}}$}%
 \measuring@false\Let@
 \attag@{#1}\tabskip\centering@\halign to\displaywidth
 {\span\preamble@@\crcr#2\crcr\black@{\wd\savealignat@}}}
\def\attag@#1{\let\Maketag@\maketag@\let\TAG@\Tag@                          
 \let\Tag@=0\let\maketag@=0
 \ifmeasuring@\def\llap@##1{\setboxz@h{##1}\hbox to\tw@\wdz@{}}%
  \def\rlap@##1{\setboxz@h{##1}\hbox to\tw@\wdz@{}}\else
  \let\llap@\llap\let\rlap@\rlap\fi                                         
 \toks@{\hfil\strut@$\m@th\displaystyle{\@lign\the\hashtoks@}$\tabskip\z@skip
  \global\advance\and@\@ne&$\m@th\displaystyle{{}\@lign\the\hashtoks@}$\hfil
  \ifxat@\tabskip\centering@\fi\global\advance\and@\@ne}
 \iftagsleft@
  \toks@@{\tabskip\centering@&\Tag@\kern-\displaywidth
   \rlap@{\@lign\maketag@\the\hashtoks@\maketag@}%
   \global\advance\and@\@ne\tabskip\displaywidth}\else
  \toks@@{\tabskip\centering@&\Tag@\llap@{\@lign\maketag@
   \the\hashtoks@\maketag@}\global\advance\and@\@ne\tabskip\z@skip}\fi      
 \atcount@#1\relax\advance\atcount@\m@ne
 \loop\ifnum\atcount@>\z@
 \toks@=\expandafter{\the\toks@&\hfil$\m@th\displaystyle{\@lign
  \the\hashtoks@}$\global\advance\and@\@ne
  \tabskip\z@skip&$\m@th\displaystyle{{}\@lign\the\hashtoks@}$\hfil\ifxat@
  \tabskip\centering@\fi\global\advance\and@\@ne}\advance\atcount@\m@ne
 \repeat                                                                    
 \xdef\preamble@{\the\toks@\the\toks@@}
 \xdef\preamble@@{\preamble@}
 \let\maketag@\Maketag@\let\Tag@\TAG@}                                      
\Invalid@\endxalignat
\def\xxalignat{\RIfMIfI@
 \DN@{\onlydmatherr@\xxalignat}\else\DN@{\csname xxalignat
  \endcsname}\fi\else
 \DN@{\onlydmatherr@\xxalignat}\fi\next@}
\expandafter\def\csname xxalignat \endcsname#1#2\endxxalignat{\inany@true
 \vspace@\allowdisplaybreak@\displaybreak@\intertext@
 \displ@y\setbox\savealignat@\hbox{$\m@th\displaystyle\Let@
 \xxattag@{#1}\vbox{\halign{\span\preamble@@\crcr#2\crcr}}$}%
 \Let@\xxattag@{#1}\tabskip\z@skip\halign to\displaywidth
 {\span\preamble@@\crcr#2\crcr\black@{\wd\savealignat@}}}
\def\xxattag@#1{\toks@{\tabskip\z@skip\hfil\strut@
 $\m@th\displaystyle{\the\hashtoks@}$&%
 $\m@th\displaystyle{{}\the\hashtoks@}$\hfil\tabskip\centering@&}%
 \atcount@#1\relax\advance\atcount@\m@ne\loop\ifnum\atcount@>\z@
 \toks@=\expandafter{\the\toks@&\hfil$\m@th\displaystyle{\the\hashtoks@}$%
  \tabskip\z@skip&$\m@th\displaystyle{{}\the\hashtoks@}$\hfil
  \tabskip\centering@}\advance\atcount@\m@ne\repeat
 \xdef\preamble@{\the\toks@\tabskip\z@skip}\xdef\preamble@@{\preamble@}}
\Invalid@\endxxalignat
\newdimen\gwidth@
\newdimen\gmaxwidth@
\def\gmeasure@#1\endgather{\gwidth@\z@\gmaxwidth@\z@\setbox@ne\vbox{\Let@
 \halign{\setboxz@h{$\m@th\displaystyle{##}$}\global\gwidth@\wdz@
 \ifdim\gwidth@>\gmaxwidth@\global\gmaxwidth@\gwidth@\fi
 &\eat@{##}\crcr#1\crcr}}}
\def\gather{\RIfMIfI@\DN@{\onlydmatherr@\gather}\else
 \ingather@true\inany@true\def\tag{&}%
 \vspace@\allowdisplaybreak@\displaybreak@\intertext@
 \displ@y\Let@
 \iftagsleft@\DN@{\csname gather \endcsname}\else
  \DN@{\csname gather \space\endcsname}\fi\fi
 \else\DN@{\onlydmatherr@\gather}\fi\next@}
\expandafter\def\csname gather \space\endcsname#1\endgather
 {\gmeasure@#1\endgather\tabskip\centering@
 \halign to\displaywidth{\hfil\strut@\setboxz@h{$\m@th\displaystyle{##}$}%
 \global\gwidth@\wdz@\boxz@\hfil&
 \setboxz@h{\strut@{\maketag@##\maketag@}}%
 \dimen@\displaywidth\advance\dimen@-\gwidth@
 \ifdim\dimen@>\tw@\wdz@\llap{\boxz@}\else
 \llap{\vtop{\normalbaselines\null\boxz@}}\fi
 \tabskip\z@skip\crcr#1\crcr\black@\gmaxwidth@}}
\newdimen\glineht@
\expandafter\def\csname gather \endcsname#1\endgather{\gmeasure@#1\endgather
 \ifdim\gmaxwidth@>\displaywidth\let\gdisplaywidth@\gmaxwidth@\else
 \let\gdisplaywidth@\displaywidth\fi\tabskip\centering@\halign to\displaywidth
 {\hfil\strut@\setboxz@h{$\m@th\displaystyle{##}$}%
 \global\gwidth@\wdz@\global\glineht@\ht\z@\boxz@\hfil&\kern-\gdisplaywidth@
 \setboxz@h{\strut@{\maketag@##\maketag@}}%
 \dimen@\displaywidth\advance\dimen@-\gwidth@
 \ifdim\dimen@>\tw@\wdz@\rlap{\boxz@}\else
 \rlap{\vbox{\normalbaselines\boxz@\vbox to\glineht@{}}}\fi
 \tabskip\gdisplaywidth@\crcr#1\crcr\black@\gmaxwidth@}}
\newif\ifctagsplit@
\def\CenteredTagsOnSplits{\global\ctagsplit@true}
\def\TopOrBottomTagsOnSplits{\global\ctagsplit@false}
\TopOrBottomTagsOnSplits
\def\split{\relax\ifinany@\let\next@\insplit@\else
 \ifmmode\ifinner\def\next@{\onlydmatherr@\split}\else
 \let\next@\outsplit@\fi\else
 \def\next@{\onlydmatherr@\split}\fi\fi\next@}
\def\insplit@{\global\setbox\z@\vbox\bgroup\vspace@\Let@\ialign\bgroup
 \hfil\strut@$\m@th\displaystyle{##}$&$\m@th\displaystyle{{}##}$\hfill\crcr}
\def\endsplit{\crcr\egroup\egroup\iftagsleft@\expandafter\lendsplit@\else
 \expandafter\rendsplit@\fi}
\def\rendsplit@{\global\setbox9 \vbox
 {\unvcopy\z@\global\setbox8 \lastbox\unskip}
 \setbox@ne\hbox{\unhcopy8 \unskip\global\setbox\tw@\lastbox
 \unskip\global\setbox\thr@@\lastbox}
 \global\setbox7 \hbox{\unhbox\tw@\unskip}
 \ifinalign@\ifctagsplit@                                                   
  \gdef\split@{\hbox to\wd\thr@@{}&
   \vcenter{\vbox{\moveleft\wd\thr@@\boxz@}}}
 \else\gdef\split@{&\vbox{\moveleft\wd\thr@@\box9}\crcr
  \box\thr@@&\box7}\fi                                                      
 \else                                                                      
  \ifctagsplit@\gdef\split@{\vcenter{\boxz@}}\else
  \gdef\split@{\box9\crcr\hbox{\box\thr@@\box7}}\fi
 \fi
 \split@}                                                                   
\def\lendsplit@{\global\setbox9\vtop{\unvcopy\z@}
 \setbox@ne\vbox{\unvcopy\z@\global\setbox8\lastbox}
 \setbox@ne\hbox{\unhcopy8\unskip\setbox\tw@\lastbox
  \unskip\global\setbox\thr@@\lastbox}
 \ifinalign@\ifctagsplit@                                                   
  \gdef\split@{\hbox to\wd\thr@@{}&
  \vcenter{\vbox{\moveleft\wd\thr@@\box9}}}
  \else                                                                     
  \gdef\split@{\hbox to\wd\thr@@{}&\vbox{\moveleft\wd\thr@@\box9}}\fi
 \else
  \ifctagsplit@\gdef\split@{\vcenter{\box9}}\else
  \gdef\split@{\box9}\fi
 \fi\split@}
\def\outsplit@#1$${\align\insplit@#1\endalign$$}
\newdimen\multlinegap@
\multlinegap@1em
\newdimen\multlinetaggap@
\multlinetaggap@1em
\def\MultlineGap#1{\global\multlinegap@#1\relax}
\def\multlinegap#1{\RIfMIfI@\onlydmatherr@\multlinegap\else
 \multlinegap@#1\relax\fi\else\onlydmatherr@\multlinegap\fi}
\def\nomultlinegap{\multlinegap{\z@}}
\def\multline{\RIfMIfI@
 \DN@{\onlydmatherr@\multline}\else
 \DN@{\multline@}\fi\else
 \DN@{\onlydmatherr@\multline}\fi\next@}
\newif\iftagin@
\def\tagin@#1{\tagin@false\in@\tag{#1}\ifin@\tagin@true\fi}
\def\multline@#1$${\inany@true\vspace@\allowdisplaybreak@\displaybreak@
 \tagin@{#1}\iftagsleft@\DN@{\multline@l#1$$}\else
 \DN@{\multline@r#1$$}\fi\next@}
\newdimen\mwidth@
\def\rmmeasure@#1\endmultline{%
 \def\shoveleft##1{##1}\def\shoveright##1{##1}
 \setbox@ne\vbox{\Let@\halign{\setboxz@h
  {$\m@th\@lign\displaystyle{}##$}\global\mwidth@\wdz@
  \crcr#1\crcr}}}
\newdimen\mlineht@
\newif\ifzerocr@
\newif\ifonecr@
\def\lmmeasure@#1\endmultline{\global\zerocr@true\global\onecr@false
 \everycr{\noalign{\ifonecr@\global\onecr@false\fi
  \ifzerocr@\global\zerocr@false\global\onecr@true\fi}}
  \def\shoveleft##1{##1}\def\shoveright##1{##1}%
 \setbox@ne\vbox{\Let@\halign{\setboxz@h
  {$\m@th\@lign\displaystyle{}##$}\ifonecr@\global\mwidth@\wdz@
  \global\mlineht@\ht\z@\fi\crcr#1\crcr}}}
\newbox\mtagbox@
\newdimen\ltwidth@
\newdimen\rtwidth@
\def\multline@l#1$${\iftagin@\DN@{\lmultline@@#1$$}\else
 \DN@{\setbox\mtagbox@\null\ltwidth@\z@\rtwidth@\z@
  \lmultline@@@#1$$}\fi\next@}
\def\lmultline@@#1\endmultline\tag#2$${%
 \setbox\mtagbox@\hbox{\maketag@#2\maketag@}
 \lmmeasure@#1\endmultline\dimen@\mwidth@\advance\dimen@\wd\mtagbox@
 \advance\dimen@\multlinetaggap@                                            
 \ifdim\dimen@>\displaywidth\ltwidth@\z@\else\ltwidth@\wd\mtagbox@\fi       
 \lmultline@@@#1\endmultline$$}
\def\lmultline@@@{\displ@y
 \def\shoveright##1{##1\hfilneg\hskip\multlinegap@}%
 \def\shoveleft##1{\setboxz@h{$\m@th\displaystyle{}##1$}%
  \setbox@ne\hbox{$\m@th\displaystyle##1$}%
  \hfilneg
  \iftagin@
   \ifdim\ltwidth@>\z@\hskip\ltwidth@\hskip\multlinetaggap@\fi
  \else\hskip\multlinegap@\fi\hskip.5\wd@ne\hskip-.5\wdz@##1}
  \halign\bgroup\Let@\hbox to\displaywidth
   {\strut@$\m@th\displaystyle\hfil{}##\hfil$}\crcr
   \hfilneg                                                                 
   \iftagin@                                                                
    \ifdim\ltwidth@>\z@                                                     
     \box\mtagbox@\hskip\multlinetaggap@                                    
    \else
     \rlap{\vbox{\normalbaselines\hbox{\strut@\box\mtagbox@}%
     \vbox to\mlineht@{}}}\fi                                               
   \else\hskip\multlinegap@\fi}                                             
\def\multline@r#1$${\iftagin@\DN@{\rmultline@@#1$$}\else
 \DN@{\setbox\mtagbox@\null\ltwidth@\z@\rtwidth@\z@
  \rmultline@@@#1$$}\fi\next@}
\def\rmultline@@#1\endmultline\tag#2$${\ltwidth@\z@
 \setbox\mtagbox@\hbox{\maketag@#2\maketag@}%
 \rmmeasure@#1\endmultline\dimen@\mwidth@\advance\dimen@\wd\mtagbox@
 \advance\dimen@\multlinetaggap@
 \ifdim\dimen@>\displaywidth\rtwidth@\z@\else\rtwidth@\wd\mtagbox@\fi
 \rmultline@@@#1\endmultline$$}
\def\rmultline@@@{\displ@y
 \def\shoveright##1{##1\hfilneg\iftagin@\ifdim\rtwidth@>\z@
  \hskip\rtwidth@\hskip\multlinetaggap@\fi\else\hskip\multlinegap@\fi}%
 \def\shoveleft##1{\setboxz@h{$\m@th\displaystyle{}##1$}%
  \setbox@ne\hbox{$\m@th\displaystyle##1$}%
  \hfilneg\hskip\multlinegap@\hskip.5\wd@ne\hskip-.5\wdz@##1}%
 \halign\bgroup\Let@\hbox to\displaywidth
  {\strut@$\m@th\displaystyle\hfil{}##\hfil$}\crcr
 \hfilneg\hskip\multlinegap@}
\def\endmultline{\iftagsleft@\expandafter\lendmultline@\else
 \expandafter\rendmultline@\fi}
\def\lendmultline@{\hfilneg\hskip\multlinegap@\crcr\egroup}
\def\rendmultline@{\iftagin@                                                
 \ifdim\rtwidth@>\z@                                                        
  \hskip\multlinetaggap@\box\mtagbox@                                       
 \else\llap{\vtop{\normalbaselines\null\hbox{\strut@\box\mtagbox@}}}\fi     
 \else\hskip\multlinegap@\fi                                                
 \hfilneg\crcr\egroup}
\def\bmod{\mskip-\medmuskip\mkern5mu\mathbin{\fam\z@ mod}\penalty900
 \mkern5mu\mskip-\medmuskip}
\def\pmod#1{\allowbreak\ifinner\mkern8mu\else\mkern18mu\fi
 ({\fam\z@ mod}\,\,#1)}
\def\pod#1{\allowbreak\ifinner\mkern8mu\else\mkern18mu\fi(#1)}
\def\mod#1{\allowbreak\ifinner\mkern12mu\else\mkern18mu\fi{\fam\z@ mod}\,\,#1}
\message{continued fractions,}
\newcount\cfraccount@
\def\cfrac{\bgroup\bgroup\advance\cfraccount@\@ne\strut
 \iffalse{\fi\def\\{\over\displaystyle}\iffalse}\fi}
\def\lcfrac{\bgroup\bgroup\advance\cfraccount@\@ne\strut
 \iffalse{\fi\def\\{\hfill\over\displaystyle}\iffalse}\fi}
\def\rcfrac{\bgroup\bgroup\advance\cfraccount@\@ne\strut\hfill
 \iffalse{\fi\def\\{\over\displaystyle}\iffalse}\fi}
\def\gloop@#1\repeat{\gdef\body{#1}\iterate}
\def\endcfrac{\gloop@\ifnum\cfraccount@>\z@\global\advance\cfraccount@\m@ne
 \egroup\hskip-\nulldelimiterspace\egroup\repeat}
\message{compound symbols,}
\def\binrel@#1{\setboxz@h{\thinmuskip0mu
  \medmuskip\m@ne mu\thickmuskip\@ne mu$#1\m@th$}%
 \setbox@ne\hbox{\thinmuskip0mu\medmuskip\m@ne mu\thickmuskip
  \@ne mu${}#1{}\m@th$}%
 \setbox\tw@\hbox{\hskip\wd@ne\hskip-\wdz@}}
\def\overset#1\to#2{\binrel@{#2}\ifdim\wd\tw@<\z@
 \mathbin{\mathop{\kern\z@#2}\limits^{#1}}\else\ifdim\wd\tw@>\z@
 \mathrel{\mathop{\kern\z@#2}\limits^{#1}}\else
 {\mathop{\kern\z@#2}\limits^{#1}}{}\fi\fi}
\def\underset#1\to#2{\binrel@{#2}\ifdim\wd\tw@<\z@
 \mathbin{\mathop{\kern\z@#2}\limits_{#1}}\else\ifdim\wd\tw@>\z@
 \mathrel{\mathop{\kern\z@#2}\limits_{#1}}\else
 {\mathop{\kern\z@#2}\limits_{#1}}{}\fi\fi}
\def\oversetbrace#1\to#2{\overbrace{#2}^{#1}}
\def\undersetbrace#1\to#2{\underbrace{#2}_{#1}}
\def\sideset#1\and#2\to#3{%
 \setbox@ne\hbox{$\dsize{\vphantom{#3}}#1{#3}\m@th$}%
 \setbox\tw@\hbox{$\dsize{#3}#2\m@th$}%
 \hskip\wd@ne\hskip-\wd\tw@\mathop{\hskip\wd\tw@\hskip-\wd@ne
  {\vphantom{#3}}#1{#3}#2}}
\def\rightarrowfill@#1{\setboxz@h{$#1-\m@th$}\ht\z@\z@
  $#1\m@th\copy\z@\mkern-6mu\cleaders
  \hbox{$#1\mkern-2mu\box\z@\mkern-2mu$}\hfill
  \mkern-6mu\mathord\rightarrow$}
\def\leftarrowfill@#1{\setboxz@h{$#1-\m@th$}\ht\z@\z@
  $#1\m@th\mathord\leftarrow\mkern-6mu\cleaders
  \hbox{$#1\mkern-2mu\copy\z@\mkern-2mu$}\hfill
  \mkern-6mu\box\z@$}
\def\leftrightarrowfill@#1{\setboxz@h{$#1-\m@th$}\ht\z@\z@
  $#1\m@th\mathord\leftarrow\mkern-6mu\cleaders
  \hbox{$#1\mkern-2mu\box\z@\mkern-2mu$}\hfill
  \mkern-6mu\mathord\rightarrow$}
\def\overrightarrow{\mathpalette\overrightarrow@}
\def\overrightarrow@#1#2{\vbox{\ialign{##\crcr\rightarrowfill@#1\crcr
 \noalign{\kern-\ex@\nointerlineskip}$\m@th\hfil#1#2\hfil$\crcr}}}

\def\overleftarrow{\mathpalette\overleftarrow@}
\def\overleftarrow@#1#2{\vbox{\ialign{##\crcr\leftarrowfill@#1\crcr
 \noalign{\kern-\ex@\nointerlineskip}$\m@th\hfil#1#2\hfil$\crcr}}}
\def\overleftrightarrow{\mathpalette\overleftrightarrow@}
\def\overleftrightarrow@#1#2{\vbox{\ialign{##\crcr\leftrightarrowfill@#1\crcr
 \noalign{\kern-\ex@\nointerlineskip}$\m@th\hfil#1#2\hfil$\crcr}}}
\def\underrightarrow{\mathpalette\underrightarrow@}
\def\underrightarrow@#1#2{\vtop{\ialign{##\crcr$\m@th\hfil#1#2\hfil$\crcr
 \noalign{\nointerlineskip}\rightarrowfill@#1\crcr}}}

\def\underleftarrow{\mathpalette\underleftarrow@}
\def\underleftarrow@#1#2{\vtop{\ialign{##\crcr$\m@th\hfil#1#2\hfil$\crcr
 \noalign{\nointerlineskip}\leftarrowfill@#1\crcr}}}
\def\underleftrightarrow{\mathpalette\underleftrightarrow@}
\def\underleftrightarrow@#1#2{\vtop{\ialign{##\crcr$\m@th\hfil#1#2\hfil$\crcr
 \noalign{\nointerlineskip}\leftrightarrowfill@#1\crcr}}}
\message{various kinds of dots,}
\let\DOTSI\relax
\let\DOTSB\relax

\newif\ifmath@
{\uccode`7=`\\ \uccode`8=`m \uccode`9=`a \uccode`0=`t \uccode`!=`h
 \uppercase{\gdef\math@#1#2#3#4#5#6\math@{\global\math@false\ifx 7#1\ifx 8#2%
 \ifx 9#3\ifx 0#4\ifx !#5\xdef\meaning@{#6}\global\math@true\fi\fi\fi\fi\fi}}}
\newif\ifmathch@
{\uccode`7=`c \uccode`8=`h \uccode`9=`\"
 \uppercase{\gdef\mathch@#1#2#3#4#5#6\mathch@{\global\mathch@false
  \ifx 7#1\ifx 8#2\ifx 9#5\global\mathch@true\xdef\meaning@{9#6}\fi\fi\fi}}}
\newcount\classnum@
\def\getmathch@#1.#2\getmathch@{\classnum@#1 \divide\classnum@4096
 \ifcase\number\classnum@\or\or\gdef\thedots@{\dotsb@}\or
 \gdef\thedots@{\dotsb@}\fi}
\newif\ifmathbin@
{\uccode`4=`b \uccode`5=`i \uccode`6=`n
 \uppercase{\gdef\mathbin@#1#2#3{\relaxnext@
  \DNii@##1\mathbin@{\ifx\space@\next\global\mathbin@true\fi}%
 \global\mathbin@false\DN@##1\mathbin@{}%
 \ifx 4#1\ifx 5#2\ifx 6#3\DN@{\FN@\nextii@}\fi\fi\fi\next@}}}
\newif\ifmathrel@
{\uccode`4=`r \uccode`5=`e \uccode`6=`l
 \uppercase{\gdef\mathrel@#1#2#3{\relaxnext@
  \DNii@##1\mathrel@{\ifx\space@\next\global\mathrel@true\fi}%
 \global\mathrel@false\DN@##1\mathrel@{}%
 \ifx 4#1\ifx 5#2\ifx 6#3\DN@{\FN@\nextii@}\fi\fi\fi\next@}}}
\newif\ifmacro@
{\uccode`5=`m \uccode`6=`a \uccode`7=`c
 \uppercase{\gdef\macro@#1#2#3#4\macro@{\global\macro@false
  \ifx 5#1\ifx 6#2\ifx 7#3\global\macro@true
  \xdef\meaning@{\macro@@#4\macro@@}\fi\fi\fi}}}
\def\macro@@#1->#2\macro@@{#2}
\newif\ifDOTS@
\newcount\DOTSCASE@
{\uccode`6=`\\ \uccode`7=`D \uccode`8=`O \uccode`9=`T \uccode`0=`S
 \uppercase{\gdef\DOTS@#1#2#3#4#5{\global\DOTS@false\DN@##1\DOTS@{}%
  \ifx 6#1\ifx 7#2\ifx 8#3\ifx 9#4\ifx 0#5\let\next@\DOTS@@\fi\fi\fi\fi\fi
  \next@}}}
{\uccode`3=`B \uccode`4=`I \uccode`5=`X
 \uppercase{\gdef\DOTS@@#1{\relaxnext@
  \DNii@##1\DOTS@{\ifx\space@\next\global\DOTS@true\fi}%
  \DN@{\FN@\nextii@}%
  \ifx 3#1\global\DOTSCASE@\z@\else
  \ifx 4#1\global\DOTSCASE@\@ne\else
  \ifx 5#1\global\DOTSCASE@\tw@\else\DN@##1\DOTS@{}%
  \fi\fi\fi\next@}}}
\newif\ifnot@
{\uccode`5=`\\ \uccode`6=`n \uccode`7=`o \uccode`8=`t
 \uppercase{\gdef\not@#1#2#3#4{\relaxnext@
  \DNii@##1\not@{\ifx\space@\next\global\not@true\fi}%
 \global\not@false\DN@##1\not@{}%
 \ifx 5#1\ifx 6#2\ifx 7#3\ifx 8#4\DN@{\FN@\nextii@}\fi\fi\fi
 \fi\next@}}}
\newif\ifkeybin@
\def\keybin@{\keybin@true
 \ifx\next+\else\ifx\next=\else\ifx\next<\else\ifx\next>\else\ifx\next-\else
 \ifx\next*\else\ifx\next:\else\keybin@false\fi\fi\fi\fi\fi\fi\fi}
\def\dots{\RIfM@\expandafter\mdots@\else\expandafter\tdots@\fi}
\def\tdots@{\unskip\relaxnext@
 \DN@{$\m@th\mathinner{\ldotp\ldotp\ldotp}\,
   \ifx\next,\,$\else\ifx\next.\,$\else\ifx\next;\,$\else\ifx\next:\,$\else
   \ifx\next?\,$\else\ifx\next!\,$\else$ \fi\fi\fi\fi\fi\fi}%
 \ \FN@\next@}
\def\mdots@{\FN@\mdots@@}
\def\mdots@@{\gdef\thedots@{\dotso@}
 \ifx\next\boldkey\gdef\thedots@\boldkey{\boldkeydots@}\else                
 \ifx\next\boldsymbol\gdef\thedots@\boldsymbol{\boldsymboldots@}\else       
 \ifx,\next\gdef\thedots@{\dotsc}
 \else\ifx\not\next\gdef\thedots@{\dotsb@}
 \else\keybin@
 \ifkeybin@\gdef\thedots@{\dotsb@}
 \else\xdef\meaning@{\meaning\next..........}\xdef\meaning@@{\meaning@}
  \expandafter\math@\meaning@\math@
  \ifmath@
   \expandafter\mathch@\meaning@\mathch@
   \ifmathch@\expandafter\getmathch@\meaning@\getmathch@\fi                 
  \else\expandafter\macro@\meaning@@\macro@                                 
  \ifmacro@                                                                
   \expandafter\not@\meaning@\not@\ifnot@\gdef\thedots@{\dotsb@}
  \else\expandafter\DOTS@\meaning@\DOTS@
  \ifDOTS@
   \ifcase\number\DOTSCASE@\gdef\thedots@{\dotsb@}%
    \or\gdef\thedots@{\dotsi}\else\fi                                      
  \else\expandafter\math@\meaning@\math@                                   
  \ifmath@\expandafter\mathbin@\meaning@\mathbin@
  \ifmathbin@\gdef\thedots@{\dotsb@}
  \else\expandafter\mathrel@\meaning@\mathrel@
  \ifmathrel@\gdef\thedots@{\dotsb@}
  \fi\fi\fi\fi\fi\fi\fi\fi\fi\fi\fi\fi
 \thedots@}
\def\plainldots@{\mathinner{\ldotp\ldotp\ldotp}}
\def\plaincdots@{\mathinner{\cdotp\cdotp\cdotp}}
\def\dotsi{\!\plaincdots@}
\let\dotsb@\plaincdots@
\newif\ifextra@
\newif\ifrightdelim@
\def\rightdelim@{\global\rightdelim@true                                    
 \ifx\next)\else                                                            
 \ifx\next]\else
 \ifx\next\rbrack\else
 \ifx\next\}\else
 \ifx\next\rbrace\else
 \ifx\next\rangle\else
 \ifx\next\rceil\else
 \ifx\next\rfloor\else
 \ifx\next\rgroup\else
 \ifx\next\rmoustache\else
 \ifx\next\right\else
 \ifx\next\bigr\else
 \ifx\next\biggr\else
 \ifx\next\Bigr\else                                                        
 \ifx\next\Biggr\else\global\rightdelim@false
 \fi\fi\fi\fi\fi\fi\fi\fi\fi\fi\fi\fi\fi\fi\fi}
\def\extra@{%
 \global\extra@false\rightdelim@\ifrightdelim@\global\extra@true            
 \else\ifx\next$\global\extra@true                                          
 \else\xdef\meaning@{\meaning\next..........}
 \expandafter\macro@\meaning@\macro@\ifmacro@                               
 \expandafter\DOTS@\meaning@\DOTS@
 \ifDOTS@
 \ifnum\DOTSCASE@=\tw@\global\extra@true                                    
 \fi\fi\fi\fi\fi}
\newif\ifbold@
\def\dotso@{\relaxnext@
 \ifbold@
  \let\next\delayed@
  \DNii@{\extra@\plainldots@\ifextra@\,\fi}%
 \else
  \DNii@{\DN@{\extra@\plainldots@\ifextra@\,\fi}\FN@\next@}%
 \fi
 \nextii@}
\def\extrap@#1{%
 \ifx\next,\DN@{#1\,}\else
 \ifx\next;\DN@{#1\,}\else
 \ifx\next.\DN@{#1\,}\else\extra@
 \ifextra@\DN@{#1\,}\else
 \let\next@#1\fi\fi\fi\fi\next@}
\def\ldots{\DN@{\extrap@\plainldots@}%
 \FN@\next@}
\def\cdots{\DN@{\extrap@\plaincdots@}%
 \FN@\next@}

\def\dotsc{\relaxnext@
 \DN@{\ifx\next;\plainldots@\,\else
  \ifx\next.\plainldots@\,\else\extra@\plainldots@
  \ifextra@\,\fi\fi\fi}%
 \FN@\next@}
\def\cdot{\mathchar"2201 }

\def\mapsto{\DOTSB\mapstochar\rightarrow}

\message{special superscripts,}
\def\dddot#1{{\mathop{#1}\limits^{\vbox to-1.4\ex@{\kern-\tw@\ex@
 \hbox{\rm...}\vss}}}}
\def\ddddot#1{{\mathop{#1}\limits^{\vbox to-1.4\ex@{\kern-\tw@\ex@
 \hbox{\rm....}\vss}}}}
\def\sphat{^{\mathchoice{}{}%
 {\,\,\botsmash{\hbox{\lower4\ex@\hbox{$\m@th\widehat{\null}$}}}}%
 {\,\botsmash{\hbox{\lower3\ex@\hbox{$\m@th\hat{\null}$}}}}}}

\def\spacute{^{\!\botsmash{\hbox{\lower\@ne ex\hbox{\'{}}}}}}
\def\spgrave{^{\mathchoice{}{}{}{\!}%
 \botsmash{\hbox{\lower\@ne ex\hbox{\`{}}}}}}
\def\spdot{^{\hbox{\raise\ex@\hbox{\rm.}}}}
\def\spddot{^{\hbox{\raise\ex@\hbox{\rm..}}}}
\def\spdddot{^{\hbox{\raise\ex@\hbox{\rm...}}}}
\def\spddddot{^{\hbox{\raise\ex@\hbox{\rm....}}}}
\def\spbreve{^{\!\botsmash{\hbox{\lower4\ex@\hbox{\u{}}}}}}

\message{\string\text,}
\def\textonlyfont@#1#2{\def#1{\RIfM@
 \Err@{Use \string#1\space only in text}\else#2\fi}}
\textonlyfont@\rm\tenrm
\textonlyfont@\it\tenit
\textonlyfont@\sl\tensl
\textonlyfont@\bf\tenbf
\def\oldnos#1{\RIfM@{\mathcode`\,="013B \fam\@ne#1}\else
 \leavevmode\hbox{$\m@th\mathcode`\,="013B \fam\@ne#1$}\fi}
\def\text{\RIfM@\expandafter\text@\else\expandafter\text@@\fi}
\def\text@@#1{\leavevmode\hbox{#1}}
\def\mathhexbox@#1#2#3{\text{$\m@th\mathchar"#1#2#3$}}
\def\dag{{\mathhexbox@279}}
\def\ddag{{\mathhexbox@27A}}
\def\S{{\mathhexbox@278}}
\def\P{{\mathhexbox@27B}}
\newif\iffirstchoice@
\firstchoice@true
\def\text@#1{\mathchoice
 {\hbox{\everymath{\displaystyle}\def\textfonti{\the\textfont\@ne}%
  \def\textfontii{\the\textfont\tw@}\textdef@@ T#1}}
 {\hbox{\firstchoice@false
  \everymath{\textstyle}\def\textfonti{\the\textfont\@ne}%
  \def\textfontii{\the\textfont\tw@}\textdef@@ T#1}}
 {\hbox{\firstchoice@false
  \everymath{\scriptstyle}\def\textfonti{\the\scriptfont\@ne}%
  \def\textfontii{\the\scriptfont\tw@}\textdef@@ S\rm#1}}
 {\hbox{\firstchoice@false
  \everymath{\scriptscriptstyle}\def\textfonti
  {\the\scriptscriptfont\@ne}%
  \def\textfontii{\the\scriptscriptfont\tw@}\textdef@@ s\rm#1}}}
\def\textdef@@#1{\textdef@#1\rm\textdef@#1\bf\textdef@#1\sl\textdef@#1\it}
\def\rmfam{0}
\def\textdef@#1#2{%
 \DN@{\csname\expandafter\eat@\string#2fam\endcsname}%
 \if S#1\edef#2{\the\scriptfont\next@\relax}%
 \else\if s#1\edef#2{\the\scriptscriptfont\next@\relax}%
 \else\edef#2{\the\textfont\next@\relax}\fi\fi}
\scriptfont\itfam\tenit \scriptscriptfont\itfam\tenit
\scriptfont\slfam\tensl \scriptscriptfont\slfam\tensl
\newif\iftopfolded@
\newif\ifbotfolded@
\def\topfoldedtext{\topfolded@true\botfolded@false\foldedtext@}
\def\botfoldedtext{\botfolded@true\topfolded@false\foldedtext@}
\def\foldedtext{\topfolded@false\botfolded@false\foldedtext@}
\Invalid@\foldedwidth
\def\foldedtext@{\relaxnext@
 \DN@{\ifx\next\foldedwidth\let\next@\nextii@\else
  \DN@{\nextii@\foldedwidth{.3\hsize}}\fi\next@}%
 \DNii@\foldedwidth##1##2{\setbox\z@\vbox
  {\normalbaselines\hsize##1\relax
  \tolerance1600 \noindent\ignorespaces##2}\ifbotfolded@\boxz@\else
  \iftopfolded@\vtop{\unvbox\z@}\else\vcenter{\boxz@}\fi\fi}%
 \FN@\next@}
\message{math font commands,}
\def\bold{\RIfM@\expandafter\bold@\else
 \expandafter\nonmatherr@\expandafter\bold\fi}
\def\bold@#1{{\bold@@{#1}}}
\def\bold@@#1{\fam\bffam\relax#1}
\def\slanted{\RIfM@\expandafter\slanted@\else
 \expandafter\nonmatherr@\expandafter\slanted\fi}
\def\slanted@#1{{\slanted@@{#1}}}
\def\slanted@@#1{\fam\slfam\relax#1}
\def\roman{\RIfM@\expandafter\roman@\else
 \expandafter\nonmatherr@\expandafter\roman\fi}
\def\roman@#1{{\roman@@{#1}}}
\def\roman@@#1{\fam\rmfam\relax#1}
\def\italic{\RIfM@\expandafter\italic@\else
 \expandafter\nonmatherr@\expandafter\italic\fi}
\def\italic@#1{{\italic@@{#1}}}
\def\italic@@#1{\fam\itfam\relax#1}
\def\Cal{\RIfM@\expandafter\Cal@\else
 \expandafter\nonmatherr@\expandafter\Cal\fi}
\def\Cal@#1{{\Cal@@{#1}}}
\def\Cal@@#1{\noaccents@\fam\tw@#1}
\mathchardef\Gamma="0000
\mathchardef\Delta="0001
\mathchardef\Theta="0002
\mathchardef\Lambda="0003
\mathchardef\Xi="0004
\mathchardef\Pi="0005
\mathchardef\Sigma="0006
\mathchardef\Upsilon="0007
\mathchardef\Phi="0008
\mathchardef\Psi="0009
\mathchardef\Omega="000A
\mathchardef\varGamma="0100
\mathchardef\varDelta="0101
\mathchardef\varTheta="0102
\mathchardef\varLambda="0103
\mathchardef\varXi="0104
\mathchardef\varPi="0105
\mathchardef\varSigma="0106
\mathchardef\varUpsilon="0107
\mathchardef\varPhi="0108
\mathchardef\varPsi="0109
\mathchardef\varOmega="010A
\let\alloc@@\alloc@
\def\hexnumber@#1{\ifcase#1 0\or 1\or 2\or 3\or 4\or 5\or 6\or 7\or 8\or
 9\or A\or B\or C\or D\or E\or F\fi}
\def\loadmsam{%
 \font@\tenmsa=msam10
 \font@\sevenmsa=msam7
 \font@\fivemsa=msam5
 \alloc@@8\fam\chardef\sixt@@n\msafam
 \textfont\msafam=\tenmsa
 \scriptfont\msafam=\sevenmsa
 \scriptscriptfont\msafam=\fivemsa
 \edef\next{\hexnumber@\msafam}%
 \mathchardef\dabar@"0\next39
 \edef\dashrightarrow{\mathrel{\dabar@\dabar@\mathchar"0\next4B}}%
 \edef\dashleftarrow{\mathrel{\mathchar"0\next4C\dabar@\dabar@}}%
 \let\dasharrow\dashrightarrow
 \edef\ulcorner{\delimiter"4\next70\next70 }%
 \edef\urcorner{\delimiter"5\next71\next71 }%
 \edef\llcorner{\delimiter"4\next78\next78 }%
 \edef\lrcorner{\delimiter"5\next79\next79 }%
 \edef\yen{{\noexpand\mathhexbox@\next55}}%
 \edef\checkmark{{\noexpand\mathhexbox@\next58}}%
 \edef\circledR{{\noexpand\mathhexbox@\next72}}%
 \edef\maltese{{\noexpand\mathhexbox@\next7A}}%
 \global\let\loadmsam\empty}%
\def\loadmsbm{%
 \font@\tenmsb=msbm10 \font@\sevenmsb=msbm7 \font@\fivemsb=msbm5
 \alloc@@8\fam\chardef\sixt@@n\msbfam
 \textfont\msbfam=\tenmsb
 \scriptfont\msbfam=\sevenmsb \scriptscriptfont\msbfam=\fivemsb
 \global\let\loadmsbm\empty
 }
\def\widehat#1{\ifx\undefined\msbfam \DN@{362}%
  \else \setboxz@h{$\m@th#1$}%
    \edef\next@{\ifdim\wdz@>\tw@ em%
        \hexnumber@\msbfam 5B%
      \else 362\fi}\fi
  \mathaccent"0\next@{#1}}
\def\widetilde#1{\ifx\undefined\msbfam \DN@{365}%
  \else \setboxz@h{$\m@th#1$}%
    \edef\next@{\ifdim\wdz@>\tw@ em%
        \hexnumber@\msbfam 5D%
      \else 365\fi}\fi
  \mathaccent"0\next@{#1}}
\message{\string\newsymbol,}
\def\newsymbol#1#2#3#4#5{\define#1{}%
  \count@#2\relax \advance\count@\m@ne 
 \ifcase\count@
   \ifx\undefined\msafam\loadmsam\fi \let\next@\msafam
 \or \ifx\undefined\msbfam\loadmsbm\fi \let\next@\msbfam
 \else  \Err@{\Invalid@@\string\newsymbol}\let\next@\tw@\fi
 \mathchardef#1="#3\hexnumber@\next@#4#5\space}

\def\Bbb{\RIfM@\expandafter\Bbb@\else
 \expandafter\nonmatherr@\expandafter\Bbb\fi}
\def\Bbb@#1{{\Bbb@@{#1}}}
\def\Bbb@@#1{\noaccents@\fam\msbfam\relax#1}
\message{bold Greek and bold symbols,}
\def\loadbold{%
 \font@\tencmmib=cmmib10 \font@\sevencmmib=cmmib7 \font@\fivecmmib=cmmib5
 \skewchar\tencmmib'177 \skewchar\sevencmmib'177 \skewchar\fivecmmib'177
 \alloc@@8\fam\chardef\sixt@@n\cmmibfam
 \textfont\cmmibfam\tencmmib
 \scriptfont\cmmibfam\sevencmmib \scriptscriptfont\cmmibfam\fivecmmib
 \font@\tencmbsy=cmbsy10 \font@\sevencmbsy=cmbsy7 \font@\fivecmbsy=cmbsy5
 \skewchar\tencmbsy'60 \skewchar\sevencmbsy'60 \skewchar\fivecmbsy'60
 \alloc@@8\fam\chardef\sixt@@n\cmbsyfam
 \textfont\cmbsyfam\tencmbsy
 \scriptfont\cmbsyfam\sevencmbsy \scriptscriptfont\cmbsyfam\fivecmbsy
 \let\loadbold\empty
}
\def\boldnotloaded#1{\Err@{\ifcase#1\or First\else Second\fi
       bold symbol font not loaded}}
\def\mathchari@#1#2#3{\ifx\undefined\cmmibfam
    \boldnotloaded@\@ne
  \else\mathchar"#1\hexnumber@\cmmibfam#2#3\space \fi}
\def\mathcharii@#1#2#3{\ifx\undefined\cmbsyfam
    \boldnotloaded\tw@
  \else \mathchar"#1\hexnumber@\cmbsyfam#2#3\space\fi}
\edef\bffam@{\hexnumber@\bffam}
\def\boldkey#1{\ifcat\noexpand#1A%
  \ifx\undefined\cmmibfam \boldnotloaded\@ne
  \else {\fam\cmmibfam#1}\fi
 \else
 \ifx#1!\mathchar"5\bffam@21 \else
 \ifx#1(\mathchar"4\bffam@28 \else\ifx#1)\mathchar"5\bffam@29 \else
 \ifx#1+\mathchar"2\bffam@2B \else\ifx#1:\mathchar"3\bffam@3A \else
 \ifx#1;\mathchar"6\bffam@3B \else\ifx#1=\mathchar"3\bffam@3D \else
 \ifx#1?\mathchar"5\bffam@3F \else\ifx#1[\mathchar"4\bffam@5B \else
 \ifx#1]\mathchar"5\bffam@5D \else
 \ifx#1,\mathchari@63B \else
 \ifx#1-\mathcharii@200 \else
 \ifx#1.\mathchari@03A \else
 \ifx#1/\mathchari@03D \else
 \ifx#1<\mathchari@33C \else
 \ifx#1>\mathchari@33E \else
 \ifx#1*\mathcharii@203 \else
 \ifx#1|\mathcharii@06A \else
 \ifx#10\bold0\else\ifx#11\bold1\else\ifx#12\bold2\else\ifx#13\bold3\else
 \ifx#14\bold4\else\ifx#15\bold5\else\ifx#16\bold6\else\ifx#17\bold7\else
 \ifx#18\bold8\else\ifx#19\bold9\else
  \Err@{\string\boldkey\space can't be used with #1}%
 \fi\fi\fi\fi\fi\fi\fi\fi\fi\fi\fi\fi\fi\fi\fi
 \fi\fi\fi\fi\fi\fi\fi\fi\fi\fi\fi\fi\fi\fi}
\def\boldsymbol#1{%
 \DN@{\Err@{You can't use \string\boldsymbol\space with \string#1}#1}%
 \ifcat\noexpand#1A%
   \let\next@\relax
   \ifx\undefined\cmmibfam \boldnotloaded\@ne
   \else {\fam\cmmibfam#1}\fi
 \else
  \xdef\meaning@{\meaning#1.........}%
  \expandafter\math@\meaning@\math@
  \ifmath@
   \expandafter\mathch@\meaning@\mathch@
   \ifmathch@
    \expandafter\boldsymbol@@\meaning@\boldsymbol@@
   \fi
  \else
   \expandafter\macro@\meaning@\macro@
   \expandafter\delim@\meaning@\delim@
   \ifdelim@
    \expandafter\delim@@\meaning@\delim@@
   \else
    \boldsymbol@{#1}%
   \fi
  \fi
 \fi
 \next@}
\def\mathhexboxii@#1#2{\ifx\undefined\cmbsyfam
    \boldnotloaded\tw@
  \else \mathhexbox@{\hexnumber@\cmbsyfam}{#1}{#2}\fi}
\def\boldsymbol@#1{\let\next@\relax\let\next#1%
 \ifx\next\cdot\mathcharii@201 \else
 \ifx\next\prime{{\null\mathcharii@030 \null}}\else
 \ifx\next\lbrack\mathchar"4\bffam@5B \else
 \ifx\next\rbrack\mathchar"5\bffam@5D \else
 \ifx\next\{\mathcharii@466 \else
 \ifx\next\lbrace\mathcharii@466 \else
 \ifx\next\}\mathcharii@567 \else
 \ifx\next\rbrace\mathcharii@567 \else
 \ifx\next\surd{{\mathcharii@170}}\else
 \ifx\next\S{{\mathhexboxii@78}}\else
 \ifx\next\P{{\mathhexboxii@7B}}\else
 \ifx\next\dag{{\mathhexboxii@79}}\else
 \ifx\next\ddag{{\mathhexboxii@7A}}\else
 \DN@{\Err@{You can't use \string\boldsymbol\space with \string#1}#1}%
 \fi\fi\fi\fi\fi\fi\fi\fi\fi\fi\fi\fi\fi}
\def\boldsymbol@@#1.#2\boldsymbol@@{\classnum@#1 \count@@@\classnum@        
 \divide\classnum@4096 \count@\classnum@                                    
 \multiply\count@4096 \advance\count@@@-\count@ \count@@\count@@@           
 \divide\count@@@\@cclvi \count@\count@@                                    
 \multiply\count@@@\@cclvi \advance\count@@-\count@@@                       
 \divide\count@@@\@cclvi                                                    
 \multiply\classnum@4096 \advance\classnum@\count@@                         
 \ifnum\count@@@=\z@                                                        
  \count@"\bffam@ \multiply\count@\@cclvi
  \advance\classnum@\count@
  \DN@{\mathchar\number\classnum@}%
 \else
  \ifnum\count@@@=\@ne                                                      
   \ifx\undefined\cmmibfam \DN@{\boldnotloaded\@ne}%
   \else \count@\cmmibfam \multiply\count@\@cclvi
     \advance\classnum@\count@
     \DN@{\mathchar\number\classnum@}\fi
  \else
   \ifnum\count@@@=\tw@                                                    
     \ifx\undefined\cmbsyfam
       \DN@{\boldnotloaded\tw@}%
     \else
       \count@\cmbsyfam \multiply\count@\@cclvi
       \advance\classnum@\count@
       \DN@{\mathchar\number\classnum@}%
     \fi
  \fi
 \fi
\fi}
\newif\ifdelim@
\newcount\delimcount@
{\uccode`6=`\\ \uccode`7=`d \uccode`8=`e \uccode`9=`l
 \uppercase{\gdef\delim@#1#2#3#4#5\delim@
  {\delim@false\ifx 6#1\ifx 7#2\ifx 8#3\ifx 9#4\delim@true
   \xdef\meaning@{#5}\fi\fi\fi\fi}}}
\def\delim@@#1"#2#3#4#5#6\delim@@{\if#32%
\let\next@\relax
 \ifx\undefined\cmbsyfam \boldnotloaded\@ne
 \else \mathcharii@#2#4#5\space \fi\fi}
\def\vert{\delimiter"026A30C }
\def\Vert{\delimiter"026B30D }
\let\|\Vert
\def\backslash{\delimiter"026E30F }
\def\boldkeydots@#1{\bold@true\let\next=#1\let\delayed@=#1\mdots@@
 \boldkey#1\bold@false}  
\def\boldsymboldots@#1{\bold@true\let\next#1\let\delayed@#1\mdots@@
 \boldsymbol#1\bold@false}
\message{Euler fonts,}

\def\frak{\mathfont@\frak}

\def\loadmathfont#1{%
   \expandafter\font@\csname ten#1\endcsname=#110
   \expandafter\font@\csname seven#1\endcsname=#17
   \expandafter\font@\csname five#1\endcsname=#15
   \edef\next{\noexpand\alloc@@8\fam\chardef\sixt@@n
     \expandafter\noexpand\csname#1fam\endcsname}%
   \next
   \textfont\csname#1fam\endcsname \csname ten#1\endcsname
   \scriptfont\csname#1fam\endcsname \csname seven#1\endcsname
   \scriptscriptfont\csname#1fam\endcsname \csname five#1\endcsname
   \expandafter\def\csname #1\expandafter\endcsname\expandafter{%
      \expandafter\mathfont@\csname#1\endcsname}%
 \expandafter\gdef\csname load#1\endcsname{}%
}
\def\mathfont@#1{\RIfM@\expandafter\mathfont@@\expandafter#1\else
  \expandafter\nonmatherr@\expandafter#1\fi}
\def\mathfont@@#1#2{{\mathfont@@@#1{#2}}}
\def\mathfont@@@#1#2{\noaccents@
   \fam\csname\expandafter\eat@\string#1fam\endcsname
   \relax#2}
\message{math accents,}
\def\accentclass@{7}
\def\noaccents@{\def\accentclass@{0}}
\def\makeacc@#1#2{\def#1{\mathaccent"\accentclass@#2 }}
\makeacc@\hat{05E}
\makeacc@\check{014}
\makeacc@\tilde{07E}
\makeacc@\acute{013}
\makeacc@\grave{012}
\makeacc@\dot{05F}
\makeacc@\ddot{07F}
\makeacc@\breve{015}
\makeacc@\bar{016}

\newcount\skewcharcount@
\newcount\familycount@
\def\theskewchar@{\familycount@\@ne
 \global\skewcharcount@\the\skewchar\textfont\@ne                           
 \ifnum\fam>\m@ne\ifnum\fam<16
  \global\familycount@\the\fam\relax
  \global\skewcharcount@\the\skewchar\textfont\the\fam\relax\fi\fi          
 \ifnum\skewcharcount@>\m@ne
  \ifnum\skewcharcount@<128
  \multiply\familycount@256
  \global\advance\skewcharcount@\familycount@
  \global\advance\skewcharcount@28672
  \mathchar\skewcharcount@\else
  \global\skewcharcount@\m@ne\fi\else
 \global\skewcharcount@\m@ne\fi}                                            
\newcount\pointcount@
\def\getpoints@#1.#2\getpoints@{\pointcount@#1 }
\newdimen\accentdimen@
\newcount\accentmu@
\def\dimentomu@{\multiply\accentdimen@ 100
 \expandafter\getpoints@\the\accentdimen@\getpoints@
 \multiply\pointcount@18
 \divide\pointcount@\@m
 \global\accentmu@\pointcount@}
\def\Makeacc@#1#2{\def#1{\RIfM@\DN@{\mathaccent@
 {"\accentclass@#2 }}\else\DN@{\nonmatherr@{#1}}\fi\next@}}
\def\unbracefonts@{\let\Cal@\Cal@@\let\roman@\roman@@\let\bold@\bold@@
 \let\slanted@\slanted@@}
\def\mathaccent@#1#2{\ifnum\fam=\m@ne\xdef\thefam@{1}\else
 \xdef\thefam@{\the\fam}\fi                                                 
 \accentdimen@\z@                                                           
 \setboxz@h{\unbracefonts@$\m@th\fam\thefam@\relax#2$}
 \ifdim\accentdimen@=\z@\DN@{\mathaccent#1{#2}}
  \setbox@ne\hbox{\unbracefonts@$\m@th\fam\thefam@\relax#2\theskewchar@$}
  \setbox\tw@\hbox{$\m@th\ifnum\skewcharcount@=\m@ne\else
   \mathchar\skewcharcount@\fi$}
  \global\accentdimen@\wd@ne\global\advance\accentdimen@-\wdz@
  \global\advance\accentdimen@-\wd\tw@                                     
  \global\multiply\accentdimen@\tw@
  \dimentomu@\global\advance\accentmu@\@ne                                 
 \else\DN@{{\mathaccent#1{#2\mkern\accentmu@ mu}%
    \mkern-\accentmu@ mu}{}}\fi                                             
 \next@}\Makeacc@\Hat{05E}
\Makeacc@\Check{014}
\Makeacc@\Tilde{07E}
\Makeacc@\Acute{013}
\Makeacc@\Grave{012}
\Makeacc@\Dot{05F}
\Makeacc@\Ddot{07F}
\Makeacc@\Breve{015}
\Makeacc@\Bar{016}
\def\Vec{\RIfM@\DN@{\mathaccent@{"017E }}\else
 \DN@{\nonmatherr@\Vec}\fi\next@}
\def\accentedsymbol#1#2{\csname newbox\expandafter\endcsname
  \csname\expandafter\eat@\string#1@box\endcsname
 \expandafter\setbox\csname\expandafter\eat@
  \string#1@box\endcsname\hbox{$\m@th#2$}\define
  #1{\copy\csname\expandafter\eat@\string#1@box\endcsname{}}}
\message{roots,}
\def\sqrt#1{\radical"270370 {#1}}
\let\underline@\underline
\let\overline@\overline
\def\underline#1{\underline@{#1}}
\def\overline#1{\overline@{#1}}
\Invalid@\leftroot
\Invalid@\uproot
\newcount\uproot@
\newcount\leftroot@
\def\root{\relaxnext@
  \DN@{\ifx\next\uproot\let\next@\nextii@\else
   \ifx\next\leftroot\let\next@\nextiii@\else
   \let\next@\plainroot@\fi\fi\next@}%
  \DNii@\uproot##1{\uproot@##1\relax\FN@\nextiv@}%
  \def\nextiv@{\ifx\next\space@\DN@. {\FN@\nextv@}\else
   \DN@.{\FN@\nextv@}\fi\next@.}%
  \def\nextv@{\ifx\next\leftroot\let\next@\nextvi@\else
   \let\next@\plainroot@\fi\next@}%
  \def\nextvi@\leftroot##1{\leftroot@##1\relax\plainroot@}%
   \def\nextiii@\leftroot##1{\leftroot@##1\relax\FN@\nextvii@}%
  \def\nextvii@{\ifx\next\space@
   \DN@. {\FN@\nextviii@}\else
   \DN@.{\FN@\nextviii@}\fi\next@.}%
  \def\nextviii@{\ifx\next\uproot\let\next@\nextix@\else
   \let\next@\plainroot@\fi\next@}%
  \def\nextix@\uproot##1{\uproot@##1\relax\plainroot@}%
  \bgroup\uproot@\z@\leftroot@\z@\FN@\next@}
\def\plainroot@#1\of#2{\setbox\rootbox\hbox{$\m@th\scriptscriptstyle{#1}$}%
 \mathchoice{\r@@t\displaystyle{#2}}{\r@@t\textstyle{#2}}
 {\r@@t\scriptstyle{#2}}{\r@@t\scriptscriptstyle{#2}}\egroup}
\def\r@@t#1#2{\setboxz@h{$\m@th#1\sqrt{#2}$}%
 \dimen@\ht\z@\advance\dimen@-\dp\z@
 \setbox@ne\hbox{$\m@th#1\mskip\uproot@ mu$}\advance\dimen@ 1.667\wd@ne
 \mkern-\leftroot@ mu\mkern5mu\raise.6\dimen@\copy\rootbox
 \mkern-10mu\mkern\leftroot@ mu\boxz@}
\def\boxed#1{\setboxz@h{$\m@th\displaystyle{#1}$}\dimen@.4\ex@
 \advance\dimen@3\ex@\advance\dimen@\dp\z@
 \hbox{\lower\dimen@\hbox{%
 \vbox{\hrule height.4\ex@
 \hbox{\vrule width.4\ex@\hskip3\ex@\vbox{\vskip3\ex@\boxz@\vskip3\ex@}%
 \hskip3\ex@\vrule width.4\ex@}\hrule height.4\ex@}%
 }}}
\message{commutative diagrams,}
\let\ampersand@\relax
\newdimen\minaw@
\minaw@11.11128\ex@
\newdimen\minCDaw@
\minCDaw@2.5pc
\def\minCDarrowwidth#1{\RIfMIfI@\onlydmatherr@\minCDarrowwidth
 \else\minCDaw@#1\relax\fi\else\onlydmatherr@\minCDarrowwidth\fi}
\newif\ifCD@
\def\CD{\bgroup\vspace@\relax\let\ampersand@&\iffalse}\fi
 \CD@true\vcenter\bgroup\Let@\tabskip\z@skip\baselineskip20\ex@
 \lineskip3\ex@\lineskiplimit3\ex@\halign\bgroup
 &\hfill$\m@th##$\hfill\crcr}
\def\endCD{\crcr\egroup\egroup\egroup}
\newdimen\bigaw@
\atdef@>#1>#2>{\ampersand@                                                  
 \setboxz@h{$\m@th\ssize\;{#1}\;\;$}
 \setbox@ne\hbox{$\m@th\ssize\;{#2}\;\;$}
 \setbox\tw@\hbox{$\m@th#2$}
 \ifCD@\global\bigaw@\minCDaw@\else\global\bigaw@\minaw@\fi                 
 \ifdim\wdz@>\bigaw@\global\bigaw@\wdz@\fi
 \ifdim\wd@ne>\bigaw@\global\bigaw@\wd@ne\fi                                
 \ifCD@\enskip\fi                                                           
 \ifdim\wd\tw@>\z@
  \mathrel{\mathop{\hbox to\bigaw@{\rightarrowfill@\displaystyle}}%
    \limits^{#1}_{#2}}
 \else\mathrel{\mathop{\hbox to\bigaw@{\rightarrowfill@\displaystyle}}%
    \limits^{#1}}\fi                                                        
 \ifCD@\enskip\fi                                                          
 \ampersand@}                                                              
\atdef@<#1<#2<{\ampersand@\setboxz@h{$\m@th\ssize\;\;{#1}\;$}%
 \setbox@ne\hbox{$\m@th\ssize\;\;{#2}\;$}\setbox\tw@\hbox{$\m@th#2$}%
 \ifCD@\global\bigaw@\minCDaw@\else\global\bigaw@\minaw@\fi
 \ifdim\wdz@>\bigaw@\global\bigaw@\wdz@\fi
 \ifdim\wd@ne>\bigaw@\global\bigaw@\wd@ne\fi
 \ifCD@\enskip\fi
 \ifdim\wd\tw@>\z@
  \mathrel{\mathop{\hbox to\bigaw@{\leftarrowfill@\displaystyle}}%
       \limits^{#1}_{#2}}\else
  \mathrel{\mathop{\hbox to\bigaw@{\leftarrowfill@\displaystyle}}%
       \limits^{#1}}\fi
 \ifCD@\enskip\fi\ampersand@}
\begingroup
 \catcode`\~=\active \lccode`\~=`\@
 \lowercase{%
  \global\atdef@)#1)#2){~>#1>#2>}
  \global\atdef@(#1(#2({~<#1<#2<}}
\endgroup
\atdef@ A#1A#2A{\llap{$\m@th\vcenter{\hbox
 {$\ssize#1$}}$}\Big\uparrow\rlap{$\m@th\vcenter{\hbox{$\ssize#2$}}$}&&}
\atdef@ V#1V#2V{\llap{$\m@th\vcenter{\hbox
 {$\ssize#1$}}$}\Big\downarrow\rlap{$\m@th\vcenter{\hbox{$\ssize#2$}}$}&&}
\atdef@={&\enskip\mathrel
 {\vbox{\hrule width\minCDaw@\vskip3\ex@\hrule width
 \minCDaw@}}\enskip&}
\atdef@|{\Big\Vert&&}
\atdef@\vert{\Big\Vert&&}
\def\pretend#1\haswidth#2{\setboxz@h{$\m@th\scriptstyle{#2}$}\hbox
 to\wdz@{\hfill$\m@th\scriptstyle{#1}$\hfill}}
\message{poor man's bold,}
\def\pmb{\RIfM@\expandafter\mathpalette\expandafter\pmb@\else
 \expandafter\pmb@@\fi}
\def\pmb@@#1{\leavevmode\setboxz@h{#1}%
   \dimen@-\wdz@
   \kern-.5\ex@\copy\z@
   \kern\dimen@\kern.25\ex@\raise.4\ex@\copy\z@
   \kern\dimen@\kern.25\ex@\box\z@
}
\def\binrel@@#1{\ifdim\wd2<\z@\mathbin{#1}\else\ifdim\wd\tw@>\z@
 \mathrel{#1}\else{#1}\fi\fi}
\newdimen\pmbraise@
\def\pmb@#1#2{\setbox\thr@@\hbox{$\m@th#1{#2}$}%
 \setbox4\hbox{$\m@th#1\mkern.5mu$}\pmbraise@\wd4\relax
 \binrel@{#2}%
 \dimen@-\wd\thr@@
   \binrel@@{%
   \mkern-.8mu\copy\thr@@
   \kern\dimen@\mkern.4mu\raise\pmbraise@\copy\thr@@
   \kern\dimen@\mkern.4mu\box\thr@@
}}
\def\documentstyle#1{\W@{}\input #1.sty\relax}
\message{syntax check,}
\font\dummyft@=dummy
\fontdimen1 \dummyft@=\z@
\fontdimen2 \dummyft@=\z@
\fontdimen3 \dummyft@=\z@
\fontdimen4 \dummyft@=\z@
\fontdimen5 \dummyft@=\z@
\fontdimen6 \dummyft@=\z@
\fontdimen7 \dummyft@=\z@
\fontdimen8 \dummyft@=\z@
\fontdimen9 \dummyft@=\z@
\fontdimen10 \dummyft@=\z@
\fontdimen11 \dummyft@=\z@
\fontdimen12 \dummyft@=\z@
\fontdimen13 \dummyft@=\z@
\fontdimen14 \dummyft@=\z@
\fontdimen15 \dummyft@=\z@
\fontdimen16 \dummyft@=\z@
\fontdimen17 \dummyft@=\z@
\fontdimen18 \dummyft@=\z@
\fontdimen19 \dummyft@=\z@
\fontdimen20 \dummyft@=\z@
\fontdimen21 \dummyft@=\z@
\fontdimen22 \dummyft@=\z@
\def\fontlist@{\\{\tenrm}\\{\sevenrm}\\{\fiverm}\\{\teni}\\{\seveni}%
 \\{\fivei}\\{\tensy}\\{\sevensy}\\{\fivesy}\\{\tenex}\\{\tenbf}\\{\sevenbf}%
 \\{\fivebf}\\{\tensl}\\{\tenit}}
\def\font@#1=#2 {\rightappend@#1\to\fontlist@\font#1=#2 }
\def\dodummy@{{\def\\##1{\global\let##1\dummyft@}\fontlist@}}
\def\nopages@{\output{\setbox\z@\box\@cclv \deadcycles\z@}%
 \alloc@5\toks\toksdef\@cclvi\output}
\let\galleys\nopages@
\newif\ifsyntax@
\newcount\countxviii@
\def\syntax{\syntax@true\dodummy@\countxviii@\count18
 \loop\ifnum\countxviii@>\m@ne\textfont\countxviii@=\dummyft@
 \scriptfont\countxviii@=\dummyft@\scriptscriptfont\countxviii@=\dummyft@
 \advance\countxviii@\m@ne\repeat                                           
 \dummyft@\tracinglostchars\z@\nopages@\frenchspacing\hbadness\@M}
\def\first@#1#2\end{#1}
\def\printoptions{\W@{Do you want S(yntax check),
  G(alleys) or P(ages)?}%
 \message{Type S, G or P, followed by <return>: }%
 \begingroup 
 \endlinechar\m@ne 
 \read\m@ne to\ans@
 \edef\ans@{\uppercase{\def\noexpand\ans@{%
   \expandafter\first@\ans@ P\end}}}%
 \expandafter\endgroup\ans@
 \if\ans@ P
 \else \if\ans@ S\syntax
 \else \if\ans@ G\galleys
 \else\message{? Unknown option: \ans@; using the `pages' option.}%
 \fi\fi\fi}
\def\alloc@#1#2#3#4#5{\global\advance\count1#1by\@ne
 \ch@ck#1#4#2\allocationnumber=\count1#1
 \global#3#5=\allocationnumber
 \ifalloc@\wlog{\string#5=\string#2\the\allocationnumber}\fi}
\def\document{\def\alloclist@{}\def\fontlist@{}}
\let\enddocument\bye

\let\proclaim\undefined
\let\footnote\undefined
\let\=\undefined
\let\>\undefined

\catcode`\@=\active
\message{... finished}

\expandafter\ifx\csname mathdefs.tex\endcsname\relax
  \expandafter\gdef\csname mathdefs.tex\endcsname{}
\else \message{Hey!  Apparently you were trying to
  \string\input{mathdefs.tex} twice.   This does not make sense.} 
\errmessage{Please edit your file (probably \jobname.tex) and remove
any duplicate ``\string\input'' lines}\endinput\fi




\catcode`\X=12\catcode`\@=11

\def\n@wcount{\alloc@0\count\countdef\insc@unt}
\def\n@wwrite{\alloc@7\write\chardef\sixt@@n}
\def\n@wread{\alloc@6\read\chardef\sixt@@n}
\def\r@s@t{\relax}\def\v@idline{\par}\def\@mputate#1/{#1}
\def\l@c@l#1X{\firstpart.#1}\def\gl@b@l#1X{#1}\def\t@d@l#1X{{}}

\def\crossrefs#1{\ifx\all#1\let\tr@ce=\all\else\def\tr@ce{#1,}\fi
   \n@wwrite\cit@tionsout\openout\cit@tionsout=\jobname.cit 
   \write\cit@tionsout{\tr@ce}\expandafter\setfl@gs\tr@ce,}
\def\setfl@gs#1,{\def\@{#1}\ifx\@\empty\let\next=\relax
   \else\let\next=\setfl@gs\expandafter\xdef
   \csname#1tr@cetrue\endcsname{}\fi\next}
\def\m@ketag#1#2{\expandafter\n@wcount\csname#2tagno\endcsname
     \csname#2tagno\endcsname=0\let\tail=\all\xdef\all{\tail#2,}
   \ifx#1\l@c@l\let\tail=\r@s@t\xdef\r@s@t{\csname#2tagno\endcsname=0\tail}\fi
   \expandafter\gdef\csname#2cite\endcsname##1{\expandafter
     \ifx\csname#2tag##1\endcsname\relax?\else\csname#2tag##1\endcsname\fi
     \expandafter\ifx\csname#2tr@cetrue\endcsname\relax\else
     \write\cit@tionsout{#2tag ##1 cited on page \folio.}\fi}
   \expandafter\gdef\csname#2page\endcsname##1{\expandafter
     \ifx\csname#2page##1\endcsname\relax?\else\csname#2page##1\endcsname\fi
     \expandafter\ifx\csname#2tr@cetrue\endcsname\relax\else
     \write\cit@tionsout{#2tag ##1 cited on page \folio.}\fi}
   \expandafter\gdef\csname#2tag\endcsname##1{\expandafter
      \ifx\csname#2check##1\endcsname\relax
      \expandafter\xdef\csname#2check##1\endcsname{}%
      \else\immediate\write16{Warning: #2tag ##1 used more than once.}\fi
      \multit@g{#1}{#2}##1/X%
      \write\t@gsout{#2tag ##1 assigned number \csname#2tag##1\endcsname\space
      on page \number\count0.}%
   \csname#2tag##1\endcsname}}

\def\multit@g#1#2#3/#4X{\def\t@mp{#4}\ifx\t@mp\empty%
      \global\advance\csname#2tagno\endcsname by 1 
      \expandafter\xdef\csname#2tag#3\endcsname
      {#1\number\csname#2tagno\endcsnameX}%
   \else\expandafter\ifx\csname#2last#3\endcsname\relax
      \expandafter\n@wcount\csname#2last#3\endcsname
      \global\advance\csname#2tagno\endcsname by 1 
      \expandafter\xdef\csname#2tag#3\endcsname
      {#1\number\csname#2tagno\endcsnameX}
      \write\t@gsout{#2tag #3 assigned number \csname#2tag#3\endcsname\space
      on page \number\count0.}\fi
   \global\advance\csname#2last#3\endcsname by 1
   \def\t@mp{\expandafter\xdef\csname#2tag#3/}%
   \expandafter\t@mp\@mputate#4\endcsname
   {\csname#2tag#3\endcsname\lastpart{\csname#2last#3\endcsname}}\fi}
\def\t@gs#1{\def\all{}\m@ketag#1e\m@ketag#1s\m@ketag\t@d@l p
\let\realscite\scite
\let\realstag\stag
   \m@ketag\gl@b@l r \n@wread\t@gsin
   \openin\t@gsin=\jobname.tgs \re@der \closein\t@gsin
   \n@wwrite\t@gsout\openout\t@gsout=\jobname.tgs }
\outer\def\localtags{\t@gs\l@c@l}
\outer\def\globaltags{\t@gs\gl@b@l}
\outer\def\newlocaltag#1{\m@ketag\l@c@l{#1}}
\outer\def\newglobaltag#1{\m@ketag\gl@b@l{#1}}

\newif\ifpr@ 
\def\m@kecs #1tag #2 assigned number #3 on page #4.%
   {\expandafter\gdef\csname#1tag#2\endcsname{#3}
   \expandafter\gdef\csname#1page#2\endcsname{#4}
   \ifpr@\expandafter\xdef\csname#1check#2\endcsname{}\fi}
\def\re@der{\ifeof\t@gsin\let\next=\relax\else
   \read\t@gsin to\t@gline\ifx\t@gline\v@idline\else
   \expandafter\m@kecs \t@gline\fi\let \next=\re@der\fi\next}
\def\pretags#1{\pr@true\pret@gs#1,,}
\def\pret@gs#1,{\def\@{#1}\ifx\@\empty\let\n@xtfile=\relax
   \else\let\n@xtfile=\pret@gs \openin\t@gsin=#1.tgs \message{#1} \re@der 
   \closein\t@gsin\fi \n@xtfile}

\newcount\sectno\sectno=0\newcount\subsectno\subsectno=0
\newif\ifultr@local \def\ultralocal{\ultr@localtrue}
\def\firstpart{\number\sectno}
\def\lastpart#1{\ifcase#1 \or a\or b\or c\or d\or e\or f\or g\or h\or 
   i\or k\or l\or m\or n\or o\or p\or q\or r\or s\or t\or u\or v\or w\or 
   x\or y\or z \fi}

\def\resetall{\global\advance\sectno by 1\subsectno=0
   \gdef\firstpart{\number\sectno}\r@s@t}
\def\resetsub{\global\advance\subsectno by 1
   \gdef\firstpart{\number\sectno.\number\subsectno}\r@s@t}
\def\newsection#1\par{\resetall\vskip0pt plus.3\vsize\penalty-250
   \vskip0pt plus-.3\vsize\bigskip\bigskip
   \message{#1}\leftline{\bf#1}\nobreak\bigskip}
\def\subsection#1\par{\ifultr@local\resetsub\fi
   \vskip0pt plus.2\vsize\penalty-250\vskip0pt plus-.2\vsize
   \bigskip\smallskip\message{#1}\leftline{\bf#1}\nobreak\medskip}


\newdimen\marginshift

\newdimen\margindelta
\newdimen\marginmax
\newdimen\marginmin

\def\margininit{       
\marginmax=3 true cm                  
				      
\margindelta=0.1 true cm              
\marginmin=0.1true cm                 
\marginshift=\marginmin
}    

\def\t@gsjj#1,{\def\@{#1}\ifx\@\empty\let\next=\relax\else\let\next=\t@gsjj
   \def\@@{p}\ifx\@\@@\else
   \expandafter\gdef\csname#1cite\endcsname##1{\citejj{##1}}
   \expandafter\gdef\csname#1page\endcsname##1{?}
   \expandafter\gdef\csname#1tag\endcsname##1{\tagjj{##1}}\fi\fi\next}
\newif\ifshowstuffinmargin
\showstuffinmarginfalse
\def\jjtags{\ifx\shlhetal\relax 
  \else
\ifx\shlhetal\undefinedcontrolseq
\else
\showstuffinmargintrue
\ifx\all\relax\else\expandafter\t@gsjj\all,\fi\fi \fi
}

\def\tagjj#1{\realstag{#1}\mginpar{\zeigen{#1}}}
\def\citejj#1{\rechnen{#1}\mginpar{\zeigen{#1}}}     

\def\rechnen#1{\expandafter\ifx\csname stag#1\endcsname\relax ??\else
                           \csname stag#1\endcsname\fi}

\newdimen\theight

\def\marginfont{\sevenrm}

\def\trymarginbox#1{\setbox0=\hbox{\marginfont\hskip\marginshift #1}%
		\global\marginshift\wd0 
		\global\advance\marginshift\margindelta}

\def \mginpar#1{%
\ifvmode\setbox0\hbox to \hsize{\hfill\rlap{\marginfont\quad#1}}%
\ht0 0cm
\dp0 0cm
\box0\vskip-\baselineskip
\else 
             \vadjust{\trymarginbox{#1}%
		\ifdim\marginshift>\marginmax \global\marginshift\marginmin
			\trymarginbox{#1}%
                \fi
             \theight=\ht0
             \advance\theight by \dp0    \advance\theight by \lineskip
             \kern -\theight \vbox to \theight{\rightline{\rlap{\box0}}%
\vss}}\fi}


\def\t@gsoff#1,{\def\@{#1}\ifx\@\empty\let\next=\relax\else\let\next=\t@gsoff
   \def\@@{p}\ifx\@\@@\else
   \expandafter\gdef\csname#1cite\endcsname##1{\zeigen{##1}}
   \expandafter\gdef\csname#1page\endcsname##1{?}
   \expandafter\gdef\csname#1tag\endcsname##1{\zeigen{##1}}\fi\fi\next}
\def\verbatimtags{\showstuffinmarginfalse
\ifx\all\relax\else\expandafter\t@gsoff\all,\fi}
\def\zeigen#1{\hbox{$\langle$}#1\hbox{$\rangle$}}

\def\margincite#1{\ifshowstuffinmargin\mginpar{\zeigen{#1}}\fi}

\def\margintag#1{\ifshowstuffinmargin\mginpar{\zeigen{#1}}\fi}

\def\(#1){\edef\dot@g{\ifmmode\ifinner(\hbox{\noexpand\etag{#1}})
   \else\noexpand\eqno(\hbox{\noexpand\etag{#1}})\fi
   \else(\noexpand\ecite{#1})\fi}\dot@g}

\newif\ifbr@ck
\def\eat#1{}
\def\[#1]{\br@cktrue[\br@cket#1'X]}
\def\br@cket#1'#2X{\def\temp{#2}\ifx\temp\empty\let\next\eat
   \else\let\next\br@cket\fi
   \ifbr@ck\br@ckfalse\br@ck@t#1,X\else\br@cktrue#1\fi\next#2X}
\def\br@ck@t#1,#2X{\def\temp{#2}\ifx\temp\empty\let\neext\eat
   \else\let\neext\br@ck@t\def\temp{,}\fi
   \def\teemp{#1}\ifx\teemp\empty\else\rcite{#1}\fi\temp\neext#2X}
\def\resetbr@cket{\gdef\[##1]{[\rtag{##1}]}}
\def\references{\resetbr@cket\newsection References\par}

\newtoks\symb@ls\newtoks\s@mb@ls\newtoks\p@gelist\n@wcount\ftn@mber
    \ftn@mber=1\newif\ifftn@mbers\ftn@mbersfalse\newif\ifbyp@ge\byp@gefalse
\def\defm@rk{\ifftn@mbers\n@mberm@rk\else\symb@lm@rk\fi}
\def\n@mberm@rk{\xdef\m@rk{{\the\ftn@mber}}%
    \global\advance\ftn@mber by 1 }
\def\rot@te#1{\let\temp=#1\global#1=\expandafter\r@t@te\the\temp,X}
\def\r@t@te#1,#2X{{#2#1}\xdef\m@rk{{#1}}}
\def\b@@st#1{{$^{#1}$}}\def\str@p#1{#1}
\def\symb@lm@rk{\ifbyp@ge\rot@te\p@gelist\ifnum\expandafter\str@p\m@rk=1 
    \s@mb@ls=\symb@ls\fi\write\f@nsout{\number\count0}\fi \rot@te\s@mb@ls}
\def\byp@ge{\byp@getrue\n@wwrite\f@nsin\openin\f@nsin=\jobname.fns 
    \n@wcount\currentp@ge\currentp@ge=0\p@gelist={0}
    \re@dfns\closein\f@nsin\rot@te\p@gelist
    \n@wread\f@nsout\openout\f@nsout=\jobname.fns }
\def\m@kelist#1X#2{{#1,#2}}
\def\re@dfns{\ifeof\f@nsin\let\next=\relax\else\read\f@nsin to \f@nline
    \ifx\f@nline\v@idline\else\let\t@mplist=\p@gelist
    \ifnum\currentp@ge=\f@nline
    \global\p@gelist=\expandafter\m@kelist\the\t@mplistX0
    \else\currentp@ge=\f@nline
    \global\p@gelist=\expandafter\m@kelist\the\t@mplistX1\fi\fi
    \let\next=\re@dfns\fi\next}
\def\symbols#1{\symb@ls={#1}\s@mb@ls=\symb@ls} 
\def\bigsymbol{\textstyle}
\symbols{\bigsymbol\ast,\dagger,\ddagger,\sharp,\flat,\natural,\star}
\def\ftnumbers{\ftn@mberstrue} \def\ftsymbols{\ftn@mbersfalse}
\def\paginal{\byp@ge} \def\resetftnumbers{\ftn@mber=1}
\def\ftnote#1{\defm@rk\expandafter\expandafter\expandafter\footnote
    \expandafter\b@@st\m@rk{#1}}

\long\def\jump#1\endjump{}
\def\ssum{\mathop{\lower .1em\hbox{$\textstyle\Sigma$}}\nolimits}

\def\qed{\nobreak\kern 1em \vrule height .5em width .5em depth 0em}
\def\newneq{\hbox{\rlap{\hbox to 1\wd9{\hss$=$\hss}}\raise .1em 
   \hbox to 1\wd9{\hss$\scriptscriptstyle/$\hss}}}
\def\subsetne{\setbox9 = \hbox{$\subset$}\mathrel{\hbox{\rlap
   {\lower .4em \newneq}\raise .13em \hbox{$\subset$}}}}
\def\supsetne{\setbox9 = \hbox{$\subset$}\mathrel{\hbox{\rlap
   {\lower .4em \newneq}\raise .13em \hbox{$\supset$}}}}

\def\vbar{\mathchoice{\vrule height6.3ptdepth-.5ptwidth.8pt\kern-.8pt}
   {\vrule height6.3ptdepth-.5ptwidth.8pt\kern-.8pt}
   {\vrule height4.1ptdepth-.35ptwidth.6pt\kern-.6pt}
   {\vrule height3.1ptdepth-.25ptwidth.5pt\kern-.5pt}}
\def\f@dge{\mathchoice{}{}{\mkern.5mu}{\mkern.8mu}}
\def\b@c#1#2{{\rm \mkern#2mu\vbar\mkern-#2mu#1}}
\def\b@b#1{{\rm I\mkern-3.5mu #1}}
\def\b@a#1#2{{\rm #1\mkern-#2mu\f@dge #1}}
\def\bb#1{{\count4=`#1 \advance\count4by-64 \ifcase\count4\or\b@a A{11.5}\or
   \b@b B\or\b@c C{5}\or\b@b D\or\b@b E\or\b@b F \or\b@c G{5}\or\b@b H\or
   \b@b I\or\b@c J{3}\or\b@b K\or\b@b L \or\b@b M\or\b@b N\or\b@c O{5} \or
   \b@b P\or\b@c Q{5}\or\b@b R\or\b@a S{8}\or\b@a T{10.5}\or\b@c U{5}\or
   \b@a V{12}\or\b@a W{16.5}\or\b@a X{11}\or\b@a Y{11.7}\or\b@a Z{7.5}\fi}}

\catcode`\X=11 \catcode`\@=12




\let\thischap\jobname

\def\partof#1{\csname returnthe#1part\endcsname}
\def\chapof#1{\csname returnthe#1chap\endcsname}

\def\setchapter#1,#2,#3.{%
  \expandafter\def\csname returnthe#1part\endcsname{#2}%
  \expandafter\def\csname returnthe#1chap\endcsname{#3}%
}

\setchapter 300a,A,I.
\setchapter 300b,A,II.
\setchapter 300c,A,III.
\setchapter 300d,A,IV.
\setchapter 300e,A,V.
\setchapter 300f,A,VI.
\setchapter 300g,A,VII.
\setchapter   88,B,I.
\setchapter  600,B,II.
\setchapter  705,B,III.

\def\cprefix#1{
\edef\theotherpart{\partof{#1}}\edef\theotherchap{\chapof{#1}}%
\ifx\theotherpart\thispart
   \ifx\theotherchap\thischap 
    \else 
     \theotherchap%
    \fi
   \else 
     \theotherpart.\theotherchap\fi}

\def\sectioncite[#1]#2{%
     \cprefix{#2}#1}

\edef\thispart{\partof{\thischap}}
\edef\thischap{\chapof{\thischap}}



\expandafter\ifx\csname citeadd.tex\endcsname\relax
\expandafter\gdef\csname citeadd.tex\endcsname{}
\else \message{Hey!  Apparently you were trying to
\string\input{citeadd.tex} twice.   This does not make sense.} 
\errmessage{Please edit your file (probably \jobname.tex) and remove
any duplicate ``\string\input'' lines}\endinput\fi

\def\sciteu{\sciteerror{undefined}}

\def\scitet{\sciteerror{ambiguous}}
\def\scitetphantom{\complainaboutcitation{ambiguous}}

\def\sciteerror#1#2{{\mathortextbf{\scite{#2}}}\complainaboutcitation{#1}{#2}}
\def\mathortextbf#1{\hbox{\bf #1}}
\def\complainaboutcitation#1#2{%
\vadjust{\line{\llap{---$\!\!>$ }\qquad scite$\{$#2$\}$ #1\hfil}}}

\sectno=-1   
\localtags
\jjtags
\NoBlackBoxes
\define\mr{\medskip\roster}
\define\sn{\smallskip\noindent}
\define\mn{\medskip\noindent}
\define\bn{\bigskip\noindent}
\define\ub{\underbar}
\define\wilog{\text{without loss of generality}}
\define\ermn{\endroster\medskip\noindent}
\define\dbca{\dsize\bigcap}
\define\dbcu{\dsize\bigcup}
\define\nl{\newline}
\magnification=\magstep 1
\documentstyle{amsppt}

{    
\catcode`@11

\ifx\alicetwothousandloaded@\relax
  \endinput\else\global\let\alicetwothousandloaded@\relax\fi

\gdef\subjclass{\let\savedef@\subjclass
 \def\subjclass##1\endsubjclass{\let\subjclass\savedef@
   \toks@{\def\usualspace{{\rm\enspace}}\eightpoint}%
   \toks@@{##1\unskip.}%
   \edef\thesubjclass@{\the\toks@
     \frills@{{\noexpand\rm2000 {\noexpand\it Mathematics Subject
       Classification}.\noexpand\enspace}}%
     \the\toks@@}}%
  \nofrillscheck\subjclass}
} 


\expandafter\ifx\csname alice2jlem.tex\endcsname\relax
  \expandafter\xdef\csname alice2jlem.tex\endcsname{\the\catcode`@}
\else \message{Hey!  Apparently you were trying to
\string\input{alice2jlem.tex}  twice.   This does not make sense.}
\errmessage{Please edit your file (probably \jobname.tex) and remove
any duplicate ``\string\input'' lines}\endinput\fi

\expandafter\ifx\csname bib4plain.tex\endcsname\relax
  \expandafter\gdef\csname bib4plain.tex\endcsname{}
\else \message{Hey!  Apparently you were trying to \string\input
  bib4plain.tex twice.   This does not make sense.}
\errmessage{Please edit your file (probably \jobname.tex) and remove
any duplicate ``\string\input'' lines}\endinput\fi

\def\renewcommand{\newcommand}	       
\edef\cite{\the\catcode`@}%
\catcode`@ = 11
\let\@oldatcatcode = \cite
\chardef\@letter = 11
\chardef\@other = 12
%
%
%
%
\def\@innerdef#1#2{\edef#1{\expandafter\noexpand\csname #2\endcsname}}%
%
%
\@innerdef\@innernewcount{newcount}%
\@innerdef\@innernewdimen{newdimen}%
\@innerdef\@innernewif{newif}%
\@innerdef\@innernewwrite{newwrite}%
%
%
%
\def\@gobble#1{}%
%
%
%
\ifx\inputlineno\@undefined
   \let\@linenumber = \empty 
\else
   \def\@linenumber{\the\inputlineno:\space}%
\fi
%
%
%
\def\@futurenonspacelet#1{\def\cs{#1}%
   \afterassignment\@stepone\let\@nexttoken=
}%
\begingroup 
\def\\{\global\let\@stoken= }%
\\ 
\endgroup
\def\@stepone{\expandafter\futurelet\cs\@steptwo}%
\def\@steptwo{\expandafter\ifx\cs\@stoken\let\@@next=\@stepthree
   \else\let\@@next=\@nexttoken\fi \@@next}%
\def\@stepthree{\afterassignment\@stepone\let\@@next= }%
%
%
%
\def\@getoptionalarg#1{%
   \let\@optionaltemp = #1%
   \let\@optionalnext = \relax
   \@futurenonspacelet\@optionalnext\@bracketcheck
}%
%
%
\def\@bracketcheck{%
   \ifx [\@optionalnext
      \expandafter\@@getoptionalarg
   \else
      \let\@optionalarg = \empty
      \expandafter\@optionaltemp
   \fi
}%
\def\@@getoptionalarg[#1]{%
   \def\@optionalarg{#1}%
   \@optionaltemp
}%
%
%
%
\def\@nnil{\@nil}%
\def\@fornoop#1\@@#2#3{}%
\def\@for#1:=#2\do#3{%
   \edef\@fortmp{#2}%
   \ifx\@fortmp\empty \else
      \expandafter\@forloop#2,\@nil,\@nil\@@#1{#3}%
   \fi
}%
\def\@forloop#1,#2,#3\@@#4#5{\def#4{#1}\ifx #4\@nnil \else
       #5\def#4{#2}\ifx #4\@nnil \else#5\@iforloop #3\@@#4{#5}\fi\fi
}%
\def\@iforloop#1,#2\@@#3#4{\def#3{#1}\ifx #3\@nnil
       \let\@nextwhile=\@fornoop \else
      #4\relax\let\@nextwhile=\@iforloop\fi\@nextwhile#2\@@#3{#4}%
}%
%
%
%
\@innernewif\if@fileexists
\def\@testfileexistence{\@getoptionalarg\@finishtestfileexistence}%
\def\@finishtestfileexistence#1{%
   \begingroup
      \def\extension{#1}%
      \immediate\openin0 =
         \ifx\@optionalarg\empty\jobname\else\@optionalarg\fi
         \ifx\extension\empty \else .#1\fi
         \space
      \ifeof 0
         \global\@fileexistsfalse
      \else
         \global\@fileexiststrue
      \fi
      \immediate\closein0
   \endgroup
}%
%
%
%
%
\def\bibliographystyle#1{%
   \@readauxfile
   \@writeaux{\string\bibstyle{#1}}%
}%
\let\bibstyle = \@gobble
%
%
\let\bblfilebasename = \jobname
\def\bibliography#1{%
   \@readauxfile
   \@writeaux{\string\bibdata{#1}}%
   \@testfileexistence[\bblfilebasename]{bbl}%
   \if@fileexists
      \nobreak
      \@readbblfile
   \fi
}%
\let\bibdata = \@gobble
%
%
\def\nocite#1{%
   \@readauxfile
   \@writeaux{\string\citation{#1}}%
}%
\@innernewif\if@notfirstcitation
%
%
\def\cite{\@getoptionalarg\@cite}%
%
%
\def\@cite#1{%
   \let\@citenotetext = \@optionalarg
   \printcitestart
   \nocite{#1}%
   \@notfirstcitationfalse
   \@for \@citation :=#1\do
   {%
      \expandafter\@onecitation\@citation\@@
   }%
   \ifx\empty\@citenotetext\else
      \printcitenote{\@citenotetext}%
   \fi
   \printcitefinish
}%
\newif\ifweareinprivate
\weareinprivatetrue
\ifx\shlhetal\undefinedcontrolseq\weareinprivatefalse\fi
\ifx\shlhetal\relax\weareinprivatefalse\fi
\def\@onecitation#1\@@{%
   \if@notfirstcitation
      \printbetweencitations
   \fi
   \expandafter \ifx \csname\@citelabel{#1}\endcsname \relax
      \if@citewarning
         \message{\@linenumber Undefined citation `#1'.}%
      \fi
     \ifweareinprivate
      \expandafter\gdef\csname\@citelabel{#1}\endcsname{%
\strut 
\vadjust{\vskip-\dp\strutbox
\vbox to 0pt{\vss\parindent0cm \leftskip=\hsize 
\advance\leftskip3mm
\advance\hsize 4cm\strut\openup-4pt 
\rightskip 0cm plus 1cm minus 0.5cm ?  #1 ?\strut}}
         {\tt
            \escapechar = -1
            \nobreak\hskip0pt\pfeilsw
            \expandafter\string\csname#1\endcsname
             \pfeilso
            \nobreak\hskip0pt
         }%
      }%
     \else  
      \expandafter\gdef\csname\@citelabel{#1}\endcsname{%
            {\tt\expandafter\string\csname#1\endcsname}
      }%
     \fi  
   \fi
   \csname\@citelabel{#1}\endcsname
   \@notfirstcitationtrue
}%
%
%
\def\@citelabel#1{b@#1}%
%
%
\def\@citedef#1#2{\expandafter\gdef\csname\@citelabel{#1}\endcsname{#2}}%
%
%
%
\def\@readbblfile{%
   \ifx\@itemnum\@undefined
      \@innernewcount\@itemnum
   \fi
   \begingroup
      \def\begin##1##2{%
         \setbox0 = \hbox{\biblabelcontents{##2}}%
         \biblabelwidth = \wd0
      }%
      \def\end##1{}
      %
      %
      \@itemnum = 0
      \def\bibitem{\@getoptionalarg\@bibitem}%
      \def\@bibitem{%
         \ifx\@optionalarg\empty
            \expandafter\@numberedbibitem
         \else
            \expandafter\@alphabibitem
         \fi
      }%
      \def\@alphabibitem##1{%
         \expandafter \xdef\csname\@citelabel{##1}\endcsname {\@optionalarg}%
         \ifx\biblabelprecontents\@undefined
            \let\biblabelprecontents = \relax
         \fi
         \ifx\biblabelpostcontents\@undefined
            \let\biblabelpostcontents = \hss
         \fi
         \@finishbibitem{##1}%
      }%
      \def\@numberedbibitem##1{%
         \advance\@itemnum by 1
         \expandafter \xdef\csname\@citelabel{##1}\endcsname{\number\@itemnum}%
         \ifx\biblabelprecontents\@undefined
            \let\biblabelprecontents = \hss
         \fi
         \ifx\biblabelpostcontents\@undefined
            \let\biblabelpostcontents = \relax
         \fi
         \@finishbibitem{##1}%
      }%
      \def\@finishbibitem##1{%
         \biblabelprint{\csname\@citelabel{##1}\endcsname}%
         \@writeaux{\string\@citedef{##1}{\csname\@citelabel{##1}\endcsname}}%
         \ignorespaces
      }%
      %
      %
      \let\em = \bblem
      \let\newblock = \bblnewblock
      \let\sc = \bblsc
      \frenchspacing
      \clubpenalty = 4000 \widowpenalty = 4000
      \tolerance = 10000 \hfuzz = .5pt
      \everypar = {\hangindent = \biblabelwidth
                      \advance\hangindent by \biblabelextraspace}%
      \bblrm
      \parskip = 1.5ex plus .5ex minus .5ex
      \biblabelextraspace = .5em
      \bblhook
      \input \bblfilebasename.bbl
   \endgroup
}%
%
%
\@innernewdimen\biblabelwidth
\@innernewdimen\biblabelextraspace
%
%
%
\def\biblabelprint#1{%
   \noindent
   \hbox to \biblabelwidth{%
      \biblabelprecontents
      \biblabelcontents{#1}%
      \biblabelpostcontents
   }%
   \kern\biblabelextraspace
}%
%
%
%
\def\biblabelcontents#1{{\bblrm [#1]}}%
%
%
\def\bblrm{\rm}%
%
%
\def\bblem{\it}%
%
%
\def\bblsc{\ifx\@scfont\@undefined
              \font\@scfont = cmcsc10
           \fi
           \@scfont
}%
%
%
\def\bblnewblock{\hskip .11em plus .33em minus .07em }%
%
%
\let\bblhook = \empty
%
%
%
\def\printcitestart{[}
\def\printcitefinish{]}
\def\printbetweencitations{, }
\def\printcitenote#1{, #1}
%
%
%
\let\citation = \@gobble
%
%
%
\@innernewcount\@numparams
%
%
\def\newcommand#1{%
   \def\@commandname{#1}%
   \@getoptionalarg\@continuenewcommand
}%
%
%
\def\@continuenewcommand{%
   \@numparams = \ifx\@optionalarg\empty 0\else\@optionalarg \fi \relax
   \@newcommand
}%
%
%
\def\@newcommand#1{%
   \def\@startdef{\expandafter\edef\@commandname}%
   \ifnum\@numparams=0
      \let\@paramdef = \empty
   \else
      \ifnum\@numparams>9
         \errmessage{\the\@numparams\space is too many parameters}%
      \else
         \ifnum\@numparams<0
            \errmessage{\the\@numparams\space is too few parameters}%
         \else
            \edef\@paramdef{%
               \ifcase\@numparams
                  \empty  No arguments.
               \or ####1%
               \or ####1####2%
               \or ####1####2####3%
               \or ####1####2####3####4%
               \or ####1####2####3####4####5%
               \or ####1####2####3####4####5####6%
               \or ####1####2####3####4####5####6####7%
               \or ####1####2####3####4####5####6####7####8%
               \or ####1####2####3####4####5####6####7####8####9%
               \fi
            }%
         \fi
      \fi
   \fi
   \expandafter\@startdef\@paramdef{#1}%
}%
%
%
%
%
\def\@readauxfile{%
   \if@auxfiledone \else 
      \global\@auxfiledonetrue
      \@testfileexistence{aux}%
      \if@fileexists
         \begingroup
            \endlinechar = -1
            \catcode`@ = 11
            \input \jobname.aux
         \endgroup
      \else
         \message{\@undefinedmessage}%
         \global\@citewarningfalse
      \fi
      \immediate\openout\@auxfile = \jobname.aux
   \fi
}%
%
%
\newif\if@auxfiledone
\ifx\noauxfile\@undefined \else \@auxfiledonetrue\fi
%
%
%
%
\@innernewwrite\@auxfile
\def\@writeaux#1{\ifx\noauxfile\@undefined \write\@auxfile{#1}\fi}%
%
%
%
\ifx\@undefinedmessage\@undefined
   \def\@undefinedmessage{No .aux file; I won't give you warnings about
                          undefined citations.}%
\fi
%
%
\@innernewif\if@citewarning
\ifx\noauxfile\@undefined \@citewarningtrue\fi
%
%
%
\catcode`@ = \@oldatcatcode

\def\pfeilso{\leavevmode
            \vrule width 1pt height9pt depth 0pt\relax
           \vrule width 1pt height8.7pt depth 0pt\relax
           \vrule width 1pt height8.3pt depth 0pt\relax
           \vrule width 1pt height8.0pt depth 0pt\relax
           \vrule width 1pt height7.7pt depth 0pt\relax
            \vrule width 1pt height7.3pt depth 0pt\relax
            \vrule width 1pt height7.0pt depth 0pt\relax
            \vrule width 1pt height6.7pt depth 0pt\relax
            \vrule width 1pt height6.3pt depth 0pt\relax
            \vrule width 1pt height6.0pt depth 0pt\relax
            \vrule width 1pt height5.7pt depth 0pt\relax
            \vrule width 1pt height5.3pt depth 0pt\relax
            \vrule width 1pt height5.0pt depth 0pt\relax
            \vrule width 1pt height4.7pt depth 0pt\relax
            \vrule width 1pt height4.3pt depth 0pt\relax
            \vrule width 1pt height4.0pt depth 0pt\relax
            \vrule width 1pt height3.7pt depth 0pt\relax
            \vrule width 1pt height3.3pt depth 0pt\relax
            \vrule width 1pt height3.0pt depth 0pt\relax
            \vrule width 1pt height2.7pt depth 0pt\relax
            \vrule width 1pt height2.3pt depth 0pt\relax
            \vrule width 1pt height2.0pt depth 0pt\relax
            \vrule width 1pt height1.7pt depth 0pt\relax
            \vrule width 1pt height1.3pt depth 0pt\relax
            \vrule width 1pt height1.0pt depth 0pt\relax
            \vrule width 1pt height0.7pt depth 0pt\relax
            \vrule width 1pt height0.3pt depth 0pt\relax}

\def\pfeilsw{ \leavevmode 
            \vrule width 1pt height0.3pt depth 0pt\relax
            \vrule width 1pt height0.7pt depth 0pt\relax
            \vrule width 1pt height1.0pt depth 0pt\relax
            \vrule width 1pt height1.3pt depth 0pt\relax
            \vrule width 1pt height1.7pt depth 0pt\relax
            \vrule width 1pt height2.0pt depth 0pt\relax
            \vrule width 1pt height2.3pt depth 0pt\relax
            \vrule width 1pt height2.7pt depth 0pt\relax
            \vrule width 1pt height3.0pt depth 0pt\relax
            \vrule width 1pt height3.3pt depth 0pt\relax
            \vrule width 1pt height3.7pt depth 0pt\relax
            \vrule width 1pt height4.0pt depth 0pt\relax
            \vrule width 1pt height4.3pt depth 0pt\relax
            \vrule width 1pt height4.7pt depth 0pt\relax
            \vrule width 1pt height5.0pt depth 0pt\relax
            \vrule width 1pt height5.3pt depth 0pt\relax
            \vrule width 1pt height5.7pt depth 0pt\relax
            \vrule width 1pt height6.0pt depth 0pt\relax
            \vrule width 1pt height6.3pt depth 0pt\relax
            \vrule width 1pt height6.7pt depth 0pt\relax
            \vrule width 1pt height7.0pt depth 0pt\relax
            \vrule width 1pt height7.3pt depth 0pt\relax
            \vrule width 1pt height7.7pt depth 0pt\relax
            \vrule width 1pt height8.0pt depth 0pt\relax
            \vrule width 1pt height8.3pt depth 0pt\relax
            \vrule width 1pt height8.7pt depth 0pt\relax
            \vrule width 1pt height9pt depth 0pt\relax
      }


\def\widestnumber#1#2{}

\def\citewarning#1{\ifx\shlhetal\relax 
    \else
    \par{#1}\par
    \fi
}

\def\rm{\fam0 \tenrm}

\def\fakesubhead#1\endsubhead{\bigskip\noindent{\bf#1}\par}



%
%
%

%

\font\textrsfs=rsfs10
\font\scriptrsfs=rsfs7
\font\scriptscriptrsfs=rsfs5

\newfam\rsfsfam
\textfont\rsfsfam=\textrsfs
\scriptfont\rsfsfam=\scriptrsfs
\scriptscriptfont\rsfsfam=\scriptscriptrsfs

\edef\oldcatcodeofat{\the\catcode`\@}
\catcode`\@11

\def\Cal@@#1{\noaccents@ \fam \rsfsfam #1}

\catcode`\@\oldcatcodeofat


\expandafter\ifx \csname margininit\endcsname \relax\else\margininit\fi

\long\def\red#1\endred{}
\long\def\green#1\endgreen{}
\long\def\blue#1\endblue{}

\def\endred{ \unmatched endred! }
\def\endgreen{ \unmatched endgreen! }
\def\endblue{ \unmatched endblue! }

\ifx\latexcolors\undefinedcs\def\latexcolors{}\fi

\def\emptycs{}
\def\evaluatelatexcolors{%
        \ifx\latexcolors\emptycs\else
        \expandafter\xxevaluate\latexcolors\xxfertig\evaluatelatexcolors\fi}
\def\xxevaluate#1,#2\xxfertig{\setupthiscolor{#1}%
        \def\latexcolors{#2}}

\font\smallfont=cmsl7
\def\rutgerscolor{\ifmmode\else\endgraf\fi\smallfont
\advance\leftskip0.5cm\relax}
\def\setupthiscolor#1{\edef\tmptmpcs{\noexpand\bgroup\noexpand\rutgerscolor
\noexpand\def\noexpand\currentcolor{#1}%
\noexpand}%
\expandafter\let\csname#1\endcsname\tmptmpcs
\def\tmptmpcs{\checkColorUnmatched{#1}\popthecolor}
\expandafter\let\csname end#1\endcsname\tmptmpcs}

\def\checkColorUnmatched#1{\def\expectcolor{#1}%
    \ifx\expectcolor\currentcolor   
    \else \edef\failhere{\noexpand\tryingToClose '\currentcolor' with end\expectcolor}\failhere\fi}

\def\currentcolor{???}

\def\popthecolor{\ifmmode\else\endgraf\fi\egroup}

\expandafter\def\csname#1\endcsname{}

\evaluatelatexcolors

 \let\outerhead\head
 \def\head{\innerhead}
 \let\innerhead\outerhead

 \let\outersubhead\subhead
 \def\subhead{\innersubhead}
 \let\innersubhead\outersubhead

 \let\outersubsubhead\subsubhead
 \def\subsubhead{\innersubsubhead}
 \let\innersubsubhead\outersubsubhead

 \def\proclaim{\innerproclaim}
 \let\innerproclaim\outerproclaim

 %
 %
 %
 %

\def\demo#1{\medskip\noindent{\it #1.\/}}
\def\enddemo{\smallskip}

\def\remark#1{\medskip\noindent{\it #1.\/}}
\def\endremark{\smallskip}

\pageheight{8.5truein}
\topmatter
\title{{\it A More General Iterable Condition \\
Ensuring $\aleph_1$ is not Collapsed}} \endtitle
\rightheadtext{General Iterable Condition}
\author {Saharon Shelah 
\thanks {\null\newline
I would like to thank Alice Leonhardt for the beautiful typing.\newline
Publication 311.} \endthanks} \endauthor
\affil{Institute of Mathematics\\
 The Hebrew University\\
 91 904 Jerusalem, Israel
 \medskip
 Rutgers University\\
 Mathematics Department\\
 New Brunswick, NJ  USA} \endaffil

\mn
\abstract   In a self-contained way, we deal with revised countable
support iterated forcing for the reals; improve theorems on
preservation of a property UP, weaker than semi proper, and hopefully
improve the presentation.  We continue \cite[Ch.X,XI]{Sh:b} (or see
\cite[Ch.X,XI]{Sh:f}), and 
Gitik Shelah \cite{GiSh:191} and \cite[Ch.XIII,XIV]{Sh:f} and 
particularly Ch.XV; concerning ``no new reals" see lately
Larson Shelah \cite{LrSh:746}.  In particular, we fulfill some promises from
\cite{Sh:f} and give a more streamlined version.  \endabstract
\endtopmatter
\document  
 
\newpage

\head {Annotated Content} \endhead  \resetall 
\bn
\S1 \ub{Preliminaries}, p.5
\mr
\item "{{}}"  [We agree that forcing notion $\Bbb P$ has actually also pure
$(\le_{pr})$ quasi-order and very pure $(<_{vpr})$ quasi-order.  For a
$\lessdot$-increasing sequence $\bar{\Bbb Q}$ of forcing notions we 
define what is a $\bar{\Bbb Q}$-named and a $\bar{\Bbb Q}$-named 
$[j,\alpha)$-ordinal.  Then we define $\kappa-Sp_e(W)$-iterations 
(revised support of size $< \kappa$, including the
case $\kappa$ inaccessible) with finite apure support, countable pure support
(the revised version) and Easton or W-Easton very pure support, similar to
\cite[XIV]{Sh:f} and prove its basic properties (this is done by simultaneous
induction).]
\endroster
\bn
\S2 \ub{Trees of Models}, p.31
\mr
\item "{{}}"  [We quote the basic definitions and theorems concerning
trees with $\omega$ levels tagged by ideal and partition theorems.]
\endroster
\bn
\S3 \ub{Ideals and Partial Orders}, p.36
\mr
\item "{{}}" [We can replace the families $\Bbb I$ of ideals by corresponding
partial orders or quasi orders (we ``ignore" the distinction).  This is
essentially equivalent (for ``some $\lambda$-complete $\Bbb I$" with ``for 
some $\lambda$-complete ${\Cal L}$")
but the ${\Cal L}$'s have better ``pullback" from forcing extensions, so we
can replace ${\underset\tilde {}\to{\Cal L}}$ in a forcing extension
of $\bold V$ by ${\Cal L}'$ in $\bold V$ 
preserving ${\underset\tilde {}\to{\Cal L}} \le_{RK} 
{\Cal L}'$ and preserving the amount of completeness we have, so a similar 
situation holds for a set of ideals; in the cases we have in mind here 
increasing those sets $\Bbb I$ or ${\Cal L}$ do not matter.]
\endroster
\bn
\S4 \ub{UP Reintroduced}, p.42
\mr
\item "{{}}" [We define when $\bar N$ is an $\Bbb I$-tagged tree of 
models, when
it is $\Bbb I$-suitable, or $(\Bbb I,\bold W)$-suitable, and when it is
strictly or $\lambda$-strictly, etc., where $\Bbb I$ is a family of
ideals.   Similarly we define $N$ is 
$\lambda$-strictly $(\Bbb I,\bold S,\bold W)$-suitable; i.e. can serve as 
$N_{\langle \rangle}$ and prove some basic properties.  Such models will
fulfill here the role that ``any countable $N \prec ({\Cal H}(\chi),
\in)$" fulfills in theorems on semi-proper iteration.  Lastly, we define 
when a forcing notion $\Bbb P$ satisfies $UP^\ell(\Bbb I,\bold W)$ for 
$\ell = 0,1,2$, and here
$UP^1,UP^2$ replace $\bold W$-proper, $\bold W$-semi-proper, where $\bold W$
is a stationary subset of $\omega_1$.  All those properties imply 
that forcing by $\Bbb P$ does not collapse $\aleph_1$, preserve 
the stationarity of $\bold W$
and even of any stationary subset of it.  They are all relatives of 
``semi-properness" for strictly $(\Bbb I,\bold W)$-suitable models, they
speak on $\Bbb I$-tagged trees of countable models.]
\endroster
\bn
\S5 \ub{An Iteration Theorem for $UP^1$}, p.53 
\mr
\item "{{}}" [We prove that satisfying $UP^1(\bold W)$, i.e. satisfying it 
for some family $\Bbb I$ of ideals, complete enough, is preserved by 
$\aleph_1-Sp_6(W)$-iterations, $\bar{\Bbb Q} = \langle \Bbb P_i,
{\underset\tilde {}\to {\Bbb Q}_i}:i < \alpha \rangle$; that is if each 
${\underset\tilde {}\to {\Bbb Q}_i}$ is like that then $P_\alpha =
\aleph_1-Sp_6(W)$-Lim$(\bar{\Bbb Q})$ is like that, provided some 
mild condition holds (say $\Bbb Q_i$ is \nl
$UP^1(\Bbb I_i,\bold W),{\underset\tilde {}\to {\Bbb I}_i}$
is ${\underset\tilde {}\to \kappa_i}$-complete, $\Bbb P_i$ satisfies the
$\kappa_i$-c.c.; we can even make ${\underset\tilde {}\to {\Bbb I}_i},
{\underset\tilde {}\to \kappa_i}$ to be just $\Bbb P_{i+1}$-names, 
see there).  
The proof is more similar to the proofs of preservation of properness and 
semi-properness \ub{than} with the proofs in \cite[XI]{Sh:b}, 
(=\cite[XI]{Sh:f}), \cite{GiSh:191}, \cite[XV]{Sh:f}, and hopefully more 
transparent.  The proof will be non-trivially shorter if we use just 
the particular case of the revised countable support (i.e.,
$\le_{\text{vpr}}$ is equality and $\le_{\text{pr}}$ is $\le$).  
We give a sufficient
condition for $\alpha$ 
not being collapsed by $\Bbb P_\alpha$ e.g. $\alpha$ is
strongly inaccessible, $\beta < \alpha \Rightarrow |\Bbb P_\beta| 
< \alpha$ and:
$W$ stationary in $\alpha$ (so $\alpha$ is Mahlo) or
$\le_{\text{vpr}}$ is equality
and the iteration is suitable enough.  Lastly, if e.g. $\alpha$ is the first
strongly inaccessible, $i < \alpha \Rightarrow |\Bbb P_i| < \alpha$ we give a
sufficient condition for $\alpha$ not being collapsed.]
\endroster
\bn
\S6 \ub{Preservation of $UP^0$}, p.72
\mr
\item "{{}}" [Here we make the condition more similar to semi-proper
iteration, that is the demand is that for suitable models $N$ (one on which
``the right trees grow") above each $p \in \Bbb Q \cap N$ there is an 
$(N,\Bbb Q)$-semi-generic $q$.  There is some price though. \nl
[?] However, if $\Bbb Q$ satisfies UP$^0$ and the $\kappa$-c.c., then 
$\Bbb Q * {\underset\tilde {}\to {\text{Levy}}}(\aleph_1,< \kappa)$ 
is appropriate in our iteration.]
\endroster  
\bn
\S7 \ub{No New Reals} - replacements for completeness, p.91
\mr
\item "{{}}"  [Here we deal with the parallel of $``\Bbb Q$ 
add no new real because it is \nl
$\bold W$-complete for some stationary $\bold W \subseteq \omega_1"$.]
\endroster
\bn
\S8 \ub{Examples}, p.98
\mr
\item "{{}}"  [We show that various forcing notions fall under our
context.  In particular (variants of) Namba forcing, shooting a club through
a stationary $S \subseteq \{\delta < \lambda:\text{cf}(\delta) = \aleph_0\}$
where $\lambda = \text{ cf}(\lambda) > \aleph_0$, and prove that the 
older condition from \cite{Sh:f} implies the present one.]
\endroster
\bn
\S9 \ub{Reflection in $[\omega_2]^{\aleph_0}$}, p.104
\mr
\item "{{}}"  [We answer a question of Jech, on the consistency of
$2^{\aleph_0} = \aleph_2 + {\Cal D}_{\omega_1}$ 
is $\aleph_2$-saturated + every
stationary subset of $[\omega_2]^{\aleph_0}$ reflects and 
there is a special projectively stationary subset of 
$[\omega_2]^{\aleph_0}$.]
\endroster
\bn
\S10  \ub{Mixing finitary norms and ideals}, p.110
\mr
\item "{{}}"  [We consider a common generalization of creature forcing
(see \cite{RoSh:470}) and relatives of Namba forcing.]
\endroster
\newpage

\head {\S1 Preliminaries} \endhead  \resetall \sectno=1
\bigskip

\definition{\stag{1.1} Definition/Notation}  1) A forcing notion 
here, $\Bbb P$, is a nonempty set (abusing notation, it too is denoted by
$\Bbb P$) and three partial orders $\le_0,\le_1,\le_2$ 
(more exactly quasi-orders and $\le^*_\ell = \le^{\Bbb P}_\ell$)
 and a minimal element $\emptyset_{\Bbb P} \in \Bbb P$ (so
$\emptyset_{\Bbb P} \le_\ell p$ for $p \in \Bbb P$) and for 
$\ell = 0,1$ we have $[p \le_\ell q \rightarrow p \le_{\ell +1} q]$.  
We call $p \in \Bbb P$ very pure
if $\emptyset_{\Bbb P} \le_0 p$ 
and we call $q$ a very pure extension of $p$ if $p \le_0 q$.
We call $p \in \Bbb P$ pure if $\emptyset_{\Bbb P} \le_1 p$ and 
we call $q$ a pure extension of $p$ if $p \le_1 q$.  Let $\le$ 
be $\le_2$, let $\le_{\text{pr}}$ be $\le_1$ and $\le_{\text{vp}}$ 
be $\le_0$.

We call $\Bbb P$ $\kappa$-vp-complete if: for 
any $<_{\text{vpr}}$-increasing sequence $\langle p_\alpha:
\alpha < \delta \rangle,\delta < \kappa$ with $p_0 = 
\emptyset_{\Bbb P}$ there is a $\le_{text{vpr}}$-upper bound 
$p$.  We define vp-$\kappa$-complete similarly waiving $p_0 =
\emptyset_{\Bbb P}$.
We define $\kappa-\le_\ell$-complete and $\le_\ell-\kappa$-complete
similarly.

The forcing relation, of course, refers to the 
partial order $\le$.  We denote
forcing notions by $\Bbb P,\Bbb Q,\Bbb R$.  
Let $\Bbb P^1 \subseteq \Bbb P^2$ mean
$p \in \Bbb P^1 \Rightarrow p \in \Bbb P^2,\le^{{\Bbb P}^1}_\ell 
= \le^{{\Bbb P}^2}_\ell \restriction
\Bbb P^1$ for $\ell = 0,1,2$ and let $\Bbb P^1 \subseteq_{\text{ic}}
\Bbb P^2$ means $\Bbb P^1 \subseteq \Bbb P^2$ and for $\ell \le 2$, if 
$p,q \in \Bbb P^1$ are $\le_\ell$-incompatible in $\Bbb P^1$, then they
are $\le_\ell$-incompatible in $\Bbb P^2$.  Let $\Bbb P^1 \lessdot
\Bbb P^2$ means $\Bbb P^1 \subseteq_{\text{ic}} \Bbb P^2 \and 
(\Bbb P^1,\le) \lessdot (\Bbb P^2,\le)$. \newline
2) $\bar{\Bbb P}$ denotes a $\lessdot$-increasing 
sequence of forcing notions.
$\bar{\Bbb Q}$ denotes a sequence of the form 
$\langle \Bbb P_i,{\underset\tilde {}\to {\Bbb Q}_i}:i <
\alpha \rangle$ such that 
$\langle \Bbb P_i:i < \alpha)$ is a $\lessdot$-increasing sequence.
Usually ${\underset\tilde {}\to {\Bbb Q}_i}$ is a $\Bbb P_i$-name, 
$\Vdash_{{\Bbb P}_i} ``\Bbb P_{i+1}/\Bbb P_i \cong 
{\underset\tilde {}\to {\Bbb Q}_i}"$. \newline
3) Convention: If $\bar{\Bbb Q} = \langle \Bbb P_i,
{\underset\tilde {}\to {\Bbb Q}_i}:i <
\alpha \rangle,\Bbb P_i$ is $\lessdot$-increasing, we 
may write $\bar{\Bbb Q}$ instead
of $\langle \Bbb P_i:i < \alpha \rangle$. \newline
4) For a forcing notion $\Bbb P$ (as in part (1)) we 
define $\hat{\Bbb P}$:
\mr
\item "{$(a)$}"   the set of elements of $\hat{\Bbb P}$ is 
$$
\align
\biggl\{A:&A \subseteq \Bbb P,\text{ and for some } p \in A \text{ (called a 
witness) we have} \\
  &(i) \quad (\forall q \in A)(\exists r)[q \le r \in A \and p 
\le_{\text{vpr}} r] \\
  &(ii) \quad \text{there is an upper bound } p^* \in \Bbb P \text{ of } A
\text{ such that } p \le_{\text{vpr}} p^* \\
  &\qquad \text{ moreover } (\forall p' \in A)(p \le_{\text{vpr}} p' 
\Rightarrow p' \le_{vpr} p^*) \biggr\}
\endalign
$$
(we call $p^*$ an outer witness for 
$A$ or for $A \in \hat{\Bbb P}$ if clause (ii) hold), and
\sn
\item "{$(b)$}"   $\hat{\Bbb P}$ is ordered by: $A \le B$ iff: 
$A=B$ or $A = \emptyset$ or for some $q \in B,
(\forall p \in A)(p \le q)$ and we call $q$ a witness to $A \le B$
\sn
\item "{$(c)$}"  we define the order $\le_\ell$ on 
$\hat{\Bbb P}$ by: $A \le_\ell B$ \underbar{iff}
$A \le B$ and $A \ne B$ implies that for every witnesses $p$ for $A$
and every witness $q$ for $B$ we have $p \in A \and q \in B \and (\forall p')
(p \le_{\text{vpr}} p' \in A \rightarrow p' \le_\ell q)$; we call
such a pair $(p,q)$ a witness for $A \le_\ell B$. [See \scite{1.8}(5)]
\sn
\item "{$(d)$}"  we stipulate sometime $\emptyset \le_\ell A$ for
every $A \in \hat{\Bbb P}$ \ub{or} $\emptyset = \emptyset_{\Bbb P}$
[Saharon].
\ermn
5) CC$(\Bbb P)$ is the minimal regular uncountable cardinal $\theta$
such that the $\Bbb P$ satisfies the $\theta$-c.c.  We may add
\mr
\item "{$(iii)$}"  if $p' \in A$ satisfies clause $(i) + (ii)$ then
there is $p'' \in A$ such that $p' \le_{text{vpr}} p'' \and p
\le_{\text{vpr}} p''$.
\endroster
\enddefinition
\bigskip

\demo{\stag{1.1A} Observation}:  1) For any forcing notion $\Bbb P$, (as in
\scite{1.1}(1), of course), also $\hat{\Bbb P}$ is a forcing notion 
(in particular $\le^{\hat{\Bbb P}_\ell}$ is a quasi order for $\ell
\le 2$) and $\Bbb P \subseteq_{\text{ic}} \hat{\Bbb P}$ and $\Bbb P 
\lessdot \hat{\Bbb P}$
and $\Bbb P$ is $\le_{\text{vpr}}$-dense in 
$\hat{\Bbb P}$ when we identify $p$ and $\{p\}$. \nl
2) If $A_i \le_\ell B$ for $i < i^*$ then $B \le_{\text{vpr}} B^+$
where $B^+ = \dbcu_{i < i^*} A_i \cup B$. \nl
3) If $\ell \in \{0,1,2\}$ and $(\Bbb P,\le_\ell)$ is $\theta$-complete (i.e.,
an increasing sequence of length $< \theta$ has an upper bound) then so 
is $(\hat{\Bbb P},\le_\ell)$.
\enddemo
\bigskip

\demo{Proof}  1) Check. \nl
2) Easy. \nl
3) If $\delta < \theta,\langle A_i:i < \delta \rangle$ is
$\le_\ell$-increasing let $(p_i,q_i)$ witness $A_i \in \hat{\Bbb P}$.  
If $\langle p_i:i < \delta \rangle$ is eventually constant then 
$\langle A_i:i < \delta \rangle$ is eventually constant and
$A_j$ for $j$ large enough can serve.  If not, without loss of generality 
$(\forall i < \delta)p_i \ne p_{i+1}$ and let $(p_i,q_i)$
witness $A_i \le_\ell A_{i+1}$.  Clearly $\langle q_i:i < \delta \rangle$ has
a $\le_\ell$-upper bound in $\Bbb P$, call it $q$.  Now $\{q\}$ is 
as required.  \hfill$\square_{\scite{1.1A}}$\margincite{1.1A}
\enddemo
\bigskip

\definition{\stag{1.2} Definition}  Let $MAC(\Bbb P)$ be the set of maximal
antichains of the forcing notion $\Bbb P$.
\enddefinition
\bigskip

\remark{\stag{1.3} Remark}  1) Note: $|MAC(\Bbb P)| \le 2^{|\Bbb
P|},\Bbb P$ satisfies the
$|\Bbb P|^+$-c.c. and if $\Bbb P$ satisfies the $\lambda$-c.c. then
$|MAC(\Bbb P)| \le |\Bbb P|^{< \lambda}$. \newline
2) Note:
\mr
\item "{$(*)$}"  if 
$\Bbb Q$ is a forcing notion, $\lambda = \lambda^{< \lambda}
> |\Bbb Q| + \aleph_0,\Vdash_\Bbb Q ``(\forall \mu < \lambda)
\mu^{\aleph_0} < \lambda"$
and $\Bbb Q' = \Bbb Q * \text{ Levy}(\aleph_1,< \lambda)$ \ub{then}
$|MAC(\Bbb Q')| = |\Bbb Q'| =
\lambda$.
\endroster
\endremark
\bigskip

\demo{\stag{1.4} Notation}  Car is the class of cardinals. \newline
IRCar is the class of (infinite) regular cardinals. \newline
RCar = IRCar $\cup \{1\}$. \newline
URCar is the class of uncountable regular cardinals. \newline
${\Cal D}^{cb}_\lambda$ is the filter of co-bounded subsets of $\lambda$.
\newline
${\Cal D}_\lambda$ is the club filter on $\lambda$ for $\lambda$
regular uncountable. \nl
$\eta^- = \eta \restriction (\ell g(\eta) - 1)$ for a finite sequence $\eta$
of length $>0$.
\enddemo
\bigskip

\demo{\stag{1.5} Notation}  1) ${\Cal H}(\chi)$ is the family of sets with
transitive closure of power $< \chi$; let $<^*_\chi$ be a well ordering of
${\Cal H}(\chi)$. \nl
2) Let $W$ be a function from the set of strongly inaccessible cardinals to
$\{0,1,\frac 12\}$; if $\alpha \notin \text{ Dom}(W)$ we understand
$W(\alpha) = 0$ and let $\alpha \in W$ means $W(\alpha) = 1$. 
\enddemo
\bigskip

\definition{\stag{1.6} Definition}  1) Assume $\bar{\Bbb P}$ is a 
$\lessdot$-increasing sequence of forcing notions. \newline
Let

$$
\align
\text{Gen}^r(\bar{\Bbb P}) =: 
\biggl\{ G:&\text{ for some (set) forcing notion }
\Bbb P^* \text{ we have } \dsize \bigwedge_{i < \alpha} \Bbb P_i \lessdot
\Bbb P^* \\
  &\text{ and for some } G^* \subseteq 
\Bbb P^* \text{ generic over } \bold V \text{ we have } \\
  &\,G = G^* \cap \dsize \bigcup_{i < \alpha} \Bbb P_i \biggr\}.
\endalign
$$
\mn
2) If $\bar{\Bbb Q} = \langle \Bbb P_i:i < \alpha \rangle$ 
or $\bar{\Bbb Q} = \langle  
\Bbb P_i,{\underset\tilde {}\to {\Bbb Q}_i}:i < \alpha \rangle$ where
$\Bbb P_i$ is a $\lessdot$-increasing we define a 
$\bar{\Bbb Q}$-name $\underset\tilde {}\to \tau$
almost as we define $(\dsize \bigcup_{i < \alpha} 
\Bbb P_i)$-names, but we do not
use maximal antichains of $\dsize \bigcup_{i < \alpha} \Bbb P_i$, that is:
\medskip
\roster
\item "{$(*)$}"  $\underset\tilde {}\to \tau$ is a function, 
Dom$(\underset\tilde {}\to \tau) \subseteq \dsize \bigcup_{i < \alpha}
\Bbb P_i$ and for every $G \in \text{ Gen}^r(\bar{\Bbb Q}),
\underset\tilde {}\to \tau[G]$ is defined iff 
Dom$(\underset\tilde {}\to \tau)
\cap G \ne \emptyset$ and then $\underset\tilde {}\to \tau[G] \in
\bold V[G]$
[from where ``every $G \ldots$" is taken?  E.g., $\bold V$ is countable, $G$ any set
from the true universe] and $\underset\tilde {}\to \tau[G]$ is definable with
parameters from $\bold V$ and the parameter $\dbcu_{i < \alpha} \Bbb P_i \cap G$ 
(so $\underset\tilde {}\to \tau$ is really a first-order formula with the 
variable $\dbcu_{i < \alpha} \Bbb P_i \cap G$ and parameters from $\bold V$).
\endroster
\medskip

\noindent
Now $\Vdash_{\bar{\Bbb Q}}$ has a natural meaning. \newline
3) For $p \in \bar{\Bbb Q}$ (i.e. $p \in \dsize \bigcup_{i < \alpha}
\Bbb P_i$) and
$\bar{\Bbb Q}$-names ${\underset\tilde {}\to \tau_0},\dotsc,
{\underset\tilde {}\to \tau_{n-1}}$ we let $\{{\underset\tilde {}\to \tau_0},
\dotsc,{\underset\tilde {}\to \tau_{n-1}}\}$ be the name for the set that
contains exactly those ${\underset\tilde {}\to \tau_i}[\bar{\Bbb Q}]$ 
that are
defined.  We let $p \Vdash ``\underset\tilde {}\to \tau = x"$ if for every
$G$ such that $p \in G \in \text{ Gen}^r(\bar{\Bbb Q})$ we have
$\underset\tilde {}\to \tau[G] = x$. 
If $\beta < \alpha$ and $G_\beta \subseteq \Bbb P_\beta$, we let
$\underset\tilde {}\to \tau[G_\beta] = x$ means that for some 
$p \in G_\beta$ we have $p \Vdash_{\bar{\Bbb Q}} 
``\underset\tilde {}\to \tau = x"$, so possibly no $p \in \dbcu_{i <
\alpha} \Bbb P_i$ forces a value to $\underset\tilde {}\to \tau$ and
no such $p$ forces $\tau$ is not definable. \newline
4) We say a $\bar{\Bbb Q}$-name $\underset\tilde {}\to x$ is full
\ub{if} $\underset\tilde {}\to x[G]$ is well defined for every $G \in
\text{ Gen} {}^r(\bar{\Bbb Q})$. \nl
5) A simple $\bar{\Bbb Q}$-named$^1$ \, $[j,\beta)$-ordinal 
$\underset\tilde {}\to \zeta$ 
is a $\bar{\Bbb Q}$-name $\underset\tilde {}\to \zeta$ such that: if
$G \in \text{ Gen}^r(\bar{\Bbb Q})$ and
$\underset\tilde {}\to \zeta[G] = \xi$ then $j \le \xi < \beta$ and
$(\exists p \in G \cap \Bbb P_{\xi \cap \alpha}) p \Vdash_{\bar{\Bbb Q}}
``\underset\tilde {}\to \zeta = \xi"$ (where $\alpha = 
\ell g(\bar{\Bbb Q})$);
however, we allowed $\underset\tilde {}\to \zeta[G]$ to be undefined.
If we omit ``$[j,\beta)$" we mean $[0,\ell g(\bar{\Bbb Q})) = 
[0,\alpha)$.  If we
omit ``simple", we mean 
replacing $(\exists p \in G \cap \Bbb P_{\xi \cap \alpha})$
by $(\exists p \in G \cap \Bbb P_{(\xi + 1) \cap \alpha})$ (this is
used in \cite[Ch.X,\S1]{Sh:f}, we shall only remark on it here). \nl
6) A simple $\bar{\Bbb Q}$-named$^2$ \, $[j,\beta)$-ordinal
$\underset\tilde {}\to \zeta$ is a simple $\bar{\Bbb Q}$-named$^2$ \, 
$[j,\beta)$-ordinal of depth $\Upsilon$ for some ordinal $\Upsilon$, where
this is defined below by induction on $\Upsilon$.  In all cases it is a
$\bar{\Bbb Q}$-name of an ordinal from the interval $[j,\beta)$ so may 
be undefined, i.e., we allow non full such names.
\enddefinition
\bn
\ub{Case 1}:  $\Upsilon = 0$.

This is an ordinal $\in [j,\beta)$, or is ``undefined" (in the full
case this is forbidden).
\bn
\ub{Case 2}:  $\Upsilon > 0$.

For some $\gamma < \ell g(\bar{\Bbb Q}) \cap \beta$ and maximal antichain
${\Cal I} = \{p_\varepsilon:\varepsilon < \varepsilon^*\}$ of 
$\Bbb P_\gamma$, for each $i < i^*$ there is a sequence
$\langle {\underset\tilde {}\to \zeta_\varepsilon}:
\varepsilon < \varepsilon^* \rangle$
such that ${\underset\tilde {}\to \zeta_\varepsilon}$ is a simple 
$\bar{\Bbb Q}$-named$^2$ \, $[\text{Max}\{j,\gamma\},\beta)$-ordinal of depth
$\Upsilon_\varepsilon < \Upsilon$ and: $\underset\tilde {}\to \zeta
[G_\xi] = \xi$ iff $\xi \ge \gamma$ and for some $\varepsilon$ we have
$p_\varepsilon \in G_\xi \cap \Bbb P_\gamma$ and 
${\underset\tilde {}\to \zeta_\varepsilon}[G_\xi] = \xi$ 
(including the case: not defined).  If we omit ``$[j,\beta)$" we mean
$[0,\ell g(\bar{\Bbb Q})) = [0,\alpha)$. \nl
7)  If we omit ``simple" in (6) we mean that in case 2,
${\underset\tilde {}\to \zeta_i}$ is a not necessarily simple
$\bar{\Bbb Q}$-name$^r$ and ${\Cal I}
\subseteq \Bbb P_{\gamma +1}$.
\nl
8) We say $\bar{\Bbb P}$ is $W$-continuous or $(\bar{\Bbb P},W)$ is
continuous when for every $\delta \in W \cap \ell g(\bar{\Bbb P})$
when if $(\forall i < \delta)$ [density $(\Bbb P_i) < \delta$, or just
$\Bbb P_i$ satisfies the cf$(\delta)$-c.c.], 
\ub{then} $\Bbb P_\delta = \dbcu_{i < \delta}\Bbb P_i$.  
We say $\Bbb P$ is $W$-smooth or $(\bar{\Bbb P},W)$ is smooth
if $\delta \in W \Rightarrow \Bbb P_\delta = \dbcu_{i < \delta} \Bbb P_i$.
We say $\underset\tilde {}\to \zeta$ is a simple$^\ell \,
(\bar{\Bbb Q},W)$-named
$[j,\beta)$-ordinal \ub{if}
it is a simple$^\ell \, \bar{\Bbb Q}$-named ordinal 
and $\delta \in W \cap (\ell g(\bar{\Bbb Q}) +1) 
\Rightarrow (\exists \alpha < \delta)(\Vdash
\underset\tilde {}\to \zeta \notin [\alpha,\delta))$. 
\bigskip

\proclaim{\stag{1.6A} Claim}  1) Assume that $\bar Q$ is $W$-continuous.
If $\underset\tilde {}\to \zeta$ is a simple
$\bar{\Bbb Q}$-named$^2 \, [0,\alpha)$-ordinal, $\gamma \in W$ is regular
and $\beta < \gamma$ implies density$(\Bbb P_\beta) < \gamma$ or just
($\Bbb P_\beta$ satisfies the ${\text{\rm cf\/}}(\alpha))$-c.c., 
\ub{then} for some $\beta < \gamma,\underset\tilde {}\to \zeta$ is 
a simple $\bar{\Bbb Q} \restriction
\beta$-named$^2 \, [0,\beta)$-ordinal. \nl
2) If $\bar{\Bbb Q}$ is $W$-continuous and $\gamma \in W \Rightarrow
\gamma$ regular and $\underset\tilde {}\to \zeta$ is a simple 
$\bar{\Bbb Q}$-named$^2 [0,\alpha)$-ordinal \ub{then} $\underset\tilde
{}\to \zeta$ is a simple $(\bar{\Bbb
Q},W)$-named$^2[0,\alpha)$-ordinal. \nl
3) If $\underset\tilde {}\to \zeta$ is a simple $\bar{\Bbb
Q}$-named$^2[0,\alpha)$-ordinal \ub{then} there is a full simple
$\bar{\Bbb Q}$-named$^2[0,\alpha)$-ordinal ${\underset\tilde {}\to
\zeta'}$ such that $\Vdash_{\bar{\Bbb Q}}$ ``if $\underset\tilde {}\to
\zeta$ is well defined then it is equal to $\zeta'$".
\endproclaim
\bigskip

\demo{Proof}  By induction on the depth of $\underset\tilde {}\to \zeta$.
\hfill$\square_{\scite{1.6A}}$\margincite{1.6A}
\enddemo
\bigskip

\remark{\stag{1.7} Remark}  1) We can restrict in the definition of
Gen$^r(\Bbb Q)$ to $\Bbb P^*$ in some class $K$, and get a $K$-variant
of our notions. \newline
2) Note: even if in \scite{1.6}(1) we ask Dom$(\underset\tilde {}\to \tau)$
to be a maximal antichain of $\dbcu_{i < \delta} \Bbb P_i$
 it will not be meaningful as in the appropriate
$\Bbb P_\delta$, we have $\dsize \bigwedge_{i < \delta} \Bbb P_i \lessdot
\Bbb P_\delta$ but not necessarily $\dbcu_{i < \delta} \Bbb P_i
\lessdot P_\delta$ hence it will not in general be a maximal antichain. \nl
3) Note that in the simple case we wrote 
$\Bbb P_{\xi \cap \alpha}$ not $\Bbb P_{(\xi + 1) \cap \alpha}$.
Compare this the remark \cite[Ch.XIV,1.1B]{Sh:f}.  Here in the main
case we use full simple $\bar{\Bbb Q}$-named$^2$ ordinals, 
though we shall remark on the affect of the non-simple case;
as a result we will not have a general associativity law, but
the definition of Sp$_3 - \text{ Lim}_\kappa
(\bar{\Bbb Q})$ will be somewhat simplified.  
As said earlier, we can interchange
decisions on this matter.  Of course, also \cite[Ch.XV]{Sh:f} can be 
represented with this iteration. \nl
4) The ``name$^1$" is necessary for the $\kappa > \aleph_1$ case, but
``name$^2$" is preferable for $\kappa = \aleph_1$, so we could have
concentrated on name$^1$ for $\kappa > \aleph_1$, name$^2$ for $\kappa
= \aleph_1$, but actually we concentrated on simple, name$^2$ for
$\kappa = \aleph_1$; see \scite{1.13}(B) below.
\endremark
\bigskip

\demo{\stag{1.8} Fact}  1) For $\bar{\Bbb P} = \langle \Bbb P_i:
i < \ell g(\bar{\Bbb P})
\rangle$, a $\lessdot$-increasing sequence of forcing notions, $\ell \in
\{1,2\}$ and simple $\bar{\Bbb P}$-named$^\ell$ \, $[j,\beta)$-ordinal 
$\underset\tilde {}\to \zeta$ and $p \in
\dsize \bigcup_{i < \alpha} \Bbb P_i$ there 
are $\xi,q$ and $q_1$ such that $p \le
q \in \dsize \bigcup_{i < \ell g(\bar{\Bbb P})}
\Bbb P_i$ and: either $q \Vdash_{\bar{\Bbb P}}
``q_1 \in \underset\tilde {}\to G",q_1 \in \Bbb P_\xi,\xi < \alpha,[p \in
\Bbb P_\xi \Rightarrow q = q_1]$ and $q_1 \Vdash_{\bar{\Bbb P}} 
``\underset\tilde {}\to \zeta 
= \xi"$ or $q \Vdash_{\bar{\Bbb P}} 
``\underset\tilde {}\to \zeta$ is not defined"
(and even $p \Vdash_{\bar{\Bbb P}} 
``\underset\tilde {}\to \zeta$ is not defined").
\newline
2) For $\bar{\Bbb P}$ and $\ell \in \{1,2\}$ as above, and simple 
$\bar{\Bbb P}$-named$^\ell$ \, $[j,\beta)$-ordinals 
$\underset\tilde {}\to \zeta,\underset\tilde {}\to \xi$, also
max$\{\underset\tilde {}\to \zeta,\underset\tilde {}\to \xi\}$ and
min$\{\underset\tilde {}\to \zeta,\underset\tilde {}\to \xi\}$ are simple
$\bar{\Bbb Q}$-named$^\ell$ \, $[j,\beta)$-ordinals (naturally
defined, so max $\{\underset\tilde {}\to \zeta,
\underset\tilde {}\to \xi\}[G]$
is defined iff a $\underset\tilde {}\to \zeta[G],
\underset\tilde {}\to \xi[G]$ are defined, and 
min$\{\underset\tilde {}\to \zeta,\underset\tilde {}\to \xi\}
[G]$ is defined iff $\underset\tilde {}\to \zeta[G]$ is defined or
$\underset\tilde {}\to \xi[G]$ is defined). 
If $\underset\tilde {}\to \zeta,\underset\tilde {}\to \xi$ are full
then so are max$\{\underset\tilde {}\to \zeta,\underset\tilde {}\to
\xi\}$ and min$\{\underset\tilde {}\to \zeta,\underset\tilde {}\to \xi\}$.
\newline
3) For $\bar{\Bbb P}$ and $\ell$ as above, $n < \omega$ and simple 
$\bar{\Bbb P}$-named$^\ell$ ordinals
${\underset\tilde {}\to \xi_1},\dotsc,{\underset\tilde {}\to \xi_n}$ and
$p \in \dsize \bigcup_{i < \ell g(\bar{\Bbb P})} \Bbb P_i$ there are 
$\zeta < \alpha$ and
$q \in \Bbb P_\zeta$ such that, first: $p\le q$ or at least $q
\Vdash_{{\Bbb P}_\zeta} 
``p \in \Bbb P_i/{\underset\tilde {}\to G_{{\Bbb P}_\zeta}}$ 
for some $i < \ell g
(\bar{\Bbb P})"$ and second: for some $\ell \in \{1,\dotsc,n\}$ we have $q 
\Vdash_{\bar{\Bbb P}} ``\zeta = 
{\underset\tilde {}\to \xi_\ell} = \text{max}\{{\underset\tilde {}\to \xi_1},
\dotsc,{\underset\tilde {}\to \xi_n}\}$" or $q \Vdash_{\bar{\Bbb P}} ``\text{ max}
\{{\underset\tilde {}\to \xi_1},\dotsc,{\underset\tilde {}\to \xi_n}\}$ not
defined".  Similarly for min. \newline
4) The same holds for simple $(\bar{\Bbb Q},W)$-names$^\ell$ and we 
can omit simple. \nl
5) If ${\underset\tilde {}\to \zeta_i}$ is a simple $\bar{\Bbb Q}$-named$^1 
[\beta_i,\gamma_i)$-ordinal for $i < i^*$ \ub{then} 
\footnote{this seems lacking for ``name$^2$".} 
$\sup\{{\underset\tilde {}\to \zeta_i}:i < i^*\}$ 
is simple $\bar{\Bbb Q}$-named$^1 \, [{\underset i {}\to {\text{Min}}} \, \beta_i,
{\underset i {}\to {\text{sup}}} \gamma_i)$-ordinal. \nl
6) Similarly for $\{{\underset \tilde {}\to \zeta_i}:i\}$ and when we
omit ``simple". \nl
7A)  A simple $\bar{\Bbb P}$-name$^\ell[j,\beta)$-ordinal
\mr
\item "{$(a)$}"  $\zeta$ is a $\bar{\Bbb P}$-named$^\ell[j,\beta)$-ordinal;
\sn
\item "{$(b)$}"  if $j_2 \le j_1 < \beta_1 \le \beta_2$ then any
[simple] $\bar{\Bbb P}$-named $[j_1,\beta_1)$-ordinal is a [simple]
$\bar{\Bbb P}$-named$^\ell [j_2,\beta_2)$-ordinal;
\sn
\item "{$(c)$}"  a [simple] $\bar{\Bbb P}$-named$^2 [j,\beta)$ ordinal
is a [simple] $\bar{\Bbb P}$-named$^1 [j,\beta)$-ordinal
\sn
\item "{$(d)$}"  if $\beta \le \alpha' \le \ell g(\bar{\Bbb P})$ then
any [simple] $\bar{\Bbb P}$-named$^\ell [j,\beta)$-ordinal is a
$(\bar{\Bbb P} \restriction \alpha')$-named$^\ell [j,\beta$)-ordinal.
\endroster
\enddemo
\bigskip

\demo{Proof}  Straight.
\enddemo
\bigskip

\demo{\stag{1.9} Discussion}  We have in defining our iteration 
several possible variants, some of our particular choices are not 
important: we can
make it like revised countable support as in \cite[Ch.X,\S1]{Sh:f} or like
$\aleph_1$-RS in \cite[Ch.XIV,\S1]{Sh:f}, or as in \cite[Ch.XIV,\S2]{Sh:f}
(as here); for most uses $\kappa = \aleph_1$ and we could restrict ourselves
to $\le_{\text{pr}} = \le_{\text{vpr}}$ as equality; but in 
\cite{GoSh:511} we need the three partial orders.

So below we have finite support for non-pure, countable for pure and Easton
for very pure.
\enddemo
\bigskip

\definition{\stag{1.10} Definition}  1) For a forcing notion $\Bbb P$ (as in 
\scite{1.1}) let $\Bbb P^{[\text{pr}]}$ 
be defined like $\Bbb P$ except that we make
$\emptyset_{\Bbb P} \le_{\text{pr}} p$ for every $p \in \Bbb P$. \newline
2) For a forcing notion $\Bbb P$ 
(as in \scite{1.1}) let $\Bbb P^{[\text{vp}]}$ be defined 
like $\Bbb P$ except that we make $\emptyset_{\Bbb P} \le_{\text{vp}} p$
(and $\emptyset \le_{\text{pr}} \Bbb P$, of course) for every $p \in \Bbb P$.
\enddefinition
\bigskip

\demo{\stag{1.11} Fact}  1) For a forcing notion $\Bbb P$ and 
$x \in \{\text{pr,vp}\},
\Bbb P^{[x]}$ is also a 
forcing notion, and they are equivalent as forcing notions. \nl
2) For $x \in \{\text{pr,vp}^r\}$ the operations $\Bbb P \mapsto
\hat{\Bbb P}$ and $\Bbb P \mapsto \Bbb P^{[x]}$ commute. \nl
3) If $(x_1,x_2,x_3) \in
\{(\text{pr,pr,pr),(vpr,pr,vpr),(vpr,vpr,pr)),(vpr,vpr,vpr)}\}$
then $\Bbb P^{[x_3]} = (\Bbb P^{[x_1][x_2]}$. \nl
4) $\theta$-completeness is preserved in the natural cases.
\enddemo
\bigskip

\demo{\stag{1.12} Discussion}  1) Why do we 
bother with $\Bbb P^{[\text{pr}]},\Bbb P^{[\text{vpr}]}$?  If
in the iteration defined below we use only $Q^{[\text{pr}]}_i,
{\underset\tilde {}\to {\Bbb Q}^{[\text{vp}]}_i}$, we get a variant 
of the definition
without the need to repeat it.  We may want that: if $\ell g(\bar{\Bbb Q}) = 
\lambda$ inaccessible and $i < \kappa \Rightarrow |\Bbb P_i| < \lambda$ then
$\dsize \bigcup_{i < \lambda} \Bbb P_i = \Bbb P_\lambda$ (here as done in 
\cite[Ch.XIV,\S2]{Sh:f} we can just impose it). \nl
Some other restrictions are for simplicity only. \nl
2) Below the case $e=6$ is the main one. [Saharon]
\enddemo
\bigskip

\definition{\stag{1.13} Definition/Claim}  We define and prove the
following by induction on $\alpha$.  Below $\kappa = \aleph_1,e = G$
so we can omit them (they are meaningful in \S11)
\mr
\item "{$(A)$}"  \,\,[Definition] $\quad \bar{\Bbb Q} = 
\langle \Bbb P_i,{\underset\tilde {}\to Q_i}:i < \alpha
\rangle$ is a $\kappa-\text{Sp}_e$-iteration for 
$W$ or $\kappa-\text{Sp}_e(W)$-iteration
(if $W$ is absent we mean $\{ \beta \le \alpha:\beta 
\text{ strongly inaccessible}\}$); $\alpha$ is called the length of
$\bar{\Bbb Q},\ell g(\bar{\Bbb Q})$.
\sn
\item "{$(B)$}"  \,\,[Definition] $\quad$ A simple 
$(\bar{\Bbb Q},W)$-named$_e$ ordinal
$\underset\tilde {}\to \zeta$ and $\underset\tilde {}\to \zeta \restriction
[\alpha,\beta)$. 
\sn
\item "{$(C)$}"  \,\,[Definition] $\quad$ A simple 
$(\bar{\Bbb Q},W)$-named$_e$ 
atomic condition $\underset\tilde {}\to q$
(or atomic $[j,\beta)$-condition where $j \le \beta \le \alpha$); also 
we define $\underset\tilde {}\to q \restriction \xi,
\underset\tilde {}\to q \restriction \{\varepsilon\},
\underset\tilde {}\to q \restriction [\xi,\zeta)$ for a simple 
$\bar{\Bbb Q}$-named$_e$ atomic condition 
$\underset\tilde {}\to q$ and ordinals $\varepsilon <
\alpha,\xi \le \zeta \le \alpha$ 
(or simple $\bar{\Bbb Q}$-named$_e$ ordinals 
$\underset\tilde {}\to \xi,\underset\tilde {}\to \zeta$ 
instead $\xi,\zeta$).  We may add pure/very pure
as adjectives to the condition.  
\sn
\item "{$(D)$}"  \,\,[Claim] $\quad$  Assume 
$\underset\tilde {}\to \zeta$ is a simple 
$(\bar{\Bbb Q},W)$-named$_e [j,\beta)$-ordinal.  \ub{Then} for any $\xi,
\underset\tilde {}\to \zeta \restriction \xi$ is a 
simple $(\bar{\Bbb Q},W)$-named
$[j,\text{min}\{\beta,\xi\})$-ordinal and $\Vdash_{\bar{\Bbb Q}}$ ``if
$\underset\tilde {}\to \zeta < \xi$ then $\underset\tilde {}\to \zeta =
\underset\tilde {}\to \zeta \restriction \xi$; if
$\underset\tilde {}\to \zeta \ge \xi$ or $\underset\tilde {}\to \zeta$ is
undefined, then $\underset\tilde {}\to \zeta \restriction \xi$ is undefined",
also $\underset\tilde {}\to \zeta \restriction \xi$ is a simple
$(\bar{\Bbb Q} \restriction \xi,W)$-named ordinal. \nl
Similarly $\underset\tilde {}\to \zeta \restriction [\xi_1,\xi_2)$.  If
${\underset\tilde {}\to \xi_1} \le {\underset\tilde {}\to \xi_2} \le \alpha$
are simple $(\bar{\Bbb Q},W)$-named $[\alpha_1,\alpha_2)$-ordinals (the
${\underset\tilde {}\to \xi_1} \le {\underset\tilde {}\to \xi_2}$ means
$\Vdash_{\bar{\Bbb Q}} ``{\underset\tilde {}\to \xi_1} \le
{\underset\tilde {}\to \xi_2}"$), \ub{then}
$\underset\tilde {}\to \zeta \restriction [{\underset\tilde {}\to \xi_1},
{\underset\tilde {}\to \xi_2})$ is a simple $(\bar{\Bbb Q},W)$-named
$[\alpha_1,\alpha_2)$-ordinal and $\Vdash_{\bar{\Bbb Q}}$ ``if
$\underset\tilde {}\to \zeta \in [{\underset\tilde {}\to \xi_1},
{\underset\tilde {}\to \xi_2})$, then $\underset\tilde {}\to \zeta =
\underset\tilde {}\to \zeta \restriction [{\underset\tilde {}\to \xi_2},
{\underset\tilde {}\to \xi_2})$ otherwise $\underset\tilde {}\to \zeta
\restriction [\xi_1,\xi_2)$ is undefined.  If in addition $\beta =
\text{ Min}\{\alpha,\alpha_2,\ell g(\bar{\Bbb Q})\}$ and $\beta \le \gamma,
\alpha'_1 \le \alpha_1$ then $\underset\tilde {}\to \zeta
\restriction[{\underset\tilde {}\to \xi_1},{\underset\tilde {}\to \xi_2})$ 
is a simple $(\bar{\Bbb Q} \restriction
\beta,W)$-named $[\alpha'_1,\gamma)$-ordinal.  Also if
$n < \omega$, for $\ell \in \{1,\dotsc,n\},{\underset\tilde {}\to \xi_\ell}$
is a simple $\bar{\Bbb Q}$-named $[\beta_1,\beta_2)$-ordinal then
Max$\{{\underset\tilde {}\to \xi_1},\dotsc,{\underset\tilde {}\to \xi_n}\}$
is a simple $\bar{\Bbb Q}$-named 
$[\beta_1,\beta_2)$-ordinal.  Similarly for Min.
\sn
\item "{$(E)$}"  \,\,[Claim] $\quad$ If 
$q$ is a simple $(\bar{\Bbb Q},W)$-named 
atomic $[j,\beta)$-condition, $\xi < \alpha$, then $\underset\tilde {}\to q 
\restriction \xi$ is a simple $(\bar{\Bbb Q} \restriction
\xi,W)$-named atomic 
$[j,\text{min}\{\beta,\xi\})$-condition and $\underset\tilde {}\to q 
\restriction \{\xi\}$ is a $\Bbb P_\xi$-name of a member of 
${\underset\tilde {}\to {\Bbb Q}_\xi}$ or undefined (and then 
it may be assigned the value 
$\emptyset_{{\underset\tilde {}\to {\Bbb Q}_\xi}}$, the minimal member of
${\underset\tilde {}\to {\Bbb Q}_\xi}$).  If $q$ is a simple 
$(\bar{\Bbb Q},W)$-named 
atomic condition, $\underset\tilde {}\to \xi \le \underset\tilde {}\to \zeta 
\le \alpha$ are simple $(\bar{\Bbb Q},W)$-named ordinals then 
$\underset\tilde {}\to q \restriction[\underset\tilde {}\to \xi,
\underset\tilde {}\to \zeta)$ is a simple $(\bar{\Bbb Q},W)$-named ordinal.  
Also $\underset\tilde {}\to q \restriction \{\zeta\} = 
\underset\tilde {}\to q
\restriction [\zeta,\zeta +1)$, and if $\underset\tilde {}\to q$ is a simple
$(\bar{\Bbb Q},W)$-named $[\zeta,\xi)$-ordinal, $\zeta' < \xi'$ and
$\Vdash_{\bar{\Bbb Q}} ``\underset\tilde {}\to \zeta \in
[\zeta',\xi')"$, then it is a simple $(\bar{\Bbb Q},W)$-named 
$[\zeta',\xi')$-ordinal.
Also ``pure" and ``very pure" are preserved by restriction.
\sn  
\item "{$(F)$}" \,\,[Definition] $\quad$  The $\kappa-\text{Sp}_e(W)$-limit of 
$\bar{\Bbb Q},\text{Sp}_e(W)$-Lim$_\kappa(\bar{\Bbb Q})$,
denoted by $\Bbb P_\alpha$ for $\bar{\Bbb Q}$ as in clause (A) in
particular of length $\alpha$, and $p \restriction
\xi$ and Dom$(p)$ for $p \in \text{ Sp}_e(W)$-Lim$_\kappa
(\bar{\Bbb Q}),\xi$ an ordinal $\le \alpha$; (similarly for a
simple $(\bar{\Bbb Q},W)$-named $[0,\ell g(\bar{\Bbb Q}))$ 
ordinal $\underset\tilde {}\to \xi$, etc.  
We also define $\Bbb P_{\underset\tilde {}\to \zeta}$ for
$\underset\tilde {}\to \zeta$ a $(\bar{\Bbb Q},W)$-named ordinal.
\sn
\item "{$(G)$}" \,\,[Theorem] $\quad \text{ Sp}_e(W)$-Lim$_\kappa(\bar{\Bbb
Q})$ is a forcing notion (in the sense of \scite{1.1}(1)).
\sn
\item "{$(H)$}" \,\,[Theorem] $\quad$ Assume $\beta < \alpha = 
\ell g(\bar{\Bbb Q})$ or more generally, $\underset\tilde {}\to \beta$ 
is a full simple $(\bar{\Bbb Q},W)$-named ordinal (see end of clause 
(F) above).  Then $\Bbb P_\beta \subseteq_{\text{ic}} 
\text{ Sp}_e(W)$-Lim$_\kappa(\bar{\Bbb Q})$ (so a submodel 
with the three partial orders, 
even compatibilities are preserved) and 
$[p \in \Bbb P_\beta \Rightarrow p \restriction \beta = p]$ and 
$[\Bbb P_\alpha \models ``p \le_\ell q" \Rightarrow \Bbb P_\beta \models
``p \restriction \beta \le_\ell q \restriction \beta"]$ (for $\ell = 0,1,2$,
of course) and $\Bbb P_\alpha \models ``p \restriction \beta \le_\ell
p"$.  Also $q \in \Bbb P_\beta,p \in 
\text{ Sp}_e(W)$-Lim$_\kappa(\bar{\Bbb Q})$ are compatible iff $q,p
\restriction \beta$ are compatible in 
$\Bbb P_\beta$.  In fact, if $q \in \Bbb P_\beta,
\Bbb P_\beta \models ``p \restriction \beta \le q"$ 
then $q \cup (p \restriction [\beta,\alpha))$ belongs 
to $\Bbb P_\alpha$ and is a least upper bound of $p,q$ 
and if $\Bbb P_\beta \models ``p \restriction \beta \le_\ell q"$ even a 
$\le_\ell$-least upper bound of $q$.
Hence $\Bbb P_\beta \lessdot (\text{Sp}_e(W)$-Lim$_\kappa(\bar q))$.
\sn
\item "{$(I)$}"  \,\,[Claim] $\quad$  The set of 
$p \in \Bbb P_\alpha$ such that for every $\beta < \alpha$ 
we have $\Vdash_{{\Bbb P}_\beta} ``p \restriction \{\beta\}$ is a singleton
or empty", is a dense subset of $\Bbb P_\alpha$.  Also we can replace
${\underset\tilde {}\to {\Bbb Q}_\beta}$ by $\hat{\Bbb Q}_\beta$ and 
the set of ``old" $p \in \Bbb P_\alpha$ is a dense subset of the new 
(but actually do not use this).
\sn
\item "{$(J)$}"  \,\, [Claim] $\quad$  If 
$\alpha$ is strongly inaccessible, $\zeta < \alpha
\Rightarrow |\Bbb P_\zeta| < \alpha$ or just $\zeta < \alpha
\Rightarrow CC(\Bbb P_\zeta) \le \alpha$
 and $\alpha \in W$, \ub{then} $\Bbb P_\alpha =
\dbcu_{\zeta < \alpha} \Bbb P_\zeta$.
\endroster
\enddefinition
\bn
\underbar{Proof and Definition}:  
\mr
\item "{$(A)$}"  $\bar{\Bbb Q} = \langle \Bbb P_i,
{\underset\tilde {}\to {\Bbb Q}_i}:i < \alpha
\rangle$ is a $\kappa-\text{Sp}_e(W)$-iteration if 
$\bar{\Bbb Q} \restriction \beta$ is a
$\kappa-\text{Sp}_e(W)$-iteration for 
every $\beta < \alpha$, and if $\alpha = 
\beta + 1$, \ub{then} $\Bbb P_\beta = 
\text{ Sp}_e(W)$-Lim$_\kappa(\bar{\Bbb Q} \restriction 
\beta)$ 
and ${\underset\tilde {}\to {\Bbb Q}_\beta}$ is a $\Bbb P_\beta$-name of a 
forcing notion as in \scite{1.1}(1) here.
\sn
\item "{$(B)$}"  We say $\underset\tilde {}\to \zeta$ is a simple
$(\bar{\Bbb Q},W)$-named$_e$ \, $[j,\beta)$-ordinal \ub{if}
$\underset\tilde {}\to \zeta$ is a simple $(\bar{\Bbb Q},W)$-named$^2$ \, 
$[j,\beta)$-ordinal.
\sn
\item "{$(C)$}"  We say $\underset\tilde {}\to q$ is a simple 
$(\bar{\Bbb Q},W)$-named atomic $[j,\beta)$-condition when: 
$\underset\tilde {}\to q$ is a $\bar{\Bbb Q}$-name, and for some 
$\underset\tilde {}\to \zeta = 
{\underset\tilde {}\to \zeta_{\underset\tilde {}\to q}}$, a simple
$(\bar{\Bbb Q},W)$-named $[j,\beta)$-ordinal, we have 
$\Vdash_{\bar{\Bbb Q}} ``\underset\tilde {}\to \zeta$
has a value iff $\underset\tilde {}\to q$ has, and if they have then $j \le
\underset\tilde {}\to \zeta < \text{ min}\{\beta,\ell g(\bar{\Bbb Q})\}$ and
$\underset\tilde {}\to q \in 
{\underset\tilde {}\to {\Bbb Q}_{\underset\tilde {}\to \zeta}}"$.  If we omit
``$[j,\beta)$" we mean ``$[0,\alpha)$".
Now $\underset\tilde {}\to q \restriction \xi$ will have a value iff
${\underset\tilde {}\to \zeta_{\underset\tilde {}\to q}}$ has a value $< \xi$
and then its value is the value of $\underset\tilde {}\to q$.  Lastly,
$\underset\tilde {}\to q \restriction \{\xi\}$ will have a value iff
${\underset\tilde {}\to \zeta_{\underset\tilde {}\to q}}$ has value $\xi$ and
then its value is the value of $\underset\tilde {}\to q$.  Similarly for
$\underset\tilde {}\to q \restriction [\zeta,\xi)$ and 
$\underset\tilde {}\to q \restriction \underset\tilde {}\to \xi,
\underset\tilde {}\to q \restriction \{ \underset\tilde {}\to \xi\}$,
$\underset\tilde {}\to q \restriction
[\underset\tilde {}\to \zeta,\underset\tilde {}\to \xi)$.  We say
$\underset\tilde {}\to q$ is pure if 
$\Vdash_{\bar{\Bbb Q}}$ ``for $\xi < \alpha$,
if ${\underset\tilde {}\to \zeta_q} = \xi$ and a
$\underset\tilde {}\to q \in {\underset\tilde {}\to {\Bbb Q}_\xi}$ then
${\underset\tilde {}\to {\Bbb Q}_\xi} \models 
\emptyset_{\underset\tilde {}\to {\Bbb Q}_\xi} 
\le_{\text{pr}} \underset\tilde {}\to q$".
We say $\underset\tilde {}\to q$ is very pure 
if $\Vdash_{\bar{\Bbb Q}}$ ``for 
$\xi < \alpha$, if $\underset\tilde {}\to q \in 
{\underset\tilde {}\to {\Bbb Q}_\xi}$, then ${\underset\tilde {}\to
{\Bbb Q}_\xi} \models \emptyset_{\underset\tilde {}\to {\Bbb Q}_i} 
\le_{\text{vpr}} \underset\tilde {}\to q$".
\sn
\item "{$(D),(E)$}"  Left to the reader.
\sn
\item "{$(F)$}"  We are defining 
$\text{Sp}_e(W)$-Lim$_\kappa(\bar{\Bbb Q})$ (where
$\bar{\Bbb Q} = \langle \Bbb P_\beta,
{\underset\tilde {}\to {\Bbb Q}_\beta}:\beta < \alpha
\rangle$, of course).  It is a quadruple $\Bbb P_\alpha = (\Bbb P_\alpha,\le,
\le_{\text{pr}},\le_{\text{vpr}})$ where
\sn
{\roster
\itemitem{ $(a)$ }  $\Bbb P_\alpha$ is the set of $p = \{
{\underset\tilde {}\to q_i}:i < i^*\}$ satisfying for some witness 
${\underset\tilde {}\to {\bar \zeta}}$:
\newline

$\qquad \quad (i) \quad$ each 
${\underset\tilde {}\to q_i}$ is a simple $(\bar{\Bbb Q},W)$-named
atomic condition,  \newline

$\qquad \qquad \quad$ and for every $\xi < \alpha$, we have \nl

$\qquad \qquad \qquad \Vdash_{{\Bbb P}_\xi} ``p \restriction
\{\xi\} =: \{{\underset\tilde {}\to q_i} 
\restriction \{\xi\}:i < i^*\} \cup 
\{\emptyset_{\underset\tilde {}\to {\Bbb Q}_\xi}\} \in 
{\underset\tilde {}\to {\hat{\Bbb Q}}_\xi}"$ \newline

$\qquad \quad (ii) \quad$ if $\alpha \in W$ is strongly inaccessible
$> CC(\Bbb P_i) + \kappa$ \nl  \newline

$\qquad \qquad \qquad$ for every $i < \alpha$, \ub{then} 
$i^* < \alpha$ \newline

$\qquad \quad (iii) \quad {\underset\tilde {}\to {\bar \zeta}} = \langle
{\underset\tilde {}\to \zeta_\varepsilon}:\varepsilon < j \rangle$
where $j < \kappa$ and

$\qquad \qquad \qquad$ each ${\underset\tilde {}\to \zeta_\varepsilon}$ is
a simple $(\bar{\Bbb Q},W)$-named $[0,\alpha)$-ordinal, \nl

$\qquad \qquad \qquad$ [the reader should think of
$\{{\underset\tilde {}\to \zeta_\varepsilon}:\varepsilon < j\}$ \nl

$\qquad \qquad \qquad$ as the non-very-pure support of $p$] \nl

$\qquad \quad (iv) \quad$ for every $\xi < \alpha$ we have (we may 
replace \nl

$\qquad \qquad \qquad 
\Vdash_{\bar{\Bbb Q}}$ by $\Vdash_{{\Bbb P}_\xi}$ as we use simple
names) \nl

$\qquad \qquad \qquad \Vdash_{\bar{\Bbb Q}}$ ``if $(\forall \varepsilon < j)
({\underset\tilde {}\to \zeta_\varepsilon}
[\underset\tilde {}\to G \cap \Bbb P_\xi]$ 
is $\ne \xi)$ (for example is \newline

$\qquad \qquad \qquad$ not well defined) \ub{then}
$\emptyset_{\underset\tilde {}\to {\Bbb Q}_\xi} \le_{\text{vpr}} 
p \restriction \{\xi\}$ in ${\underset\tilde {}\to {\hat{\Bbb Q}}_\xi}$"
\nl

$\qquad \quad (v)$ if $\beta < \alpha$ then $p \restriction \beta =:
\{q_i \restriction \beta:i < i^*\}$ belongs to $\Bbb P_\beta$
\nl

$\qquad \quad (vi)$ if $\alpha \in W$ is strongly inaccessible 
$> CC(\Bbb P_i) + \kappa$ \nl

$\qquad \qquad \qquad$ for every $i < \alpha$ \ub{then} for some $\beta <
\alpha$ \nl

$\qquad \qquad \qquad$ every ${\underset\tilde {}\to \zeta_i}$ is a simple
$(\bar{\Bbb Q},W)$-named $[0,\beta)$-ordinal; \nl

$\qquad \qquad \qquad$ needed, e.g., 
in \scite{si6.11} (note: this demand follows by \scite{1.6A}) \nl

$\qquad \quad (vii)$ 
\sn
\itemitem{ $(b)$ }  for $p \in 
\text{ Sp}_e(W)$-Lim$_\kappa(\bar{\Bbb Q})$ and $\xi < \ell g
(\bar{\Bbb Q})$ we let:

$$
p \restriction \xi =: \{r \restriction \xi:r \in p\}
$$

$$
p \restriction \{\xi\} =: \{r \restriction \{\xi\}:r \in p\}
$$
\newline

\noindent
$\quad$ we define similarly $p \restriction [\zeta,\xi),p \restriction \{
\underset\tilde {}\to \zeta\}, p \restriction [\underset\tilde {}\to \zeta,
\underset\tilde {}\to \xi)$.
\sn
\itemitem{ $(c)$ }  $\Bbb P_\alpha 
\models ``p^1 \le_{\text{vpr}} p^2"$ iff for every
$\xi < \alpha$ we have (letting \nl

$\quad p^\ell = \{ q^\ell_i:i < i^\ell(*)\}$ for $\ell = 1,2$): 
$$
p^2 \restriction \xi \Vdash_{{\Bbb P}_\xi}
``{\underset\tilde {}\to {\hat{\Bbb Q}}_\xi} 
\models p^1 \restriction \{ \xi \} 
\le_{\text{vpr}} p^2 \restriction \{ \xi \}"
$$
\sn
\itemitem{ $(d)$ }  $\Bbb P_\alpha 
\models ``p^1 \le_{\text{pr}} p^2"$ \underbar{iff}
\newline
\smallskip

$\qquad \quad (i) \quad$ for every 
$\xi < \ell g(\bar{\Bbb Q})$, we have \footnote{recall
\scite{1.1}(4)(d) we can omit ``$p^1 \restriction \{xi\} \ne \emptyset
\Rightarrow$"} \newline

$\qquad \qquad \qquad
p^2 \restriction \xi \Vdash_{{\Bbb P}_\xi}$ ``then $\hat{\Bbb Q}_\xi
\models p^1 \restriction \{\xi\} \le_{\text{pr}}
p^2 \restriction \{\xi\}$" \nl
\smallskip

$\qquad \quad (ii) \quad$ for some ordinal $j < \kappa$ and simple
$(\bar{\Bbb Q},W)$-named \nl
 
$\qquad \qquad \qquad [0,\alpha)$-ordinals 
${\underset\tilde {}\to \zeta_\varepsilon}$ for 
$\varepsilon < j$, for every $\xi < \ell g(\bar{\Bbb Q})$ we have: \newline

$\qquad \qquad \qquad
p^2 \restriction \xi 
\Vdash_{{\Bbb P}_\xi}$ ``if for no $\varepsilon < j$ do we have
${\underset\tilde {}\to \zeta_\varepsilon}
[{\underset\tilde {}\to G_{{\Bbb P}_\xi}}] = \xi$, then \newline

$\qquad \qquad \qquad
p^1 \restriction \{\xi\} \le_{\text{vpr}} p^2 \restriction \{\xi\}$ in
${\underset\tilde {}\to {\hat Q}_\xi}$"; we call
$\langle {\underset\tilde {}\to \zeta_\varepsilon}:\varepsilon < j \rangle$
a witness.
\sn
\itemitem{ $(e)$ }  $\Bbb P_\alpha \models p^1 \le p^2$ iff \newline
\smallskip

$\qquad \quad (i) \quad$ for every $\xi < \ell g(\bar{\Bbb Q})$ we have:
$$
p^2 \Vdash_{{\Bbb P}_\xi} ``{\underset\tilde {}\to {\hat{\Bbb Q}}_\xi} 
\models p^1 \restriction \{\xi\} \le
p^2 \restriction \{\xi\}"
$$

$\qquad \quad (ii) \quad$ as in the definition of $\le_{pr}$; 
clause (ii) \newline
\smallskip

$\qquad \quad (iii) \quad$ for some 
$n < \omega$ and simple $(\bar{\Bbb Q},W)$-named ordinals 
${\underset\tilde {}\to \xi_1},\dotsc,{\underset\tilde {}\to \xi_n}$
\nl 

$\qquad \qquad \qquad$ we have: \nl

$\qquad \qquad \qquad$ for each 
$\xi < \ell g(\bar{\Bbb Q})$ we have \nl

$\qquad \qquad \qquad \,\,p_2 \restriction \xi \Vdash_{{\Bbb P}_\xi}$ ``if
$\xi \ne {\underset\tilde {}\to \xi_\ell}
[{\underset\tilde {}\to G_{{\Bbb P}_\xi}}]$ for \nl

$\qquad \qquad \qquad \qquad \qquad \qquad \quad \ell = 1,\dotsc,n$ and \nl

$\qquad \qquad \qquad \qquad \qquad \qquad \quad 
\neg(\emptyset_{\underset\tilde {}\to {\Bbb Q}_\xi} 
\le_{\text{vp}} p^1 \restriction \{\xi\})]$  \nl

$\qquad \qquad \qquad \qquad \qquad \qquad \quad$ then: 
${\underset\tilde {}\to {\hat{\Bbb Q}}_\xi} 
\models p^1 \restriction
\{\xi\} \le_{\text{pr}} p^2 \restriction \{\xi\}"$ \nl

$\qquad \quad$ note that the truth value of $\zeta = 
{\underset\tilde {}\to \xi_\ell}$ is a $\Bbb P_\zeta$-name so this is \nl

$\qquad \quad$ well defined.  \nl

$\qquad \quad$ We then (i.e. if $(i) + (ii) + (iii)$) say: $p_1 \le p_2$ over
$\{{\underset\tilde {}\to \xi_1},\dotsc,{\underset\tilde {}\to \xi_n}\}$.
\sn
\itemitem{ $(f)$ }  Lastly, for $p \in \Bbb P_\alpha$ we let
Dom$(p) = \text{ Dom}_{\text{vp}}(p) =
\{{\underset\tilde {}\to \zeta_{\underset\tilde {}\to q}}:
\underset\tilde {}\to q \in p\}$ and Dom$_{\text{pr}}(p) = 
\{{\underset\tilde {}\to \zeta_\varepsilon}:\varepsilon < j\}$ where
$\bar \zeta = \langle {\underset\tilde {}\to \zeta_\varepsilon}:\varepsilon
< j \rangle$ is as in clause (F)(a) above (we can make it part of $p$).
\sn
\itemitem{ $(g)$ }  We still have to 
define $\Bbb P_{\underset\tilde {}\to \beta}$
for $\underset\tilde {}\to \beta$ a full simple 
$(\bar{\Bbb Q},W)$-named ordinal, it
is $\{p:p = \{{\underset\tilde {}\to q_i}:i < i^*\}$ and 
$\Vdash_{\bar{\Bbb Q}} 
``\zeta_{\underset\tilde {}\to q_i} < \underset\tilde {}\to \beta"$
that is if $\xi < \alpha$ then $\Vdash_{{\Bbb P}_\xi}$ ``if
$\zeta_{\underset\tilde {}\to q_i}[G_{{\Bbb P}_\xi}] = \xi$ and
$\underset\tilde {}\to \beta[G_{{\Bbb P}_\xi}]$ is well defined then
it is $> \xi"\}$.
\sn
\itemitem{ $(h)$ }  We call $p \in \Bbb P_\alpha$ full if it has a
witness $\langle {\underset\tilde {}\to \zeta_\varepsilon}:\varepsilon
< j \rangle$ with each ${\underset\tilde {}\to \zeta_\varepsilon}$ full.
\endroster}
\sn
\item "{$(G)$}"   Let us check Definition \scite{1.1}(1) for
$\Bbb P_\alpha =: \text{ Sp}_e(W)$-Lim$_\kappa(\bar{\Bbb Q})$: \nl
\ub{Proof of} $\le^{{\Bbb P}_\alpha}$ is a partial order. \nl
Suppose $p_0 \le p_1 \le p_2$.  Let $n^\ell,
{\underset\tilde {}\to \xi^\ell_1},\dotsc,
{\underset\tilde {}\to \xi^\ell_{n^\ell}}$ and $j^\ell \,(< \kappa)
\,\,{\underset\tilde {}\to \zeta^\ell_\varepsilon}$ 
(for $\varepsilon < j^\ell$) appear in the
definition of $p_\ell \le p_{\ell + 1}$.  Let $n = n^0 + n^1$, and

$$
{\underset\tilde {}\to \xi_i} = \cases
{\underset\tilde {}\to \xi^0_i} \quad &\text{if} \quad 1 \le i \le n^1 \\
{\underset\tilde {}\to \xi^1_{i-n}} &\text{if} \quad n^1 < i \le n^1 + n^2.
\endcases
$$

Let $j = j^0 + j^1$ and

$$
{\underset\tilde {}\to \zeta_\varepsilon} = \cases
{\underset\tilde {}\to \zeta^0_\varepsilon} \quad &\text{\underbar{if}} \quad 
\varepsilon < j^0 \\
{\underset\tilde {}\to \zeta^1_{\varepsilon-j^0}} &\text{\underbar{if}} \quad 
\varepsilon \in [j^0,j^0 + j^1).
\endcases
$$
\medskip

\noindent
Let us check the three clauses of (e) of part (D).
\bn
\ub{Clause (i)}: \nl
Let $\xi < \ell g(\bar{\Bbb Q})$ so for $\ell = 0,1$

$$
p_{\ell +1} \restriction \xi \Vdash_{{\Bbb P}_\xi} ``p_\ell \restriction \{\xi\} 
\le p_{\ell +1} \restriction \{\xi\} \text{ in } \hat{\Bbb Q}_\xi".
$$
\mn
As $\Bbb P_\xi \models 
``p_1 \restriction \xi \le p_2 \restriction \xi"$ (by the
induction hypothesis, clause (H)) clearly
$p_2 \restriction \xi$ forces both assertions.  As
${\underset\tilde {}\to {\hat{\Bbb Q}}_\xi}$ is a 
partial order (under $\le$) the conclusion follows.
\bn
\ub{Clause (ii)}: \nl
Let $\xi < \ell g(\bar{\Bbb Q})$, so similarly 
$p^2 \restriction \xi \Vdash_{{\Bbb P}_\xi}$ ``if $\xi \ne
{\underset\tilde {}\to \xi^\ell_m}$ for $m = 1,\dotsc,n^\ell$ 
and $\ell = 0,1$ (i.e., $\xi^\ell_m[{\underset\tilde {}\to G_{{\Bbb P}_\xi}}]$ 
is $\ne \xi$ or is not well defined), \ub{then} $p_0 \restriction
\{\xi\} \le_{\text{pr}} p_1 \restriction 
\{\xi\}$ in ${\underset\tilde {}\to {\hat{\Bbb Q}}_\varepsilon}$ and
$p_1 \restriction \{\xi\} \le_{\text{pr}} p_2 \restriction \{\xi\}$ in
${\underset\tilde {}\to {\hat{\Bbb Q}}_\xi}"$ from which the result
follows.
\bn
\ub{Clause (iii)}: \nl
Lastly, for $\xi < \alpha$ we have $p_2 \restriction \xi
\Vdash_{{\Bbb P}_\xi}$ ``if $\xi \notin 
\{{\underset\tilde {}\to \zeta^\ell_\varepsilon}[G_{{\Bbb P}_\zeta}]:
{\underset\tilde {}\to \zeta^\ell_\varepsilon}[G_{{\Bbb P}_\zeta}]$ well defined,
$\varepsilon < j^\ell$ and $\ell \in \{0,1\}\}$ then 
$p_0 \restriction \{ \xi \} \le_{\text{vpr}} p_1
\restriction \{ \xi\} \le_{\text{vpr}} p_2 \restriction \{ \xi\}$ hence $p_0
\restriction \{\xi\} \le_{\text{vpr}} p_2 \restriction \{\xi\}$". 
\nl
\mn
So we have proved the three conditions needed for $p_0 \le p_2$ by 
the definition above so 
really $p_0 \le p_2$ holds, so $\le$ is a partial order. \nl
\sn
\ub{Proof of} $\le_{\text{pr}}$ is a partial order. \nl

Similar proof.
\sn
\ub{Proof of} $\le_{\text{vpr}}$ is a partial order.  \nl

Similar proof just easier. \nl
\sn
\ub{Proof of} $p \le_{\text{pr}} q \Rightarrow p \le q$: 

By the definition; easy. \nl
\sn
\ub{Proof of} $p \le_{\text{vpr}} q \Rightarrow p \le_{\text{pr}} q$:

By the definition, check. \nl
So in \scite{1.1}(1) all the requirements on $\Bbb P_\alpha$ holds.
\sn
\item "{$(H),(I),(J)$}"  We leave the checking to the reader 
(actually we prove $(I)$ in the proof of \scite{1.17} below). 
\hfill$\square_{\scite{1.13}}$\margincite{1.13}
\endroster
\bigskip

\demo{\stag{1.14} Fact}  1) If $\bar{\Bbb Q}$ is a 
$\kappa-\text{Sp}_e$-iteration and for each $i < \ell g(\bar{\Bbb Q})$ 
we have it is forced (i.e., $\Vdash_{{\Bbb P}_i}$) that
$\le^{{\Bbb Q}_i}_{\text{pr}} = \le^{{\Bbb Q}_i}$ and 
$\le^{{\Bbb Q}_i}_{\text{vpr}}$ is 
equality, \ub{then} $\bar{\Bbb Q}$ is a variant of $\kappa$-RS iteration 
(as in \cite[Ch.XIV,\S1]{Sh:f}), i.e. they are the same if we use there 
only simple $\bar{\Bbb Q}$-named ordinals (or allow here non-simple 
ones so the version here is exactly as in \cite[Ch.XIV,2.6]{Sh:f}). 
\enddemo
\bigskip

\demo{Proof}  Straightforward.
\enddemo
\bigskip

\proclaim{\stag{1.15} Claim}  1) In \scite{1.13} in Definition (F),
clause (a)(iii) we can demand that each ${\underset\tilde {}\to
\zeta_\varepsilon}$ is full (simple $(\bar{\Bbb Q},W)$-named
$[0,\alpha)$-ordinal) and similarly in (d)(ii), (=(e)(ii)) and
(e)(iii). \nl
2) Suppose $\bar{\Bbb Q} = \langle \Bbb P_i,
{\underset\tilde {}\to {\Bbb Q}_i}:i < \alpha 
\rangle$ is a $\kappa-{\text{\rm Sp\/}}_e(W)$-iteration
(so $\Bbb P_\alpha = { \text{\rm Sp\/}}_e(W)$-{\text{\rm
Lim\/}}$_\kappa(\bar{\Bbb Q})$).  If $p \le q$ in 
$\Bbb P_\alpha$, \ub{then} there are $r,n < \omega$ and 
$\xi_1 < \ldots < \xi_n < \alpha$ such that:
\medskip
\roster
\item "{$(a)$}"  $r \in \Bbb P_\alpha$
\sn
\item "{$(b)$}"  $q \le r$
\sn
\item "{$(c)$}"  $p \le r$ above $\{\xi_1,\dotsc,\xi_n\}$.
\ermn
3)  If $p^1,p^2 \in \Bbb P_\alpha$ and \scite{1.13}(F)(e) holds but we 
allow $\underset\tilde {}\to n$ to be a full simple $(\bar{\Bbb
Q},W)$-named ordinal, \ub{then} $p^2 \Vdash_{{\Bbb P}_\alpha} ``p_1
\in {\underset\tilde {}\to G_{{\Bbb P}_\alpha}}"$. 
\endproclaim
\bigskip

\remark{Remark}  In fact, in \scite{1.15}(1) we can have 
$r \restriction [\xi_n,\alpha) = q \restriction [\xi_n,\alpha)$.
\endremark
\bigskip

\demo{Proof}  1) As increasing those sets $(\{{\underset\tilde {}\to
\zeta_\varepsilon}:\varepsilon < \zeta\},\{{\underset\tilde {}\to
\xi_1},\dotsc,{\underset\tilde {}\to \xi_n}\}$, respectively) cause no
harm. \nl
2) We prove this by induction on $\alpha$.
\enddemo
\bn
\underbar{Case 1}:  $\alpha = 0$. \newline
Trivial.
\bn
\underbar{Case 2}:  $\alpha = \beta + 1$.

Apply the induction 
hypothesis to $\bar{\Bbb Q} \restriction \beta,p \restriction
\beta,q \restriction \beta$ (clearly $\bar{\Bbb Q} \restriction 
\beta$ is an \nl
$\kappa-\text{Sp}_e(W)$-iteration, 
$p \restriction \beta \in \Bbb P_\beta,q \restriction 
\beta \in \Bbb P_\beta$ and $\Bbb P_\beta \models ``p \le q"$, 
by \scite{1.13}, clause (H)).  So we can find $r',m < \omega$ 
and $\{\xi'_1,\dotsc,\xi'_m\}$ such that:
\medskip
\roster
\item "{$(a)'$}"  $r' \in \Bbb P_\beta$
\sn
\item "{$(b)'$}"  $\Bbb P_\beta \models q \restriction \beta \le r'$
\sn
\item "{$(c)'$}"  $p \restriction \beta \le r'$ (in $\Bbb P_\beta$) above
$\{\xi_1,\dotsc,\xi_m\}$
\sn
\item "{$(d)'$}"  $\xi'_1 < \ldots < \xi'_m$.
\endroster
\medskip

Let $n =: m+1$ and

$$
\xi_\ell = \cases \xi'_\ell \quad &\text{if} \quad \ell \in \{1,\dotsc,m\} \\
\beta \quad &\text{if} \quad \ell = n
\endcases
$$
\mn
and lastly $r = r' \cup (q \restriction \{\beta\})$.
\bn
\underbar{Case 3}:  $\alpha$ is a limit ordinal. \nl
Let $p \le q$ (in $\Bbb P_\alpha$) above 
$\{{\underset\tilde {}\to \xi_1},\dotsc,
{\underset\tilde {}\to \xi_n}\}$.  We choose by 
induction on $\ell \le n$, the objects
$r_\ell,\beta_\ell,\xi^*_\ell$ such that:
\mr
\item "{$(\alpha)$}"  $r_\ell \in \Bbb P_{\beta_\ell}$
\sn
\item "{$(\beta)$}"  $r_\ell \le r_{\ell +1}$
\sn
\item "{$(\gamma)$}"  $q \restriction \beta_\ell \le r_\ell$
\sn
\item "{$(\delta)$}"  $\beta_\ell \le \beta_{\ell +1} < \alpha$
\sn
\item "{$(\varepsilon)$}"  $\beta_0 = 0,r_0 = \emptyset_{{\Bbb P}_0}$
\sn
\item "{$(\zeta)$}"  for $\ell \in \{1,\dotsc,n\}$ we have: either
$r_\ell \Vdash_{{\Bbb P}_{\beta_\ell}} ``{\underset\tilde {}\to \xi_\ell} = 
\xi^*_\ell"$ and $\xi^*_\ell \le \beta_\ell$ or \nl
$\beta_\ell = \beta_{\ell-1} \and r_\ell = 
r_{\ell-1}$ and $r_\ell \cup (q \restriction [\beta,\alpha))
\Vdash_{{\Bbb P}_\alpha}
``{\underset\tilde {}\to \xi_\ell}$ is not defined or is $\ge \alpha"$.
\endroster
\medskip

Carrying the definition is straight: for $i=0$ use clause $(\varepsilon)$.
For $\ell + 1 \le n$ when the second possibility of clause $(\zeta)$ fails
there is $r'$, such that $r_\ell \cup (q \restriction [\beta,\alpha)) \le 
r' \in 
\Bbb P_\alpha$, and $r' \Vdash_{{\Bbb P}_\alpha} ``{\underset\tilde {}\to 
\xi_{\ell+1}}$ is defined and is $< \alpha$", so there are 
$r'',\xi^*_{\ell +1} < \alpha$ such that
$r' \le r'' \in \Bbb P_\alpha$ and $r'' \Vdash ``{\underset\tilde {}\to 
\xi_{\ell+1}} = \xi^*_{\ell +1}$" so as ``$\xi^*_{\ell+1}$ is a simple
$\bar{\Bbb Q}$-named ordinal" 
we know that $\xi^*_{\ell +1} < \alpha$ and $r''
\restriction \xi^*_{\ell+1} \Vdash_{{\Bbb P}_{\xi^*_{\ell+1}}} ``
{\underset\tilde {}\to \xi_{\ell+1}} = \xi^*_{\ell +1}$".  Let
$\beta_{\ell+1} =: \text{ max}\{\beta_\ell,\xi^*_{\ell+1}\}$, and
$r_{\ell +1} =: r'' \restriction \beta_{\ell+1}$.  So we have carried the
induction.

We now apply 
the induction hypothesis to $\bar{\Bbb Q} \restriction \beta_n,p \restriction
\beta_n,r_n$; it is applicable as $\beta_n < \alpha$, and $\Bbb P_{\beta_n}
\models ``p \restriction \beta_n \le q \restriction \beta_n \le r_n"$.  So
there are $m < \omega,\xi_1 < \ldots < \xi_m < \beta_n$ and $r^*$ such that
$\Bbb P_{\beta_n} \models ``r_n \le r^*"$ and 
$p \le r^*$ (in $\Bbb P_{\beta_n}$)
above $\{\xi_1,\dotsc,\xi_m\}$.  Now let $r =: r^* \cup (q \restriction
[\beta_n,\alpha))$, clearly $q \le r$ and $p \le r$ above $\{\xi_1,\dotsc,
\xi_m,\beta_n\}$. \nl
3) So $p^1,p^2,\{\xi_\ell:\ell =1,\dotsc,n\}$ are given.  Assume that
$p^2 \le q_1 \in \Bbb P_\alpha$.  We can find $q_2$ such that $q_1 \le
_2 \in \Bbb P_\alpha$ and $q_2$ forces a value to $\underset\tilde
{}\to n$ say $n(*)$.  Next we can find $q_3$ such that $q_2 \le q_3
\in \Bbb P_\alpha$ and $q_3$ forces values to ${\underset\tilde {}\to
\xi_\ell} \, (\ell =1,\dotsc,n(*)\}$, say
$\xi^0_1,\dotsc,\xi^0_{n(*)}$.  Now repeating the proof of
``$\le^{{\Bbb P}_\alpha}$ is a partial order" (in clause (H) of
\scite{1.13}) with $p^1,p^2,q_3$ here standing for $p_0,p_1,p_2$ there
we let $\Bbb P_\alpha \models p^1 \le q_3$.  As we have assumed just
$p^2 \le q_1 \in \Bbb P_\alpha$ and $q_1 \le q_3$ we are done.
\hfill$\square_{\scite{1.15}}$\margincite{1.15}  
\bigskip

\citewarning{\noindent \llap{---$\!\!>$} MARTIN WARNS: Label 1.16 on next line is also used somewhere else (Perhaps should have used scite instead of stag?)}
\proclaim{\stag{1.16} Claim}  Let 
$\bar{\Bbb Q}$ be a $\kappa-{\text{\rm Sp\/}}_e(W)$-iteration
of length $\alpha$. \nl
0) If $\underset\tilde {}\to \zeta$ 
is a simple $(\bar{\Bbb Q},W)$-named ordinal
\ub{then} for some ordinal $\gamma$ we have:
$\underset\tilde {}\to \zeta$ is a simple $(\bar{\Bbb Q},W)$-named
$[0,\gamma)$-ordinal. \newline
1) 
\mr
\item "{$(i)$}"  If 
$\beta < \alpha$ and $\underset\tilde {}\to \zeta$ is a $\Bbb P_\beta$-name
of a [full] simple 
$(\bar{\Bbb Q},W)$-named $[\beta,\alpha)$-ordinal \ub{then} for some [full]
simple $(\bar{\Bbb Q},W)$-named $[\beta,\alpha)$-ordinal 
$\underset\tilde {}\to \xi$
we have
$$
\Vdash_{\bar{\Bbb Q}} 
``\underset\tilde {}\to \zeta = \underset\tilde {}\to \xi"
$$
\sn
\item "{$(ii)$}"  if $\beta \le \alpha,\beta \le \gamma_1 \le \gamma_2$
and $\underset\tilde {}\to \zeta$ is a $\Bbb P_\beta$-name of a [full]
simple $(\bar{\Bbb Q},W)$-named $[\gamma_1,\gamma_2)$-ordinal
\ub{then} for some [full] simple
$(\bar{\Bbb Q},W)$-named 
$[\gamma_1,\gamma_2)$-ordinal $\underset\tilde {}\to \xi$
we have $\Vdash_{\bar{\Bbb Q}} ``\underset\tilde {}\to \zeta =
\underset\tilde {}\to \xi"$.
\ermn
2) The same holds if we replace ``ordinal" by ``atomic condition"
(so in (ii) we should demand $\gamma_2 \le \alpha$). \newline
3) If $\alpha \le \gamma$ and $\underset\tilde {}\to \beta$ is a full simple 
$(\bar{\Bbb Q},W)$-named $[0,\alpha)$-ordinal, and for each 
$\beta < \alpha,{\underset\tilde {}\to \zeta_\beta}$ is a [full] simple
$(\bar{\Bbb Q},W)$-named $[\beta,\gamma)$-ordinal \ub{then} for some
[full] simple $(\bar{\Bbb Q},W)$-named $[0,\gamma)$-ordinal 
$\underset\tilde {}\to \xi$ we have

$$
\Vdash_{\bar{\Bbb Q}} ``\text{if } \underset\tilde {}\to \beta
[\underset\tilde {}\to G] = \beta (\text{so } \beta < \alpha) \text{ then }
\underset\tilde {}\to \xi[\underset\tilde {}\to G] = \zeta_\beta
[\underset\tilde {}\to G]".
$$
\mn
4) If $\underset\tilde {}\to \beta$ is a full
simple $(\bar{\Bbb Q},W)$-named $[0,\alpha)$-ordinal and 
for each $\beta < \alpha,{\underset\tilde {}\to p_\beta}$ is a 
$(\bar{\Bbb Q},W)$-named $[\beta,\alpha)$-atomic 
condition \underbar{then} for 
some $(\bar{\Bbb Q},W)$-named atomic condition $p$ we have

$$
\Vdash_{\bar{\Bbb Q}} ``\text{if } \underset\tilde {}\to \beta
[\underset\tilde {}\to G] = \beta \text{ then }
\underset\tilde {}\to p[\underset\tilde {}\to G] = 
{\underset\tilde {}\to p_\beta}[\underset\tilde {}\to G]".
$$
\endproclaim
\bigskip

\demo{Proof}  Easy.  \hfill$\square_{\scitet{1.16}}$\scitetphantom{1.16}\margincite{1.16}
\enddemo
\bigskip

\proclaim{\stag{1.16A} Claim}  1) Suppose $F$ is a function, \ub{then} 
for every ordinal $\alpha$ there is \nl
$\kappa-{\text{\rm Sp\/}}_e(W)$-iteration $\bar{\Bbb Q} = \langle \Bbb P_i,
{\underset\tilde {}\to {\Bbb Q}_i}:i < \alpha^\dag \rangle$, such that:
\mr
\item "{$(a)$}"  for every $i$ we have ${\underset\tilde {}\to {\Bbb Q}_i} = 
F(\bar{\Bbb Q} \restriction i)$,
\sn
\item "{$(b)$}"  $\alpha^\dag \le \alpha$
\sn
\item "{$(c)$}"  either $\alpha^\dag = \alpha$ or $F(\bar{\Bbb Q})$ is not an
${\text{\rm Sp\/}}_e(W)-{\text{\rm Lim\/}}_\kappa
(\bar{\Bbb Q})$-name of a forcing notion.
\ermn
2) Suppose $\beta < \alpha,G_\beta \subseteq \Bbb P_\beta$ is generic
over $\bold V$, \ub{then} in $\bold V[G_\beta],\bar{\Bbb Q}/G_\beta 
= \langle \Bbb P_i/G_\beta,{\underset\tilde {}\to {\Bbb Q}_i}:\beta 
\le i < \alpha \rangle$ is an $\kappa-{\text{\rm Sp\/}}_e$-iteration 
and $\kappa-{\text{\rm Sp\/}}_e(W)-{\text{\rm Lim\/}}(\bar{\Bbb Q}) = 
\Bbb P_\beta * ({\text{\rm Lim\/}}(\bar{\Bbb Q})/
{\underset\tilde {}\to G_\beta}$) (essentially; more exactly
up to equivalence) where, of course, $\Bbb P_i/G_\beta = \{p \in \Bbb P_i:p
\restriction \beta \in G_\beta\}$. \newline
3) If $\bar{\Bbb Q}$ is an $\kappa-{\text{\rm Sp\/}}_e(W)$-iteration, 
$p \in { \text{\rm Sp\/}}_e(W)-{\text{\rm Lim\/}}_\kappa
(\bar{\Bbb Q}),\Bbb P'_i = \{q \in \Bbb P_i:q \ge p 
\restriction i\}$, \nl
${\underset\tilde {}\to {\Bbb Q}'_i} = \{ p \in
{\underset\tilde {}\to {\Bbb Q}_i}:p \ge p \restriction \{i\}\}$ and
$\emptyset_{\underset\tilde {}\to {\Bbb Q}'_i} = p \restriction \{i\},
\le^{\underset\tilde {}\to {\Bbb Q}'_i}_\ell =
\le^{\underset\tilde {}\to {\Bbb Q}_i}_\ell 
\restriction {\underset\tilde {}\to {\Bbb Q}'_i}$ 
for $\ell = 0,1,2$ \ub{then} 
$\bar{\Bbb Q} = \langle \Bbb P'_i,\Bbb Q'_i:
i < \ell g(\bar{\Bbb Q}) \rangle$ is
(essentially) a $\kappa-{\text{\rm Sp\/}}_e(W)$-iteration 
(and ${\text{\rm Sp\/}}_e(W)-{\text{\rm Lim\/}}_\kappa(\bar{\Bbb Q}')$ is 
$P'_{\ell g(\bar{\Bbb Q})}$).
\endproclaim
\bigskip

\demo{Proof}  Should be clear.
\enddemo
\bigskip

\proclaim{\stag{1.17} Claim}  Suppose 
\mr
\item "{$(a)$}"  $\bar{\Bbb Q} = \langle \Bbb P_i,{\underset\tilde
{}\to {\Bbb Q}_i}:i < \alpha
\rangle$ is a $\kappa-{\text{\rm Sp\/}}_e(W)$-iteration 
(and $\Bbb P_\alpha = { \text{\rm Sp\/}}_e(W)-{\text{\rm
Lim\/}}_\kappa(\bar{\Bbb Q}))$
\sn
\item "{$(b)$}"  $\ell(*) \in \{0,1\}$
\sn
\item "{$(c)$}"  $\Vdash_{{\Bbb P}_i} ``({\underset\tilde {}\to {\Bbb
Q}_i},\le_{\ell(*)})$ is a $\theta$-complete" for each $i < \alpha$. 
\ermn
\ub{Then}: \nl
1) $(\Bbb P_\alpha,\le_{\ell(*)})$ is $\theta$-complete, i.e. 
if $\delta < \theta,\langle p_i:i < \delta \rangle$ is 
$\le_{\ell(*)}$-increasing \ub{then} it has an 
$\le_{\ell(*)}$-upper bound provided that:
\nl
$\theta \le \kappa$ \ub{or} $\ell(*) = 0 \and 
\theta \le { \text{\rm Min\/}}\{\delta:
\delta \in W \text{ is strongly inaccessible and } (\forall \beta < \delta)
(|\Bbb P_\beta| < \delta)\}$. \newline
2) Moreover for $\beta < \alpha$ we have 
$(\Bbb P_\alpha/\Bbb P_\beta,\le_{\ell(*)})$ is 
$\theta$-complete. \newline
3) In fact, we can get $\le_{\ell(*)}$-lub (provided that 
there are such lub's for each ${\underset\tilde {}\to {\Bbb Q}_i}$.
\endproclaim
\bigskip

\remark{Remark}   We deal with $\theta$-complete rather than strategically
$\theta$-complete (here and later) just for simplicity presentation, as it
does not matter much by \cite[CH.XIV,2.4]{Sh:f}.
\endremark
\bigskip

\demo{Proof}  Straightforward but we elaborate. \newline
1) So assume $\delta < \theta$ and $p_i \in \Bbb P_\alpha$ for $i < 
\delta$ and $[i < j < \delta \Rightarrow p_i \le_{\ell(*)} p_j]$.  
Now it is enough to find $p \in \Bbb P_\alpha$ such that

$$
i < \delta \Rightarrow p_i \le_{\ell(*)} p.
$$
\mn
Let $p_i = \{{\underset\tilde {}\to q^i_\gamma}:\gamma < \gamma_i\}$ 
and for each $\gamma < \gamma_i,{\underset\tilde {}\to q^i_\gamma}$ is 
a simple $(\bar{\Bbb Q},W)$-named atomic condition, say
$\Vdash_{\bar{\Bbb Q}} ``{\underset\tilde {}\to q^i_\gamma} \in 
{\underset\tilde {}\to {\Bbb Q}_{\underset\tilde {}\to
\zeta^i_\gamma}}"$, where 
${\underset\tilde {}\to \zeta^i_\gamma}$ is a simple
$(\bar{\Bbb Q},W)$-named ordinal 
(which is 
${\underset\tilde {}\to \zeta_{\underset\tilde {}\to q^i_\gamma}}$).
Now for each $\beta < \alpha$ let ${\underset\tilde {}\to <^*_\beta}$ 
be a $\Bbb P_\beta$-name of a well ordering of
${\underset\tilde {}\to {\Bbb Q}_\beta}$.  For 
each $i(*) < \delta,\gamma(*) <
\gamma_i$ let ${\underset\tilde {}\to r^{i(*)}_{\gamma(*)}}$ be the following
simple $(\bar{\Bbb Q},W)$-named atomic condition:

Let $\zeta < \alpha,G_\zeta \subseteq \Bbb P_\zeta$ generic over
$\bold V$ and
${\underset\tilde {}\to \zeta^{i(*)}_{\gamma(*)}}[G_\zeta] = \zeta$, now work
in $\bold V[G_\zeta]$, let $w_\zeta = \{i < \delta:\text{for some } \gamma <
\gamma_i \text{ we have } {\underset\tilde {}\to \zeta^i_\gamma}[G_\zeta] = 
\zeta\}$.
We let $u^\zeta_i = \{\gamma < \gamma_i:
{\underset\tilde {}\to \zeta^i_\gamma}
[G_\zeta] = \zeta\}$ for each $i \in w_\zeta$.
(As $p_i$ is $\le_{\ell(*)}$-increasing, $w_\zeta$ is an end segment of 
$\delta$ and 
$i(*) \in w_\zeta,\gamma(*) \in u^\zeta_{i(*)}$).  For $i \in w_\zeta$ let
$q^*_{i,\zeta} = (p_i \restriction \{\zeta\})[G_\zeta]$.
Now define ${\underset\tilde {}\to  r^{i(*)}_{\gamma(*)}}[G_\zeta]$ as
follows.
\mn
\underbar{Case 1}:  For some $j < \delta$ the sequence 
$\langle q^*_{i,\zeta}:i \in w_\zeta \backslash j \rangle$ is constant.

Let ${\underset\tilde {}\to r^{i(*)}_{\gamma(*)}}[G_\zeta] =
q_{\text{Min}(w_\zeta \backslash j),i}$.
\mn
\underbar{Case 2}:  Not Case 1, but for some $\ell < 2,j < \delta$
the sequence $\langle q_{i,\zeta}:i \in w_\zeta \backslash j \rangle$ 
is $\le_\ell$-increasing,
without loss of generality $\ell$ minimal (on all possible $j$) and then
${\underset\tilde {}\to r^{i(*)}_{\gamma(*)}}[G_\zeta] \in
{\underset\tilde {}\to {\Bbb Q}_\zeta}[G_\zeta]$ is the
$<^*_\zeta$-first $\le_\ell$-upper bound of
$\{{\underset\tilde {}\to q^*_{i,\zeta}}:i \in w_\zeta \backslash j\}$ where
$<^*_\zeta = {\underset\tilde {}\to <^*_\zeta}[G_\zeta]$.
It exists by \scite{1.1A}(2).  
\mn
\ub{Case 3}:  Neither Case 1 nor Case 2 ${\underset\tilde {}\to
r^{i(*)}_{\gamma(*)}}[G_\zeta]$ is $\emptyset_{\underset\tilde {}\to
{\Bbb Q}_\zeta}$. \nl
Let $p = \{ {\underset\tilde {}\to  r^{i(*)}_{\gamma(*)}}:i(*) < \delta
\text{ and } \gamma < \gamma_{i(*)}\}$.

If $\ell(*) = 0$ (i.e., $p_i$ is $\le_{\text{vpr}}$-increasing) and 
$\langle {\underset\tilde {}\to \zeta_j}:j < j^* \rangle$
witness $p_0 \in \Bbb P_\alpha$, then it witnesses 
$p \in \Bbb P_\alpha$ and easily
$i < \delta \Rightarrow p_i \le_0 p$. \nl
So assume $\ell(*)=1$, that is $p_i$ is $\le_{\text{pr}}$-increasing.  
For $i_0 < i_2 < \delta$ let 
$\{ {\underset\tilde {}\to \zeta^{i_0,i_1}_j}:j < j_{\ell_0,\ell_1}\}$ be a
witness to $p_{i_0} \le_{pr} p_{i_1}$ and let 
$\{{\underset\tilde {}\to \zeta^i_j}:j < j_i\}$  witness $p_i \in \Bbb P$.  
Letting $\kappa^-$ be the maximal cardinal $< \kappa$, clearly
$\{{\underset\tilde {}\to \zeta^{i_0,i_1}_j}:i_0 < i_1 < \delta \text{ and }
j < j_{i_0,i_1}\}$ has cardinality $\kappa^-$, so we can order it as  
$\{ {\underset\tilde {}\to \zeta^{i_0(\varepsilon),i_1(\varepsilon)}
_{j(\varepsilon)}}:\varepsilon < \kappa^-\}$ and some
$\{ \zeta^{i(\varepsilon)}_{j'(\varepsilon)}:\varepsilon < \kappa^-\}$ list
$\{ {\underset\tilde {}\to \zeta^i_j}:i < \delta,j < j_i\}$.  Now
$\{ \zeta^{i(\varepsilon)}_{j'(\varepsilon)}:
\varepsilon < \kappa^-\}$ witness $p \in \Bbb P_\alpha$ and
$\{ \zeta^{i_0(\varepsilon),i_1(\varepsilon)}_{j(\varepsilon)}:\varepsilon
< \kappa^-\}$ witness $p_i \le_1 p_\delta$ for every $i < \delta$.
\mn
2), 3)  Similar proof. \hfill$\square_{\scite{1.17}}$\margincite{1.17}
\enddemo
\bigskip

\definition{\stag{1.18} Definition}  Let $\bar{\Bbb Q} = \langle \Bbb P_i,
{\underset\tilde {}\to {\Bbb Q}_i}:i < \alpha \rangle$ be an 
$\kappa-\text{Sp}_e(W)$-iteration. \newline
1) We say $\underset\tilde {}\to y$ is a 
$(\bar{\Bbb Q},W,\underset\tilde {}\to \zeta)$-name if: $\underset\tilde {}\to y$ is
a $\Bbb P_\alpha$-name, $\underset\tilde {}\to \zeta$ is a simple $\bar{\Bbb Q}$-named 
$[0,\alpha)$-ordinal, and: if $\beta < \alpha,G_{{\Bbb P}_\alpha} \subseteq
\Bbb P_\alpha$ is generic over $\bold V$ and for some 
$r \in G_{{\Bbb P}_\alpha} \cap \Bbb P_\beta$
we have $r \Vdash_{\bar{\Bbb Q}} ``\underset\tilde {}\to \zeta = \beta"$, then
$\underset\tilde {}\to y[G_{{\Bbb P}_\alpha}] \in \bold V[G_{{\Bbb
P}_\beta}]$ is well defined and
depends only on $G_{{\Bbb P}_\alpha} \cap \Bbb P_\beta$ so we write
$\underset\tilde {}\to y[G_{{\Bbb P}_\alpha} \cap \Bbb P_\beta]$; and if
$G_{{\Bbb P}_\alpha} \subseteq \Bbb P_\alpha$ is generic over $\bold V$ and
$\underset\tilde {}\to \zeta[G_{{\Bbb P}_\alpha}]$ not well defined then 
$\underset\tilde {}\to y[G_{{\Bbb P}_\alpha}]$
is not well defined (do not arise if $\underset\tilde {}\to \zeta$ is
full). \nl
2) If $p \in \Bbb P_\alpha$ and $G_{{\Bbb P}_\alpha} \subseteq 
\Bbb P_\alpha$ is generic over $\bold V$, (or
just in Gen$^r(\bar{\Bbb Q})$), \ub{then} 
$p[G_{{\Bbb P}_\alpha}]$ is the following 
function, Dom$(p[{\underset\tilde {}\to G_{{\Bbb P}_\alpha}}]) = \{
{\underset\tilde {}\to \zeta_{\underset\tilde {}\to q}}[G_{{\Bbb P}_\alpha}]:
\underset\tilde {}\to q \in p\}$ and $(p[G_{{\Bbb P}_\alpha}])(\varepsilon) =
\{\underset\tilde {}\to q[G_{{\Bbb P}_\alpha}]:\underset\tilde {}\to q
\in p$ and ${\underset\tilde {}\to \zeta_{\underset\tilde {}\to q}}
[G_{{\Bbb P}_\alpha}] = \varepsilon\}$. 
\enddefinition
\bigskip

\proclaim{\stag{1.19} Claim}  Suppose
\mr
\item "{$(a)$}"  $\bar{\Bbb Q} = \langle \Bbb P_i,
{\underset\tilde {}\to {\Bbb Q}_i}:i < \alpha
\rangle$ is a $\kappa-{\text{\rm Sp\/}}_e(W)$-iteration.
\sn
\item "{$(b)$}"  $p \in \Bbb P_\alpha$ and $\underset\tilde {}\to \zeta$ is a
simple $(\bar{\Bbb Q},W)$-named $[0,\alpha)$-ordinal
\sn
\item "{$(c)$}"  $\underset\tilde {}\to r$ is a 
$(\bar{\Bbb Q},W,\underset\tilde {}\to \zeta)$-named member of 
$\Bbb P_\alpha/\Bbb P_{\underset\tilde {}\to \zeta}$.
\ermn
\ub{Then}: \newline
1)  There is $q \in \Bbb P_\alpha$ satisfying $p \le q$ such that:
\mr
\item "{$(*)$}"  if $\xi < \alpha,G_\xi \subseteq \Bbb P_\xi$ generic over
$\bold V$, then
{\roster
\itemitem{ $(\alpha)$ }  $\underset\tilde {}\to \zeta[G_\xi] = \xi$ \nl
implies $(p \restriction \xi)[G_\xi] = (q \restriction \xi)[G_\xi]$
and
\sn
\itemitem{ $(\beta)$ }  $\underset\tilde {}\to \zeta[G_\xi] = \xi$
implies \nl
$(q \restriction [\xi,\alpha))[G_\xi] = \underset\tilde {}\to r[G_\xi]$.
\endroster}
\ermn
2) If in addition (for any $\ell < 3)$ clause (c)$^+$ below,
\ub{then} we can in $(*)$ add $p \le_\ell q$
\mr
\item "{$(c)^+$}"  $\underset\tilde {}\to r$ is a 
$(\bar{\Bbb Q},W,\underset\tilde {}\to r)$-named member of $\Bbb P_\alpha/
\Bbb P_{\underset\tilde {}\to \zeta}$ which is $\le_\ell$-above $p 
\restriction [\underset\tilde {}\to \zeta,\alpha)$.
\endroster
\endproclaim
\bigskip

\demo{Proof}  Straightforward.
\enddemo
\bn
Central here is pure decidability.
\definition{\stag{1.25} Definition}  1) A forcing notion $\Bbb Q$ has pure
$(\theta_1,\theta_2)$-decidability if: for every $p \in \Bbb Q$ and
$\Bbb Q$-name $\underset\tilde {}\to \gamma < \theta_1$, there are 
$a \subseteq \theta_1,|a| < \theta_2$ (but $|a| > 0$) and $r \in 
\Bbb Q$ such that $p \le_{\text{pr}} r$ and 
$r \Vdash_{\Bbb Q} ``\underset\tilde {}\to \gamma \in a"$ (for $\theta_1=2$,
alternatively, $\underset\tilde {}\to \gamma$ is a truth value).  If we write
$``\le \theta_2"$ we mean $|a| \le |\theta_2|$. \nl
2) A forcing notion $\Bbb Q$ has pure $\theta$-decidability 
where $\theta$ is an
ordinal \ub{if}: for every $p \in \Bbb Q$ and 
${\underset\tilde {}\to {\Bbb Q}}$-name
$\underset\tilde {}\to \gamma < \theta$ 
there are $\gamma < \theta$ and $r \in
\Bbb Q$ such that $p \le_{\text{pr}} r$ and 
$r \Vdash_{\Bbb Q}$ ``if $\theta < \omega$
then $\underset\tilde {}\to \gamma = \gamma$ \ub{and} if 
$\theta \ge \omega$ is
a limit ordinal then $\underset\tilde {}\to \gamma < \gamma$".
\enddefinition
\bigskip

\demo{\stag{1.26} Observation}  1) If $\aleph_0 > \theta_2 > 2$ \ub{then} pure
$(\theta_2,2)$-decidability is equivalent to pure $(2,2)$-decidability. \nl
2) If $\Bbb Q$ is purely semi-proper (see \cite[X]{Sh:f} or here xxx) 
or just $\Bbb Q$ 
satisfies UP$^0(\Bbb I,\bold W)$
(see \S5) \ub{then} $\Bbb Q$ has pure $(\aleph_1,\aleph_1)$-decidability. \nl
3) If $\Bbb Q$ is purely proper, \ub{then} $\Bbb Q$ has 
$(\lambda,\aleph_1)$-decidability for every $\lambda$. \nl
4) If $\Bbb Q$ has the c.c.c. (and we let 
$\le_{\text{pr}}$ be equality if not defined), 
\ub{then} $\Bbb Q$ is purely proper.
\nl
5) If $\le_{\text{pr}} = \le$, \ub{then} $\Bbb Q$ has pure 
$(\lambda,2)$-decidability for every $\lambda$.
\enddemo
\bigskip

\demo{Proof}  Think of the definitions.
\enddemo
\bigskip

\definition{\stag{1.26A} Definition}  1) A forcingnotion $\Bbb P$ is
purely proper for $\chi$ large enough (e.g., ${\Cal P}(\Bbb P) \in
{\Cal H}(\chi)$ is enough) and $N$ is an elementary submodel of
$(H(\chi),\in)$ to which $\Bbb P$ belongs and $p \in N \cap \Bbb P$
\ub{then} three is $(N,\Bbb P)$-semi-generic satisfying $p
\le_{\text{pr}} q \in \Bbb P$, see below. \nl
2) $q$ is $(N,\Bbb P)$-generic if $q \Vdash_{\Bbb P}$ ``if
$\underset\tilde {}\to \tau$ which belongs to $N$ is a $\Bbb P$-name
of an ordinal then $\underset\tilde {}\to \tau[G_{\Bbb P}] \in N \cap
\text{ Ord}$.  \nl
3) A forcing notion $\Bbb P$ is purely semi-proper if in part 4) we
replace $(N,\Bbb P)$-generic by $(N,\Bbb P)$-semi-generic.
\nl
4) $q$ is $(N,\Bbb P)$-semi-generic if $q \Vdash_{\Bbb P}$ ``if
$\underset\tilde {}\to \tau$ which belongs to $N$, is a $\Bbb P$-name
of a countable ordinal then $\underset\tilde {}\to \tau[G_{\Bbb P}]
\in N \cap \omega_1$".
\enddefinition
\bigskip

\proclaim{\stag{1.27} Claim}  Let $\bar{\Bbb Q}$ 
be a $\kappa-{\text{\rm Sp\/}}_e(W)$-iteration.
\nl
1) The property ``$\Bbb Q$ has pure $\delta^*$-decidability and pure 
$(2,2)$-decidability" is preserved by
$\aleph_1-{\text{\rm Sp\/}}_e(W)$-iterations \ub{if} $\delta^*$ is a limit
ordinal. \nl
2) The property ``$\Bbb Q$ has pure $(2,2)$-decidablity" is preserved by
$\aleph_1-{\text{\rm Sp\/}}_e(W)$-iterations.
\endproclaim
\bigskip

\remark{\stag{1.27A} Remark}  1) This is like \cite[Ch.XIV,2.13]{Sh:f} and is
reasonable for iterations not adding reals.  For getting rid of pure
$(2,2)$-decidability at the expense of others, natural demands, see
\S5. \nl
2) Is this not suitable for name$^1$ ordinals only?  By UP help. \nl
3) See proof of \scite{si.7} for use of $\dbcu_n q_n \cup p^*$
and more cases phrase as a subclaim?
\endremark
\bigskip

\demo{Proof}  Let $\alpha = \ell g(\bar{\Bbb Q})$ and let $\le^*_\chi$ be a
well ordering of ${\Cal H}(\chi)$, let
${\underset\tilde {}\to \le^*_{\chi,\Bbb P_\beta}}$ 
be a $\Bbb P_\beta$-name of a well ordering of 
${\Cal H}(\chi)^{\bold V^{{\Bbb P}_\beta}}$.
Let $p \in \Bbb P_\alpha$ and $\underset\tilde {}\to \tau$ be
a $\Bbb P_\alpha$-name of an ordinal $< \theta,\theta \in 
\{2,\delta^*\}$ and let
${\underset\tilde {}\to {\bar \zeta}^0} = \langle
{\underset\tilde {}\to \zeta^0_\varepsilon}:\varepsilon < j \rangle$ be a
witness for $p$ (see \scite{1.13}(F), clause (a)) \wilog \, each
${\underset\tilde {}\to \zeta^0_\varepsilon}$ is full.  
For each $\varepsilon < j$ and $\xi < \alpha$ below we shall
${\underset\tilde {}\to r^0_{\varepsilon,\xi}},
{\underset\tilde {}\to \gamma^0_{\varepsilon,\xi}}$ and
${\underset\tilde {}\to {\bold t}_{\varepsilon,\xi}}$ 
such that ${\underset\tilde {}\to r^0_{\varepsilon,\xi}}$ is a 
$\Bbb P_{\xi +1}$-name of a condition with domain $\subseteq (\xi +1,\alpha), 
{\underset\tilde {}\to \gamma^0_{\varepsilon,\xi}}$ is a
$\Bbb P_{\xi +1}$-name of an ordinal $< \theta$ and 
${\underset\tilde {}\to {\bold t}_{\varepsilon,\xi}}$ is a 
$\Bbb P_{\xi +1}$-name of a truth value, satisfying the following.  
Let $\xi < \alpha,G_{{\Bbb P}_{\xi +1}} 
\subseteq \Bbb P_{\xi +1}$ be generic over
$\bold V$ and ${\underset\tilde {}\to \zeta^0_\varepsilon}[G_{\xi+1}]= \xi$:
\mr
\item "{$(*)_1$}"  if there are $r \in \Bbb P_\alpha/G_{{\Bbb
P}_{\xi+1}}$ and 
$\gamma < \theta$ satisfying (a) + (b) below then $\gamma =
{\underset\tilde {}\to \gamma^0_{\varepsilon,\xi}}
[G_{{\Bbb P}_{\xi +1}}],r =
{\underset\tilde {}\to r^0_{\varepsilon,\xi}}[G_{{\Bbb P}_{\xi+1}}]$ are
such objects, first by the fixed well ordering $<^*_\chi$ and 
${\underset\tilde {}\to {\bold t}_{\varepsilon,\zeta}}
[G_{{\Bbb P}_{\xi+1}}]$ is
truth, and does not depend on $\varepsilon$; if there are no such $\gamma,r$
we let ${\underset\tilde {}\to r^0_{\varepsilon,\zeta}}
[G_{{\Bbb P}_{\xi+1}}]$ be
the empty condition, ${\underset\tilde {}\to \gamma^0_{\varepsilon,\zeta}}
[G_{{\Bbb P}_{\xi+1}}] = 0$ and 
${\underset\tilde {}\to {\bold t}_{\varepsilon,\zeta}}
[G_{{\Bbb P}_{\xi+1}}]$ is false.
{\roster
\itemitem{ $(a)$ }  $p \le_{\text{pr}} 
r \in \Bbb P_\alpha/G_{{\Bbb P}_{\xi+1}},p \restriction
(\xi + 1) = r \restriction (\xi +1)$ and we have
$r \Vdash_{{\Bbb P}_\alpha/G_{{\Bbb P}_{\xi+1}}}$ ``if $\theta =2$ then 
$\underset\tilde {}\to \tau = \gamma$ and if $\theta \ge \aleph_0$ then
$\underset\tilde {}\to \tau < \gamma$"
\sn
\itemitem{ $(b)$ }  if $\xi_1 < \xi$ then no $r'$ satisfies (a) with $\xi_1,
G_{{\Bbb P}_{\xi+1}} \cap \Bbb P_{\xi_1 +1}$.
\endroster}
\ermn
So by \scitet{1.16}(4), there is in $\Bbb P_\alpha$ a condition 
${\underset\tilde {}\to r^0_\varepsilon}$ which is
${\underset\tilde {}\to r^0_{\varepsilon,\xi}}$ if
${\underset\tilde {}\to \zeta^0_\varepsilon}[G] = \xi,\bold t_{\varepsilon,
\zeta}[G] =$ truth, is well defined (the $\underset\tilde {}\to \beta$ there
is ${\underset\tilde {}\to \zeta^0_\varepsilon}$ here!, hence it is
full).  
So easily
$p_1 = p \cup \{{\underset\tilde {}\to r^0_\varepsilon}:\varepsilon < j\}$
belongs to $\Bbb P_\alpha$ and is a pure extension of $p$ (using
$<^*_\chi$, noting that for each $\xi + 1 \le \alpha$ and
$G_{P_{\xi+1}} \subseteq \Bbb P_{\xi+1}$ generic over $\bold V$, if
${\underset\tilde {}\to \zeta^0_{\xi+1}}[G_{P_{\xi+1}}] = \xi =
{\underset\tilde {}\to \zeta^0_{\varepsilon_1}}[G_{P_{\xi+1}}]$ then
${\underset\tilde {}\to r^0_{\varepsilon_1}}[G_{{\Bbb P}_{\xi+1}}] =
{\underset\tilde {}\to r^0_{\varepsilon_2}}[G_{{\Bbb P}_{\xi+1}}])$.
\mn
We now define $p_2 = p_1 \cup \{{\underset\tilde {}\to r^1_\varepsilon}:
\varepsilon < j\}$ where ${\underset\tilde {}\to r^1_\varepsilon}$ is an
atomic $\bar{\Bbb Q}$-named condition with
$\zeta_{\underset\tilde {}\to r^1_\varepsilon} =
{\underset\tilde {}\to \zeta_\varepsilon}$ defined as follows
\mr
\item "{$(*)$}"  if $\beta < \alpha,G_{{\Bbb P}_\beta} \subseteq \Bbb P_\beta$ generic
over $\bold V$ and ${\underset\tilde {}\to \zeta_\varepsilon}[G_{{\Bbb P}_\beta}] =
\beta$ then in $\bold V[G_{{\Bbb P}_\beta}]$ we 
have $r \in {\underset\tilde {}\to {\Bbb Q}_\beta}
[G]$ is ${\underset\tilde {}\to \le^*_{\chi,\Bbb P_\beta}}[G_{{\Bbb P}_\beta}]$-minimal
such that
{\roster 
\itemitem{ $(i)$ }  $\hat{\Bbb Q}_\beta[G_{{\Bbb P}_\beta}] \models ``p_2 \restriction
\{\beta\} \le_{\text{pr}}\{r\} \in \hat{\Bbb Q}_\beta[G_{{\Bbb P}_\beta}]"$
\sn
\itemitem{ $(ii)$ }  for some $r_1 \in \Bbb P_{\beta +1}$ and $\gamma < \theta$
we have: $r_1 \restriction \beta \in G_{{\Bbb P}_\beta}$ and $r_1 \restriction \beta
=r$  and: 
$r_1 \Vdash_{{\Bbb P}_{\beta +1}} ``\theta = 2 \and 
{\underset\tilde {}\to \gamma^0_{\varepsilon,\beta}} = \gamma$ or 
$\theta \ge \aleph_0 \and {\underset\tilde {}\to \gamma^0_{\varepsilon,\xi}} 
< \gamma$ and 
${\underset\tilde {}\to {\bold t}_{\varepsilon,\beta}} =$ truth"
or $r_1$ forces ($\Vdash_{{\Bbb P}_{\beta +1}}$) that
${\bold t}_{\varepsilon,\beta} =$ false.
\endroster}
\ermn
Let us choose now $\beta \le \alpha$ and $r_1$ with $\beta$ minimal such that
\mr
\item "{$\otimes$}"   $r_1 \in \Bbb P_\beta$ and there are $q \in \Bbb P_\alpha$ and 
$\gamma < \theta$ such
that $p_2 \le_{\text{pr}} q$ and $q \restriction \beta \le r_1$ and $r_1
\cup (q \restriction [\beta,\alpha)) \Vdash ``\theta = 2,
\underset\tilde {}\to \tau = \gamma$ or $\theta \ne 2,
\underset\tilde {}\to \tau < \gamma"$.
\ermn
There is such $\beta$ as $\beta = \alpha (= \ell g(\bar \theta))$ is O.K.
\enddemo
\bn
\ub{Case 1}:  $\beta =0$.

We are done.
\bn
\ub{Case 2}:  $\beta$ is limit.

Without loss of generality, by \scite{1.15} for some $n < \omega$ and $\xi_1
< \ldots < \xi_n < \beta$ we have: $p_2 \restriction \beta \le r_1$ above
$\{\xi_1,\dotsc,\xi_n\}$.  If $n=0$ we are done (as $\beta =0$) so
assume $n > 0$.  Let $\beta' = \xi_n +1,r' = r_1 \restriction
(\xi_n + 1)$ and there is $q'$ defined by $q' \restriction (\xi_n +1) = q
\restriction (\xi_n +1),q' \restriction (\xi_n +1,\beta)$ is $r_1 \restriction
(\xi_n +1,\beta)$ if $r_1 \restriction (\xi_n +1) \in
{\underset\tilde {}\to G_{{\Bbb P}_{\xi_n+1}}}$ and is
$r_1 \restriction (\xi_n +1,\beta)$ otherwise and lastly
$q' \restriction [\beta,\alpha) = q \restriction [\beta,\alpha)$.  Now
$\beta',r',q'$ satisfies: $r' \in \Bbb P_{\xi_n +1} = \Bbb P_{\beta'},p_2 
\le_{\text{pr}} q',q' \restriction \beta' \le r'$ and $r' \cup q \restriction
[\beta',\alpha) \Vdash_{{\Bbb P}_\alpha} ``\theta = 2,\underset\tilde {}\to \tau =
\gamma$ or $\theta \ne 2,\underset\tilde {}\to \tau < \gamma"$ and
$\beta' < \beta$.
So we get a contradiction to the choice of $\beta$. 
\bn
\ub{Case 3}:  $r_1 \Vdash ``\beta_0 \notin \{
{\underset\tilde {}\to \zeta^0_\varepsilon}:\varepsilon < j\}"$ where
$\beta = \beta_0 +1$.

The proof is similar to the one of case 2 using $\beta' = \beta_0$.
\bn
\ub{Case 4}:  None of the above.  

So by ``neither case 1 nor case 2" we have
$\beta = \beta_0 +1$, and as we can increase $r_1$ \wilog \, $r_1$ forces
$\beta_0 \in \{{\underset\tilde {}\to \zeta^0_\varepsilon}:\varepsilon 
< j\}$, so \wilog \, $r_1 \Vdash ``\beta_0 = 
{\underset\tilde {}\to \zeta^0_\varepsilon}"$ where $\varepsilon < j$.

Let $r_1 \in G_\beta \subseteq \Bbb P_\beta$ with $G_\beta$ generic
over $\bold V$; let
$G_{\beta'} = G_\beta \cap \Bbb P_{\beta'}$ for $\beta' \le \beta$.
So ${\underset\tilde {}\to \zeta^0_\varepsilon}[G_{\beta_0}] = \beta_0$. \nl
We first ask: is there $\varepsilon_1 < j$ such that
${\underset\tilde {}\to \zeta^0_{\varepsilon_1}}[G_{\beta_0}]$ is well 
defined so call it $\zeta$, (so necessarily $\zeta \le \beta_0$) and
${\underset\tilde {}\to {\bold t}_{\varepsilon_1,\zeta}}[G_\beta]$ 
is truth? \nl
If yes, then we get a contradiction to the minimality of $\beta$ as $\zeta$
can serve by the choice of $\beta_2$, so assume not.  Now considering 
${\underset\tilde {}\to r^0_{
\varepsilon,\xi}},{\underset\tilde {}\to \gamma^0_{\varepsilon,\xi}}$,
clause (b) holds and $r_1,q$ exemplifies
${\underset\tilde {}\to {\bold t}_{\varepsilon,\xi}}[G_p] =$ truth. 
\hfill$\square_{\scite{1.27}}$\margincite{1.27}
\bigskip

\remark{\stag{1.27B} Remark}  1)  You may ask why we do not use the
${\underset\tilde {}\to \zeta^*}$ defined by
${\underset\tilde {}\to \zeta^*}
[G_{\xi +1}] = \xi +1$ if ${\underset\tilde {}\to {\bold t}_{\varepsilon,
\zeta}}[G_{\xi +1}] =$ truth for some $\varepsilon < j$?  The reason is that
(as for $e=6,\kappa = \aleph_1$) this seems not to be a simple 
$\bar{\Bbb Q}$-named ordinal.
\nl
2) By the proof, if $\bar{\Bbb Q}$ is a $\kappa$-Sp$_e(W)$-iteration,
$\alpha \le \ell g(\bar{\Bbb Q}),p \in \Bbb P_\alpha$ and for each
$\beta < \alpha,{\underset\tilde {}\to {\bold t}_\beta}$ is a
$\Bbb P_\beta$-name of a truth value, ${\underset\tilde {}\to p_\beta}$ a
$\Bbb P_\beta$-name of some $p \in \Bbb P_\alpha/
{\underset\tilde {}\to G_{{\Bbb P}_\beta}}$ such that $\Vdash_{{\Bbb P}_\beta}
``p \le_{pr} {\underset\tilde {}\to p_\beta},p \restriction  \beta =
{\underset\tilde {}\to p_\beta} \restriction \beta"$ \ub{then} we can find
$q$ such that $p \le_{\text{pr}} q \in \Bbb P_\alpha$ and: if $G_\beta \subseteq
\Bbb P_\beta$ is generic over $V,\beta < \alpha,
{\underset\tilde {}\to {\bold t}_\beta}[G_\beta] =$ truth and $\gamma < \beta
\Rightarrow ``{\underset\tilde {}\to {\bold t}_\gamma}[G_\beta \cap
\Bbb P_\gamma] =$ false then $q \Vdash_{{\Bbb P}_\alpha/G} ``
{\underset\tilde {}\to {\Bbb P}_\beta} 
\in {\underset\tilde {}\to G_{{\Bbb P}_\beta}}"$.
\endremark  
\bn
We now consider some variants of the $\lambda$-c.c.
\definition{\stag{1.28} Definition}  1) We say $\Bbb P$ satisfies the local
$\underset\tilde {}\to \partial$-c.c. \ub{if} $\underset\tilde {}\to \kappa$
is a $\Bbb P$-name and $\{p \in \Bbb P:\Bbb P \restriction \{q:p \le q 
\in \Bbb P\} \text{ satisfies the }
\partial' \text{-c.c. and } p 
\Vdash_{\Bbb P} ``\underset\tilde {}\to \partial = 
\partial'" \text{ for some } \kappa'\}$ is dense in $\Bbb P$. \nl
2) We say $\Bbb P$ satisfies the local $\underset\tilde {}\to 
\partial$-c.c. purely if the set above is dense in $(\Bbb P,
\le_{\text{pr}})$. \nl
3) We say $\Bbb P$ satisfies lc.pr. $\underset\tilde {}\to \partial$-c.c. if:
\mr
\item "{$(a)$}"  $\underset\tilde {}\to \kappa$ is a 
$(\Bbb P,\le_{\text{pr}})$-name, usually of a regular cardinal of
$\bold V$ \nl
(could use just a partial function from $\Bbb P$ to cardinals such that
$\underset\tilde {}\to \kappa(p) = \kappa \wedge p \le_{\text{pr}} q
\Rightarrow \underset\tilde {}\to \kappa(q) = \kappa$, but abusing
notation we write $q \Vdash ``\underset\tilde {}\to \kappa = \kappa"$
if $\underset\tilde {}\to \kappa(q) = \kappa$)
\sn
\item "{$(b)$}"  for every $p \in \Bbb P$ for some $q,\partial$ we have
$p \le_{\text{pr}} q,q \Vdash_{({\Bbb P},\le_{\text{pr}})}
``\underset\tilde {}\to \partial = \partial"$ 
and $\Bbb P_{\ge q}$ satisfies the $\partial$-c.c. \nl
(we could use: if $\underset\tilde {}\to \kappa(p) = \kappa$ then
$\Bbb P_{\ge p}$ i.e. $(\{q:p \le q \in \Bbb P\},
\le^{\Bbb P})$ satisfies the $\kappa$-c.c.)
\ermn
4) If $\Bbb P$ satisfies the lc.pr. $\underset\tilde {}\to \partial$-c.c.
and $q \in \Bbb P$ let $\kappa^{\text{mcc}}_{\underset\tilde {}\to
\partial}(q,\Bbb P) = \partial$ means $q 
\Vdash_{({\Bbb P},\le_{\text{pr}})}$ ``$\underset\tilde {}\to \kappa =
\kappa$" and $\Bbb P_{\ge q}$ satisfies the $\kappa$-c.c. \nl
5) Let $\underset\tilde {}\to \partial^{\text{mcc}}
(\Bbb P)$ be minimal such that $\Bbb P$ satisfies
the lc.pr.  $\underset\tilde {}\to \kappa$-c.c.; that is
${\underset\tilde {}\to \partial'}(q,\Bbb P) 
= \text{ Min}\{\kappa:\Bbb P_{\ge q}$
satisfies that $\kappa$-c.c.$\}$ and 
$\underset\tilde {}\to \partial(q,\Bbb P) = \kappa$ 
if $(\forall r')(q \le_{\text{pr}} 
r \rightarrow {\underset\tilde {}\to \partial'}
(r,\Bbb P) = \kappa$) (see below) and let
$\partial^{\text{mcc}}(q,\Bbb P) = \partial$ 
mean $\kappa^{\text{mcc}}_{\underset\tilde {}\to \partial({\Bbb P})}
(q,\Bbb P) = \partial$ where $\underset\tilde {}\to \kappa = 
\underset\tilde {}\to \kappa_{\text{mcc}}(\Bbb P)$.
\enddefinition
\bigskip

\proclaim{\stag{1.29} Claim}  1) For a forcing notion 
$\Bbb P$ (as in \scite{1.1}) 
the $(\Bbb P,\le_{\text{pr}})$-name 
$\underset\tilde {}\to \partial^{\text{mcc}}(\Bbb P)$ is well
defined, so \nl
2) If $\Bbb P$ satisfies the lc.pr. 
$\underset\tilde {}\to \partial$-c.c. and $p \in \Bbb P$ 
\ub{then} for some $q$ we have $p \le_{\text{pr}} q$ and
$\partial_{\text{mcc}}(q,\Bbb P)$ is well defined.
\endproclaim
\bigskip

\demo{Proof}  Straight.
\enddemo
\bigskip

\definition{\stag{1.30} Definition}  1) We say $\Bbb Q$ has strong pr.
$({\underset\tilde {}\to \partial_1},
{\underset\tilde {}\to \partial_2})$-decidability when 
$\underset\tilde {}\to \kappa,{\underset\tilde {}\to \kappa_2}$ are
$(\Bbb Q,\le_{\text{pr}})$-names of regular cardinals of $\bold V$ 
and \ub{if} $p \in \Bbb Q,
p \Vdash_{({\Bbb Q},\le_{\text{pr}})} 
``{\underset\tilde {}\to \partial_1} = \theta_1$
and ${\underset\tilde {}\to \kappa_2} = \theta_2"$ and
${\underset\tilde {}\to \zeta_\varepsilon}$ is a $\Bbb Q$-name of an ordinal
$< \theta_1$ for $\varepsilon < \varepsilon^* < \theta_2$ \ub{then} for some
$a \subseteq \theta_1$ of cardinality $< \theta_2$ and $q$ such that
$p \le_{\text{pr}} q \in \Bbb Q$ we have $q \Vdash_{\Bbb Q} 
``{\underset\tilde {}\to \zeta_\varepsilon} \in a$ for $\varepsilon <
\varepsilon^*"$. \nl
2) We say $\Bbb Q$ has strong pr. 
$\underset\tilde {}\to \partial$-decidability
\ub{if} for any $\theta$ it has 
pr. $(\theta,\underset\tilde {}\to \kappa)$-decidability 
(i.e. each ${\underset\tilde {}\to \zeta_\varepsilon}$ is a $\Bbb Q$-name 
of an ordinal $< \theta$). \nl
3) We use ``weak" instead of ``strong" in parts (1), (2) if above we restrict
ourselves to the case $\varepsilon^*=1$. \nl
4) We let ${\underset\tilde {}\to \partial^w_\otimes}(\Bbb P)$ is the minimal
$(\Bbb P,\le_{\text{pr}})$-name $\underset\tilde {}\to \kappa$ of a regular
cardinal from $\bold V$ such that $\Bbb P$ has weak pr.
$\underset\tilde {}\to \kappa$-decidability.  Similarly
${\underset\tilde {}\to \partial^{\text{St}}_\otimes}(\Bbb P)$ for strong
pr. $\underset\tilde {}\to \kappa$-decidability.
\enddefinition
\bigskip

\demo{\stag{1.31} Note/Observation}  If $\bar{\Bbb Q}$ is 
an $\kappa-\text{Sp}_e(W)$-iteration, \ub{then} 
\mr
\item "{$(a)$}"  $\langle (\Bbb P_i,\le^{{\Bbb P}_i}_{\text{pr}}):
i \le \ell g(\bar{\Bbb Q})
\rangle$ is a $\lessdot$-increasing sequence
\sn
\item "{$(b)$}"  if $i < j \le \ell g(\bar{\Bbb Q}),q \in \Bbb P_j,
q \restriction i \le_{\text{pr}} p \in \Bbb P_i$ then $p,q$ has a 
$\le_{\text{pr}}$-lub, $p \cup q \restriction [i,j)$.
\endroster
\enddemo
\bigskip

\demo{Proof}  Check.
\enddemo
\bigskip

\proclaim{\stag{1.32} Claim}  1) If $\Bbb P$ satisfies the lc.pr.
$\underset\tilde {}\to \partial$-c.c. and $\partial =
\partial^{\text{mcc}}_{\underset\tilde {}\to \partial}(p,\Bbb P) \Rightarrow 
\Vdash_{\Bbb P} ``\partial$ is regular"
\ub{then} 
$\Bbb P$ has a strong pr $\underset\tilde {}\to \partial$-decidability. \nl
2) Let $\bar{\Bbb Q}$ be a $\kappa_1-{\text{\rm Sp\/}}_e(W)$-iteration
$e \in \{4\}$.  If $\delta \le \ell g(\bar{\Bbb Q})$ is a limit ordinal, $u$ an unbounded subset of
$\delta$, for $i \in u$ we have $\Bbb P_i$ has the strong pr.
${\underset\tilde {}\to \partial_i}$-decidability \ub{then} letting
$\underset\tilde {}\to \partial = \text{ Min}\{\partial:
\partial \text{ a regular cardinal in } \bold V$ and $\partial 
\ge {\underset\tilde {}\to \partial_i}$ for $i \in u\}$ we have
\mr
\item "{$(i)$}"  $\underset\tilde {}\to \partial$ is a $(\Bbb P_\delta,
\le_{\text{pr}})$-name of a regular cardinal of $\bold V$
\sn
\item "{$(ii)$}"  $\Bbb P_\delta$ has weak pr. 
$\underset\tilde {}\to \partial$-decidability.
\ermn
3) Similarly for 
$\Bbb P_\delta/G_{{\Bbb P}_\alpha}$ when $\alpha < \delta$ and even
$\Bbb P_\delta/\Bbb P_\alpha$ where $\underset\tilde {}\to \alpha$ is 
a simple $\bar{\Bbb Q}$-named $[0,\delta)$-ordinal. \nl
4) If $\Bbb P$ 
satisfies the strong pr $\underset\tilde {}\to \partial$-decidability
and $p \in \Bbb P,\kappa = \kappa_{\underset\tilde {}\to
\kappa}(p,\Bbb P)$ \ub{then} $p \Vdash ``\partial$ is a regular cardinal.
\endproclaim
\bigskip

\demo{Proof}  1) Trivial. \nl
2),3) Similar to the proof of \scite{1.27} [Saharon!] \nl
4) Trivial.  \hfill$\square_{\scite{1.32}}$\margincite{1.32}
\enddemo
\bn
\margintag{n.0}\ub{\stag{n.0} Convention}:  Let $\bar{\Bbb Q} = \langle \Bbb P_j,
{\underset\tilde {}\to {\Bbb Q}_i}:j \le \alpha,i < \alpha \rangle$ be a
$\kappa$-Sp$_e(W)$-iteration.
\bigskip

\proclaim{\stag{n.1} Claim}  Assume $(\bar{\Bbb Q},W)$ is smooth, (see
Definition \scite{1.6}(8)). \nl
1) If $p^* \in \Bbb P_j$ \ub{then} $\langle \Bbb P'_j,
{\underset\tilde {}\to {\Bbb Q}'_i}:j \le
\alpha,i < \alpha \rangle$ is a $\kappa_1$-Sp$_e(W)$-iteration where \nl
$\Bbb P'_j = \{p \in \Bbb P_j:p^* \restriction j \le p\}$, \nl
$\Bbb Q'_j = \{p \in {\underset\tilde {}\to {\Bbb Q}_j}:
\Vdash_{{\Bbb P}_j} ``p^* \restriction
\{j\} \le p \text{ in } \hat{\Bbb Q}_j"$.
\mn
2) If $\gamma < \alpha$ and 
$G_j \subseteq \Bbb P_\gamma$ is generic over $\bold V$ \ub{then}
$\langle \Bbb P_{\gamma +j}/G_\gamma,
{\underset\tilde {}\to {\Bbb Q}_{\gamma +i}}:j \le
\alpha - \gamma$ and $i < \alpha - \gamma \rangle$ is a
$\aleph_1-{\text{\rm Sp\/}}_e(W)$-iteration.
\endproclaim
\bigskip

\demo{Proof}  Straight.
\enddemo
\newpage

\head {\S2 Tree of Models} \endhead  \resetall \sectno=2
\bn
We present here needed information on trees and tagged trees.
On partition of tagged trees, see Rubin, Shelah
\cite{RuSh:117} and, \cite{Sh:b}, 
\cite[XI,3.5,3.5A,3.7,XV,2.6,2.6A,2.6B]{Sh:f} and
\cite[2.4,2.5,p.111-113]{Sh:136}; or the representation in
\cite[AP,\S1]{Sh:e}; on history see \cite{RuSh:117} and \cite{Sh:f}. 
\bigskip

\definition{\stag{2.1} Definition}  A tagged tree is a pair $(T,\bold I)$ 
such that:
\mr
\item  $T$ is a $\omega$-tree, which here means a nonempty set of finite
sequences of ordinals such that if $\eta \in T$ then any initial segment of
$\eta$ belongs to $T$.  $T$ is ordered by initial segments; i.e., $\eta 
\triangleleft \nu$ iff $\eta$ is a proper initial segment of $\nu$.
\sn
\item  $\bold I$ is a partial function such that for every $\eta \in
T \cap \text{ Dom}(\bold I)$: if $\bold I(\eta) = \bold I_\eta$ is 
defined then $\bold I(\eta)$ is an
ideal of subsets of some set called the domain of $\bold I_\eta$,
Dom$(\bold I_\eta)$, and \footnote{in this section it is not
unreasonable to demand equality but this is very problematic in
Definition \scite{4.1}(1), clause (d)}

$$
\text{Suc}_T(\eta) =: \{\nu:\nu \text{ is an immediate successor of } \eta
\text{ in } T\} \subseteq \text{ Dom}(\bold I_\eta),
$$
\mn
Usually $\bold I_\eta$ is $\aleph_2$-complete.
\sn
\item  For every $\eta \in T$ we have Suc$_T(\eta) \ne \emptyset$. 
\endroster
\enddefinition
\bigskip

\demo{\stag{2.2} Convention} For any tagged tree $(T,\bold I)$ we can define
$\bold I^\dag$ by:

$$
\text{Dom}(\bold I^\dag) = \{\eta:\eta \in T \cap \text{ Dom}(\bold I)
\text{ and Suc}_T(\eta) \subseteq \bold I_\eta
\text{ and  Suc}_T(\eta) \notin \bold I_\eta\}
$$
\mn
and

$$
\bold I^\dag_\eta = \bigl\{ \{\alpha:\eta \char 94 \langle \alpha \rangle \in
A\}:A \in \bold I_\eta \bigr\};
$$
\mn
we sometimes, in an abuse of notation, do not distinguish between $\bold I$
and $\bold I^\dag$; e.g. if $\bold I^\dag_\eta$ is constantly $I^*$, we write
$I^*$ instead of $\bold I$.  Also, if $\bold I = \bold I_x$, we may write
$\bold I^x_\eta$ for $\bold I_x(\eta)$.
\enddemo
\bigskip

\definition{\stag{2.3} Definition}  1) Let $\eta$ be called a splitting point
of $(T,\bold I)$ if $\bold I_\eta$ is well defined and Suc$_T(\eta) \notin
\bold I_\eta$ (normally this follows but we may ``forget" 
to decrease the domain of $\bold I$).  Let split$(T,\bold I)$ be 
the set of splitting points of
$(T,\bold I)$.  We usually consider trees where each $\omega$-branch meets 
split$(T,\bold I)$ infinitely often (see Definition \scite{2.5}(6),
i.e. are normal). \nl
2) For $\eta \in T$ let $T^{[\eta]} =: \{\nu \in T:\nu = \eta \text{ or }
\nu \triangleleft \eta \text{ or } \eta \triangleleft \nu\}$. \nl
3) For a tree $T$, let lim$(T)$ be the set of branches of $T$; i.e. all
$\omega$-sequences of ordinals, such that every finite initial segment 
of them is a member of $T:\text{lim}(T) =$ \nl
$\{s \in {}^\omega\text{Ord}:(\forall n)s \restriction n \in T\}$. 
We call them also $\omega$-branches. \newline
4) A subset $Z$ of a tree $T$ is a front if: $\eta \ne \nu \in Z$ implies
none of them is an initial segment of the other, and every $\eta \in
\text{ lim}(T)$ has an initial segment which is a member of $Z$. 
\enddefinition
\bigskip

\definition{\stag{2.4} Definition}  We now define orders between tagged trees:
\mr
\item "{$(a)$}"  $(T_1,\bold I_1) \le (T_2,\bold I_2)$ \ub{if} $T_2 \subseteq
T_1$, and split$(T_2,\bold I_2) \subseteq \text{ split}(T_1,\bold I_1)$, and
\nl
for every $\eta \in \text{ split}(T_2,\bold I_2)$ we have 
$\bold I_2(\eta) \restriction \text{ Suc}_{T_2}(\eta) = \bold I_1(\eta) 
\restriction \text{ Suc}_{T_2}(\eta)$ \nl
(where $I \restriction A = \{B:B \subseteq A \text{ and } B \in I\})$. \nl
(So every splitting point of $T_2$ is a splitting point of $T_1$, and
$\bold I_2$ is completely determined by $\bold I_1$ and 
split$(T_2,\bold I_2)$ and $T_2$.)
\sn
\item "{$(b)$}"  $(T_1,\bold I_1) \le^* (T_2,\bold I_2)$ \ub{if} 
$(T_1,\bold I_1)
\le (T_2,\bold I_2)$ and \nl
split$(T_2,\bold I_2) = \text{ split}(T_1,\bold I_1) \cap T_2$.
\sn
\item "{$(c)$}"  For any set 
$A,(T_1,\bold I_1) \le^\otimes_A (T_2,\bold I_2)$
if $T_1 \supseteq T_2$ and $\eta \in A \cap T_2 \Rightarrow \text{ Suc}_{T_2}
(\eta) = \text{ Suc}_{T_1}(\eta)$ and $\eta \in T_2 \cap
\text{ split}(T_1,\bold I_1)
\backslash A \Rightarrow \text{ Suc}_{T_2}(\eta) \ne \emptyset \text{ mod }
\bold I^1_\eta$.
\sn
\item "{$(d)$}"  In (c) we may omit the subscript $A$ when $A = T_2
\backslash \text{ split}(T_1,\bold I_1)$.
\sn
\item "{$(e)$}"  $(T_1,\bold I_1) \le^\boxtimes_\kappa (T_2,\bold I_2)$
\ub{if} $T_2 \subseteq T_1$ is a subtree and if $\nu \triangleleft
\eta \in \text{ lim}(T_2),\bold I_\nu$ is $\kappa$-complete, then for some
$k \ge \ell g(\nu),\bold I^1_\nu \le_{RK} \bold I^2_{\eta \restriction k},
\bold I^2_{\eta \restriction k}$ is 
$\kappa$-complete and $\eta \restriction k
\in \text{ split}(T_2,\bold I_2)$.  We can replace $\kappa$ and 
$\kappa$-complete by $\varphi$ and ``satisfying $\varphi$", e.g. $\in
\Bbb I$.
\sn
\item "{$(f)$}"  $(T_1,\bold I_1) \le^\kappa (T_2,\bold I_2)$
\ub{when} $T_2 \subseteq T_1$, and if $\eta \in T_2 \cap$split$(T_1,\bold
I_1)$ and $\bold I^1_\eta$ is $\kappa$-complete then 
$\eta \in \text{ split}(T_2,\bold I_2)$ and $\bold I^2_\eta = \bold
I^1_\eta$ and every $\eta \in \text{ split}(T_2,\bold I_2)$ is like
that.
\endroster
\enddefinition
\bigskip

\definition{\stag{2.5} Definition}  1) For a set $\Bbb I$ of ideals, a 
tagged tree $(T,\bold I)$ is an $\Bbb I$-tree if for every splitting point 
$\eta \in T$ we have $\bold I_\eta \in \Bbb I$ (up to an isomorphism). \nl
2) For a tagged tree $(T,\bold I)$ and set $\Bbb I$ of ideals let
$\bold I \restriction \Bbb I = \bold I \restriction \{\eta \in \text{ Dom}
(\bold I)$ and $\bold I_\eta \in \Bbb I\}$. \nl
3) Let in (2), $(T,\bold I) \restriction \Bbb I = (T,\bold I \restriction
\Bbb I)$. \nl
4) A tagged tree 
$(T,\bold I)$ is standard \ub{if} for every non-splitting point $\eta \in T,
|\text{Suc}_T(\eta)| = 1$. \newline
5) A tagged tree $(T,\bold I)$ is full \ub{if} 
every $\eta \in T$ is a splitting point. \nl
6) A tagged tree $(T,\bold I)$ is normal \ub{if} 
for every $\eta \in \text{ lim}(T)$ for
infinitely many $k < \omega$ we have $\eta \restriction k \in \text{
split}(T,\bold I)$.
\enddefinition
\bigskip

\remark{\stag{2.6} Remark}  1) Of course, the set 
lim$(T)$ is not absolute; i.e., if
$\bold V_1 \subseteq \bold V_2$ are two universes of set 
theory then in general
(lim$(T))^{{\bold V}_1}$ will be a proper 
subset of (lim$(T))^{{\bold V}_2}$. \newline
2) However, the notion of being a front is absolute: if $\bold V_1 
\models ``Z$
contains a front in $T$", \ub{then} there is a depth function $\rho:T 
\rightarrow \text{ Ord}$ 
satisfying $\eta \triangleleft \nu \and \forall k \le \ell g(\eta)
[\eta \restriction k \notin Z] \rightarrow \rho(\eta) > \rho(\nu)$.  This
function will also witness in $\bold V_2$ that $Z$ is a front. \newline
3) $Z \subseteq T$ contains a front \ub{iff} $Z$ 
meets every branch of $T$.  So if
$Z \subseteq T$ contains a front of $T$ and $T' \subseteq T$ is a subtree, 
\ub{then} $Z \cap T'$ contains a front of $T'$.  This is absolute, too.
\endremark
\bigskip

\definition{\stag{2.7} Definition}  1) An ideal $I$ is
$\lambda$-complete  \ub{if}
any union of less than $\lambda$ members of $I$ is still a member of $I$.
\newline
2) A tagged tree $(T,\bold I)$ is $\lambda$-complete if for each $\eta \in
T \cap \text{ Dom}(\bold I)$ or just $\eta \in \text{ split}(T,\bold I),$
the ideal $\bold I_\eta$ is $\lambda$-complete.
\newline
3) A family $\Bbb I$ of ideals is $\lambda$-complete \ub{if} 
each $I \in \Bbb I$ is
$\lambda$-complete.  We will only consider $\aleph_2$-complete families 
$\Bbb I$. \newline
4) A family $\Bbb I$ is called restriction-closed \ub{if}: 
$I \in \Bbb I,A \subseteq \text{ Dom}(I),A \notin I$ 
implies $I \restriction A = \{B \in I:B \subseteq 
A\}$ belongs to $\Bbb I$. \newline
5) The restriction closure of $\Bbb I$, res-c$\ell(\Bbb I)$ is 
$\{I \restriction A:I \in \Bbb I,
A \subseteq \text{ Dom}(I),A \notin I\}$. \nl
6) $I$ is $\lambda$-indecomposable \ub{if}
for every $A \subseteq \text{ Dom}(I),
A \notin I$ and $h:A \rightarrow \lambda$ there is $Y \subseteq \lambda,|Y|
< \lambda$ such that $h^{-1}(Y) \notin I$.  We say $\bold I$ 
(or we say $\Bbb I$) is
$\lambda$-indecomposable \ub{if} each 
$\bold I_\eta$ where $\eta \in \text{ split}(T,
\bold I)$ (or $I \in \Bbb I$) is $\lambda$-indecomposable. \newline
7) $I$ is strongly $\lambda$-indecomposable \ub{if} for any
$A_i \in I \,\,(i < \lambda)$ and $A \subseteq \text{ Dom}(I),
A \notin I$ we can find $B \subseteq
A$ of cardinality $< \lambda$ such that for no $i < \lambda$ does $A_i$ 
include $B$. \newline
8) Let $\Bbb I^{[\kappa]} = 
\{I \in \Bbb I:I \text{ is } \kappa \text{-complete}\}$.
\enddefinition
\bigskip

\demo{\stag{2.9} Fact}  If $\lambda$ is a regular cardinal and $I$ 
is a strongly $\lambda$-indecomposable,
\ub{then} $I$ is $\lambda$-indecomposable.
\enddemo
\bigskip

\demo{Proof}  Given 
$A,h$ as in \scite{2.7}(6), let $A_i = h^{-1}(\{j:j < i\})$ and $A'_i
= h^{-1}(\{i\})$; if
for some $i < \lambda,A_i \notin I$ we are done, 
otherwise by Definition \scite{2.7}(7)
there is $B \subseteq A,|B| < \lambda$ such that: $i < \lambda \Rightarrow
B \nsubseteq A_i$.  But as $\lambda$ is regular, $B \subseteq A,|B| <
\lambda$ and $\langle A_i:i < \lambda \rangle$ is a
$\subseteq$-increasing sequence of sets with union $A$, clearly
for some $j < \lambda,B \subseteq A_j$, contradiction.
\hfill$\square_{\scite{2.9}}$\margincite{2.9}
\enddemo
\bigskip

\definition{\stag{2.10} Definition}  For a subset of $A$ of 
(an $\omega$-tree) $T$ we define by induction on the length 
of a sequence $\eta$, res$_T(\eta,A)$
for each $\eta \in T$.  Let res$_T(\langle \rangle,A) = \langle \rangle$.
Assume res$_T(\eta,A)$ is already defined and we define res$_T(\eta \char 94
\langle \alpha \rangle,A)$ for all members $\eta \char 94 \langle \alpha
\rangle$ of Suc$_T(\eta)$.  If $\eta \in A$ then res$_T(\eta \char 94
\langle \alpha \rangle,A) = \text{ res}_T(\eta,A) \char 94 \langle 
\alpha \rangle$, and if $\eta \notin A$ then res$_T(\eta \char 94 \langle 
\alpha \rangle,A) = \text{ res}_T(\eta,A) \char 94 \langle 0 \rangle$.
If $\eta \in \text{ lim}(T)$, we let \newline
res$(\eta,A) = \dsize \bigcup_{k \in \omega} \text{ res}
(\eta \restriction k,A)$.
\enddefinition
\bigskip

\demo{\stag{2.11} Explanation}  Thus res$(T,A) =: \{\text{ res}_T(\eta,A):\eta
\in T\}$ is a tree obtained by projecting $T$; i.e., glueing 
together all members of Suc$_T(\nu)$ whenever $\nu \notin A$.
\enddemo
\bigskip

We state now (Lemma \scite{2.12} is \cite[Ch.XI,5.3;p.559]{Sh:f} and Lemma
\scite{2.13} is \nl
\cite[XV,2.6;p.738]{Sh:f} and Lemma \scite{2.14} is
\cite[XI,3.5(2); p.546-7]{Sh:f}.
\proclaim{\stag{2.12} Lemma}  Let $\lambda,\mu$ be uncountable cardinals
satisfying $\lambda^{< \mu} = \lambda$ and let $(T,\bold I)$ be a tagged
tree in which for each $\eta \in T$ either 
$|{\text{\rm Suc\/}}_T(\eta)| < \mu$ or
$\bold I(\eta)$ is $\lambda^+$-complete.  \ub{Then} for every function $H:T
\rightarrow \lambda$ there exists a subtree $T'$ of $T$ 
satisfying $(T,\bold I) \le^* (T',
\bold I)$ such that for $\eta^1,\eta^2 \in T'$ satisfying
${\text{\rm res\/}}_T(\eta^1,A) = { \text{\rm res\/}}_T(\eta^2,A)$ we
have:
\mr
\widestnumber\item {$(iii)$}
\item "{$(i)$}"  $H(\eta^1) = H(\eta^2)$
\sn
\item "{$(ii)$}"  $\eta^1 \in A \Leftrightarrow \eta^2 \in A$,
\sn
\item "{$(iii)$}"  $\text{if } \eta \in T' \cap A$, then 
${\text{\rm Suc\/}}_T(\eta) = { \text{\rm Suc\/}}_{T'}(\eta)$.
\endroster
\endproclaim
\bigskip

\demo{Proof}  See \cite[XV,2.6;p.738]{Sh:f}.
\enddemo
\bigskip

\proclaim{\stag{2.13} Lemma}  Let $\theta$ be an uncountable regular cardinal
(the main case here is $\theta = \aleph_1$).  Assume
\mr
\item "{$(\alpha)$}"   $\Bbb I$ be a family of $\theta^+$-complete
ideals,
\sn
\item "{$(\beta)$}"   $(T_0,\bold I)$ a tagged tree,
\sn
\item "{$(\gamma)$}"  $A =: \{\eta \in T:|\text{Suc}_{T_0}(\eta)| 
\le \theta\}$,
\sn
\item "{$(\delta)$}"  $[\eta \in T_0 \backslash A \Rightarrow
\bold I_\eta \in \Bbb I \and { \text{\rm Suc\/}}_{T_0}(\eta) 
\notin \bold I_\eta]$ 
\sn
\item "{$(\varepsilon$}"  $[\eta \in A \Rightarrow { \text{\rm
Suc\/}}_{T_0}(\eta) \subseteq \{ \eta \char 94
\langle i \rangle:i < \theta\}]$ 
\sn
\item "{$(\zeta)$}"   $H:T_0 \rightarrow \theta$ and
\sn
\item "{$(\eta)$}"  $\bar e = 
\langle e_\eta:\eta \in A,|{\text{\rm Suc\/}}_{T_0}(\eta)| = \theta \rangle$ 
is such that $e_\eta$ is a club of $\theta$.  
\ermn
\ub{Then} there is a club $C$ of $\theta$ 
such that: for each $\delta \in C$ there is $T_\delta \subseteq T_0$ 
satisfying:
\mr
\item "{$(a)$}"  $T_\delta$ is a tree
\sn
\item "{$(b)$}"  if $\eta \in T_\delta,|{\text{\rm
Suc\/}}_{T_0}(\eta)| < \theta$, then ${\text{\rm
Suc\/}}_{T_\delta}(\eta) 
= { \text{\rm Suc\/}}_{T_0}(\eta)$ and if ${\text{\rm Suc\/}}(\eta) = 
\theta$, then ${\text{\rm Suc\/}}_{T_\delta}(\eta) = \{ \eta \char 94 
\langle i \rangle:i < \delta\}$ and $\delta \in e_\eta$
\sn
\item "{$(c)$}"  $\eta \in T_\delta \backslash A$ implies
${\text{\rm Suc\/}}_{T_\delta}(\eta) \notin \bold I_\eta$
\sn
\item "{$(d)$}"  for every $\eta \in T_\delta$ we have $H(\eta) < \delta$.
\endroster
\endproclaim
\bigskip

\demo{Proof}  See \cite[XV,2.6]{Sh:f}.
\enddemo
\bigskip

\proclaim{\stag{2.14} Lemma}   Let $(T,\bold I)$ be a $\theta^+$-complete
$\bold I$-tree, and assume $\theta = { \text{\rm cf\/}}(\theta)$.  
If ${\text{\rm lim\/}}(T) = \dsize \bigcup_{i < \theta} B_i$
where $B_i$ is a Borel subset of ${\text{\rm lim\/}}(T)$, 
\underbar{then} for some $i <
\theta$ and $(T',\bold I')$ we have $(T,\bold I) \le^* (T',\bold I')$ and
${\text{\rm lim\/}}(T') \subseteq B_i$.
\endproclaim
\bigskip

\demo{Proof}  By \cite[XI,3.5(2);p.546-7]{Sh:f}.
\enddemo
\newpage

\head {\S3 Ideals and Partial Orders} \endhead  \resetall \sectno=3
\bigskip

\definition{\stag{3.1} Definition}  1) We call an ideal $J$ non-atomic if
$\{x\} \in J$ for every \newline
$x \in \text{ Dom}(J)$. \newline
2) We call the ideal with domain $\{0\}$, which is $\{\emptyset\}$, the
trivial ideal.
\enddefinition
\bigskip

\definition{\stag{3.2} Definition}  1) For ideals $J_1,J_2$ we say

$$
J_1 \le_{\text{RK}} J_2
$$
\mn
$\text{ Dom}(J_1)$ is such that

$$
\text{for every } A \subseteq \text{ Dom}(J_2) \text{ we have}:
A \ne \emptyset \text{ mod }
J_2 \Rightarrow h''(A) \ne \emptyset \text{ mod } J_1
$$
\mn
or equivalently,

$$
J_2 \supseteq \{h^{-1}(A):A \in J_1\}.
$$
\mn
2) For families $\Bbb I_1,\Bbb I_2$ of ideals we say $\Bbb I_1 
\le_{\text{RK}} \Bbb I_2$ if there is a function $H$ witnessing it; i.e.,
\mr
\item "{$(i)$}"  $H$ is a function from $\Bbb I_1$ into $\Bbb I_2$
\sn
\item "{$(ii)$}"  for every 
$J \in \Bbb I_1$ we have $J \le_{\text{RK}} H(J)$.
\ermn
3) For families $\Bbb I_1,\Bbb I_2$ of ideals, $\Bbb I_1 
\equiv_{\text{RK}} \Bbb I_2$ \ub{if} $\Bbb I_1 \le_{\text{RK}} \Bbb I_2 \and 
\Bbb I_2 \le_{\text{RK}} \Bbb I_1$.
\enddefinition
\bigskip

\proclaim{\stag{3.3} Claim}  1) If an ideal $J$ is not non-atomic \ub{then} 
$J \le_{\text{RK}}$ ``the trivial ideal". \newline
2) $\le_{\text{RK}}$ is a partial quasi-order (among ideals and also among
families of ideals).
\endproclaim
\bigskip

\demo{Proof}  Easy.
\enddemo
\bigskip

\definition{\stag{3.4} Definition}  1) For an (upward) directed quasi 
order \footnote{no real difference if we ask partial order or just quasi
orders; i.e., partial orders satisfy 
$x \le y \and y \le x \nRightarrow x = y$, quasi order not
necessarily; note that in the case
there is $x^* \in L$ such that $(\forall y \in L)(y \le x^*)$ gives an ideal
which is not non-atomic} $L = (B,<)$ we define an ideal id$_L$:

$$
\text{id}_L = \bigl\{ A \subseteq B:\text{for some } y \in L \text{ we have }
A \subseteq \{x \in L:\neg y \le x\} \bigr\}.
$$
\mn
Equivalently, the dual filter fil$_L$ is generated by the ``cones",
where the cone of $L$ defined by $y \in L$ is 

$$
L_y =: \{ x \in L:y \le x\}.
$$
\mn
We call such an ideal a partial order ideal.  We let Dom$(L) = 
\text{ Dom}(\text{id}_L)(= B)$, but we may use $L$ instead of Dom$(L)$
(like $\forall x \in L$). \newline
2) For a partial order $L$ let dens$(L) = \min\{|X|:X \subseteq \text{ Dom}
(L)$ is dense \footnote{also called cofinal in this context}; 
i.e. $(\forall a \in \text{ Dom}(L))(\exists b \in X)
[a \le b])\}$ (this applies also to ideals considered as $(I,\subseteq)$).
\newline
3) For a family ${\Cal L}$ of directed quasi orders let 
id$_{\Cal L} = \{\text{id}_L:L \in {\Cal L}\}$.
\enddefinition
\bigskip

\demo{\stag{3.5} Fact}  1)  id$_L$ is $\lambda$-complete \ub{iff} $L$ is 
$\lambda$-directed. \newline
2) dens$(L) = \text{ dens(id}_{(L,<)},\subseteq)$. \newline
3) If $h:L_1 \rightarrow L_2$ preserves order (i.e. $\forall x,y \in L,(x \le
y \Rightarrow h(x) \le h(y))$) and has cofinal (= dense) range (i.e. $(\forall x \in
L_2)(\exists y \in L_1)[x \le h(y)]$) \ub{then} id$_{L_2} \le_{\text{RK}}
\text{ id}_{L_1}$. \newline
4) Let $h:L_1 \rightarrow L_2$.  Now $h$ exemplifies 
id$_{L_2} \le_{\text{RK}} \text{ id}
_{L_1}$ \ub{iff} for every $x_2 \in L_2$ there is $x_1 \in L_1$ such that:
$(\forall y)[y \in L_1 \and x_1 \le_{L_1} y \Rightarrow x_2 \le_{L_2} h(y)]$;
(equivalently for every 
$y \in L_1:\neg(x_2 \le_{L_2} h(y)) \Rightarrow \neg(x_1 \le_{L_1}
y)$; note that $h$ is not necessarily order preserving). \newline
5) If $L_1 \subseteq L_2$ and
$L_1$ is a dense in $L_2$, \ub{then} id$_{L_1} \equiv_{RK} 
\text{ id}_{L_2}$. \nl
6) If $\lambda = \text{ density}(L)$ \ub{then} id$_L$ is $\lambda$-based; i.e.
$X \subseteq \text{ Dom}(\text{id}_L),X \notin \text{ id}_L \Rightarrow
(\exists Y \subseteq X)[|Y| \le \lambda \and Y \notin \text{ id}_L]$.
\enddemo
\bigskip

\demo{Proof}  Straight.  E.g. \newline
3) If $A \subseteq L_1$ and $A \notin \text{ id}_{L_1}$, then $(\forall x \in
L_1)(\exists y)(x \le_{L_1} y \in A)$ hence $(\forall x \in L_2)(\exists y,
z \in L_1)(x \le_{L_2} h(y) \and y \le_{L_1} z \in A)$ hence
$(\forall x \in L_2)(\exists z)(x \le_{L_2} z \in h''(A))$ hence
$h''(A) \notin \text{ id}_{L_2}$ (and trivially $h''(A) \subseteq L_1)$. \nl
By Defintion \scite{3.2}(1) this 
shows id$_{L_2} \le_{RK} \text{ id}_{L_1}$. \nl
4) Note: $h$ exemplifies id$_{L_2} \le \text{ id}_{L_1}$ 
\sn
iff

$$
(\forall A \subseteq L_1)(A \ne \emptyset \text{ mod id}_{L_1} 
\rightarrow h''(A) \ne \emptyset \text{ mod id}_{L_2})
$$ 
\mn
\underbar{iff}

$$
(\forall A \subseteq L_1)(\forall x_2 \in L_2)[A \ne \emptyset 
\text{ mod id}_{L_1} \rightarrow h''(A) \cap \{y \in L_2:x_2 \le_{L_2} y\}
\ne \emptyset]
$$
\medskip

\noindent
iff

$$
(\forall x_2 \in L_2)(\forall A \subseteq L_1)(A \ne \emptyset 
\text{ mod id}_{L_1} \rightarrow h''(A) \cap \{y \in L_2:x_2 \le_{L_2} y\}
\ne \emptyset)
$$
\medskip

\noindent
iff

$$
(\forall x_2 \in L_2)(\{y \in L_1:\neg(x_2 \le_{L_2} h(y))\} = \emptyset 
\text{ mod id}_{L_1}]
$$
\mn
iff

$$
(\forall x_2 \in L_2)(\exists x_1 \in L_1)(\forall y \in L_1)[\neg(x_2 
\le_{L_2} h(y)) \rightarrow \neg x_1 \le_{L_1} y]
$$
\mn
iff

$$
(\forall x_2 \in L_2)(\exists x_1 \in L_1)(\forall y \in L_1)
[x_1 \le_{L_1} y \rightarrow x_2 \le_{L_2} h(y)].
$$
\mn
5) Letting $h_1$ be the identity map on $L_1$ by part (3) we get
id$_{L_2} \le_{\text{RK}} \text{ id}_{L_1}$; choose $h_2:L_2
\rightarrow L_1$ which extends $h_1$, by part (4) we get id$_{L_1}
\le_{\text{RK}} \text{ id}_{L_2}$, together we are done. \nl
6) Easily.  \hfill$\square_{\scite{3.5}}$\margincite{3.5}
\enddemo
\bigskip

\demo{\stag{3.6} Fact}  1) For any ideal $J$ (such that (Dom$(J)) \notin J$),
letting $J_1 = \text{ id}_{(J,\subseteq)}$, we have
\mr
\widestnumber\item{(iii)}
\item "{$(i)$}"  $J_1$ is a partial order ideal
\sn
\item "{$(ii)$}"  $|\text{Dom}(J_1)| = |J| \le 2^{|\text{Dom}(J)|}$
\sn
\item "{$(iii)$}"  $J \le_{\text{RK}} J_1$
\sn
\item "{$(iv)$}"  if $J$ is $\lambda$-complete, then $(J,\subseteq)$ is
$\lambda$-directed hence $J_1$ is $\lambda$-complete
\sn
\item "{$(v)$}"  dens$(J,\subseteq) = \text{ dens}(J_1,\subseteq)$
\sn
\item "{$(vi)$}"  dens$(J,\subseteq) \le |J| \le 2^{|\text{Dom}(J)|}$.
\ermn
2) For every ideal $J$ and 
dense $X \subseteq J$ we can use id$_{(X,\subseteq)}$ and get
the same conclusions. \newline
3) For every ideal $J$ there is a directed order $L$ such that:

$$
J \le_{\text{RK}} \text{ id}_L, \text{ dens}(J) = \text{ dens}(L)
$$
\medskip

\noindent
and: for every $\lambda$ if $J$ is $\lambda$-complete, then so is
id$_L$.
\enddemo
\bigskip

\demo{Proof}  Least trivial is (1)(iii), let $h:J \rightarrow \text{ Dom}(J)$
be such that $h(A) \in (\text{Dom}(J)) \backslash A$ (exists as (Dom$(J))
\notin J$).  Let $J_1 = \text{ id}_{(J,\subseteq)}$, so $h$ is a
function from Dom$(J_1)$ into Dom$(h)$ and we shall prove that $h$
exemplifies the desired conclusion $J \le_{\text{RK}} J_1$ by
Definition \scite{3.2}(1)

Assume toward contradicion that 
$X \subseteq \text{ Dom}(J_1) = J,X \notin J_1$ and $A =: h''(X)$ belongs
to $J$.  So
$Y =: \{B \in J:\neg (A \subseteq B)\} \in \text{ id}_{(J,\subseteq)} =
J_1$ (by the definition of $J_1 = \text{ id}_{(J,\subseteq)}$) 
hence (as $X \notin J_1$)
for some $B \in X$ we have $B \notin Y$ hence by definition
$A \subseteq B$, so 
$h(B) \in h''(X) = A$ contradicting the
choice of $h(B)$ (as $A \subseteq B$). \hfill$\square_{\scite{3.6}}$\margincite{3.6}
\enddemo
\bigskip

\remark{\stag{3.7} Remark}  So we can replace a family of ideals by a
family of directed quasi orders without changing much the relevant 
invariants 
such as completeness or density as long as we do not mind adding ``larger"
ones in the appropriate sense.
\endremark
\bigskip

\demo{\stag{3.8} Conclusion}  For any family of ideals $\Bbb I$ there is 
a family of ${\Cal L}$ of quasi order such that:
\mr
\item "{$(i)$}"  $\Bbb I \le_{\text{RK}} \{\text{id}_{(L,<)}:(L,<) \in
{\Cal L}\}$
\sn
\item "{$(ii)$}"  $|{\Cal L}| \le |\Bbb I|$
\sn
\item "{$(iii)$}"  $\sup\{|L|:(L,<) \in {\Cal L}\} = \sup\{|J|:J \in 
\Bbb I\} \le \sup\{2^{\text{Dom}(J)}:J \in \Bbb I\}$
\sn
\item "{$(iv)$}"  $\sup\{\text{dens}(L,<):(L,<) \in {\Cal L}\} = \sup
\{\text{dens}(J,\subseteq):J \in \Bbb I\}$
\sn
\item "{$(v)$}"  $\Bbb I$ is $\lambda$-complete \ub{iff} every 
$(L,<) \in {\Cal L}$ is $\lambda$-directed.
\endroster
\enddemo
\bigskip

\demo{Proof}  Easy.
\enddemo
\bigskip

\definition{\stag{3.9} Definition}  For a 
forcing notion $\Bbb Q$, satisfying the
$\kappa$-c.c., a $\Bbb Q$-name $\underset\tilde {}\to L$ of a 
quasi order with Dom$(\underset\tilde {}\to L) \in \bold V$ 
for notational simplicity given (i.e., is not just a $\Bbb Q$-name); 
let $L^* = \text{ ap}_\kappa(\underset\tilde {}\to L) = \text{
ap}_\kappa(\underset\tilde {}\to L,\Bbb Q)$ be the following quasi order

$$
\text{Dom}(L^*) = \{a:a \subseteq \text{ Dom}(\underset\tilde {}\to L)
\text{ and } |a| < \kappa\}
$$

$$
a \le^* b \text{ iff } \Vdash_{\Bbb Q} ``(\forall y \in a)(\exists x \in b)
[\underset\tilde {}\to L \models y < x]"
$$
\mn
(this is a quasi-order only, e.g. maybe $a \le^* b \le^* a$ but
$a \ne b$).
\enddefinition
\bigskip

\proclaim{\stag{3.10} Claim}  1) For a forcing notion $\Bbb Q$ satisfying the 
$\kappa$-c.c. and a 
$\Bbb Q$-name $\underset\tilde {}\to L$ of a $\lambda$-directed
quasi order (with ${\text{\rm Dom\/}}(\underset\tilde {}\to L) 
\in \bold V$ given,
not just a $\Bbb Q$-name, for simplicity) such
that $\lambda \ge \kappa$ we have:
\mr
\widestnumber\item{$(iii)$}
\item "{$(i)$}"  ${\text{\rm ap\/}}_\kappa(\underset\tilde {}\to L)$ 
is $\lambda$-directed
quasi order (in $\bold V$ and also in $\bold V^{\Bbb Q}$)
\sn
\item "{$(ii)$}"  $|{\text{\rm ap\/}}_\kappa(\underset\tilde {}\to L)| \le
|{\text{\rm Dom\/}}(\underset\tilde {}\to L)|^{< \kappa}$ 
\sn
\item "{$(iii)$}"  $\Vdash_{\Bbb Q} 
``{\text{\rm id\/}}_{{\underset\tilde {}\to L}[G]}
\le_{\text{RK}} { \text{\rm id\/}}_{\text{ap}_\kappa({\underset\tilde
{}\to L})}"$
\ermn
2) For a forcing notion $\Bbb Q$ satisfying the $\mu$-c.c., the local
$\underset\tilde {}\to \kappa$-c.c. and a $\Bbb Q$-name $\underset\tilde {}\to L$
of a $\lambda$-directed quasi order such that $\lambda \ge \kappa$ we
can find ${\Cal L} \in \bold V$ such that
\mr
\widestnumber\item{$(iii)$}
\item "{$(i)$}"  ${\Cal L}$ is a family of $< \mu \, \lambda$-directed quasi
orders (in $\bold V$ and also in $\bold V^{\Bbb Q}$)
\sn
\item "{$(ii)$}"  for each $L \in {\Cal L}$ for some $\theta$ we have
$|L| \le \theta^{< \kappa}$ and $\nVdash 
``|{\text{\rm Dom\/}}(\underset\tilde {}\to L)| \ne \mu,
\underset\tilde {}\to \kappa = \kappa"$
\sn
\item "{$(iii)$}"  $\Vdash_{\Bbb Q} 
``{\text{\rm id\/}}_{\underset\tilde {}\to L[G]}
\le_{\text{RK}} { \text{\rm id\/}}_L"$ for some $L \in {\Cal L}$.
\ermn
3) In (2), of course, we can replace $\Bbb I$ by a $\lambda$-complete
$\le_{\text{RK}}$-upper bound.
\endproclaim
\bigskip

\demo{Proof}  1) We leave (i), (ii) to the reader.  We check (iii).  Let
$G \subseteq \Bbb Q$ be generic over $\bold V$, and in $\bold V[G]$ 
we define a function $h$
from ap$_\kappa(\underset\tilde {}\to L)$ 
to Dom$(\underset\tilde {}\to L[G])$ as follows:
$h(a)$ will be an element of Dom$(\underset\tilde {}\to L[G])$ such that

$$
(\forall x \in a)[\underset\tilde {}\to L[G] \models ``x \le h(a)"].
$$
\mn
It exists by the ``$\lambda$-directed", ``$\lambda \ge \kappa$" assumptions.
We can now easily verify the condition in \scite{3.5}(4). \nl
2) Let $\{p_i:i < i^*\}$ be a maximal antichain of $\Bbb Q$ such that:
$p_i \Vdash ``\text{Dom}(\underset\tilde {}\to L)$ has cardinality $\mu_i,
\underset\tilde {}\to \kappa 
= \kappa_i"$ and $\Bbb Q_{\ge p_i}$ satisfies the
$\kappa_i$-c.c and \wilog \, $p_i \Vdash ``\text{Dom}(\underset\tilde {}\to
L) = \mu_i"$.  
Let $\Bbb Q_i = \Bbb Q_{\ge p_i}$ and ${\underset\tilde {}\to L_i}$ be
$\underset\tilde {}\to L$ restricted to $\Bbb Q_i$ and lastly let
${\Cal L} = \{\text{ap}_{\kappa_i}({\underset\tilde {}\to L_i},\Bbb
Q_i):i < i^*\}$, by part (1) we have $p_i \Vdash_{\Bbb Q}
``\text{id}_{\underset\tilde {}\to L[{\underset\tilde {}\to G_{\Bbb
Q}}]} \le_{\text{RK}}$; ap$_{\kappa_i}
({\underset\tilde {}\to L_i},{\Bbb Q}_i)"$ by $\{p_i:i < i^*\}$ is a maixmal antichain of $\Bbb
Q$ hence $\Vdash_{\Bbb Q} ``\text{id}_{\underset\tilde {}\to
L[{\underset\tilde {}\to G_{\Bbb Q}}]} \le_{\text{RK}}
\{\text{ap}_{\kappa_i}({\underset\tilde {}\to L_i},\Bbb Q_i):i <
i^*\}$. \nl
A base $|\text{ap}_{\kappa_i}({\underset\tilde {}\to L_i},\Bbb Q_i)|
\le |\text{Dom}({\underset\tilde {}\to L_i})|^{< \kappa_i} = \mu^{<
\kappa_i}_i$ and ap$_{\kappa_i}({\underset\tilde {}\to L_i},\Bbb Q_i)$
is $\lambda$-directed.  Together 
we are done as necessarily $i^* < \mu$. \nl
3) Easy by \scite{3.13}(3),(4) below.  \hfill$\square_{\scite{3.10}}$\margincite{3.10}
\enddemo
\bigskip

\demo{\stag{3.11} Conclusion}  1) 
Suppose $\Bbb Q$ is a forcing notion satisfying
the local $\kappa$-c.c., ${\underset\tilde {}\to {\Bbb I}_1}$ 
a $\Bbb Q$-name of 
a family of $\lambda$-complete filters and $\lambda \ge \kappa$. 
\underbar{Then} there is, (in $\bold V$), a family $\Bbb I_2$ of 
$\lambda$-complete filters such that:
\medskip
\roster
\widestnumber\item{$(iii)$}
\item "{$(i)$}"  $\Vdash_{\Bbb Q} ``{\underset\tilde {}\to {\Bbb I}_1} 
\le_{\text{RK}} \Bbb I_2"$
\sn
\item "{$(ii)$}"  if $\Bbb Q$ satisfies the $\mu$-c.c. then
$|\Bbb I_2| = |{\underset\tilde {}\to {\Bbb I}_1}| + \mu$ i.e.
$\Vdash_{\Bbb Q} ``|\Bbb I_2| = |\Bbb I_1|^V + \mu"$
\sn
\item "{$(iii)$}"  $\sup\{|\text{Dom}(J)|:J \in \Bbb I_2\} = 
\sup\{(2^\mu)^{< \kappa}:\text{some } q \in \Bbb Q$ \newline

$\qquad \qquad \qquad \qquad \qquad \qquad \qquad \qquad$ forces that some 
$J \in \Bbb I_1$ 
\newline

$\qquad \qquad \qquad \qquad \qquad \qquad \qquad \qquad$ has domain of 
power $\mu\}$.
\ermn
2) If in $\bold V^{\Bbb Q}$ the set $\Bbb I_1$ 
has the form $\{\text{id}_{(L,<)}:(L,<) \in {\Cal L}_1\}$; 
i.e., is a family of quasi order ideals, \underbar{then} in (iii) we can have
\mr
\item "{$(iii)'$}"   $\sup\{|\text{Dom}(J)|:J \in \Bbb I_2\} = 
\sup\{\mu^{< \kappa}:\text{some } q \in \Bbb Q$ \newline

$\qquad \qquad \qquad \qquad \qquad \qquad \qquad \quad$ 
force some $(L,<) \in {\Cal L}$ has power $\mu\}$
\sn
\item "{$(iv)$}"  $\Bbb I_2$ is a family of quasi order ideals.
\endroster
\enddemo
\medskip

\demo{Proof}  Easy.
\enddemo
\bigskip

\remark{\stag{3.12} Remark}  The aim of \scite{3.10}, \scite{3.11} is the
following:  we will consider iterations $\langle \Bbb P_i,
{\underset\tilde {}\to {\Bbb Q}_i}:i < \alpha \rangle$ where
$\Vdash_{{\Bbb P}_i} ``{\underset\tilde {}\to {\Bbb Q}_i}$ satisfies 
UP$({\underset\tilde {}\to {\Bbb I}_i})"$, but 
${\underset\tilde {}\to {\Bbb I}_i}$ may 
not be a subset of the ground model $\bold V$.  
Now \scite{3.11} gives us a good
$\le_{\text{RK}}$-bound $\Bbb I'_i$ of $\Bbb I_i$ in $\bold V$, and we can 
prove (under suitable assumptions) that $\Bbb P_\alpha$ will satisfy the 
UP$(\dsize \bigcup_{i < \alpha} \Bbb I'_i)$.
\endremark
\bigskip

\definition{\stag{3.13} Definition}  1) For a family ${\Cal L}$ of
directed quasi
order let the $(< \kappa)$-closure $c \ell_\kappa({\Cal L})$ be

$$
{\Cal L} \cup \biggl\{ \dsize \prod_{i < \alpha} L_i:L_i \in {\Cal L}
\text{ for } i < \alpha \text{ and } \alpha < \kappa \biggr\}
$$
\mn
(the partial order on $\dsize \prod_{i < \alpha} L_i$ is natural). \newline
2) ${\Cal L}$ is $(< \kappa)$-closed if for any $\alpha < \kappa$ and $L_i
\in {\Cal L}$ for $i < \alpha$ \underbar{there is} $L \in {\Cal L}$ such
that $i < \alpha \Rightarrow L_i \le_{\text{RK}} L$.
\enddefinition
\bigskip

\proclaim{\stag{3.14} Claim}  1) If ${\Cal L}$ is $\lambda$-complete,
\underbar{then} $c \ell_\kappa({\Cal L})$ is $\lambda$-complete. \newline
2) $L_j \le_{\text{RK}} \dsize \prod_{i < \alpha} L_i$ for $j < \alpha$.
\newline
3) If $\kappa$ is regular, \ub{then} 
$c \ell_\kappa({\Cal L})$ is $\kappa$-directed
under $\le_{\text{RK}}$ and is $(< \kappa)$-closed.
\endproclaim
\bigskip

\demo{Proof}  Use \scite{3.5}(4).
\enddemo
\newpage

\head {\S4 UP reintroduced} \endhead  \resetall \sectno=4
\bn
The reader may concentrate on the case $\bold S = \{\aleph_1\}$.
\demo{\stag{4.0} Convention}  $\Bbb I$ will be a set of quasi order ideals,
i.e. $\Bbb I = \Bbb I_{\Cal L}$.
\enddemo
\bigskip

\definition{\stag{4.1} Definition}  Fix $\Bbb I,\bold S,\bold W$ assuming
\medskip
\roster
\item "{$(*)(a)$}"  $\Bbb I$ a family of $\aleph_2$-complete quasi-order
ideals
\sn
\item "{$(b)$}"  $\bold S$ is a set of regular cardinals to which $\aleph_1$
belongs
\sn
\item "{$(c)$}"  $\bold W \subseteq \omega_1$ is stationary.
\endroster
\medskip

\noindent
1) We say $\bar N$ is an $\Bbb I$-tagged tree of models (for $\chi$ or for
$(\chi,x)$) \ub{if} there is an $\Bbb I$-tagged tree $(T,\bold I)$ such that
$\bar N = \langle N_\eta:\eta \in (T,\bold I) \rangle$ satisfies 
the following:
\medskip
\roster
\item "{$(a)$}"  for $\eta \in T$ we have $N_\eta \prec ({\Cal H}(\chi),\in,
<^*_\chi)$ is a countable model
\sn
\item "{$(b)$}"  $\Bbb I \in N_{\langle \rangle}$ and $x \in
N_{\langle \rangle}$ (if $x$ is present, we can use $\Bbb I$ as a predicate)
\sn
\item "{$(c)$}"  $\eta \triangleleft \nu \in T$ implies $N_\eta \prec
N_\nu$
\sn
\item "{$(d)$}"  for $\eta \in T$ we have $\eta \in N_\eta$ and $\bold I_\eta
\in N_\eta \cap \Bbb I$.
\endroster
\medskip

\noindent
Whenever we have such an $\Bbb I$-tagged tree $\bar N$ of models, we write
$N_\eta = \dsize \bigcup_{k < \omega} N_{\eta \restriction k}$ for each
$\eta \in \text{ lim}(T)$.  
 \newline
1A)  $\bar N = \langle N_\eta:\eta \in (T,\bold I) \rangle$ is a tagged 
tree of models if this occurs for some $\aleph_2$-complete $\Bbb I$.  
If $x$ is not mentioned, we assume it contains all necessary information, 
in particular $\Bbb I,\bold S,\bold W$. \newline
1B) In part (1) we say ``weak $\Bbb I$-tagged tree of models" 
\ub{if} we replace
clause (d) by
\mr
\item "{$(d)^-(i)$}"  for $\eta \in T \cap \text{ Dom}(\bold I)$ 
we have $\bold I_\eta \in N_\eta \cap \Bbb I$
\sn
\item "{$(ii)$}"  if $\eta \in \text{ Split}(T,\Bbb I)$ then for some
one-to-one function $f \in N_\eta$ with domain $\supseteq \text{
Suc}_T(\eta)$ and range $\subseteq \text{ Dom}(\bold I_\eta)$ we have 
$\nu \in \text{ Suc}_T(\eta) \Rightarrow f(\nu) \in N_\nu$.
\ermn
We say weaker if we omit clause $(d)^-(ii)$.  In all the definitions
below we can use this version (i.e., adding weak/weaker and replacing
$(d)$ by $(d)^-/(d)^-(i)$
\nl
2) We say $\bar N$ is truely $\Bbb I$-suitable (tagged tree of models)
\ub{if} clauses $(a)-(d)$ and:
\medskip
\roster
\item "{$(e)$}"  $\Bbb I \in N_{\langle \rangle}$ and if $\eta \in T$ and
$I \in \Bbb I \cap N_\eta$ \ub{then} the set
\endroster

$$
\{\nu \in T^{[\eta]}:\nu \in \text{ split}(T) \text{ and } \bold I_\nu =
I \in N_\eta\}
$$
\medskip

\noindent
contains a front of $T^{[\eta]}$. 
So ``$\bar N$ is truly $\Bbb I$-suitable tree (of models)" does not imply
``$\bar N$ is an $\Bbb I$-tagged tree of models" as possibly $\Bbb I
\notin N_{<>}$.  \newline
3) We say $\bar N$ is $\Bbb I$-suitable (a tagged tree of models)
\ub{if} clauses $(a)-(d)$ and:
\medskip
\roster
\item "{$(e)^-$}"  if $\eta \in T$ and $I \in \Bbb I \cap N$, \ub{then} 
the set
$$
\{\nu \in T^{[\nu]}:\nu \in \text{ split}(T) \text{ and } I
\le_{\text{RK}} \bold I_\nu \in N\}
$$
contains a front of $T^{[\eta]}$.
\ermn
4) We say $\bar N$ is $\lambda$-strictly $(\Bbb I,\bold W)$-suitable \ub{if}
$\bar N$ is $\Bbb I$-suitable and in addition
\mr
\item "{$(b)^+$}"  $\Bbb I \in N_{\langle \rangle},\bold W \in N$
and $x \in N_{\langle \rangle}$
\sn
\item "{$(f)$}"  one of the following cases holds:
{\roster
\itemitem{ $(i)$ }   $\lambda = \aleph_2$ and for some $\delta \in
\bold W$, for all $\eta \in T$ we have: $N_\eta \cap \omega_1 = \delta$
\sn
\itemitem{ $(ii)$ }   $\lambda = \text{ cf}(\lambda) > \aleph_1,
N_{<>} \cap \omega_1 \in \bold W$ and $\eta
\triangleleft \nu \in T \Rightarrow N_\eta <_\lambda N_\nu$ (i.e.,
$N_\eta \cap \lambda$ is an initial segment of $N_\nu \cap \lambda$ \nl
(note: if $\lambda = \mu^+$ this means $N_\eta \cap \mu = N_\nu \cap \mu$
so no contradiction with the case $\lambda = \aleph_2$)
\sn
\itemitem{ $(iii)$ }  $\lambda = \aleph_1$ and we demand nothing. 
\endroster}
If we omit $\lambda$, we mean $\lambda = \aleph_2$.  So $\lambda = 
\text{ cf}(\lambda) > \aleph_1$ implies $\eta \in T \Rightarrow N_\eta 
\cap \omega_1 = N_{<>} \cap \omega_1 \in \bold W$.  If we omit ``strictly" we
demand only $(b)^+$.
\endroster
\medskip
\noindent
5) We say $\bar N$ is $\bold S$-strictly $(\Bbb I,\bold W)$-suitable, \ub{if}
in addition to clauses (a) - (d),(e)$^-$ we have:
\medskip
\roster
\item "{$(b)^+$}"  $\Bbb I \in N_{\langle \rangle},\bold W \in
N_{\langle \rangle},\bold S \in N_{\langle \rangle}$ and $x \in
N_{\langle \rangle}$
\sn
\item "{$(g)$}"  for all $\nu \in T$ and $\lambda \in \bold S \cap N_\delta$
there is $\delta_\lambda,\nu < \lambda$ such that \nl
$(\forall \eta \in T)[\nu \trianglelefteq \eta \in T
\Rightarrow \sup(N_\eta \cap \lambda) = \delta_\lambda]$.
\ermn
6) We say $\bar N$ is $\lambda^+$-uniformly $(\Bbb I,\bold
W)$-suitable \ub{if} it is $(\Bbb I,\bold W)$-suitable and
\mr
\item "{$(b)^+$}"  $\Bbb I \in N_{\langle \rangle},\bold W \in 
N_{\langle \rangle},\lambda \in N_{\langle \rangle}$ and $x \in
N_{\langle \rangle}$
\sn
\item "{$(f)'$}"  for some $a \in [\lambda]^{\aleph_0}$ for every
$\eta \in \text{ lim}(T)$ we have
$$
N_\eta \cap \lambda = a
$$
and
$$
\text{if } \lambda^+ = \aleph_2 \text{ then } a \in \bold W.
$$
\ermn
6A) We say $\bar N$ is $\bold S$-uniformly $(\Bbb I,\bold W)$-suitable
\ub{if} it is $(\Bbb I,\bold W)$-suitable and
\mr
\item "{$(b)^+$}"  $\Bbb I \in N_{\langle \rangle},\bold W \in
N_{\langle \rangle},\bold S \in N_{\langle \rangle}$ and $x \in
N_{\langle \rangle}$
\sn
\item "{$(g)'$}"  for every $\eta \in T,\lambda \in \bold S \cap N_\eta$
for some $a \in [\lambda]^{\aleph_0}$ for every $\nu$ satisfying $\eta 
\triangleleft \nu \in \text{ lim}(T)$ we have
$$
N_\nu \cap \lambda = a.
$$
\ermn
7) In 4), 5) if we add "truely" if (e)$^-$ is replaced by (e). \newline
8) If ${\underset\tilde {}\to {\bold S}}$ is a $\Bbb P$-name then 
in the clauses
above we mean $\bold S^* = \{\lambda:\,\,\Vdash_{\Bbb P} ``\lambda \notin 
{\underset\tilde {}\to {\bold S}}"\}$. \nl
9) $\Bbb I^{[\lambda]} = 
\{I \in \Bbb I:I \text{ is } \lambda \text{-complete}\}$.
\enddefinition
\bigskip

\definition{\stag{4.2} Definition}  1) We say $N \prec ({\Cal H}(\chi),\in,
<^*_\chi)$ is ${\underset\tilde {}\to {\bold S}}$-strictly or
$\lambda$-strictly $(\Bbb I,\bold W)$-suitable model if there is an
${\underset\tilde {}\to {\bold S}}$-strictly or a $\lambda$-strictly $(\Bbb I,
\bold W)$-suitable $\bar N$ such that $N = N_{\langle \rangle}$, (see 
\scite{4.1}(5),(9), and see \scite{4.1}(8), it applies).  
We can add ``truely". \nl
2) We say $N$ is $(\Bbb I,\bold W)$-suitable if it is strictly
$(\Bbb I,\bold W)$-suitable, that is $\aleph_2$-strictly $(\Bbb
I,\bold W)$-suitable.
\enddefinition
\bigskip

\definition{\stag{4.3} Definition}   In Definitions \scite{4.1},
\scite{4.2} we may omit $\bold W$ when it is 
$\omega_1$, and may omit $\bold S$ when $\bold S = \{\aleph_2\}$,
for \scite{4.1}(5) and \scite{4.1}(6),
\scite{4.1}(6A).  We may replace $\bold S$ by $*$ if $\bold S = U$Reg.
Let $\eta \in (T,\bold I)$ means $\eta \in T$ and we write $T$ 
when $\bold I$ is clear.
\enddefinition
\bigskip

\proclaim{\stag{4.8} Claim}  Assume
\mr
\item "{$(i)$}"  ${\underset\tilde {}\to {\bold S}},
\Bbb P,\Bbb I,\bold W,x \in
{\Cal H}(\chi)$ and $\bold S^*$ is as in \scite{4.1}(8) \nl
(i.e., $\bold S^* = \{\theta:\theta = \aleph_1 \text{ or } \aleph_1 \le
\theta = \text{ cf}(\theta) \le |\Bbb P|$ and $\nVdash_{\Bbb P} ``\theta \in
{\underset\tilde {}\to {\bold S}}"\} \cup \{\aleph_1\}$, \nl
e.g., ${\underset\tilde {}\to {\bold S}} = \{\aleph_1\}$), and
\sn
\item "{$(ii)$}"  $\Bbb I$ 
is a $\aleph_2$-complete or for each $I \in \Bbb I,
\kappa \in \bold S,I$ is $\kappa$-indecomposable.
\ermn
\ub{Then} there is an $\bold S$-strictly, truely $(\Bbb I,\bold W)$-suitable 
tree $\bar N$ with $x \in N_{\langle \rangle}$.
\endproclaim
\bigskip

\demo{Proof}  We will construct this tree in three steps: first we find a
suitable tree, then we thin it out to be a $\bold S^*$-uniformly suitable 
tree, then we blow up the models to make it $\bold S^*$-strict.  For 
notational simplicity let $\bold S^* = \{\aleph_1\}$ 
(the reader can check the others).
\mn
\underbar{First Step}:  An easy bookkeeping argument (to ensure 
\scite{4.1}(1)(e)) yields a truly
$(\Bbb I \cup \{J^{\text{bd}}_{\omega_1}\})$-suitable tree
$\langle N_\eta:\eta \in (T,\bold I) \rangle$ satisfying $\nu
\triangleleft \eta \Rightarrow \sup(N_\nu \cap \omega_1) < \sup(N_\eta
\cap \omega_1)$ such that $N_\eta \prec ({\Cal H}(\chi),\in,<^*_\chi)$
and $x \in N_{<>}$; so for $\eta \in \text{ lim}
(T)$ we let $N_\eta = \dsize \bigcup_{\ell < \omega} N_{\eta \restriction
\ell}$.  

Moreover we can we get that 
for all $\eta \in \text{ lim}(T)$, for each $I \in (\Bbb I 
\cap N_\eta) \cup \{J^{\text{bd}}_{\omega_1}\}$, there are 
infinitely many $k$ such 
that $\eta \restriction k \in \text{ split}(T,\bold I)$ 
and $\bold I_{\eta \restriction k} = I$ and Suc$_T(\eta \restriction k) =
\{\eta \char 94 \langle x \rangle:x \in \text{ Dom}(I)\}$.
\medskip

\noindent
\underbar{Second Step}:  Define $H:T \rightarrow \omega_1$ by $H(\eta) = \sup
(N_\eta \cap \omega_1) < \omega_1$.  Apply \scite{2.13} to get a subtree
$T'$ and a limit ordinal $\delta \in \bold W \subseteq \omega_1$ such that 
clauses (a) - (d) of \scite{2.13} hold for $\delta$.  By clause (d) of \scite{2.13}, for 
all $\eta \in T',N_\eta \cap \omega_1 \subseteq \delta$.  Let $\delta_0 < 
\delta_1 < \ldots \dsize \bigcup_n \delta_n = \delta$, and let

$$
\align
T'' = \biggl\{ \eta \in T':&\text{ for each } 
\forall k < \ell g(\eta),\text{ if Suc}_T(\eta
\restriction k) = \{\eta \restriction k \char 94 \langle \alpha \rangle:
\alpha < \omega_1\} \\
  &\text{ so Suc}_{T'}(\eta \restriction k) = \{\eta \restriction k \char 94
\langle \alpha \rangle:\alpha < \delta\}) \\
  &\text{ then } \eta(k) = \delta_k \biggr\}.
\endalign
$$
\medskip

\noindent
Clearly $T''$ will be $\aleph_1$-uniformly suitable; i.e. \newline
$\eta \in \text{ lim}(T) \Rightarrow N_{\eta,\ell} \cap \omega_1 = \delta$.
\mn
\underbar{Third Step}:  For $\eta \in T_2$, 
let $N'_\eta =$ the Skolem Hull of
$N_\eta \cup \delta$ in $({\Cal H}(\chi),\in,<^*_\chi)$.  
So $N'_\eta \cap \omega_1 \supseteq \delta$.
Conversely, let $\nu \in \text{ lim}(T_2),\eta \triangleleft \nu$, then
$N_\eta \cup \delta \subseteq N_\nu$, so $N'_\eta \subseteq N_\nu$ hence
$N'_\eta \cap \omega \subseteq \delta$.  So $N'_\eta \cap \omega_1 = \delta$,
i.e. $\langle N'_\eta:\eta \in T \rangle$ is an $\aleph_1$-strict, $(\Bbb I,
{\underset\tilde {}\to {\bold S}},\bold W)$-tree of models (see Definition
\scite{4.1}(4)).

We claim that this tree is still truly suitable.  
Indeed, assume $\eta \in T_2,\nu
\in \text{ lim}(T_2),\eta \trianglelefteq \nu$ and $I \in \Bbb I 
\cap N'_\eta$.
Then for some $\alpha < \delta,I$ is in the Skolem hull of $N_\eta \cup
\alpha$.  Let $k < \omega$ be such that $\alpha \in N_{\nu \restriction k}
\cap \omega_1$ and $k \le \ell g(\eta)$.  
Then since $\langle N_\eta:\eta \in T_2
\rangle$ was suitable, there is $\ell \ge k$ such that 
$\bold I_{\nu \restriction \ell} = I$.  So $\langle N'_\eta:\eta \in T_2
\rangle$ is also $(\Bbb I,{\underset\tilde {}\to {\bold S}},
\bold W)$-suitable. \hfill$\square_{\scite{4.8}}$\margincite{4.8}
\enddemo
\bigskip

\demo{\stag{4.4} Fact}  Assume $\Bbb I' \le_{\text{RK}} \Bbb I$, where 
$\Bbb I,\Bbb I'$ are families of ideals. \newline
1) If $\langle N_\eta:\eta \in (T,\bold I) \rangle$ is a $\Bbb I$-suitable 
tree and $\Bbb I' \in N_{\langle \rangle}$, \ub{then}
$\langle N_\eta:\eta \in (T,\bold I) \rangle$ is also $\Bbb I'$-suitable. \nl
2) If $\langle N_\eta:\eta \in (T,\bold I) \rangle$ is $\Bbb I$-suitable and 
$\Bbb I' \in N_{\langle \rangle}$, \underbar{then} there is a tree $(T',
\bold I')$ satisfying the following for some $T'',f$:
\medskip
\roster
\item "{$(a)$}"  $T'' \subseteq T$ and $f$ is an isomorphism from
$T''$ onto $T'$, (i.e., is one to one onto preserving length and
$\triangleleft$ and its negation) and $\eta\in T'' \Rightarrow \bold
I''_\eta \le_{\text{RK}} \bold I_{f(\eta)},\bold I''_\eta \ne \bold I'_\eta$
\sn
\item "{$(b)$}"  $\langle N_\eta:\eta \in (T',\bold I') \rangle$ is truely
$\Bbb I'$-suitable when we let $I''_\eta = \bold I'_{f(\eta)}$
\sn
\item "{$(c)$}"  split$(T'',\bold I'') = 
T'' \cap \text{ split}(T,\bold I)$ if $\Bbb I' \equiv_{\text{RK}} \Bbb I$.
\sn
\ermn
3) We can weaken the hypothesis to $\Bbb I' \le_{\text{RK}} 
\Bbb I \cup \{\text{the trivial ideal}\}$.  The same holds in 
similar situations. \newline
4) In Definition \scite{4.3}, if $\bold S = \{\theta:\aleph_1 \le \theta =
\text{ cf}(\theta) \le \lambda\},\lambda = \text{ cf}(\lambda)$, 
\underbar{then} clause (f) of \scite{4.1}(4) (i.e., $\lambda^+$-strictly) and
clause (g) of \scite{4.1}(5) (i.e. the demand concerning $\bold S$,
i.e., $\bold S$-strictly) are equivalent.  Similarly \scite{4.1}(6),
\scite{4.1}(6A) are equivalent. \nl
5) In part (2), $\langle N_\eta:\eta \in (T'',\bold I) \rangle$ is a
weak $\Bbb I'$-tagged tree, truely $\Bbb I'$-suitable; moreover it is
enough to assume $\langle N_\eta:\eta \in (T,\bold I) \rangle$ is a
weak $\Bbb I$-suitable tree (see \scite{5.2}). \nl
6) If $\langle N_\eta:\eta \in (T,\bold I) \rangle$ is a weak $\Bbb
I$-tagged tree \ub{then} for some tree $T'$ and tree isomorphic $f$ from
$T'$ onto $T$ letting $\bold I' = \langle \bold I_{f(\eta)}:\eta \in
T' \rangle$ we have $\langle N_{f(\eta)}:\eta \in (T',\bold I')
\rangle$ is a $\Bbb I$-tagged tree.  All relevant properties are
preserved.  [Check, see \scite{5.2}.]
\enddemo
\bigskip

\demo{Proof}  1) Should be clear, as $\le_{\text{RK}}$ is transitive (as a
relation among ideals and also among families of ideals). \newline
2) For every $\eta \in Y =: \{\eta \in T:(\exists I' \in \Bbb I')
(I' \le \bold I_\eta)$
and $\eta \in \text{ split}(T,\bold I)\}$ pick an ideal $\bold I''_\eta
\in \Bbb I \cap N_\eta,\bold I''_\eta \le_{\text{RK}} \bold I_\eta$
such that:
for every $\nu \in T$, for every 
$I' \in \Bbb I \cap N_\nu$ the set $\{\eta \in
T^{[\nu]}:I' = \bold I''_\eta \text{ and } \nu \triangleleft \eta$ and
$\eta \in Y\}$ contains a front of $T^{[\nu]}$.  This can
be done using a bookkeeping argument.

Now define $T'$ as follows.  We choose by induction on $n$, a function
$f_n$ with domain $\subseteq T \cap {}^n\text{Ord}$, such that $\eta
\in \text{ Dom}(f_\eta) \Rightarrow f(\eta) \in N_\eta \cap
{}^n\text{Ord}$.  Let $f_0$ be the identity on $\{<>\}$.  Assume $f_n$
has been defined and $\eta \in \text{ Dom}(f_n)$, and we shall define
$f_{n+1} \restriction \text{ Suc}_T(\eta)$.  If 
$\eta \in T \backslash Y$, then Dom$(f_{n+1}) \cap \text{ Suc}_T(\eta) \subseteq
\text{ Suc}_T(\eta)$ is a singleton $\{\nu_\eta\}$ and let
$f_{n+1}(\nu_\eta) = f_n(\eta) \char 94 \langle (\nu_\eta(n) \rangle$.  
If $\eta \in T \cap Y$, then $\bold I'_\eta$ is already defined and it 
belongs to $N_\eta$.  Let $g_\eta$
be a witness for $\bold I'_\eta \le_{\text{RK}} \bold I_\eta$ and
stipulate Dom$(\bold I_\eta) \supseteq \{x:\eta \char 94 <x> \in
\text{ Suc}_T(\eta)\}$.  Now $g_\eta$
introduces an equivalence relation on Dom$(\bold I_\eta)$.  Let $A_\eta$ be a
selector set for this equivalence relation; i.e. $g_\eta \restriction A_\eta$
is 1-1 and has the same range as $g_\eta$.  Note that we can choose $g_\eta$
in $N_\eta$ as $\bold I_\eta,\bold I'_\eta \in N_\eta$ 
(whereas Suc$_T(\eta)$
does not necessarily belong to $N_\eta$) and then choose $A_\eta$ and
let $A'_\eta =: \{x \in A_\eta:(\exists y \in \text{
Suc}_T(\eta))[g_\eta(x)= g_\eta(y)]\}$ (so possibly
$A_\eta \notin N_\eta$).  Now for $x \in A'_n$ so $\eta \char 94
<x> \in \text{ Suc}_T(\eta)$ we let $f_{n+1}(\eta \char 94 <x>) =
f_n(\eta) \char 94 \langle g_n(x) \rangle$ so Dom$(f_{n+1}) \cap
\text{ Suc}_T(\eta) = \{\eta \char 94 <x>:x \in A'_\eta\}$.  Lastly, let $T' =
\cup\{\text{Rang}(f_n):n < \omega\},T'' = \cup\{\text{Dom}(f_n):n <
\omega\}$ and for $\eta \in \text{ Dom}(f_n),n < \omega$ let $\bold
I'_{f_n(\eta)} = \bold I''_\eta$ and $f = \cup\{f_n:n < \omega\}$.
Now check. \nl
3), 4)  Left to the reader.   \hfill$\square_{\scite{4.4}}$\margincite{4.4}
\enddemo
\bigskip

\definition{\stag{4.5A} Definition}  1) Let $\chi > \aleph_0,
\Bbb I \in {\Cal H}(\chi)$ a set of ideals and $\bold S \in {\Cal H}(\chi)$
a set of regular cardinals (or just limit ordinals) such that $\aleph_1 \in
\bold S$.  For $N$ a countable elementary submodel of ${\frak B} =
({\Cal H}(\chi),\in)$ (or ${\frak B}$ an expansion of $({\Cal
H}(\chi),\in)$ with countable
vocabulary) such that $\Bbb I,\bold S \in N$ we define Dp$(N) = 
\text{ Dp}_{\Bbb I}(N) = \text{ Dp}_{\Bbb I}(N,{\frak B}) = 
\text{ Dp}_{\Bbb I}(N,\bold S,{\frak B}) \in \text{ Ord } \cup \{\infty\}$,
by defining when Dp$(N) \ge \alpha$ for an ordinal $\alpha$, by induction
on $\alpha$:

$$
\align
\text{Dp}(N) \ge \alpha \text{ \ub{iff} } &N \text{ is as above and for
every } I \in \Bbb I \cap N \\
  &\text{ and for every } \beta < \alpha \text{ and } X \in I 
\text{ there is } M \text{ satisfying}: \\
  &(i) \quad \text{ Dp}(M) \ge \beta ( \text{hence } M \prec {\frak B}
\text{ is countable})  \\
  &(ii) \quad N \prec M \\
  &(iii) \quad \text{ sup}(M \cap \omega_1) = \text{ sup}(N \cap \omega_1)
\text{ moreover} \\
  &\qquad \quad \theta \in \bold S \cap N \Rightarrow \text{ sup}(N \cap
\theta) = \text{ sup}(M \cap \theta) \\
  &(iv) \quad M \cap \text{ Dom}(I) \backslash X \ne \emptyset.
\endalign
$$
\mn
2) We 
define Dp$'(N) = \text{ Dp}'_{\Bbb I}(N)$ by defining: $\text{Dp}'(N) \ge
\alpha$ iff $N$ is as above and for every $J \in \Bbb I \cap N$ and $\beta <
\alpha$ for some $I \in N \cap \Bbb I$ we have $J \le_{\text{RK}} I$ 
and for every $X \in I$ there is $M$ satisfying (i)-(iv) above.
\enddefinition
\bigskip

\proclaim{\stag{4.5B} Claim}  1)  In Definition \scite{4.5A}:
\mr
\item "{$(a)$}"  ${\text{\rm Dp\/}}(N) \in { \text{\rm Ord\/}} 
\cup \{\infty\}$ if well defined
\sn
\item "{$(b)$}"  if ${\text{\rm Dp\/}}(N) = 
\infty,I \in \Bbb I \cap N$ then we can find 
$Y \in I^+$ (i.e., $Y \subseteq { \text{\rm Dom\/}}(I),
Y \notin I$) and $\bar N =
\langle N_t:t \in Y \rangle$ such that:
{\roster
\itemitem{ $(i)$ }  ${\text{\rm Dp\/}}(N_t) = \infty$
\sn
\itemitem{ $(ii)$ }  $N \prec N_t$
\sn
\itemitem{ $(iii)$ }  ${\text{\rm sup\/}}(N \cap \omega_1) = 
{ \text{\rm sup\/}}(N_t \cap
\omega_1)$ moreover $\theta \in \bold S \cap N \Rightarrow \sup
(N \cap \theta) = \sup(N_t \cap \theta)$.
\endroster}
\ermn
2) If $\Bbb I_1 \le_{\text{RK}} \Bbb I_2$ and $\{\Bbb I_1,\Bbb I_2\} \in N$
and Dp$_{{\Bbb I}_2}(N) \ge \alpha$ then Dp$_{{\Bbb I}_1}(N) \ge \alpha$.
\nl
3) Dp$_{\Bbb I}(N) = \text{ Dp}'_{\Bbb I}(N)$.
\endproclaim
\bigskip

\demo{Proof}  Straightforward.
\enddemo
\bigskip

\proclaim{\stag{4.5C} Claim}  1) Let $N \prec ({\Cal H}(\chi),\in)$ be
countable, $\Bbb I \in N, N \cap \omega_1 \in {\Cal W}$.  \ub{Then}
\mr
\item "{$(a)$}"   $N$ is strictly 
$(\Bbb I,{\Cal W})$-suitable \ub{iff} ${\text{\rm Dp\/}}_{\Bbb I}(N) = \infty$
\sn
\item "{$(b)$}"  $N$ is $(\Bbb I,{\Cal W})$-suitable \ub{iff}
${\text{\rm Dp\/}}'_{\Bbb I}(N) = \infty$.
\ermn
2) Similarly with $\bold S$.
\endproclaim
\bigskip

\demo{Proof}  Easy.
\enddemo
\bigskip

\definition{\stag{4.6} Definition}  1) The forcing notion $\Bbb P$ satisfies
UP$^\ell_{\underset\tilde {}\to \lambda}
(\Bbb I,{\underset\tilde {}\to {\bold S}},\bold W)$ 
(note if $\ell \ne 2$ we may omit $\underset\tilde {}\to \lambda$)
(adopting the
conventions of \scite{4.1}(8); $\underset\tilde {}\to \lambda$ is a
purely decidable $\Bbb P$-name of a $\bold V$-cardinal) when: 
$\ell \in \{0,1,2\}$ and if $\chi$ is 
large enough and $\bar N$ is ${\underset\tilde {}\to {\bold S}}$-strictly 
$(\Bbb I,\bold W)$-suitable and $p \in N_{\langle \rangle} \cap \Bbb P$ and 
$P \in N_{\langle \rangle}$, of course, \ub{then} 
there is $q \in \Bbb P$ such that $p \le_{\text{pr}} q \in \Bbb P$ and:
\mr
\item "{$(a)$}"  if $\ell = 0$ then $q \Vdash ``N_{\langle \rangle} 
[{\underset\tilde {}\to G_{\Bbb P}}] \cap \omega_1 = N_{\langle \rangle} 
\cap \omega_1$ and sup$(N_{<>}[{\underset\tilde {}\to G_{\Bbb P}}] 
\cap \theta) = \text{ sup}(N_{<>} \cap \theta)$ 
if $\theta \in \underset\tilde {}\to S"$
\sn
\item "{$(b)$}"   if $\ell = 1$ then $q \Vdash$ ``for some $\eta \in
\text{ lim}(T)$ we have $N_\eta
[{\underset\tilde {}\to G_{\Bbb P}}] \cap \omega_1 = 
N_{\langle \rangle} \cap \omega_1$ and for every $ \theta \in
{\underset\tilde {}\to {\bold S}}$ we have
sup$(N_\eta[{\underset\tilde {}\to G_{\Bbb P}}] \cap \theta) =
\sup(N_\eta \cap \theta)$ where $N_\eta = \cup\{N_{\eta \restriction
\ell}:\ell < \omega\}$ and $\eta$ is not necessarily from $\bold V$"
\sn
\item "{$(c)$}"  if $\ell = 2$ then $q \Vdash ``N_{\langle \rangle} 
[{\underset\tilde {}\to G_{\Bbb P}}]$ is $({\underset\tilde {}\to {\bold S}}
\backslash {\underset\tilde {}\to \lambda'})$-strictly 
$(\Bbb I^{[\underset\tilde {}\to \lambda]},\bold W)$-suitable and \nl
sup$(N_{<>}[{\underset\tilde {}\to G_{\Bbb P}}] 
\cap \theta) = \sup(N_{<>} \cap
\theta)$ for every $\theta \in \underset\tilde {}\to S$ (in particular
$\aleph^{\bold V}_1$)" where
${\underset\tilde {}\to \lambda'} = 
\aleph_2^{{\bold V}[{\underset\tilde {}\to G_{\Bbb P}}]}$.
\ermn
2) If we omit $\ell$ we mean $\ell=0$. \nl
If $\bold W = \omega_1$ we may omit it.  We write $*$ instead of
${\underset\tilde {}\to {\bold S}}$ if \newline
${\underset\tilde {}\to {\bold S}} =
\{\lambda:\bold V[{\underset\tilde {}\to G_{\Bbb P}}] \models \lambda \in 
\text{ UReg}^{\bold V[G_{\Bbb P}]}\}$.  
If we omit ${\underset\tilde {}\to {\bold S}}$ we mean
$\{\aleph^{\bold V}_1\}$. \nl
3)  The forcing notion $\Bbb P$ 
satisfies UP$^4_{\underset\tilde {}\to \kappa,
\underset\tilde {}\to \lambda}(\Bbb I,{\underset\tilde {}\to {\bold S}},
\bold W)$ when $\underset\tilde {}\to \kappa,\underset\tilde {}\to \lambda$ 
are $(\Bbb P,\le_{\text{pr}})$-names of regular $\bold V$-cardinals and for 
some $x$ we have: \underbar{if} $\chi$ is large enough and 
$\bar N = \langle N_\eta:\eta \in (T,\bold I) \rangle$ is an
${\underset\tilde {}\to {\bold S}}$-strictly $(\Bbb I,\bold W)$-suitable 
tree of models for $(\chi,x)$ and $p \in \Bbb P \cap N_{\langle \rangle}$ and
$\langle \underset\tilde {}\to \kappa,\underset\tilde {}\to
\lambda,\Bbb P,\Bbb I \rangle 
\in N_{\langle \rangle}$, \underbar{then} for some 
$q,\underset\tilde {}\to T$ we have:
\mr
\item "{$(a)$}"  $p \le_{\text{pr}} q \in \Bbb P$
\sn
\item "{$(b)$}"  $q$ is $(\bar N \restriction \underset\tilde {}\to T,
\underset\tilde {}\to \kappa,\underset\tilde {}\to \lambda,\Bbb P)$-semi$_4$ 
generic (see below). 
\ermn
3A) If $\underset\tilde {}\to \kappa = \infty$ we can replace
$\underset\tilde {}\to T$ by $\underset\tilde {}\to \eta$ such that
$\underset\tilde {}\to T = \{\underset\tilde {}\to \eta \restriction n:n <
\omega\}$ (see below) so $\Vdash ``\underset\tilde {}\to \eta \in
\text{ lim}(T)"$ and then in clause (b) write 
$N_{\underset\tilde {}\to \eta}$ instead of $\bar N \restriction
\underset\tilde {}\to T$.  We then may omit $\kappa$. 
We may $\underset\tilde {}\to \lambda$ if $\underset\tilde {}\to
\lambda = \underset\tilde {}\to \kappa(\Bbb P)$, 
see Definition \scite{1.28}(5).
\nl
3B) We say that $q$ is $(\bar N \restriction \underset\tilde {}\to T,
\underset\tilde {}\to \kappa,\underset\tilde {}\to \lambda,\Bbb P)$-semi$_4$ 
generic  where $\bar N,\underset\tilde {}\to T,\underset\tilde {}\to \kappa,
\underset\tilde {}\to \lambda$ is as above \ub{if}:

$$
\align
q \Vdash_{\Bbb P} ``&(i) \quad \underset\tilde {}\to T \text{ is a subtree of }
T \text{ (so } \underset\tilde {}\to T \subseteq T \text{ is closed under} \\
  &\qquad \text{initial segments } <> \in \underset\tilde {}\to T,\eta \in
\underset\tilde {}\to T \Rightarrow \text{ Suc}_{\underset\tilde {}\to T}
(\eta) \ne \emptyset) \\
  &(ii) \quad N_\eta[{\underset\tilde {}\to G_{\Bbb P}}] \cap \omega_1 = 
N_{\langle \rangle} \cap \omega_1 \text{ for } \eta \in
\underset\tilde {}\to T \\
  &(iii) \quad \bar N[{\underset\tilde {}\to G_{\Bbb P}}] \restriction
\underset\tilde {}\to T \text{ has }
(\underset\tilde {}\to \lambda,\underset\tilde {}\to \lambda) \text{-covering
which means: if} \\
  &\qquad \qquad \underset\tilde {}\to \eta \text{ is an } \omega
\text{-branch of } \underset\tilde {}\to T \text{ and }
y \in N_{\underset\tilde {}\to \eta}[{\underset\tilde {}\to G_{\Bbb P}}] \\
  &\qquad \qquad \text{ is a set of } < 
\underset\tilde {}\to \lambda[{\underset\tilde {}\to G_{\Bbb P}}] \text{ ordinals} \\
  &\qquad \qquad 
\text{(if } \underset\tilde {}\to \lambda[{\underset\tilde {}\to
G_{\Bbb P}}] 
\text{ is not a cardinal, this means } \le |\underset\tilde {}\to \lambda
[{\underset\tilde {}\to G_{\Bbb P}}]|) \\
  &\qquad \qquad 
\text{\ub{then} for some } A \in \bold V \cap \dbcu_{\ell < \omega} 
N_{{\underset\tilde {}\to \eta} \restriction \ell} \text{ we have} \\
  &\qquad \qquad |A|^{\bold V} < 
\underset\tilde {}\to \lambda [{\underset\tilde {}\to G_{\Bbb P}}]
\text{ and } y \subseteq A \\
  &(iv) \quad \langle N_\eta[{\underset\tilde {}\to G_{\Bbb P}}]:
\eta \in (\underset\tilde {}\to T,\bold I) \rangle
\text{ is a strictly} \\
  &\qquad \qquad \Bbb I^{[\underset\tilde {}\to \kappa[
{\underset\tilde {}\to G_{\Bbb P}}]]} \text{-suitable tree of models}".
\endalign
$$
\mn
4)  We define UP$^3_{\underset\tilde {}\to \kappa,
\underset\tilde {}\to \lambda}(\Bbb I,{\underset\tilde {}\to {\bold S}},
\bold W)$ similarly replacing clause (b) by
the weaker
\mr
\item "{$(b)^-$}"  $q$ is $(N_{\underset\tilde {}\to \eta},
\underset\tilde {}\to \kappa,\underset\tilde {}\to \lambda,\Bbb P)$-semi$_3$
-generic, (see below).
\ermn
4A) If $\underset\tilde {}\to \kappa = \infty$ we can replace
$\underset\tilde {}\to T$ by $\underset\tilde {}\to \eta$ such that
$\underset\tilde {}\to T = \{\underset\tilde {}\to \eta \restriction n:n <
\omega\}$ so $\Vdash_{\Bbb P} 
``\underset\tilde {}\to \eta \in \text{ lim}(T)"$
and replace $\bar N \restriction \underset\tilde {}\to T$ by
$N_{\underset\tilde {}\to \eta}$.  We may omit $\underset\tilde {}\to
\lambda$ if it is $\infty$ (but see \scite{4.7}(0)).
We then may omit $\underset\tilde {}\to \kappa$ if it is $\infty$, too. \nl
4B) $q$ is $(N,\underset\tilde {}\to \kappa,
\underset\tilde {}\to \lambda,\Bbb P)$-semi$_3$-generic is defined
as in (4) only replacing clause (iii) in $\boxtimes$ of (3B) by
\mr
\item "{$(iii)^-$}"  $\bar N[{\underset\tilde {}\to G_{\Bbb P}}]
\restriction \underset\tilde {}\to T$ has $(\underset\tilde {}\to \lambda,
1)$-covering which means for every 
$y \in \bold V \cap N_{\underset\tilde {}\to \eta}
[{\underset\tilde {}\to G_{\Bbb P}}]$
for some $A \in \bold V \cap 
N_{\underset\tilde {}\to \eta}$ we have
$|A|^V < \underset\tilde {}\to \lambda
[{\underset\tilde {}\to G_{\Bbb P}}]$ and $y \in A$, recalling
$N_{\underset\tilde {}\to \eta}[{\underset\tilde {}\to G_{\Bbb P}}] =
\cup\{N_{\underset\tilde {}\to \eta \restriction
\ell}[{\underset\tilde {}\to G_{\Bbb P}}]:\ell < \omega\}"$.
\ermn
5) We allow to use ${\underset\tilde {}\to {\Bbb I}}$, a $\Bbb P$-name of an 
element of $\bold V$ as above \ub{if}:
\mr
\item "{$(a)$}"  it is decidable purely
\sn
\item "{$(b)$}"  if $q \in \Bbb P$ 
decides ${\underset\tilde {}\to {\Bbb I}} =
\Bbb I$ then $\Bbb P_{\ge q}$ satisfies 
UP$^\ell(\Bbb I,{\underset\tilde {}\to {\bold S}},\bold W)$.
\ermn
6) We say that $\Bbb P$ satisfies the UP$^5_{\underset\tilde {}\to \kappa,
\underset\tilde {}\to \lambda}(\Bbb I,{\underset\tilde {}\to {\bold S}},
\bold W)$ \ub{iff}
\mr
\item "{$(a)$}"  $\underset\tilde {}\to \kappa$ and 
$\underset\tilde {}\to \lambda$ are $\Bbb P$-names of $\bold
V$-cardinals such that $\Bbb P$
satisfies the $\underset\tilde {}\to \kappa$-c.c. purely locally
\sn
\item "{$(b)$}"  ${\underset\tilde {}\to {\Bbb I}}$ is a $\Bbb P$-name 
of a set which belongs to $\bold V$, it is a set of ideals and
${\underset\tilde {}\to {\Bbb I}}$ is decidable purely
\sn
\item "{$(c)$}"  is $\kappa$-complete if 
$\nVdash_{\Bbb P} ``\neg(\underset\tilde {}\to \kappa = \kappa \and I \in
{\underset\tilde {}\to {\Bbb I}})"$,
\sn
\item "{$(d)$}"  if $p \in \Bbb P$ forces ${\underset\tilde {}\to {\Bbb I}} =
\Bbb I,\Bbb I \subseteq \Bbb I',
\underset\tilde {}\to \lambda = \lambda,\underset\tilde {}\to \kappa =
\kappa$ and $\Bbb I' \backslash \Bbb I$ is a set of $\lambda$-complete ideals
\ub{then} for some $x$ 
{\roster
\itemitem{ $\boxtimes$ }  if $\bar N = \langle N_\eta:\eta \in (T,\bold I)
\rangle$ is strictly $\Bbb I'$-suitable, $x \in N_{<>}$, then for some
$q,\underset\tilde {}\to T$ we have: $p \le_{\text{pr}} q \in \Bbb P$
and $q$ is $(\bar N \restriction \underset\tilde {}\to T,\underset\tilde
{}\to \kappa,\underset\tilde {}\to \lambda,\Bbb P)$-semi$_5$-generic,
(see below).
\endroster}
\ermn
7) We define $\Bbb P$ satisfies UP$^5_{\underset\tilde {}\to \kappa}
({\underset\tilde {}\to {\Bbb I}},{\underset\tilde {}\to {\bold
S}},\bold W)$ as in part (6) but restrict
ourselves to $\Bbb I' = \Bbb I$.  [others?] \nl
8) Assume $\bar N = \langle N_\eta:\eta \in (T,\bold I) \rangle$ and
$\Bbb P \in N_{<>}$ satisfies UP$^\ell_{\underset\tilde {}\to \kappa}
({\underset\tilde {}\to {\Bbb I}},\bold W)$ and $\langle \Bbb P,
{\underset\tilde {}\to {\Bbb I}},\bold W,\ldots \rangle \in N_{<>}$.  We say
$q$ is $(\bar N,\underset\tilde {}\to \kappa,\Bbb P)$-semi$_5$-generic (for
$\bar N$ when not understood from the context) if:

$$
\align
q \Vdash ``&\text{for some } \underset\tilde {}\to T \text{ we have }
(T,{\underset\tilde {}\to {\Bbb I}}) \le^{\underset\tilde {}\to \kappa}
({\underset\tilde {}\to T'},\bold I), \text{ see \scite{2.4}, 
clause (f) and} \\
  &\eta \in {\underset\tilde {}\to T'} \Rightarrow N_\eta
[{\underset\tilde {}\to G_{\Bbb P}}] \cap \omega_1 
= N_{<>} \cap \omega_1 \text{ and} \\
  &\eta \in T' \and \mu \in {\underset\tilde {}\to {\bold S}}
\Rightarrow \sup(N_{\underset\tilde {}\to \eta}[G_{\Bbb P}] \cap \mu) = \sup
(N_{\underset\tilde {}\to \eta} \cap \mu)".
\endalign
$$
\mn
9) We write UP$^3({\underset\tilde {}\to {\Bbb I}},\bold S,\bold W)$ for
UP$^3_{\underset\tilde {}\to \kappa}
({\underset\tilde {}\to {\Bbb I}},\bold S, \bold W)$ where
$\underset\tilde {}\to \kappa$ is $\underset\tilde {}\to \kappa(\Bbb
P)$ see \scite{1.28}(5), \scite{1.29}(1).
\enddefinition
\bigskip

\definition{\stag{4.6A} Definition}[?]  We call $\Bbb I$ to be a name if
it is a name of an old family of ideals purely decidable.
\enddefinition
\bigskip

\proclaim{\stag{4.7} Claim}  0) In Definition \scite{4.6}(3)-(6), if
$\underset\tilde {}\to \lambda \ge \underset\tilde {}\to \kappa(\Bbb
P)$, \ub{then} the demand concerning $\underset\tilde {}\to \lambda$
(i.e., clause (iii) of \scite{4.6}(3B) holds trivially (as increasing
$p$ purely, $p \Vdash ``\underset\tilde {}\to \lambda = \lambda"$ and
$\Bbb P_{\ge p}$ satisfies the $\lambda$-c.c). \nl
1) If $\Bbb Q$ satisfies 
UP$^4_{\underset\tilde {}\to \kappa,\underset\tilde {}\to \lambda}(\Bbb I,
{\underset\tilde {}\to {\bold S}},\bold W)$, \ub{then} it satisfies
${\text{\rm UP\/}}^3_{\underset\tilde {}\to \kappa,
\underset\tilde {}\to \lambda}(\Bbb I,
{\underset\tilde {}\to {\bold S}},\bold W)$.  If $\Bbb Q$ satisfies
${\text{\rm UP\/}}^3_{\underset\tilde {}\to \kappa,
\underset\tilde {}\to \lambda}(\Bbb I,
{\underset\tilde {}\to {\bold S}},\bold W)$, \ub{then} $\Bbb Q$
satisfies ${\text{\rm UP\/}}^1
(\Bbb I,{\underset\tilde {}\to {\bold S}},\bold W)$ and
${\text{\rm UP\/}}^2_{\underset\tilde {}\to \kappa}(\Bbb I,\bold W)$.
If $\Bbb Q$ satisfies ${\text{\rm UP\/}}^2_{\underset\tilde {}\to \kappa}(\Bbb
I,{\underset\tilde {}\to {\bold S}},\bold W)$ or ${\text{\rm UP\/}}^1
(\Bbb I,\bold S,\bold W)$ then it satisfies ${\text{\rm UP\/}}^0
(\Bbb I,\bold S,\bold W)$.  
If $\ell \in \{3,4\}$ and $\Bbb Q$ satisfies 
${\text{\rm UP\/}}^\ell_{\underset\tilde {}\to \kappa,
\underset\tilde {}\to \lambda}
(\Bbb I,{\underset\tilde {}\to {\Bbb I}},\bold W)$ and
$\underset\tilde {}\to \kappa_1 \ge \underset\tilde {}\to \kappa,
{\underset\tilde {}\to \lambda_1} \ge \underset\tilde {}\to \lambda,[\ell =
4 \Rightarrow {\underset\tilde {}\to \lambda_1} =
\underset\tilde {}\to \lambda]$, \ub{then} $\Bbb Q$ satisfies
${\text{\rm UP\/}}^\ell_{{\underset\tilde {}\to \kappa_1},
{\underset\tilde {}\to \lambda_1}}(\Bbb I,
{\underset\tilde {}\to {\bold S}},\bold W)$. \nl
1A) If $\Bbb Q$ satisfies ${\text{\rm UP\/}}^\ell(\Bbb I,
{\underset\tilde {}\to {\bold S}},\bold W)$, \ub{then} it satisfies
${\text{\rm UP\/}}^0(\Bbb I,{\underset\tilde {}\to {\bold S}},\bold W)$
which implies ``$\Bbb Q$ has pure $\aleph_1$-decidability", see
Definition \scite{1.25}(2). \newline
2) The forcing notion $\Bbb Q$ satisfies ${\text{\rm UP\/}}^\ell(\Bbb I,
{\underset\tilde {}\to {\bold S}},\bold W)$ \ub{iff} its completion 
(i.e., $\Bbb Q^1$, or equivalently its completion
to a complete Boolean algebra) satisfies it assuming $\le_{\text{pr}} 
= \le$. \nl
3) If $\Bbb Q$ satisfies ${\text{\rm UP\/}}(\Bbb I,*,\bold W)$, 
(i.e., see \scite{4.6}(2)) and $\Bbb I$ is
$\mu^+$-complete (e.g., $\Bbb I = \emptyset$), \underbar{then} any ``new"
countable set of ordinals $< \mu$ is included in an ``old" countable set of
ordinals; i.e., one from $\bold V$. \newline
4) $\Bbb Q$ satisfies ${\text{\rm UP\/}}(\emptyset,*)$ 
\ub{iff} $\Bbb Q$ is purely proper (see Definition \scite{1.26A}(1)). \nl
5) $\Bbb Q$ satisfies ${\text{\rm UP\/}}(\emptyset,\{\aleph_1\})$ \ub{iff} 
$\Bbb Q$ is purely semiproper (see Definition \scite{1.26A}(2))
\footnote{that is: if $\Bbb Q \in N
\prec ({\Cal H}(\chi),\in),N$ is countable, $p \in \Bbb Q \cap N$,
\ub{then} for some $q$ we have $p \le_{pr} q$ and $q \Vdash_{\Bbb Q}
``\underset\tilde {}\to \tau \in N"$ for every $\Bbb Q$-name
$\underset\tilde {}\to \tau \in N$ of a countable ordinal}. 
\newline
6) If $\Bbb Q$ satisfies ${\text{\rm UP\/}}
(\Bbb I,{\underset\tilde {}\to {\bold S}},\bold W)$
and $\Bbb I \subseteq \Bbb I_1,{\underset\tilde {}\to {\bold S}_1} \subseteq
{\underset\tilde {}\to {\bold S}}$ and $\bold W_1 \subseteq \bold W$ 
\ub{then} $\Bbb Q$ satisfies ${\text{\rm UP\/}}^\ell
(\Bbb I_1,{\underset\tilde {}\to {\bold S}_1},
\bold W_1)$. \newline
7) In Definition \scite{4.1}, if $\Bbb P$ satisfies the $\kappa$-c.c. (e.g.
$\kappa = |\Bbb P|^+$) \underbar{then}:
\mr
\item "{$(a)$}"   we can replace
${\underset\tilde {}\to {\bold S}}$ by any set
${\underset\tilde {}\to {\bold S}'}$ of uncountable regular cardinals
of $\bold V$,
such that $\Vdash_{\Bbb P} ``{\underset\tilde {}\to {\bold S}} \cap \kappa =
{\underset\tilde {}\to {\bold S}'} \cap \kappa"$. 
\ermn
8) In Definition \scite{4.6} (in all the variants), \ub{if} we 
demand ``for $\chi$ large enough, for
some $x \in {\Cal H}(\chi)$, for every $\bar N$ such that $x \in
N_{\langle \rangle}$ and $\ldots$" we get an equivalent
definition. \nl
9) In Definition \scite{4.6} we can use weak $\Bbb I$-tagged trees,
i.e. we get with this an equivalent definition.
\endproclaim
\bigskip

\demo{Proof}  1), 2)  Trivial. \newline
3)  Straightforward. \newline  
4) Use \scite{4.8} below. \newline
5),6)  If $\Bbb I = \emptyset$, then every 
$N \prec ({\Cal H}(\chi),\in,<^*_\chi)$ is a $\Bbb I$-model. \newline
7) Easy. \newline
8) Check the Definition. \nl
9)  As in \cite{Sh:f}. \nl
10) The ``weak" version allows more trees of models so apparently is a
stronger condition, but by \scite{4.5B}(4) it is equivalent.
${{}}$ \hfill$\square_{\scite{4.7}}$\margincite{4.7}
\enddemo
\bigskip

\demo{\stag{4.9} Conclusion}  If $\Bbb P$ satisfies UP$^\ell(\Bbb I,
{\underset\tilde {}\to {\bold S}},\bold W)$ and 
${\underset\tilde {}\to {\bold S}}$ is as in \scite{4.1}$(*)(b)$ (or
${\underset\tilde {}\to {\bold S}} = \{\aleph_1\})$ (recall that this 
notation
implies $\Bbb I$ is $\aleph_2$-complete, $\aleph_1 \in
{\underset\tilde {}\to {\bold S}}$ and $\bold W \subseteq \omega_1$ 
stationary)
\underbar{then} $\Vdash_P ``\bold W$ is stationary".  Moreover, if $\bold W'
\subseteq \bold W$ is stationary then also $\Vdash_P ``\bold W'$ is a
stationary subset of $\omega_1"$.
\enddemo
\bigskip

\demo{Proof}  The ``moreover" fact is by \scite{4.7}(7) (i.e., monotonicity
in $\bold W$).

Assume that $p \Vdash ``\underset\tilde {}\to C$ is a club of $\omega_1$ and
$\underset\tilde {}\to C \cap \bold W = \emptyset"$.  By \scite{4.8} we can
find an $\aleph_1$-strictly $(\Bbb I,{\underset\tilde {}\to {\bold S}},
\bold W)$-suitable tree of models
$\langle N_\eta:\eta \in (T,\bold I)\rangle$ with $\underset\tilde {}\to C,p
\in N_{\langle \rangle}$.  Let $\delta = N_{<>} \cap \omega_1$, so $\delta \in
\bold W$.  By UP$^\ell(\Bbb I,{\underset\tilde {}\to {\bold S}},\bold W)$ we 
can find
a condition $q$ as in Definition \scite{4.1} in particular 
$p \le_{\text{pr}} q$.
Clearly $q \Vdash ``N_{\langle \rangle}[G] \cap \omega_1 = \delta"$ and,
trivially $p \Vdash_P ``\underset\tilde {}\to C$ is unbounded in
$N_{\langle \rangle}[G] \cap \omega_1"$ hence $p \Vdash ``N_{\langle \rangle}
[G] \cap \omega_1 \in \underset\tilde {}\to C"$.  So $q \Vdash_Q``\delta \in 
\underset\tilde {}\to C \cap \bold W"$, contradiction.
\hfill$\square_{\scite{4.9}}$\margincite{4.9}
\enddemo
\bigskip

\remark{\stag{4.10} Remark}  Usually we assume $\Bbb I,\bold S$ satisfies
\scite{4.1}(*)(a) + (c), ${\underset\tilde {}\to {\bold S}} = \{\aleph_1\}$
is the main case.
\endremark
\bigskip

\remark{\stag{4.11} Remark}  1) From the proof of \scite{4.8} we can 
conclude that in \scite{4.1}; we can replace ``$\bold S$-strictly 
$(\Bbb I,\bold W)$-suitable, $N_\eta \cap \omega_1 = \delta \in \bold W"$ by
$``(\Bbb I,\bold S,\bold W)$-suitable", and then the condition $q$ will be
$N_{\langle \rangle}$-semi-generic. \nl
2) As at present ${\underset\tilde {}\to {\bold S}} = \{\aleph_1\}$ seem to
suffice, we shall use only it for notational simplicity.
\endremark
\newpage

\head {\S5 An iteration theorem for UP} \endhead  \resetall \sectno=5
\bigskip

\proclaim{\stag{5.1} Claim}  1) If $\bar N = \langle N_\eta:\eta \in (T,
\bold I) \rangle$ is a tagged tree of models for $(\chi,\langle x,
\Bbb P \rangle),\Bbb P$ a forcing notion and $\Bbb P \in N_{\langle \rangle}$, \ub{then}
$\Vdash_{\Bbb P} ``\langle N_\eta[{\underset\tilde {}\to G_{\Bbb P}}]:
\eta \in (T,\bold I) \rangle$ is a tagged tree
of models for $(\chi,\langle x,\Bbb P,G \rangle)"$. \newline
2) If in addition $\Bbb P$ satisfies the $\kappa$-c.c. and 
$\Bbb I \in N_{\langle \rangle}$ is $\kappa$-closed (see Definition 
\scite{3.13}(2)) and $\bar N$ is $\Bbb I$-suitable, \ub{then}
\mr
\item "{$(*)$}" $\Vdash_{\Bbb P} ``\langle N_\eta[{\underset\tilde
{}\to G_{\Bbb P}}]:\eta \in (T,\bold I) \rangle \text{ is } \Bbb I \text{-suitable}"$.
\ermn
3) If in part (2) assume in addition that 
$\underset\tilde {}\to \kappa,{\underset\sim {}\to {\Bbb I}}$ 
to be $\Bbb P$-names of objects from $\bold V$ such that
${\underset\tilde {}\to {\Bbb I}}$ is purely decidable and 
$\Bbb P$ satisfies the local $\underset\tilde {}\to \kappa$-c.c. purely,  
\ub{then} for every $p \in N_{<>} \cap \Bbb P$ 
there is $q,p \le_{\text{pr}} q \in N_{<>} \cap \Bbb P$ 
forcing $(*)$ above.  If ${\underset\tilde {}\to {\Bbb I}}$ is purely
decidable and $\Bbb P$ is locally $\underset\tilde {}\to \kappa$-c.c. 
purely then we
can find $q$ satisfying $p \le_{\text{pr}} q \in N_{<>} \cap \Bbb P$ and forcing $(*)$.
\endproclaim
\bigskip

\demo{Proof}  1) Straight. \newline
2) So $\Bbb P$ satisfies the $\kappa$-c.c. and let
$G \subseteq \Bbb P$ be generic over $\bold V$.   Now from Definition 
\scite{4.1}(1), clearly $\langle N_\eta[G]:\eta \in (T,\bold I) \rangle$
satisfies clauses (a) - (d), so it is enough to check clause (e)$^-$ of
Definition \scite{4.1}(3).  So let $I \in \Bbb I \cap
N_\eta[G]$ where $\eta \in T$.  Hence there is $\underset\tilde {}\to I \in
N_\eta$ such that $\underset\tilde {}\to I$ is a $\Bbb P$-name and
$\underset\tilde {}\to I[G] = I$.  Let $\Bbb I' = \{J \in \Bbb I:
\text{for some } p \in \Bbb P$ we have $p \Vdash_{\Bbb P} ``\underset\tilde {}\to I = J"\}$.
So $\Bbb I'$ belongs to $\bold V$ and 
is a subset of $\Bbb I$ of cardinality $< \kappa$ and
$\Bbb I' \in N_\eta$ hence there is $I^* \in \Bbb I$ such that
$(\forall J)(J \in \Bbb I' \Rightarrow J \le_{RK} I^*)$, so \wilog \,
$I^* \in N_\eta$, hence as $N$ is $\Bbb I$-suitable clearly 
$\{\nu:\eta \triangleleft \nu \in T \text{ and } 
I^* \le_{RK} \bold I_\nu\}$ contains a front of $T[\eta]$.  
Hence in $\bold V[G]$, the set $\{\nu \in T:I \le_{\text{RK}} \bold I_\nu\}$
contains a front of $T^{[\eta]}$ as required.  \nl
3) Left to the reader.  \hfill$\square_{\scite{5.1}}$\margincite{5.1}
\enddemo
\bn
The point of the following claim is to get more from some UP$^\ell$ 
than seems on the surface; our aim is to help iterating.

\proclaim{\stag{5.2} Claim}  1) Assume 
$\ell \in \{3,4\}$ and the forcing notion
$\Bbb Q$ satisfies ${\text{\rm UP\/}}^\ell_{\underset\tilde {}\to \kappa,
\underset\tilde {}\to \lambda}({\underset\tilde {}\to {\Bbb I}},\bold W)$ and
$\underset\tilde {}\to \kappa,\underset\tilde {}\to \lambda,
{\underset\tilde {}\to {\Bbb I}},{\underset\tilde {}\to {\Bbb I}^+}$ are
$\Bbb Q$-names with pure decidability and 
$\Vdash ``{\underset\tilde {}\to {\Bbb I}^+} \backslash
{\underset\tilde {}\to {\Bbb I}}$ is 
$\underset\tilde {}\to \lambda$-complete".
\ub{Then} the forcing notion $\Bbb Q$ satisfies the 
${\text{\rm UP\/}}^\ell_{\underset\tilde {}\to \kappa,
\underset\tilde {}\to \lambda}
({\underset\tilde {}\to {\Bbb I}^+},\bold W)$. \nl 
2)  Suppose 
\mr
\item "{$(a)$}"  $\Bbb I_0,\Bbb I_1,\Bbb I_2,\Bbb I_3$ are sets of quasi-order
ideals, $\Bbb I_1 \subseteq \Bbb I_0 \subseteq \Bbb I_2,\Bbb I_3 = \Bbb I_1
\cup (\Bbb I_2 \backslash \Bbb I_0),\Bbb I_2 \backslash \Bbb I_0 =
\Bbb I_3 \backslash I_1$ is
$\kappa$-closed, $\kappa \le \lambda$ and $\Bbb I_2
\backslash \Bbb I_0$ is $\lambda$-complete
\sn
\item "{$(b)$}"  $\bar N = \langle N_\eta:\eta \in (T^*,\bold I^*) \rangle$
is a strict truely $\Bbb I_2$-suitable tree of models (for $\chi$ and $x =
\langle \Bbb Q,\Bbb I_0,\Bbb I_1,\Bbb I_2,\kappa,\lambda \rangle$
\sn
\item "{$(c)$}"  $p \in N_{<>},
\ell \in \{3,4\}$ and $\Bbb Q_{\ge p}$ satisfies the $\lambda$-c.c.
\sn
\item "{$(d)$}"  $\varphi(-)$ is a property with 
$\bar N,G_{\Bbb Q}$ as parameters (and possibly others)
\sn
\item "{$(e)$}"  for any $T'$ a subtree of $T^*$ such that $\langle N_\eta:
\eta \in (T',\bold I) \rangle$ is 
a truely $\Bbb I_0$-suitable tree of models 
there are $q = q_{T'},\underset\tilde {}\to T = 
\underset\tilde {}\to T(T')$ such that
{\roster
\itemitem{ $(i)$ }  $p \le_{\text{pr}} q \in \Bbb Q$
\sn
\itemitem{ $(ii)$ }  $(T',\bold I^*) \le (\underset\tilde {}\to T,\bold I^*)$
\sn
\itemitem{ $(iii)$ }  $q \Vdash_{\Bbb Q}$ ``$\langle N_\eta
[{\underset\tilde {}\to G_{\Bbb Q}}]:\eta \in (\underset\tilde {}\to T,\bold I)
\rangle$ is a truely $\Bbb I_1$-suitable tree of models and for every $\eta \in 
\underset\tilde {}\to T$ we have
\sn
\itemitem{ ${{}}$ }  $(\alpha) \quad 
N_\eta[{\underset\tilde {}\to G_{\Bbb Q}}] \cap \omega_1 = N_{<>} \cap \omega_1$ and
\sn
\itemitem{ ${{}}$ }  $(\beta) \quad \varphi[\eta]$
\sn
\itemitem{ ${{}}$ }  $(\gamma) \quad$ if $\ell = 3,y \in
N_\eta[{\underset\tilde {}\to G_{\Bbb Q}}]$ is a member of $\bold V$ then \nl

$\qquad \quad \{\nu:\eta
\triangleleft \nu \in \underset\tilde {}\to T$, and for some $A \in N_\nu,A$
a set of cardinality $< \lambda$ \nl

$\qquad \quad$ we have $y \in A\}$ contains a front of
$T^{[\eta]}$
\sn
\itemitem{ ${{}}$ }  $(\delta) \quad$ if $\ell = 4$ and $y \in N_\eta
[{\underset\tilde {}\to G_{\Bbb Q}}]$ is a set of $< \lambda$ members
of $\bold V$ then
\nl

$\qquad \quad \{\nu:\eta \triangleleft \nu \in \underset\tilde {}\to T$ 
and for some $A \in N_\eta$ a set of cardinality $< \lambda$ \nl

$\qquad \quad$ we have $y \subseteq A\}$
contains a front of $\underset\tilde {}\to T^{[\eta]}$.
\endroster}
\ermn
\ub{Then} there are $q,\underset\tilde {}\to T$ such that:
\mr
\widestnumber\item{$(iii)$}
\item "{$\circledast(i)$}"  $p \le_{\text{pr}} q \in \Bbb Q$
\sn
\item "{$(ii)$}"  $(T^*,\bold I^*) \le (\underset\tilde {}\to T,\bold I^*)$
\sn
\item "{$(iii)$}"  $q \Vdash_{\Bbb Q} ``\langle 
N_\eta[{\underset\tilde {}\to G_{\Bbb Q}}]:\eta \in (\underset\tilde {}\to T,
\bold I)\rangle$ is a truely $\Bbb I_3$-suitable tree of models and
for every $\eta \in \underset\tilde {}\to T$ we have
\sn
\item "{${{}}$}"  $(\alpha) \quad \quad N_\eta
[{\underset\tilde {}\to G_{\Bbb Q}}] \cap \omega_1 
= N_{<>} \cap \omega$
\sn
\item "{${{}}$}"   $(\beta) \quad \quad \models \varphi(\eta)$
\sn
\item "{${{}}$}"  $(\gamma),(\delta) \quad$  as in clause (iii)
of (e) above.
\ermn
3) In part (2), if $\Bbb Q$ satisfies the $\lambda$-c.c., 
\ub{then} we can omit $(\gamma),(\delta)$ in their two
appearances as they follow. \nl
4) In part (2), we can replace ``truely $\Bbb I_\ell$-suitable" by
``weakly $\Bbb I_\ell$-suitable".
\endproclaim
\bigskip

\remark{\stag{5.2a} Remark}  In part (2) clause (e) we can restrict $T'$ to
those needed.
\endremark
\bigskip

\demo{Proof}  1) As in the definition of UP$^\ell_{\underset\tilde
{}\to \kappa,\underset\tilde {}\to \lambda}(\Bbb I^+,\bold W)$ let
$\bar N = \langle N_\eta:\eta \in (T,\Bbb I) \rangle$ be a
strict $\Bbb I^+$-suitable tree of models for $\chi$ and $x \equiv
\langle \Bbb Q,
{\underset\tilde {}\to {\Bbb I}},{\underset\tilde {}\to {\Bbb I}^+},
\bold W,\underset\tilde {}\to \kappa,
\underset\tilde {}\to \lambda \rangle$ and $p \in N_{<>}$.  
We can find $p' \in N_{<>},p \le_{\text{pr}} p' \in \Bbb Q$ which forces 
$\underset\tilde {}\to \kappa,
\underset\tilde {}\to \lambda,{\underset\tilde {}\to {\Bbb I}},
{\underset\tilde {}\to {\Bbb I}^+}$ to be say $\kappa,\lambda,\Bbb I,\Bbb I^+$
respectively.
Now we can apply part (2) of the claim
with $\bar N,\Bbb I,\Bbb I^{[\kappa]},\Bbb I^+,\Bbb I^{[\kappa]} \cup
(\Bbb I^+ \backslash \Bbb I),x=x$ here standing for $\bar N,\Bbb I_0,\Bbb I_1,
\Bbb I_2,\Bbb I_3,\varphi$ there. \nl
2) Let ${\Cal T} = \{T:T$ is a subtree of $T^*$ such that $\bar N \restriction
T$ is a truely $\Bbb I_0$-suitable tree of models$\}$ or just 
${\Cal T} = \{T:T$ a subtree of
$T^*$ such that ($<> \in T,\eta \in T \Rightarrow \emptyset \ne
\text{ Suc}_T(\eta) \subseteq \text{ Suc}_{T^*}(\eta),T$ closed under initial
segments and): if $\eta \in \text{ split}(T^*,\bold I^*) \and \bold I_\eta
\in \Bbb I_0$ then Suc$_T(\eta) \in \bold I^+_\eta$ and if $\eta \notin
\text{ split}(T^*,\bold I^*) \vee \bold I_\eta \notin \Bbb I_0$ then
$|\text{Suc}_T(\eta)| = 1\}$. \nl
For any $T \in {\Cal T}$ by assumption (e) there are $q_T,
{\underset\tilde {}\to T^{\Bbb Q}}[T]$ as required there.
\sn
We shall show that some such $q_T$ is as required.  We define a $\Bbb Q$-name
${\underset\tilde {}\to T^\oplus}$ as follows:

for $G_{\Bbb Q} \subseteq \Bbb Q$ generic over $\bold V$ we let

$$
\align
T^\otimes = {\underset\tilde {}\to T^\otimes}[G_{\Bbb Q}] = \{\eta \in T^*:&N_\eta
[G] \cap \omega_1 = N_{<>} \cap \omega_1 \\
  &\text{and } \ell \le \ell g(\eta) \Rightarrow \varphi(\eta \restriction
\ell)\}.
\endalign
$$
\mn
Clearly
\mr
\item "{$(*)_1$}"  $q_T 
\Vdash_{\Bbb Q} ``{\underset\tilde {}\to T^{\Bbb Q}}[T]
\subseteq T^\otimes"$ for every $T \in {\Cal T}$.
\ermn
Working in $\bold V[G_{\Bbb Q}]$ we define a depth function Dp function from
$T^\otimes$ to Ord $\cup \{\infty\}$ by defining for any ordinal $\alpha$
when $Dp(\eta) \ge \alpha$ as follows:
\mr
\item "{$\boxtimes$}"  Dp$(\eta) \ge \alpha$ \ub{iff} the following
conditions hold:
{\roster
\itemitem{ $(\alpha)$ }  $\eta \in T^\otimes$
\sn
\itemitem{ $(\beta)$ }  for every $\beta < \alpha$ there is $\nu \in
\text{ Suc}_{T^\otimes}(\eta)$ such that $Dp(\nu) \ge \beta$
\sn
\itemitem{ $(\gamma)$ }  if $\beta < \alpha$ and $\ell \le \ell g(\eta)$
and $I \in N_\eta[G] \cap \Bbb I_3$ \ub{then} for some $\nu$ with 
Dp$(\nu) \ge \beta$ we have $\eta \trianglelefteq \nu \in T^\otimes \and 
\bold I = \bold I_\nu \in \Bbb I_3$ and \nl
$\{ \nu \in \text{ Suc}_{T^\otimes}(\nu):Dp(\nu) \ge \beta\} \ne \emptyset 
\text{ mod } \bold I_\nu$
\sn
\itemitem{ $(\delta)$ }  if $\beta < \alpha$ and $\underset\tilde {}\to y \in
N_\eta$ is a $\Bbb Q$-name of a member of $\bold V$ when $\ell=3$ and is a
set of cardinality $< \lambda$ when $\ell = 4$ \ub{then} for some 
$\eta',\eta \trianglelefteq \eta'$ with Dp$(\eta') \ge \beta$ and $A \in
N_{\eta'}$ of cardinality $< \lambda$ we have $[\ell = 3 \Rightarrow
\underset\tilde {}\to y[{\underset\tilde {}\to G_{\Bbb Q}}] \in A]$ and $[\ell=4
\Rightarrow \underset\tilde {}\to y[G_{\Bbb Q}] \subseteq A]$.
\endroster}
\ermn
Clearly it is enough to show in $\bold V$ that for some $T \in {\Cal T}$ we have
$q_T \Vdash ``\text{Dp}(<>) = \infty"$.  Note
\mr
\item "{$(*)_2$}"  if $\eta \triangleleft \nu \in T^\otimes$ then
Dp$(\eta) \ge \text{ Dp}(\nu)$.
\ermn
Clearly in $\bold V^{\Bbb Q}$ we have
\mr
\item "{$(*)_3$}"  if $\eta \in 
T^\otimes$, Dp$(\eta) < \infty$ and $\bold I_\eta
\in \Bbb I_3$ then $\{\nu \in \text{ Suc}_{T^\otimes}(\eta):\text{Dp}(\eta) <
\text{ Dp}(\nu)\} = \emptyset \text{ mod } \bold I_\eta$
\sn
\item "{$(*)_4$}"  if $\eta \in T^\otimes$, Dp$(\eta) < \infty$ 
and $\bold I_\eta \notin \Bbb I_3,\nu \in \text{ Suc}_{T^\otimes}
(\eta)$ then Dp$(\eta) \ge \text{ Dp}(\nu)$.
\ermn
For each $\eta \in T^\otimes$ such that $\bold I_\eta \in \Bbb I_2
\backslash \Bbb I_0$ define the set $A_\eta$ as follows:
\sn
First, ${\Cal A}_\eta$ is the minimal family of sets satisfying
\mr
\widestnumber\item{$(iii)$}
\item "{$(i)$}"   if $\ell =3$ and $\underset\tilde {}\to y \in N_\eta$ is a
$\Bbb Q$-name of a member of $\bold V$ \ub{then} the set
$A^3_{\eta,\underset\tilde {}\to y} =:
\{\rho \in \text{ Suc}_{T^*}(\eta):
\text{ for some set } A \in N_\rho 
\subseteq \bold V$ of cardinality $< \kappa,
\underset\tilde {}\to y \in A\}$ belongs to ${\Cal A}_\eta$
\sn
\item "{$(ii)$}"  if $\ell =4$, parallely (using
$A^4_{\eta,\underset\tilde {}\to y}$), i.e., if $y \in N_\eta$ is a
$\Bbb Q$-name of a family of $< \kappa$ members of $\bold V$ then the
set $A^4_{\eta,\underset\tilde {}\to y} =: \{\rho \in \text{
Suc}_{T^*}(\eta)$: for some set $A \in N_\rho \subseteq \bold V$ of
cardinality $, \kappa$ we have $\underset\tilde {}\to y \subseteq A\}$
belongs to ${\Cal A}_\eta$
\sn
\item "{$(iii)$}"  if $\ell \le \ell g(\eta),\beta = \text{ Dp}(\eta
\restriction \ell)$ then the set $A^*_{\eta,\ell} = 
\{\rho \in \text{ Suc}_{T^*}(\eta):\text{Dp}(\rho)
\ge \beta\}$ belongs to ${\Cal A}_\eta$.
\ermn
Let ${\Cal A}'_\eta = \{A \in {\Cal A}_\eta:A = \emptyset$ mod
$\bold I_\eta\}$, note that ${\Cal A}'_\eta$ is a 
countable family of members of $\bold I_\eta$ (more exactly the ideal
$\bold I^{{\bold V}^{\Bbb Q}}_\eta$, which $\bold I_\eta$ generates in
$\bold V^{\Bbb Q}$), and
so actually we have defined a $\Bbb Q$-name ${\underset\tilde {}\to
A_\eta} = \cup\{A:A \in {\Cal A}'_\eta\}$.
\nl
Now we can define in $\bold V$ a sequence $\langle B^*_\eta:\eta \in T^* \rangle$
such that
\mr
\item "{$(*)_5$}"  $B_\eta = \emptyset \text{ mod } \bold I_\eta$ for
$\eta \in T^*$
\sn
\item "{$(*)_6$}" $(i) \quad$ if Suc$_{T^*}(\eta) = \emptyset \text{ mod }
\bold I_\eta$ or $\bold I_\eta \notin \Bbb I_2 \backslash \Bbb I_0$ then
$B^*_\eta = \emptyset$
\sn
\item "{${{}}$}" $(ii) \quad$ if Suc$_{T^*}(\eta) \ne \emptyset \text{ mod }
\bold I_\eta$ and $\bold I_\eta \in \Bbb I_2 \backslash \Bbb I_0$ 
\ub{then} \nl

$\qquad \quad 
p \Vdash_{\Bbb Q} ``{\underset\tilde {}\to A_\eta} \subseteq B^*_\eta$ if
$\eta \in {\underset\tilde {}\to T^\otimes}"$. 
\ermn
This is possible as $\Bbb Q_{\ge p}$ satisfies the $\lambda$-c.c. and each
$I \in \Bbb I_2 \backslash \Bbb I_0$ is $\lambda$-complete.
\sn
Now we define

$$
T_0 = \{\eta \in T^*:\text{ for every } \ell < \ell g(\eta) \text{ we have }
\eta \restriction (\ell +1) \notin B^*_{\eta \restriction \ell}\}.
$$
\mn
Clearly we can find $T_1 \in {\Cal T}$ such that $T_1 \subseteq T_0$
(in particular by $(*)_6$(i)).
So if $q_{T_1} \Vdash_{\Bbb Q} ``\text{Dp}(<>) = 
\infty"$ then we are done so toward
contradiction we assume that this fails.  Hence, there is a subset
$G_{\Bbb Q}$ of $\Bbb Q$ generic over $\bold V$ 
such that $q_{T_1} \in G_{\Bbb Q}$ and $\bold V[G_{\Bbb Q}] \models
``\text{Dp}(<>) < \infty"$.  
As $q_{T_1} \in G_{\Bbb Q}$ clearly $\underset\tilde {}\to T
(T)[G] \subseteq T^\otimes[G_{\Bbb Q}]$. \nl
By the 
choice of $G_{\Bbb Q}$ we have $Dp(<>) < \infty$ hence by $(*)_2$ we have
$\eta \in \underset\tilde {}\to T(T_1)[G_{\Bbb Q}] \Rightarrow \eta \in 
{\underset\tilde {}\to T^\otimes}[G_{\Bbb Q}] \Rightarrow Dp(\eta) < \infty$.

Now we shall prove by induction on $\alpha \in \text{ Ord}$ that $\eta \in
\underset\tilde {}\to T(T_1)[G_{\Bbb Q}] \Rightarrow \text{ Dp}(\eta) \ge
\alpha$.  For
$\alpha = 0,\alpha$ limit we have no problem, so let $\alpha = \beta +1$,
assume toward contradiction Dp$(\eta) = \beta$ for some $\eta \in
\underset\tilde {}\to T(T_1)[G_{\Bbb Q}]$, hence by $(*)_2$ 
and the induction hypothesis we have
\mr
\item "{$\boxtimes_2$}"  $\eta 
\trianglelefteq \nu \in \underset\tilde {}\to T
(T_1)[G_{\Bbb Q}] \Rightarrow \text{ Dp}(\nu) = \beta$
\sn
\item "{$\boxtimes_3$}"  if $I \in (\Bbb I_2 \backslash \Bbb I_0) \cap
N_\eta[G_{\Bbb Q}]$ then for some $\rho \in (\underset\tilde {}\to
T(T_1)[G_{\Bbb Q}],
\eta \triangleleft \rho$ (in fact ``many $\rho$'s) and $J \in \Bbb I_3 \cap
N_\rho$ we have $I = \bold I_\rho$ \nl
[why?  by clause (e)(iii) of the assumption, i.e. choice of
$q_{T_1},\underset\tilde {}\to T(T_1)$.]
\sn
\item "{$\boxtimes_4$}"  if $I \in \Bbb I_3 \cap N_\eta[G_{\Bbb Q}]$ then there is
$\nu$ satisfying $\eta \triangleleft \nu \in \underset\tilde {}\to T(T_1)
[G_{\Bbb Q}]$ such that $I = \bold I_\nu$ \nl
[why? if $I \in \Bbb I_2 \backslash \Bbb I_1$ by $\boxtimes_3$, if 
$I \in \Bbb I_0$ use the
choice of $\underset\tilde {}\to T(T_1)$.]
\ermn
Now
\mr
\item "{$\boxtimes_5$}"  if $I \in N_\eta \cap \Bbb I_3$
\ub{then} for some $\nu$ satisfying
$\eta \triangleleft \nu \in \underset\tilde {}\to T
(T_1)[G_{\Bbb Q}]$ we have $\bold I_\ell = \bold I_\nu$ and \nl
$\{\rho \in \text{ Suc}_{T^\otimes}(\nu):\text{Dp}
(\rho) \ge \beta\} \ne \emptyset \text{ mod } \bold I_\nu$.
\ermn
[Why true?  We can choose $\nu$ such that $I = \bold I_\nu$ and 
$\eta \triangleleft \nu \in \underset\tilde {}\to T
(T)[G_{\Bbb Q}]$ and Suc$_{{\underset\tilde {}\to T}(T_1)
[G_{\Bbb Q}]}(\nu) \ne \emptyset$
mod $I_\nu$ and choose $\rho' \in 
\text{ Suc}_{{\underset\tilde {}\to T}(T_1)[G_{\Bbb Q}]}(\nu) \subseteq
{\underset\tilde {}\to T^\otimes}[G_{\Bbb Q}]$.  First assume $\bold
I_\nu \in \Bbb I_3 \backslash \Bbb I_0$. \nl
Now easily 

$$
{\underset\tilde {}\to A^*_{\nu,\ell g(\eta)}}[G_{\Bbb Q}] 
= \emptyset \text{ mod } \bold I_\nu
\Rightarrow {\underset\tilde {}\to A^*_{\nu,\ell g(\eta)}}
[G_{\Bbb Q}] \subseteq B^*_\nu
\Rightarrow \rho' \notin {\underset\tilde {}\to A_\nu}[G_{\Bbb Q}]
$$
\mn
by the definition of ${\underset\tilde {}\to A^*_{\nu,\ell
g(\eta)}}[G_{\Bbb Q}]$ and $\beta$ we get

$$
\rho' \in \{\rho \in \text{ Suc}_{T^\otimes[G_{\Bbb Q}]}(\nu):
\text{Dp}(\rho) \ge \text{ Dp}(\nu)(= \beta)\}
$$
\mn
easy contradiction.
\sn
Next assume $\bold I \in \Bbb I_1$.  Now
Suc$_{\underset\tilde {}\to T(T_1)[G_{\Bbb Q}]}(\nu)$ is a set witnessing the
requirement.]  \nl
Now for $\eta \in \underset\tilde {}\to T(G_1)
[G_{\Bbb Q}]$ we check the definition
of Dp$(\eta) \ge \beta + 1$: clause $(\alpha)$ holds as
$\underset\tilde {}\to T(T_1)[G_{\Bbb Q}] \subseteq \underset\tilde {}\to T^\otimes
[G_{\Bbb Q}]$, clause $(\beta)$ holds by the induction hypothesis, 
clause $(\gamma)$ holds by
$\boxtimes_3 + \boxtimes_4 + \boxtimes_5$ and clause $(\delta)$ by the choice
of $q_{T_1},\underset\tilde {}\to T(T_1)$.  So Dp$(\eta) \ge \beta +1$
contradiction. \nl
3),4)  Similar to the proof of part (2). \hfill$\square_{\scite{5.2}}$\margincite{5.2}
\enddemo
\bn
We now can deduce more implications between the UP$^\ell$-s.
\demo{\stag{5.3} Conclusion}  Assume that 
$\underset\tilde {}\to \kappa,\underset\tilde
{}\to \lambda$ are purely decided $\Bbb Q$-names and
$\Bbb I^{[\underset\tilde {}\to \kappa]}$ is 
$(< \underset\tilde {}\to \lambda)$-closed, (which just means: if
$r \Vdash ``\underset\tilde {}\to \kappa = \kappa,\underset\tilde
{}\to \lambda = \lambda$ and
${\underset\tilde {}\to {\Bbb I}} = \Bbb I"$ where $r \in \Bbb P$ then
$\Bbb I^{[\kappa]}$ is $\lambda$-closed) and $\Bbb Q$ satisfies the local
$\underset\tilde {}\to \lambda$-c.c. purely. \nl
1) If $\Bbb Q$ satisfies UP$^1_{\underset\tilde {}\to
\kappa,\underset\tilde {}\to \lambda}(\Bbb I,\bold W)$ 
\underbar{then} $\Bbb Q$ satisfies UP$^2_{\underset\tilde {}\to
\kappa}({\underset\tilde {}\to {\Bbb I}},
\bold W)$ and UP$^4_{\underset\tilde {}\to \kappa,\underset\tilde
{}\to \lambda}({\underset\tilde {}\to {\Bbb I}},\bold W)$,
UP$^3_{\underset\tilde {}\to \kappa,\underset\tilde {}\to
\lambda}({\underset\tilde {}\to {\Bbb I}},\bold W)$. \nl
2) UP$^1({\underset\tilde {}\to {\Bbb I}},\bold W)$ implies
UP$^5_{\underset\tilde {}\to \kappa,\underset\tilde {}\to
\lambda}({\underset\tilde {}\to {\Bbb I}},\bold W)$ (for $\Bbb Q$)
\ub{if} $\underset\tilde {}\to \kappa = \underset\tilde {}\to \lambda$
is a $\Bbb Q$-name decidably pure, $\Bbb Q$ satisfies the local
$\underset\tilde {}\to \kappa$-c.c. \nl
3) 4-5 fill!!!
\enddemo
\bigskip
 
\demo{Proof}  Why?  By Definition \scite{4.6} and Lemma \scite{5.2}.
\hfill$\square_{\scite{5.3}}$\margincite{5.3}
\enddemo
\bigskip

\definition{\stag{5.4} Definition}  1) We say $\bar \Bbb Q = \langle
\Bbb P_i,{\underset\tilde {}\to {\Bbb Q}_i},
{\underset\tilde {}\to {\Bbb I}_i},
{\underset\tilde {}\to \kappa_i},{\underset\tilde {}\to {\bold S}_i}:
i < \alpha \rangle$ is UP$^{4,e}(\bold W,W)$-suitable iteration \ub{if}:
\medskip
\roster
\item "{$(a)$}"  $\langle \Bbb P_i,
{\underset\tilde {}\to {\Bbb Q}_j}:i < \alpha,j <
\alpha \rangle$ is an $\aleph_1-\text{Sp}_e(W)$-iteration
\footnote{the reader can fix $W$ as the class of strongly inaccessible
cardinals}
\sn
\item "{$(b)$}"  ${\underset\tilde {}\to {\Bbb I}_i}$ is a 
$\Bbb P_i$-name of 
a set of quasi order ideals with domain a 
cardinal in $\bold V^{{\Bbb P}_i}$ for 
notational simplicity \ub{or} even just a $\Bbb P_{i+1}$-name of such objects
(i.e., $\Vdash_{{\Bbb P}_{i+1}} ``{\underset\tilde {}\to {\Bbb I}} \in
\bold V^{{\Bbb P}_i}"$) such that in $\bold V^{{\Bbb P}_i},
{\underset\tilde {}\to {\Bbb I}_i}/
{\underset\tilde {}\to G_{{\Bbb P}_i}}$ which is a 
${\underset\tilde {}\to {\Bbb Q}_i}$-name, is purely decidable
\sn
\item "{$(c)$}"  $\bold W \subseteq \omega_1$ is stationary
\sn
\item "{$(d)$}"  for each $i < \alpha$, we have: ${\underset\tilde {}\to 
\kappa_i}$ is a $\Bbb P_i$-name of a regular uncountable cardinal of
$\bold V^{{\Bbb P}_i}$, purely decidable
\sn
\item "{$(e)$}"  $\Vdash_{{\Bbb P}_i} 
``{\underset\tilde {}\to {\Bbb Q}_i}$ satisfies
UP$^4_{{\underset\tilde {}\to \kappa_i},{\underset\tilde {}\to \kappa_{i+1}}}
({\underset\tilde {}\to {\Bbb I}_i},
{\underset\tilde {}\to {\bold S}_i},\bold W)$ and 
${\underset\tilde {}\to {\Bbb I}_i}$ is
${\underset\tilde {}\to \kappa_i}$-complete set of partial order ideals
(from $\bold V$) and $({\underset\tilde {}\to {\Bbb Q}_i},
\le_{\text{vpr}})$ is $\aleph_1$-complete (see \scite{1.1})"
\sn
\item "{$(f)$}"  or $i<j$ we have $\Vdash_{P_j} ``
{\underset\tilde {}\to \kappa_i} \le {\underset\tilde {}\to \kappa_j}"$
\sn
\item "{$(g)$}"  $\Bbb P_{i +1}$ satisfies the ${\underset\tilde {}\to
\kappa_i}$-c.c., or just for every $p \in \Bbb P_{i+1}$ 
there are $\kappa',q$ such that 
{\roster
\itemitem{ $(\alpha)$ }   $p \le_{\text{pr}} q \in \Bbb P_{i+1}$
\sn
\itemitem{ $(\beta)$ }   $q \Vdash_{({\Bbb P}_{i+1},
\le_{\text{pr}})} ``{\underset\tilde {}\to \kappa_{i+1}} = \kappa'"$
\sn
\itemitem{ $(\gamma)$ }   $\kappa \le \kappa' \in$ UReg and ${\underset\tilde {}\to 
\zeta_\varepsilon}$ is a $\Bbb P_i$-name of an ordinal for $\varepsilon <
\varepsilon^* < \kappa'$ and $p',q \restriction i \le p' \in \Bbb P_i$ 
\ub{then} there is
$q',p' \le_{\text{pr}} q' \in \Bbb P_i$ and set $a$ of $< \kappa'$
ordinals such that $q' \Vdash_P ``{\underset\tilde {}\to \zeta_\varepsilon} \in a$ for
$\varepsilon < \varepsilon^*"$ 
\endroster}
\item "{$(h)$}"  [?] for any $i < \alpha$ for some $n < \omega$, if $i < j <
\alpha$ and $p \in \Bbb P_{j+1}$ and $\kappa^\otimes_{\underset\tilde {}\to \kappa_j}(p,
\Bbb P_{j+1}) = \kappa$ and ${\underset\tilde {}\to \zeta_\varepsilon}$ is a
$(\Bbb P_i)_{\ge p \restriction i}$-name of an ordinal for $\varepsilon <
\varepsilon^* < \kappa$ \ub{then} 
for some $q$ and $a$ we have: $p \restriction i
\le_{\text{pr}} q \in \Bbb P_i$ and $q \Vdash_{{\Bbb P}_i} 
``{\underset\tilde {}\to \zeta_\varepsilon} \in a$ for
$\varepsilon < \varepsilon^*"$ where $a$ is a set of $< \kappa$ ordinals.
\ermn
2) We may write UP$^{4,e}$ instead UP$^{4,e}(\omega_1$, the class of
strongly inaccessibles).  If we omit ${\underset\tilde {}\to {\bold S}_i}$ we mean
$\{\aleph_1\}$.  \nl
If we omit ${\underset\tilde {}\to {\Bbb I}_i}$, we mean ``some
${\underset\tilde {}\to {\Bbb I}}$ as required" (note that the requirements
on ${\underset\tilde {}\to {\Bbb I}_i}$ are actually on each member so the
family of candidates to being ${\underset\tilde {}\to {\Bbb I}}$ is closed
under union).  If we omit $e$ we mean $e=6$.  We may omit $W$ if it is the
class of strongly inaccessible cardinals.
\nl
If we omit ${\underset\tilde {}\to \kappa_i}$ we mean some such 
$\Bbb P_{i+1}$-name.  (Can we eliminate names?  
Well if we use the iteration as in \cite[Ch.X,\S1]{Sh:f},
(RCS) no, but if we waive associativity as done here, we can). \nl
3) We defined UP$^{3,e}(\bold W,W)$-suitable iterations similarly but
replace clause (e) by:
\mr
\item "{$(e)_a$}"  as above replacing 
UP$^4_{{\underset\tilde {}\to \kappa_i},
{\underset\tilde {}\to \kappa_{i+1}}}$
by UP$^3_{{\underset\tilde {}\to \kappa_i},
{\underset\tilde {}\to \kappa_{i+1}}}$.
\ermn
4) We define a UP$^\ell(\Bbb I,\bold W,W)$-iterations as above but with
${\underset\tilde {}\to {\Bbb I}_i} = 
{\underset\tilde {}\to {\Bbb I}^{[{\underset\tilde {}\to \kappa_i}]}}$
[Saharon like \S6, straight in successor, in limit work it out, Question:
${\underset\tilde {}\to \kappa_i}$ pure name?]]
\enddefinition
\bn
We can also deal with strong preservation
\definition{\stag{5.4B} Definition}  We say that $\bar{\Bbb Q} =
\langle \Bbb P_i,{\underset\tilde {}\to {\Bbb Q}_i},
{\underset\tilde {}\to \kappa_j}:
i < \alpha,j \le \alpha \rangle$ is a weak UP$^4(\bold W,W)$-iteration if
\mr
\item "{$(a)$}"  $\langle \Bbb P_i,
{\underset\tilde {}\to {\Bbb Q}_i}:i < \alpha \rangle$
is a $\aleph_1-\text{Sp}_6(W)$-iteration, 
${\underset\tilde {}\to \kappa_i}$ is
$(\Bbb P_i,\le^{{\Bbb P}_i}_{\text{pr}})$-name such that $j < i \Rightarrow
{\underset\tilde {}\to \kappa_j} \le {\underset\tilde {}\to \kappa_i}$
\sn
\item "{$(b)$}"  if $i < j \le \alpha$ and $i$ is non-limit, $p \in
\Bbb P_i$ and $p \Vdash_{({\Bbb P}_i,\le_{\text{pr}})} 
``{\underset\tilde {}\to \kappa_i} =
\kappa_i"$ then $p \Vdash_{{\Bbb P}_i} ``\Bbb P_j/G_{\Bbb P_i}$ satisfies
UP$^4_{{\underset\tilde {}\to \kappa_i},{\underset\tilde {}\to \kappa_j}}
({\underset\tilde {}\to {\Bbb I}},\bold W)$ for some $\Bbb P_i$-name
${\underset\tilde {}\to {\Bbb I}}$ of ${\underset\tilde {}\to 
\kappa_i}$-complete ideals".
\endroster
\enddefinition
\bigskip

\proclaim{\stag{5.5} Lemma}  Assume that $\bold W \subseteq \omega_1$ is 
stationary and $\bar{\Bbb Q} = 
\langle \Bbb P_i,{\underset\tilde {}\to {\Bbb Q}_i},
{\underset\tilde {}\to {\Bbb I}_i},{\underset\tilde {}\to \kappa_i}:
i < \alpha \rangle$ is a UP$^4(\bold W,W)$-suitable iteration, and 
$\Bbb P_\alpha = { \text{\rm Sp\/}}(W)-{\text{\rm Lim\/}}
(\bar{\Bbb Q})$ be the limit. \nl
For $j \le \alpha$ we define ${\underset\tilde {}\to \kappa^*_j}$, a
$(\Bbb P_j,\le_{\text{pr}})$-name of a member of 
${\text{\rm UReg\/}}^{\bold V}:
{\underset\tilde {}\to \kappa^*_j} =$ \nl
$\text{Min}\{\kappa \in { \text{\rm UReg\/}}^{\bold V}:
i < j \Rightarrow {\underset\tilde {}\to \kappa_i} \le \kappa \text{ and }
\kappa \ge \aleph_2 \text{ and } \kappa \ge j\}$.
\sn
For simplicity we can assume
\mr
\item "{$\boxtimes$}"  there are $\langle \kappa'_i:i < \alpha
\rangle,\Bbb P_i$
satisfies the $\kappa'_i$-c.c. and for each $i,i + \omega \le \alpha$ for
some $\Vdash_{(\Bbb P_i,\le_{\text{pr}})} ``{\underset\tilde {}\to \kappa_i} \le
\kappa'_i",\Vdash_{(\Bbb P_{i+1},\le_{\text{pr}})} ``\kappa'_i \le
{\underset\tilde {}\to \kappa'_{i+1}}"$ (and $\kappa^*_\beta = \beta
\Rightarrow \Bbb P_\beta = \dbcu_{i < \beta}\Bbb P_i$ (i.e. $\kappa^*_\beta$ is
strongly inaccessible $> |\Bbb P_i|$ for $i < \beta$ \ub{or} change support?)
\ermn
1) For each $\beta \le \alpha,\Bbb P_\beta$ satisfies 
UP$^5_{{\underset\tilde {}\to \kappa^*_\beta}}(\Bbb I'_\beta,\bold W)$ 
for some $\kappa^*$-complete $\Bbb I'_\beta \in \bold V$. \newline
1A) If $\gamma \le \beta \le \alpha,p \Vdash_{{\Bbb P}_\gamma}
``{\underset\tilde {}\to \kappa^*_\gamma} = \kappa_\gamma",p \in
G_{{\Bbb P}_\gamma}$ \ub{then} in $\bold V[G_{{\Bbb P}_\gamma}]$ 
the forcing notions $\Bbb P_\beta/
G_{{\Bbb P}_\gamma}$ satisfies 
${\text{\rm UP\/}}^5_{{\underset\tilde {}\to \kappa_\gamma},
{\underset\tilde {}\to \kappa_\beta}}
({\underset\tilde {}\to {\Bbb I}'_{\beta,
\delta}},\bold W)$ for some $\kappa_j$-complete
${\underset\tilde {}\to {\Bbb I}_{\beta,\gamma}}
\in \bold V^{{\Bbb P}_\gamma}$ (so this justifies the weak in
\scite{5.4B}). 
\nl
2) In fact each $I \in \Bbb I'_\beta$ has domain of cardinality 
$\le (\underset{\gamma < \beta} \to \sup\{ \lambda^\kappa:\,\, 
\nVdash_{{\Bbb P}_\gamma} ``\neg(\exists I \in \dbcu_{i<j}
{\underset\tilde {}\to {\Bbb I}_i})(\lambda = 
|{\text{\rm Dom\/}}(I)|)$ and 
${\underset\tilde {}\to \kappa_i} > \kappa$ for some
$i \le \gamma"\}$ and
$|\Bbb I'_\beta| \le 
\dsize \sum_{\gamma < \beta} (\aleph_0 + |\Bbb P_\gamma| +
\text{ sup}\{\lambda^\kappa:\,\, \nVdash_{{\Bbb P}_\gamma} 
``\neg(|\dbcu_{i \le \gamma}{\underset\tilde {}\to {\Bbb I}_i}| 
\le \lambda$ and $\kappa_i > \kappa$ for some $i \le \gamma)"\}$.
Similarly for ${\underset\tilde {}\to {\Bbb I}_{\gamma,\beta}}$. \nl
3) Similarly for weak ${\text{\rm UP\/}}^5(\bold W,W)$-iterations.
\endproclaim
\bigskip

\remark{\stag{5.5A} Remark}  1)  We can also get 
the preservation version of this Lemma. \nl
2) The reader can concentrate on the case that
${\underset\tilde {}\to \kappa'_\ell}$'s are objects and not names.
\endremark
\bigskip

\demo{Proof}  1) We prove this by induction on $\beta$, so \wilog \,
$\beta = \alpha$.  For each $\gamma < \alpha$ let ${\Cal J}_\gamma =:
\{q \in \Bbb P_{\gamma +1}:q \text{ forces a value to } 
{\underset\tilde {}\to \kappa_\gamma}$, called 
$\kappa_{\gamma,q}$ and $q$ forces 
${\underset\tilde {}\to {\Bbb I}_\gamma}$ to be equal to a $\Bbb P_\gamma$-name 
called ${\underset\tilde {}\to {\Bbb I}_{\gamma,q}}$ and
$q \restriction \gamma$ forces that 
$|{\underset\tilde {}\to {\Bbb I}_\gamma}|$ is 
$\le \mu_{\gamma,q}$ but no $q'$ such that $q \restriction \gamma 
\le_{\text{pr}} q' \in \Bbb P_\gamma$ forces a smaller bound$\}$.
Let $\mu_\gamma = \underset{q \in {\Cal J}_\gamma} \to \sup \mu_{\gamma,q}$. 
\newline
Let $q \Vdash_{{\Bbb P}_\gamma} ``{\underset\tilde {}\to {\Bbb I}_\gamma} = \{
{\underset\tilde {}\to I_{\gamma,\zeta}}:\zeta < 
{\underset\tilde {}\to \zeta_{\gamma,q}} \le \mu_{\gamma,q}\}"$ for
$q \in {\Cal J}_\gamma$ and let
${\Cal J}_{\gamma,\zeta} = \{q \in {\Cal J}_\gamma:\mu_{\gamma,q} > \zeta
\text{ and } q \Vdash ``\text{Dom}({\underset\tilde {}\to I_{\gamma,\zeta}})
\text{ is } \le \lambda_{\gamma,q,\zeta}"$ 
and no $q'$ such that $q \restriction \gamma \le_{\text{pr}} q' \in
\Bbb P_\gamma$ forces a smaller bound$\}$ and let
${\underset\tilde {}\to I_{\gamma,\zeta}} 
\text{ be id}_{{\underset\tilde {}\to L_{\gamma,q,\zeta}}}$, 
so ${\underset\tilde {}\to L_{\gamma,q,\zeta}}$ is a
$\Bbb P_\gamma$-name of a ${\underset\tilde {}\to \kappa_{\gamma,q}}$-directed 
quasi order on some $\lambda' \le \lambda_{\gamma,q,\zeta}$ 
(but $\Vdash_{{\Bbb P}_\gamma}$ ``if
$|{\underset\tilde {}\to {\Bbb I}_\gamma}| \le \zeta < \mu_\gamma$ then let
$L_{\gamma,\zeta}$ be trivial").  We can assume $L_{\gamma,q,\zeta}$ is a
quasi order on $\lambda_{\gamma,q,\zeta}$ (putting every $\beta \in
\lambda_{\gamma,q,\zeta} \backslash \text{ Dom}(I_{\gamma,q,\zeta})$ at the
bottom.

For $q \in {\Cal J}_\gamma$ let 
$L^*_{\gamma,q,\zeta}$ be ap$_{\kappa_{\gamma,q,q}}(
{\underset\tilde {}\to L_{\gamma,\zeta}})$ for the forcing notion \nl
$\Bbb P_\gamma^{[q]} = \{p \in \Bbb P_\gamma:q \restriction \gamma 
\le^{{\Bbb P}_\gamma}_{pr} p\}$ from Definition \scite{3.9}, so it is
defined in $\bold V$.  So by Claim \scite{3.10}
\medskip
\roster
\widestnumber\item{$(iii)$}
\item "{$(i)$}"   $L^*_{\gamma,q,\zeta}$ is $\kappa_{\gamma,q}$-directed
partial order on $[\lambda_{\gamma,q,\zeta}]^{< \kappa_{\gamma,q}}$
\sn
\item "{$(ii)$}"  $|L^*_{\gamma,q,\zeta}| \le (\lambda_{\gamma,q,\zeta})
^{< \kappa_{\gamma,q}}$
\sn
\item "{$(iii)$}"  $q \restriction \gamma \Vdash_{{\Bbb P}_\gamma} 
``{\underset\tilde {}\to I_{\gamma,\zeta}} = 
\text{ id}_{\underset\tilde {}\to L_{\gamma,q,\zeta}} 
\le_{\text{RK}} \text{ id}_{L^*_{\gamma,q,\zeta}}"$.
\endroster
\medskip

\noindent
Let $\kappa_\beta = \sup\{\kappa_{\gamma,q}:\gamma < \beta \text{ and } q \in
{\Cal J}_\gamma\}$. 
\sn
Let $\Bbb I^*_\beta$ be the $(< \kappa_\beta)$-closure of $\{
\text{id}_{L^*_{\gamma,q,\zeta}}:\gamma < \beta,q \in {\Cal J}_\gamma,
\zeta < \mu_{\gamma,q}\}$ (see Definition \scite{3.13}(1)).
\smallskip

Let $\bar N = \langle N_\eta:\eta \in (T^*,\bold I) \rangle$ be a
strict truely
$(\Bbb I^*_\alpha,\bold W)$-suitable tree of models for $(\chi,x),x$ coding 
enough information (so $\bar{\Bbb Q},
\Bbb I^*_\alpha,\bold W,W \in N_{\langle \rangle}$); why
truely? see \sciteu{4.x}.

Let ${\Cal T}_{\bar N}$ be the set of quadruples 
$(\gamma,q,\underset\tilde {}\to \nu,\underset\tilde {}\to T)$ such that:
\medskip
\roster
\item "{$\bigotimes_1$}"  $\gamma \le \alpha,q \in \Bbb P_\gamma,
\underset\tilde {}\to T$ is a $\Bbb P_\gamma$-name of a subtree of $T^*$, 
\newline
$q \Vdash_{({\Bbb P}_\gamma,\le_{\text{pr}})} 
``{\underset\tilde {}\to \kappa^*_\gamma} = \kappa_\gamma"$ and \nl
$q \Vdash_{{\Bbb P}_\gamma} ``\langle N_\eta[{\underset\tilde {}\to
G_{{\Bbb P}_\gamma}}]:
\eta \in (\underset\tilde {}\to T,\bold I \restriction 
\underset\tilde {}\to T) \rangle$ \nl

$\qquad \qquad$ is strictly 
$((\Bbb I^*_\alpha)^{[\kappa_\gamma]},\bold W)$-suitable tree, \nl

$\qquad \qquad N_{\langle \rangle}
[{\underset\tilde {}\to G_{{\Bbb P}_\gamma}}] \cap \omega_1 = N_{\langle \rangle}
\cap \omega_1$ and $\underset\tilde {}\to \nu \in \underset\tilde {}\to T$ 
\newline

$\qquad \qquad$ and $\gamma,\kappa \in N_{\underset\tilde {}\to \nu}
[{\underset\tilde {}\to G_{{\Bbb P}_\gamma}}]$ \newline

$\qquad \qquad$ and $\bar N[{\underset\tilde {}\to G_{{\Bbb P}_\gamma}}]$ has
$(\kappa)$-covering".
\endroster
\medskip
\noindent
Now ${\Cal T}'_{\bar N}$ is defined similarly as the set of quadruples
$(\underset\tilde {}\to \gamma,q,\underset\tilde {}\to \nu,
\underset\tilde {}\to T)$ such that:
$\underset\tilde {}\to \gamma$ is a simple $(\bar Q,W)$-named 
$[0,\alpha)$-ordinal, $q \in P_{\underset\tilde {}\to
\gamma},\underset\tilde {}\to \nu$ a $\Bbb P_{\underset\tilde {}\to \gamma}$-name and 
$\underset\tilde {}\to \gamma \in
N_{\underset\tilde {}\to \nu}[G_{{\Bbb P}_{\underset\tilde {}\to \gamma}}]$.  
(I.e. if $\zeta < \beta,G_{{\Bbb P}_\zeta} \subseteq \Bbb P_\zeta$ is generic over 
$\bold V$ and $\zeta = \underset\tilde {}\to \gamma
[G_{{\Bbb P}_\zeta}]$ then $r \in q \Rightarrow
{\underset\tilde {}\to \zeta_r}[G_{{\Bbb P}_\zeta}] < \zeta$, i.e. is well defined
$< \zeta$ \ub{or} is forced ($\Vdash_{{\Bbb P}_\alpha/G_{{\Bbb P}_\zeta}}$) to be not 
well defined), and $q \Vdash_{{\Bbb P}_{\underset\tilde {}\to \gamma}} 
``\underset\tilde {}\to \nu \in \text{ lim}(T)"$.
\medskip

We consider the statements, for $\gamma \le \beta < \alpha$
\medskip
\roster
\item "{$\boxtimes_{\gamma,\beta}$}"  for any 
$(\gamma,q,\underset\tilde {}\to \eta,\underset\tilde {}\to T) \in 
{\Cal T}_{\bar N}$ and $\underset\tilde {}\to \rho$ a $\Bbb P_\gamma$-name 
such that \newline
$q \Vdash_{P_\gamma} ``\underset\tilde {}\to \eta \triangleleft 
\underset\tilde {}\to \rho \in \underset\tilde {}\to T$ and $\gamma \in
N_{\underset\tilde {}\to \rho}[{\underset\tilde {}\to G_\gamma}]"$ \nl
and ${\underset\tilde {}\to p'}$ a $\Bbb P_\gamma$-name such that 
$q \Vdash_{{\Bbb P}_\gamma} ``{\underset\tilde {}\to p'}[G_{{\Bbb P}_\gamma}] \in
N_{\underset\tilde {}\to \rho}[{\underset\tilde {}\to G_{{\Bbb P}_\gamma}}] 
\cap \Bbb P_\beta/G_{{\Bbb P}_\gamma}$ and
$({\underset\tilde {}\to p'}[G_{{\Bbb P}_\gamma}]) \restriction \gamma \le_{pr} q"$
and ${\underset\tilde {}\to p'}[G_{{\Bbb P}_\gamma}]$ forces
$(\Vdash_{(P_\beta,\le_{\text{pr}})})$ a value $\kappa^*_\gamma$ to
${\underset\tilde {}\to \kappa_\gamma}$ (usually redundant) 
\ub{there is} $(\beta,q',
\underset\tilde {}\to \rho,{\underset\tilde {}\to T'})
\in {\Cal T}_{\bar N}$ such that 
${\underset\tilde {}\to p'} \le_{pr} q'$ \nl
(i.e., $p \Vdash_{{\Bbb P}_\gamma} ``{\underset\tilde {}\to p'}
[{\underset\tilde {}\to G_{{\Bbb P}_\gamma}}] \le_{pr} q"$) and 
$q' \restriction \gamma = q$ and $q' \Vdash_{{\Bbb P}_\beta} 
``\underset\tilde {}\to \rho \in {\underset\tilde {}\to T'} \subseteq
{\underset\tilde {}\to T}$".
\ermn
For simple $(\bar{\Bbb Q},W)$-names $[0,\alpha)$-ordinals
$\underset\tilde {}\to \gamma \le \underset\tilde {}\to \beta$ we define
\mr
\item "{$\boxtimes_{\underset\tilde {}\to \gamma,\underset\tilde {}\to
\beta}$}"
similarly $(\forall \beta < \beta^*)\forall \gamma \le
\beta(\square_{\gamma,\beta})$ and $\underset\tilde {}\to
\gamma,\underset\tilde {}\to \beta$.
\endroster
\bn
\ub{Observation}:  If $\forall < \beta^*,\forall \gamma \le
\beta(\boxtimes_{\gamma,\beta})$
and $\underset\tilde {}\to \gamma,\underset\tilde {}\to \beta$ are simple
$\bar Q$-named $[0,\beta)$-ordinals $\Vdash \underset\tilde {}\to \gamma
\le \underset\tilde {}\to \beta < \beta^*$ then
$\boxtimes_{{\underset\tilde {}\to \gamma^*},
{\underset\tilde {}\to \beta^*}}$
(defined naturally).
\enddemo
\bigskip

\demo{Proof}  By induction on the depth of $\underset\tilde {}\to \beta$
(see \scite{si.7}, fact A).

We prove by induction on $\beta \le \alpha$ that
\mr
\item "{$(a)$}"  $\Bbb P_\beta$ has pure 
${\underset\tilde {}\to \kappa_\beta}$-covering; i.e. if
$\underset\tilde {}\to \tau$ is a $\Bbb P_\beta$-name of an ordinal
$< {\underset\tilde {}\to \kappa_\beta}$ and $p \in \Bbb P_\beta$
\ub{then} for some $q$ and $a$ we have: $p \le_{\text{pr}} q \in \Bbb P_\beta,
a \in \bold V$ is a set of ordinals and $q \Vdash ``|a| <
{\underset\tilde {}\to \kappa_\beta} \and \underset\tilde {}\to \tau \le a"$
(even over $\Bbb P_\gamma$)
\sn
\item "{$(b)$}"  $\Bbb P_\beta$ has pure $(\aleph_1,\aleph_1)$-decidability
\sn
\item "{$(c)$}"  for every $\gamma \le \beta$ we have 
$\boxtimes_{\gamma,\beta}$ (but for \scite{5.5}(3) we have to
restrict ourselves to non-limit $\gamma$).
\ermn
Note that for $\gamma = \beta$ the statement in clause (c)
is trivial hence we shall consider only $\gamma < \beta$.
\bigskip

\noindent
\underbar{Case 1}:  $\beta = 0$.

Trivial.
\bigskip

\noindent
\underbar{Case 2}:  $\beta$ a successor ordinal.

Clauses (a), (b) follows easily from clause (c) so let us 
concentrate on clause (c).
As trivially $\boxtimes_{\gamma_0,\gamma_1} \and \boxtimes_{\gamma_1,\gamma_2}
\Rightarrow \boxtimes_{\gamma_0,\gamma_2}$, clearly without loss of
generality \nl
$\beta = \gamma + 1$.

Let $G_{{\Bbb P}_\gamma}$ be such that $q \in G_{{\Bbb P}_\gamma} 
\subseteq \Bbb P_\gamma$ and
$G_{{\Bbb P}_\gamma}$ generic over $\bold V$.

Let $T' = \{\nu:\rho \char 94 \nu \in \underset\tilde {}\to T[G_{{\Bbb P}_\gamma}]\},
\bar N' = \langle N'_\nu:\nu \in
(T',\bold I') \rangle$ where $N'_\nu = N_{\rho \char 94 \nu}[G_{{\Bbb P}_\gamma}],
\bold I'_\nu = \bold I^*_{\rho \char 94 \nu}$. \newline
By \scite{5.2} applied to $\bar N'$ we can find ${\underset\tilde {}\to p'},
{\underset\tilde {}\to T''}$ as required.
\bigskip

\noindent
\underbar{Case 3}:  $\beta$ a limit ordinal. \nl
[Saharon note: it would be if
$\underset\tilde {}\to \kappa$ is a $(\Bbb P_\beta,\le_{\text{pr}^+})$-name,
$\gamma < \beta,G_\gamma \subseteq \Bbb P_\gamma$ generic we should define
$\underset\tilde {}\to \kappa/G_\gamma$ anyhow we can use real names but
then ${\underset\tilde {}\to \kappa^*_\gamma}$ is just a
$(\Bbb P_\gamma,\le_{\text{pr}})$-name.  But if 
${\underset\tilde {}\to \kappa_\beta}$ are real cardinals no problem.  But see
clause (g) of definition of UP".]
\enddemo
\bigskip

\demo{Proof of Clause (a)}  If we use the c.c. version: easier, hardest 
case is $\kappa^*_\beta = \beta$, so $\beta$ strongly inaccessible. \nl
Note that we have to prove the weak version.  If the property fails for
$\beta,p,\underset\tilde {}\to \tau$ (so $p \Vdash_{P_\beta}
``\underset\tilde {}\to \tau$ an ordinal", etc.), \ub{then} by the
induction hypothesis
$\beta$ is minimal so definable in $({\Cal H}(\chi),\in)$ from 
$\bar{\Bbb Q}$ hence necessarily
$\beta \in N_{<>}$ and \wilog \, $p,\underset\tilde {}\to \tau \in N_{<>}$.
Also \wilog \, $p \Vdash_{({\Bbb P}_\beta,\le_{\text{pr}})}
``{\underset\tilde {}\to \kappa^*_\beta} = \kappa^*_\beta"$ for some
$\kappa^*_\beta \in$ UReg V.

We shall now choose $p_1$ as in the proof of \scite{1.27}.  Let
$\langle {\underset\tilde {}\to \zeta_\varepsilon}:\varepsilon < j \rangle$
be a witness for $p$ (see \scite{1.13} clause (F), in particular
$(a)(iii)$ of (F)), so \wilog \, $j \le \omega$.  For each 
$\varepsilon < j$ and $\xi < \beta$, let
${\underset\tilde {}\to a_{\varepsilon,\xi}}$ be a
$\Bbb P_{\xi +1}$-name of a set of $< \kappa^*_\beta$ ordinals and let
${\underset\tilde {}\to r_{\varepsilon,\xi}}$ be
a $\Bbb P_{\xi +1}$-name of a member of 
$\Bbb P_\beta/G_{\xi +1}$ with domain $\subseteq [\xi +1,\beta)$ 
such that if
$G_{\xi +1} \subseteq \Bbb P_{\xi +1}$ is generic over $\bold V$ and
${\underset\tilde {}\to \zeta_\varepsilon}[G_{\xi +1} \cap
\Bbb P_\xi] = \xi$, then $r = {\underset\tilde {}\to r_\varepsilon}
[{\underset\tilde {}\to G_{\xi +1}}]$ satisfies: 
\mr
\item "{$(a)$}"  if possible 
$p \le_{\text{pr}}(p \restriction (\xi +1)) \cup r \in \Bbb P_\beta$ and
$(p \restriction (\xi +1)) \cup r \Vdash_{{\Bbb P}_\beta} 
``\underset\tilde {}\to \tau \in 
a"$ where $a = {\underset\tilde {}\to a_{\varepsilon,\zeta}}[G_{\xi +1}] \in
\bold V[G_{\xi +1}]$ is a set of $< \kappa^*_\beta$ ordinals and
${\underset\tilde {} \to {\bold t}_{\varepsilon,\zeta}}[G_{\xi +1}] =$ true
\sn
\item "{$(b)$}"  if not possible $r = 
{\underset\tilde {}\to r_{\varepsilon,\xi}}[G_{\xi +1}]$ is the empty function
and ${\underset\tilde {}\to a_{\varepsilon,\xi}}
[G_{\xi +1}]  = \emptyset$ and ${\underset\tilde {}\to {\bold
t}_{\varepsilon,\xi}}[G_{\xi +1}] =$ false.
\nl
Also we can demand that:
\sn
\item "{$(c)$}"  ${\underset\tilde {}\to \zeta_{\varepsilon_1}}[G_\xi] = \xi 
= {\underset\tilde {}\to \zeta_{\varepsilon_2}}[G_\xi]$ then
${\underset\tilde {}\to r_{\varepsilon_1,\xi}}[G_{\xi +1} \cap \Bbb P_\xi] =
{\underset\tilde {}\to r_{\varepsilon_2,\xi}}[G_\xi]$. 
\ermn
Let ${\underset\tilde {}\to r_\varepsilon}[G_\beta] = r$
iff for some $\xi$ 
we have ${\underset\tilde {}\to \zeta_\varepsilon}[G_\beta \cap \Bbb P_\xi] = \xi$
and $r = {\underset\tilde {}\to r_{\varepsilon,\xi}}[G_\beta]$, similarly we
define ${\underset\tilde {}\to a_\varepsilon}$.  Let $p_1 = p \cup \bigcup
\{{\underset\tilde {}\to r_\varepsilon}:\varepsilon < j\}$.  Clearly
$p \le_{\text{pr}} p_1 \in \Bbb P_\beta$.
\mn
Next define $p_2$ as in \scite{1.27}: (recall (g) + (h) of Definition
\scite{5.3}).  I.e. for each $\varepsilon < j,\xi < \beta,G_\xi \subseteq
\Bbb P_\xi$ generic over $\bold V,p \restriction \xi \in G_\xi,
{\underset\tilde {}\to r_\varepsilon}[G_\xi] = \xi$ there are
${\underset\tilde {}\to r'_{\varepsilon,\xi}},a'_{\varepsilon,\xi} \in
\bold V[G_\beta]$ such that 
${\underset\tilde {}\to {\hat{\Bbb Q}}_\xi}[G_\xi] \models ``p_1
\restriction \{\xi\} \le_{\text{pr}} r^1_{\varepsilon,\xi}"$ and
$r'_{\varepsilon,\xi} \Vdash_{\underset\tilde {}\to {\hat{\Bbb Q}}_\xi}
``{\underset\tilde {}\to a_{\varepsilon,\xi}} \subseteq 
a'_{\varepsilon,\xi}"$.
So really \wilog \, we have $\Bbb P_\xi$-names
${\underset\tilde {}\to r'_{\varepsilon,\xi}},
{\underset\tilde {}\to a'_{\varepsilon,\xi}}$ such that
${\underset\tilde {}\to \zeta_{\varepsilon_1}}[G_\xi] = \xi =
{\underset\tilde {}\to \zeta_{\varepsilon_2}}[G_\xi]$ implies
${\underset\tilde {}\to r'_{\varepsilon_1,\xi}}[G_\xi] =
r'_{\varepsilon_2,\xi}[G_\xi]$ and
${\underset\tilde {}\to a'_{\varepsilon_1,\xi}}[G_\xi] =
{\underset\tilde {}\to a'_{\varepsilon_1,\xi}}[G_\xi]$.  We define
${\underset\tilde {}\to r'_\varepsilon},
{\underset\tilde {}\to a'_\varepsilon}$ by:
${\underset\tilde {}\to r'_\varepsilon}[G_\beta] = r$ iff for some $\xi <
\beta$ we have ${\underset\tilde {}\to \zeta_\varepsilon}[G_\beta \cap
\Bbb P_\xi] = \xi$ and 
$r = {\underset\tilde {}\to r_{\varepsilon,\xi}}[G_\beta
\cap \Bbb P_\xi]$ and similarly ${\underset\tilde {}\to
a'_\varepsilon}$.  
Now let $p_2 = p_1 \cup
\{{\underset\tilde {}\to r'_\varepsilon}:\varepsilon < j\}$ so $p_1 
\le_{\text{pr}} p_2$.  We can finish as in the proof of \scite{1.27} [fill]!!!
\bn
\ub{Clause (b)}:  

As in the proof of clause (a) \wilog \, we have $\beta,p,\underset\tilde {}\to
\zeta \in N_{<>}$. We define also $p_1$ as in the proof of clause (a) 
trying to force a countable
bound for $\underset\tilde {}\to \zeta$.  Let $\langle
{\underset\tilde {}\to \zeta^*_\varepsilon}:\varepsilon < \omega \rangle$ be
a witness for $p \in \Bbb P_\beta$ (see Definition \scite{1.13},
clause (F) in particular (a)(iii) of (F) and \wilog \, $\langle
{\underset\tilde {}\to \zeta^*_\varepsilon}:n < \omega \rangle$ belong
to $N_{<>}$).  We now choose
by induction on $n$ a quadruple $({\underset\tilde {}\to \gamma_n},
{\underset\tilde {}\to q_n},{\underset\tilde {}\to \nu_n},
{\underset\tilde {}\to t_n})$ such that
\mr
\widestnumber\item{$(iii)$}
\item "{$(i)$}"  ${\underset\tilde {}\to \gamma_n}$ is a simple
$(\bar{\Bbb Q} \restriction \beta,W)$-named $[0,\beta)$-ordinal
\sn
\item "{$(ii)$}"  ${\underset\tilde {}\to \gamma_0} = 0,
{\underset\tilde {}\to \gamma_n} < {\underset\tilde {}\to \gamma_{n+1}}$
\sn
\item "{$(iii)$}"  $\zeta_{\underset\tilde {}\to \zeta^*_n} + 1 \le
{\underset\tilde {}\to \gamma_{n+1}}$
\sn
\item "{$(iv)$}"  $({\underset\tilde {}\to \gamma_n},
{\underset\tilde {}\to q_n},{\underset\tilde {}\to \nu_n},
{\underset\tilde {}\to T_n}) \in
{\Cal T}'_{\bar N}$
\sn
\item "{$(v)$}"  ${\underset\tilde {}\to q_n} =
{\underset\tilde {}\to q_{n+1}} \restriction {\underset\tilde {}\to \gamma_n}$
\sn
\item "{$(vi)$}"  $p \restriction {\underset\tilde {}\to \gamma_n}
\le_{\text{pr}} {\underset\tilde {}\to q_n}$.
\ermn
No problem to carry it by the observation above.
\mn
Let $p_2 = \dbcu_{n < \omega} {\underset\tilde {}\to q_n} \cup p$, clearly
$p \le_{\text{pr}} p_2,{\underset\tilde {}\to q_n} \le p_2$, and so it
is enough to prove $p_2 \Vdash ``\underset\tilde {}\to \zeta < N_{<>} \cap
\omega_1"$.  So let $p_2 \in G_\beta$ with $G_\beta$ a subset of $\Bbb P_\beta$ generic
over $\bold V$.  So there is $p^+_2 \in G_\beta$ satisfying $p \le p_2
\le p^+_2$ such that 
$p^+_2 \Vdash ``\underset\tilde {}\to \zeta = \gamma^* < \omega_1"$ and so
$p \le p^+_2$ so \wilog \, $p \le p^+_2$ above $\xi_0 < \ldots <
\xi_{m-1}$, using \scite{1.15}.
There is $n$ such that $[{\underset\tilde {}\to \gamma_n}[G_\beta],
\dbcu_{\ell < \omega} {\underset\tilde {}\to \gamma_\ell}[G_\beta])$ is
disjoint to $\{\xi_0,\dotsc,\xi_{m-1}\}$, hence for some $\varepsilon <
\omega$, letting $\xi = {\underset\tilde {}\to \zeta_\varepsilon}[G_\beta]$
defining ${\underset\tilde {}\to r^0_{\varepsilon,\xi}}[G_\beta \cap
\Bbb P_{\xi +1}]$ we get ${\underset\tilde {}\to {\bold t}^0_{\varepsilon,\xi}} =$
truth (see? \scite{1.27}) and $\xi < {\underset\tilde {}\to \gamma_n}
[G_\beta]$.

Now consider $N' = N_{<>}[G_\beta \cap \Bbb P_{\xi +1}]$ we know $N' \cap \omega_1
= N_{<>} \cap \omega_2$ (by clause (iv) above), and in it we have $p
\restriction (\xi +1) \cup {\underset\tilde {}\to r^0_{\varepsilon,\xi}}
[G_\beta \cap p_{\xi +1}]$ forces a bound to $\underset\tilde {}\to \tau$,
but the condition is $\le_{\text{pr}} p_1 \le_{\text{pr}} p_2$, and it
belongs to $N'$, so the value is $< N_{<>} \cap \omega_1$ and
we are done.
\enddemo
\bn
\ub{Clause (c)}:  

By \scite{5.2} it suffices to prove
\mr
\item "{$\bigotimes_2$}"  there are $r,\underset\tilde {}\to \eta$ such that:
\nl
$\underset\tilde {}\to \eta$ is a $\Bbb P_\beta$-name, $r \in \Bbb P_\beta,
r \restriction \gamma = q \restriction \gamma,p' \le_{pr} q$ and
$r \Vdash_{{\Bbb P}_\beta} ``\underset\tilde {}\to \eta \in \text{ lim}
(\underset\tilde {}\to T)$ and
$N_{{\underset\tilde {}\to \eta} \restriction \ell}
[{\underset\tilde {}\to G_{{\Bbb P}_\beta}}] \cap \omega_1 = N_{\langle \rangle}
\cap \omega_1$ and for every $y \in \bold V \cap \dbcu_{\ell < \omega}
N_{{\underset\tilde {}\to \eta} \restriction \ell}[G_{{\Bbb P}_\beta}]$ for some
$A \in \bold V \cap \dbcu_{\ell < \omega} N_{{\underset\tilde {}\to \eta}
\restriction \ell}$ we have $|A|^{\bold V} < {\underset\tilde {}\to \kappa^*_\beta}
[G_{{\Bbb P}_\beta}]$ and $y \in A"$,
\ermn

We shall choose by induction on $n < \omega,{\underset\tilde {}\to \gamma_n},
q_n,{\underset\tilde {}\to \rho_n},{\underset\tilde {}\to T_n},
{\underset\tilde {}\to k_n},{\underset\tilde {}\to p_n}$ such that:
\medskip
\roster
\item "{$(a)$}"  $({\underset\tilde {}\to \gamma_n},q_n,
{\underset\tilde {}\to \rho_n},{\underset\tilde {}\to T_n}) \in 
{\Cal T}'_{\bar N}$ \newline
(so ${\underset\tilde {}\to \gamma_n}$ is a simple $\bar{\Bbb Q}$-named ordinal)
\sn
\item "{$(b)$}"  ${\underset\tilde {}\to k_n}$ is a
$\Bbb P_{\underset\tilde {}\to \gamma_n}$-name of a natural number
\sn
\item "{$(c)$}"  ${\underset\tilde {}\to \rho_n}$ is a 
$\Bbb P_{\underset\tilde {}\to \gamma_n}$-name
\sn
\item "{$(d)$}"  $q_n \Vdash_{{\Bbb P}_{\underset\tilde {}\to \gamma_n}}
``{\underset\tilde {}\to \rho_n} \in {\underset\tilde {}\to T_n} \cap
{}^{\underset\tilde {}\to k_n}\text{Ord}"$
\sn
\item "{$(e)$}"  ${\underset\tilde {}\to \gamma_0} = \gamma$ and
$\Vdash_{\bar{\Bbb Q}}``{\underset\tilde {}\to \gamma_n} < 
{\underset\tilde {}\to \gamma_{n+1}} < \beta$ and
${\underset\tilde {}\to \gamma_{n+1}} \text{ non-limt}"$ \newline
i.e., if $\zeta < \beta$ and $G_{{\Bbb P}_\gamma}
\subseteq \Bbb P_\zeta$ is generic over $\bold V$ and $\zeta = {\underset\tilde {}\to 
\gamma_n}[G_{{\Bbb P}_\zeta}]$ then \nl
$r \in q_n \Rightarrow {\underset\tilde {}\to \zeta_n}[G_\zeta] < \zeta$ 
(i.e., is well defined $< \zeta$ or is forced to be not well defined),
\sn
\item "{$(f)$}"  $q_{n+1} \restriction {\underset\tilde {}\to \gamma_n} =
q_n$
\sn
\item "{$(g)$}"  $q_{n+1} 
\Vdash_{{\Bbb P}_{\underset\tilde {}\to \gamma_{n+1}}} 
``{\underset\tilde {}\to \rho_n} \triangleleft
{\underset\tilde {}\to \rho_{n+1}}$, so 
${\underset\tilde {}\to k_n} < {\underset\tilde {}\to k_{n+1}}$ and
${\underset\tilde {}\to T_{n+1}} \subseteq {\underset\tilde {}\to T_n}"$
\sn
\item "{$(h)$}"  ${\underset\tilde {}\to p_n}$ is a
$\Bbb P_{\underset\tilde {}\to \gamma_n}$-name, $p_0 = p,
{\underset\tilde {}\to p_n} \restriction {\underset\tilde {}\to \gamma_n}
\le_{\text{pr}} q_n$ and \newline
$q_n \Vdash_{{\Bbb P}_{\underset\tilde {}\to \gamma_n}} ``
{\underset\tilde {}\to p_n} \in N_{\underset\tilde {}\to \rho_n}
[G_{{\Bbb P}_{\underset\tilde {}\to \gamma_n}}] \cap \Bbb P_\beta$ and 
$p_n \restriction {\underset\tilde {}\to \gamma_n} \in G_{{\Bbb P}_
{\underset\tilde {}\to \gamma_n}}"$
\sn
\item "{$(i)$}"  $q_n \Vdash_{{\Bbb P}_{\underset\tilde {}\to \gamma_n}}
``{\underset\tilde {}\to p_n} \le^{{\Bbb P}_\beta}_{\text{pr}}
{\underset\tilde {}\to p_{n+1}} \in
N_{{\underset\tilde {}\to \rho_{n+1}}}[G_{{\Bbb P}_{\gamma_n}}] 
\cap \Bbb P_\beta"$
\sn
\item "{$(j)$}"  letting $\langle {\underset\tilde {}\to \tau_{\nu,\ell}}:
\ell < \omega \rangle$ list the $\Bbb P_\beta$-names of ordinals from $N_\nu$: 
for $m,\ell \le n$ we have:
\endroster
\medskip
$$
\align
q_n \Vdash_{{\Bbb P}_{\underset\tilde {}\to \gamma_n}} ``{\underset\tilde {}\to 
p_{n+1}} \text{ force that: } &a) \,\,\,\text{if } 
{\underset\tilde {}\to \tau_{{\underset\tilde {}\to \rho_n} 
\restriction m,\ell}} 
\text{ is a countable ordinal}, m \le {\underset\tilde {}\to k_n} \\
  &\quad \text{ then it is smaller than some }  
{\underset\tilde {}\to \tau'_{{\underset\tilde {}\to \rho_n} \restriction
m,\ell}} \in N_{\underset\tilde {}\to \rho_{n+1}}[
G_{{\Bbb P}_{\underset\tilde {}\to \gamma_{n+1}}}], \\
  &\quad \text{ a } 
P_{\underset\tilde {}\to \gamma_n} \text{-name of a countable ordinal} \\
  &b) \,\,\, \text{for some } 
A \in \bold V \cap N_{{\underset\tilde {}\to \eta_{n+1}} 
\restriction {\underset\tilde {}\to k_{n+1}}}, \\
  &\quad \, |A|^{\bold V} < \kappa^*_\beta \text{ and }
{\underset\tilde {}\to \tau_{{\underset\tilde {}\to \rho_n} \restriction
m,\ell}} \in A".
\endalign
$$
\medskip

\noindent
The induction is straight (later we shall show that 
$\dsize \bigcup_{n < \omega} q_n$ and $\underset\tilde {}\to \eta = 
\dsize \bigcup_{n < \omega}
{\underset\tilde {}\to \rho_n}$ are as required in $\bigotimes_2$) 
by clauses (a) + (b) proved above.
\bn
\centerline {$* \qquad * \qquad *$}
\bn
Because we need and have $(*)_1$ or $(*)_2 + (*)_3$ below:
\medskip
\roster
\item "{$(*)_1$}"  Assume $\le_{pr},\le_{vpr}$ are equal to $\le$ \nl
(i.e., $\Vdash_{{\Bbb P}_\beta} 
``\le^{\underset\tilde {}\to {\Bbb Q}_\beta}_{vpr}$ is
$\le^{\underset\tilde {}\to {\Bbb Q}_\beta}"$ for each $\beta < \alpha$),
if $p \in \Bbb P_\beta,\gamma < \beta,\underset\tilde {}\to \tau$
a $\Bbb P_\beta$-name of an ordinal \ub{then} there are $p',
{\underset\tilde {}\to \tau'}$ such that:
{\roster
\itemitem{ (i) }  ${\underset\tilde {}\to \tau'}$ 
is a $\Bbb P_\gamma$-name of an ordinal
\sn
\itemitem{ (ii) }  $p \le_{\text{pr}} 
p' \in \Bbb P_\beta$ and $p \restriction \gamma =
p' \restriction \gamma$
\sn
\itemitem{ (iii) }  $p' \Vdash_{{\Bbb P}_\beta} ``
\underset\tilde {}\to \tau = {\underset\tilde {}\to \tau'}"$. \nl
[why?  straight by \scitet{1.16}].
\endroster}
\sn
[Saharon: maybe below we are stuck with $\underset\tilde {}\to \zeta \in
[\gamma,\beta)$, but this suffices - need to change?]
\sn
\item "{$(*)_2$}"  \ub{old proof of clause (b)}:  if 
$p \in \Bbb P_\beta,\gamma < \beta,
\underset\tilde {}\to \tau$ is a $\Bbb P_\beta$-name of a countable ordinal,
\ub{then} there are $p',\tau'$ such that
{\roster
\itemitem{ (i) }  $\tau'$ is a $\Bbb P_\gamma$-name of a countable ordinal
\sn
\itemitem{ (ii) }  $p \le_{pr} p' \in \Bbb P_\beta$ and $p \restriction \gamma =
p' \restriction \gamma$
\sn
\itemitem{ (iii) }  $p' \Vdash_{{\Bbb P}_\beta} ``\underset\tilde {}\to \tau \le
\tau'"$, 
\endroster}
\sn
[why $(*)_2$?  let $\underset\tilde {}\to \zeta$ be the following simple
$\bar{\Bbb Q}$-named $[\gamma,\beta)$-ordinal: \nl
for $G_\zeta \subseteq \Bbb P_\zeta$ is generic over $\bold V$ for 
$\zeta \in [\gamma,\beta)$ we let
$\underset\tilde {}\to \zeta[G_\zeta] = \zeta$ if 
{\roster
\itemitem{ $(a)$ }  $p \restriction \zeta \notin G_\zeta$ \ub{or}: for some 
$p' \in \Bbb P_\beta$ we have $p' \restriction \zeta = p \restriction \zeta$ 
and \nl
$\Bbb P \models p \le_{\text{pr}} p'$ and 
$p' \Vdash_{{\Bbb P}_\beta/G_\zeta} ``\underset\tilde {}\to \tau < \tau^*"$ 
for some countable ordinal $\tau^*$ 
\sn
\itemitem{ $(b)$ }  for no $\xi \in [\gamma,\zeta)$ does clause (a) hold for
$\xi,G_\zeta \cap \Bbb P_\xi$. 
\endroster}
Now if for some $\gamma \in [\alpha,\beta)$ we have 
$p \Vdash_{{\Bbb P}_\alpha} ``\underset\tilde {}\to \zeta = \gamma"$ we are
done.  Also $\Vdash_{{\Bbb P}_\alpha} ``\underset\tilde {}\to \zeta
[{\underset\tilde {}\to G_{{\Bbb P}_\alpha}}]$ is well defined" as if $p \in
G_\alpha \subseteq \Bbb P_\alpha$ and $G_\alpha$ is generic over
$\bold V$, then for
some $q \in G_\alpha$ and countable ordinal $\tau^*$ we have
$q \Vdash ``\underset\tilde {}\to \tau = \tau^*$.  By the definition of
$\aleph_1-\text{Sp}_e(W)$-iteration for some $\zeta \in
[\gamma,\beta)$ we have
$\xi \in [\zeta,\beta) \Rightarrow [p \restriction \{\xi\}
\le^{\underset\tilde {}\to {\Bbb Q}_\xi}_{pr} q \restriction \{\xi\}$ or
$e = 4 \and p \restriction \{\xi\}$ not defined[?]]. \nl
Define $p'$ by: $p' \restriction \underset\tilde {}\to \zeta = 
p \restriction \underset\tilde {}\to \zeta$, and for
$\xi \in [\underset\tilde {}\to \zeta,\beta)$ we let 
$p' \restriction \{\xi\}$ be $q \restriction \{\xi\}$ if: 
$p \restriction \{\xi\} \le^{\underset\tilde {}\to {\Bbb Q}_\xi}_{pr}
q \restriction \{\xi\}$ or $e = 4 \and p \restriction \{\xi\}$ not defined.
We shall show that $p'$ is as required, hence
really $\Vdash_{{\Bbb P}_\alpha} ``\underset\tilde {}\to \zeta \in [\gamma,\beta)$ is
well defined".
So there is a $\Bbb P_{\underset\tilde {}\to \zeta}$-name of
${\underset\tilde {}\to p'}$ as appearing in the definition of
$\underset\tilde {}\to \zeta$ and it is, essentially, a member of
$\Bbb P_\beta$.
Now as we have finite apure support, the proof of
``$\underset\tilde {}\to \zeta[{\underset\tilde {}\to G_{{\Bbb P}_\alpha}}]$ is well
defined" gives $\Vdash_{{\Bbb P}_\alpha} ``\underset\tilde {}\to \zeta$ is not a
limit ordinal $> \alpha$".  Lastly $\Vdash_{{\Bbb P}_\alpha} ``\underset\tilde {}\to 
\zeta$ is not a successor ordinal $> \gamma"$ is proved by the property of
each ${\underset\tilde {}\to {\Bbb Q}_\xi}$.]
\sn
\item "{$(*)_3$}"  \ub{old proof of clause (a)}:  if 
$p \in \Bbb P_\beta,\gamma < \beta,\underset\tilde {}\to 
\tau$ a $\Bbb P_\beta$-name of a set of $< \kappa$ ordinals, \ub{then} there are
$p',{\underset\tilde {}\to \tau'}$ such that:
\sn
{\roster
\itemitem { (i) }  ${\underset\tilde {}\to \tau'}$ is a $\Bbb P_\gamma$-name
of a set of $< \kappa$ ordinals
\sn
\itemitem{ (ii) }  $p \le_{pr} p' \in \Bbb P_\beta$
\sn
\itemitem{ (iii) }  $p' \Vdash_{{\Bbb P}_\beta} ``\underset\tilde {}\to \tau
\subseteq {\underset\tilde {}\to \tau'}$.
\endroster}
\sn
[Why?  Similar to the proof of $(*)_2$; note that it is automatic if the
${\underset\tilde {}\to \kappa_i}$'s increase fast enough].
\ermn
Finishing the induction we let
$\underset\tilde {}\to \eta = \dsize \bigcup_{n < \omega}
{\underset\tilde {}\to \rho_n}$ and we define $q_\omega \restriction
{\underset\tilde {}\to \gamma_n} = q_n$, \nl
$q_\omega \restriction [\dsize \bigcup_{n < \omega} 
{\underset\tilde {}\to \gamma_n},\beta)$ is defined as $\le_{vpr}$-upper 
bound of $\langle {\underset\tilde {}\to p_m} \restriction
[\dsize \bigcup_{n < \omega} \gamma_n,\beta):m < \omega \rangle$.
\medskip

More formally, let $\gamma^*,\beta$ and $G_{\gamma^*} \subseteq
\Bbb P_{\gamma^*}$ be such that: $G_{\gamma^*}$ is generic over $\bold V,\gamma^* =
\dsize \bigcup_{n < \omega} \gamma^*_n,\gamma^*_n = 
{\underset\tilde {}\to \gamma_n}[G_{\gamma^*}]$, let $p'_n =
{\underset\tilde {}\to p_n}[G_{\gamma^*} \cap \Bbb P_{\gamma^*_n}]$, let
$p'_n \restriction [\gamma^*,\alpha) = \{{\underset\tilde {}\to r^n_\zeta}:
\zeta < \zeta^*_n\}$ where ${\underset\tilde {}\to r^n_\zeta}$ is a simple
$[\gamma^*,\alpha)$-named atomic condition. \newline
Now we define ${\underset\tilde {}\to s^n_\zeta}$, a simple
$[\gamma^*,\alpha)$-named atomic condition as follows:
\medskip
\roster
\item "{$(a)$}"  $\zeta_{\underset\tilde {}\to s^n_\zeta} =
\zeta_{\underset\tilde {}\to r^n_\zeta}$
\sn
\item "{$(b)$}"  if $\zeta \in [\gamma^*,\alpha),G_{\gamma^*} \subseteq
G_\gamma \subseteq \Bbb P_\gamma,G_\gamma$ generic over $V,
\zeta_{\underset\tilde {}\to r^n_\zeta}[G_\zeta] = \zeta$ then
${\underset\tilde {}\to s^n_\zeta}[G_\zeta]$ is the $<^{*V[G_\zeta]}$-first
elements of ${\underset\tilde {}\to {\Bbb Q}_\zeta}[G_\zeta]$ which satisfies the
following:
{\roster
\itemitem{ $(*)(\alpha)$ }  
$({\underset\tilde {}\to p_n} \restriction \{\zeta\}) \le_{vpr} s$
\sn
\itemitem{ $(\beta)$ }  if 
$\langle \emptyset_{{\underset\tilde {}\to {\Bbb Q}_\zeta}[G_\zeta]} \rangle
\char 94 \langle ({\underset\tilde {}\to p_m} \restriction \{\xi\})
[G_\zeta]:m \in (n,\omega) \rangle$ has a $\le_{vpr}$-upper bound then
${\underset\tilde {}\to s^n_\zeta}[G_\zeta]$ is such upper bound.
\endroster}
\ermn
Now actually such $\le_{\text{vpr}}$-upper 
bound actually exists, and $q_\omega$
is as required.  \hfill$\square_{\scite{5.5}}$\margincite{5.5}
\bigskip

\proclaim{\stag{5.6} Claim}  Suppose 
$\bold W \subseteq \omega_1$ is stationary
and $\bar{\Bbb Q} = \langle \Bbb P_i,
{\underset\tilde {}\to {\Bbb Q}_i},\Bbb I_i,
{\underset\tilde {}\to \kappa_i}:i < \alpha \rangle$ is 
a ${\text{\rm UP\/}}^\ell(\bold W,W)$-suitable iteration 
(where $\ell \in \{4,5\}$), and 
$\Bbb P_\alpha = { \text{\rm Sp\/}}_6(W)-{\text{\rm Lim\/}}(\bar{\Bbb Q})$.

Each of the following is a sufficient condition for ``$\dsize \bigcup_{\beta <
\alpha} \Bbb P_\beta$ is a dense subset of \nl
$\Bbb P_\alpha$":
\medskip
\roster
\item "{$(A)_\alpha$}"  $\beta \in W$ is 
\footnote{remember $W$ is a parameter in the definition of $\kappa-Sp_e(W)$-
iteration} strongly inaccessible and
$\dsize \bigwedge_{\beta < \alpha}$ density $(\Bbb P_\beta) < \alpha$
\sn
\item "{$(B)_\alpha$}"  for every $i,
\le^{\underset\tilde {}\to {\Bbb Q}_i}_{\text{vpr}}$
is equality and $\bold V \models ``\text{cf}(\alpha) = \aleph_1"$ or at least
for some $\beta < \alpha$ we have 
$\Vdash_{{\Bbb P}_\beta} ``\text{cf}(\alpha) = \aleph_1"$.
\endroster
\endproclaim
\bigskip

\demo{Proof} \newline
\underbar{Case 1}:  $(A)_\alpha$ \newline
Straight by the definition of $\kappa-\text{Sp}_6(W)$-iteration 
(see \scite{1.13}).
\mn
\underbar{Case 2}:  $(B)_\alpha$ \newline
Follows by \scite{5.5}.  \hfill$\square_{\scite{5.6}}$\margincite{5.6}
\enddemo
\bigskip

\demo{\stag{5.7} Conclusion}  Assume $\bold W \subseteq 
\omega_1$ is stationary
and $\bar{\Bbb Q} = \langle \Bbb P_i,{\underset\tilde {}\to {\Bbb Q}_i},
{\underset\tilde {}\to {\Bbb I}_i},{\underset\tilde {}\to \kappa_i}:
i < \alpha \rangle$ is a
UP$^4(\bold W,W)$-iteration with Sp$_e(W)$-limit $\Bbb P_\alpha$.

If $\{\beta < \alpha:\bar{\Bbb Q} 
\restriction \beta$ satisfies $(A)_\beta \vee
(B)_\beta$ from \scite{5.6}$\}$ is a stationary subset of $\alpha$ and
$\beta < \alpha \Rightarrow \Bbb P_\beta$ satisfies the
cf$(\alpha)$-c.c. (e.g.
has cardinality or at least density $< \text{ cf}(\alpha)$), \ub{then}
$\Bbb P_\alpha$ satisfies the cf$(\alpha)$-c.c.
\enddemo
\bigskip

\demo{Proof}  Straight.

We may like to iterate up to e.g. the first inaccessible (we may below
weaken $|\Bbb P_\beta| < \alpha$ to 
$\Bbb P_\beta$ satisfies the $\alpha$-c.c. if
$\Bbb P_\alpha = \dbcu_{\beta < \alpha} \Bbb P_\beta$).
\enddemo
\bigskip

\proclaim{\stag{5.8} Claim}  [See \scite{si6.11}]??  Assume
\medskip
\roster
\item "{$(a)$}"  $\bold W \subseteq \omega_1$ is stationary
and $\bar{\Bbb Q} = \langle \Bbb P_i,{\underset\tilde {}\to {\Bbb Q}_i},
{\underset\tilde {}\to {\Bbb I}_i},{\underset\tilde {}\to \kappa_i}:
i < \alpha \rangle$ is a ${\text{\rm UP\/}}^4(\bold W,W)$-iteration 
$(y \in \{ a,b\})$ with ${\text{\rm Sp\/}}_6(W)$-limit $\Bbb P_\alpha$
\sn
\item "{$(b)$}"  $\Vdash_{{\Bbb P}_\beta}$ ``if in 
${\underset\tilde {}\to {\Bbb Q}_\beta},
\emptyset_{\underset\tilde {}\to {\Bbb Q}_\beta} \ne p
\le_{\text{vpr}} q$, then
$p = q$ and $\emptyset^+_{\underset\tilde {}\to {\Bbb Q}_\beta} \in
\Bbb Q_\beta$ is
$\ne \emptyset_{\underset\tilde {}\to {\Bbb Q}_\beta}$ but
$\emptyset_{\underset\tilde {}\to {\Bbb Q}_\beta} \ne p \in 
{\underset\tilde {}\to {\Bbb Q}_\beta} \Rightarrow
\emptyset^+_{\underset\tilde {}\to {\Bbb Q}_\beta} \le_{\text{pr}} p')"$ [?]
\sn
\item "{$(c)$}"  $S \subseteq \{\delta < \alpha:{\text{\rm cf\/}}(\delta) =
\aleph_1\}$ is stationary; for a club $E$ of $\kappa,\delta \in E \cap S \and
\delta \le \beta$ implies $\Vdash_{{\Bbb P}_\beta} ``(\{r \in
{\underset\tilde {}\to {\Bbb Q}_\beta}:
\emptyset_{\underset\tilde {}\to {\Bbb Q}_\beta} \le_{\text{vpr}} r\},
\le_{\text{vpr}})$ is $\delta^+$-directed (question directed above
$p,\emptyset <_{\text{vpr}} p$
\sn
\item "{$(d)$}"  $\alpha \notin W$ is strongly inaccessible and:
$\beta < \alpha \Rightarrow \Bbb P_\beta < \alpha$.
\ermn
\underbar{Then}:
\mr
\item "{$(\alpha)$}"  forcing with $\Bbb P_\alpha$ 
does not collapse $\alpha$
\sn
\item "{$(\beta)$}"  any function from any $\alpha(*) < \alpha$ to
ordinals in $\bold V^{{\Bbb P}_\alpha}$ belongs to some $\bold
V^{{\Bbb P}_\beta}$.
\endroster
\endproclaim
\bigskip

\demo{Proof}  Clearly clause $(\alpha)$ follows from clause $(\beta)$, so
we shall prove just clause $(\beta)$.  If $W \cap \alpha$ is stationary,
then by \scite{5.7} we are done, so assume not and let $E$ be a club of
$\alpha$ disjoint to $W$, \wilog \, $\beta < \delta \in E \Rightarrow$
density$(\Bbb P_\beta) < \delta$.
Suppose $p \in \Bbb P_\alpha$ and $p \Vdash_{{\Bbb P}_\alpha} ``\underset\tilde {}\to f:
\alpha(*) \rightarrow \text{ Ord}$ is not in any 
$\bold V^{{\Bbb P}_\beta}$ for $\beta < \alpha"$ where $\alpha(*) < \alpha$. \newline
We choose by induction on $\zeta < \alpha$ the tuple 
$(p_\zeta,\alpha_\zeta,\gamma_\zeta,\beta_\zeta,q_\zeta)$ such that:
\medskip
\roster
\item "{$(a)$}"  $\beta_\zeta < \alpha$ is increasing continuous in $\zeta$
and $\beta_{\zeta +1} > \text{ Min}(E \backslash B_\zeta)$
\sn
\item "{$(b)$}"  [?] $p_\zeta \in \Bbb P_{\beta_\zeta}$ is such that 
$\beta(*) \le \beta < \alpha_\zeta \Rightarrow p_\zeta \restriction \beta 
\Vdash_{{\Bbb P}_\beta} ``p_\zeta \restriction \{\beta\} \ne 
\emptyset_{\underset\tilde {}\to {\Bbb Q}_\beta}$ or
${\underset\tilde {}\to {\Bbb Q}_\beta} = 
\{\emptyset_{\underset\tilde {}\to \beta}\}"$
\sn
\item "{$(c)$}"  for $\xi < \zeta,p_\zeta \restriction \beta_\xi = p_\xi$
\sn
\item "{$(d)$}"  $p_\zeta \le q_\zeta \in \Bbb P_\alpha$
\sn
\item "{$(e)$}"  $q_\zeta \Vdash_{{\Bbb P}_\alpha} ``\underset\tilde {}\to f
(\alpha_\zeta) = \gamma_\zeta"$ but there is no $q' \in \Bbb P_{\beta_\zeta}$
compatible with $p_\zeta$ which forces this and $\alpha_\zeta < \alpha(*)$,
of course
\sn
\item "{$(f)$}"  if cf$(\zeta) = \aleph_1$, then for some 
$\beta'_\zeta < \beta_\zeta$, \newline
$\gamma \in [\beta'_\zeta,\beta_\zeta) \Rightarrow q_\zeta \restriction \gamma
\Vdash_{{\Bbb P}_\gamma} ``\emptyset_{\underset\tilde {}\to {\Bbb Q}_\gamma} 
\le_{\text{vpr}} = q_\zeta \restriction \{\gamma\}"$
\sn
\item "{$(g)$}"  $q_\zeta \in \Bbb P_{\beta_{\zeta +1}}$ and for every $\beta \in
[\beta_\zeta,\alpha)$ we have
$$
\align
p_{\zeta +1} \restriction \beta \Vdash_{{\Bbb P}_\beta} ``&\text{if }
\emptyset_{\underset\tilde {}\to {\Bbb Q}_\beta} <_{\text{vpr}} q_\zeta \restriction
\{\beta\} \text{ in } {\underset\tilde {}\to {\hat{\Bbb Q}}_\beta} \text{ then} \\
  &p_{\zeta +1} \restriction \{\beta\} = q_\zeta \restriction \{\beta\}; \\
  &\text{if not [?] then } r \in {\underset\tilde {}\to {\Bbb Q}_\beta} \text{ such
that } \emptyset_{\underset\tilde {}\to {\Bbb Q}_\beta} <_{\text{vpr}} r \text{ and }
\emptyset_{\underset\tilde {}\to {\Bbb Q}_\beta} \text{ if there is none}"
\endalign
$$
\noindent
(see clause (b) in the assumptions of \scite{5.8}).
\endroster
\medskip

\noindent
Having carried the definition, for some stationary $W' \subseteq S
\subseteq \{\delta < \alpha:\text{cf}(\delta) = \aleph_1\}$ and 
$\gamma^* < \alpha$ and $\beta'$ we have:
$\zeta \in W' \Rightarrow \gamma_\zeta = \gamma^* \and \beta'_\zeta \le
\beta' < \zeta$. \nl
As $|\Bbb P_{\gamma^*}| < \alpha = \text{ cf}(\alpha)$, \wilog \,
$\zeta \in W' \Rightarrow q_\zeta \restriction
\gamma_\zeta = q^*$.  Now choose $\xi < \zeta$ in $W'$ then $q_\zeta,q_\xi$
are compatible and an easy contradiction to clause (c) (with $q_\xi$ here
playing the role of $q'$ there). \hfill$\square_{\scite{5.8}}$\margincite{5.8}
\enddemo
\bn
Now we can refine \scite{1.16A} to the iteration theorem of this section.
\proclaim{\stag{5.9} Claim}  1) Suppose $\bold W \subseteq \omega_1$ 
be stationary, $F$ is a function, \ub{then} for every ordinal 
$\alpha$ there is ${\text{\rm UP\/}}^4(\bold W)$-iteration 
$\bar{\Bbb Q} = \langle \Bbb P_i,{\underset\tilde {}\to {\Bbb Q}_i},
\kappa_i:i < \alpha^\dag \rangle$, such that:
\mr
\item "{$(a)$}"  for every $i$ we have ${\underset\tilde {}\to {\Bbb Q}_i} = 
F(\bar{\Bbb Q} \restriction i)$ and $\kappa_i = \kappa_{cc}(\Bbb P_i *
\Bbb Q_i)$
\sn
\item "{$(b)$}"  $\alpha^\dag \le \alpha$
\sn
\item "{$(c)$}"  either $\alpha^\dag = \alpha$ or the following fails:
{\roster
\itemitem{ $(*)$ }  $F(\bar{\Bbb Q})$ is an (${\text{\rm
Sp\/}}_e(W)-{\text{\rm Lim\/}}(\bar{\Bbb Q}))$-name of a 
forcing notion forced to satisfy ${\text{\rm UP\/}}^4(\Bbb I,\bold W)$ 
for some $\Bbb I$
$\kappa$-complete set of ideals, where $\kappa$ is minimal such that
${\text{\rm Sp\/}}_e(W)-{\text{\rm Lim\/}}(\bar{\Bbb Q})$ satisfies 
the $\kappa$-c.c
\endroster}
\item "{$(d)$}"  ${\text{\rm Sp\/}}_e(W)-{\text{\rm Lim\/}}
(\bar{\Bbb Q})$ does not collapse $\aleph_1$ and
preserve stationary subsets of $\bold W$ (in fact it satisfies
${\text{\rm UP\/}}^1(\bold W)$.
\endroster
\endproclaim
\bigskip

\demo{Proof}  Straight.
\enddemo
\newpage


\head {\S6 Preservation of UP$^0$} \endhead  \resetall \sectno=6
\bn
Here we present alternatives to \S5, i.e. to UP$^4(\bold W)$-iterations. 
In UP$^6$-iteration $\langle \Bbb P_j,{\underset\tilde {}\to {\Bbb
Q}_i},
\kappa_j:
k \le \alpha,i < \alpha \rangle$ the demands are weak but the
$\kappa_{i+j}$ may be large and it is quite similar to 
semi-properness (but for
fewer models).  In UP$^6$-iteration we carry with us trees.
\nl
Recall [?]
\definition{\stag{si.1} Definition}  1) We say $q$ is 
$(N,\kappa,\Bbb Q)$-semi-generic \ub{if}
$q$ is $(N,\Bbb Q)$-semi-generic and 
$q \Vdash ``\text{if } y \in
N[\underset\tilde {}\to G] \cap \bold V$ 
then for some $A \in N,|A|^{\bold V} < \kappa$ and $y \in A"$. \nl
2) Similarly with $\underset\tilde {}\to \kappa$ instead of $\kappa$.
\enddefinition
\bigskip

\remark{Remark}  In 
\scite{si.1}(1) \wilog \, $y \in N[\underset\tilde {}\to G]
\cap \text{Ord}$.
\endremark
\bigskip

\proclaim{\stag{si.2} Lemma}  Suppose
\mr
\item "{$(A)$}"  $\Bbb Q$ is a forcing notion
\sn
\item "{$(B)$}"  $\Bbb I \in N$ is a family of ideals, $\Bbb I$ is 
$\kappa$-closed, $\lambda \ge \aleph_2$ (and it is
natural but not needed to assume that $\Bbb I$ is $\lambda$-complete)
\sn
\item "{$(C)$}"  $N$ is $\lambda$-strictly $(\Bbb I,\bold W)$-suitable
for $(\chi,x)$ as witnessed by \newline
$\bar N = \langle N_\eta:\eta \in (T,\bold I)\rangle,
\Bbb Q \in N$ (so $N \prec ({\Cal H}(\chi),\in,<^*_\chi)$ is countable,
$N = N_{<>}$ and, of
course, $x$ codes $\langle \Bbb Q,\Bbb I,\bold W \rangle)$, see Definition
\scite{4.1})
\sn
\item "{$(D)$}"  $q$ is $(N,\kappa,\Bbb Q)$-semi-generic
\sn
\item "{$(E)$}"  at least one of the following holds:
{\roster
\itemitem{ $(\alpha)$ } $|{\text{\rm MAC\/}}(\Bbb Q)| < \lambda$ 
\sn
\itemitem{ $(\beta)$ }  $q$ is $(\bar N,\kappa,\Bbb Q)$-semi-generic, 
i.e. $(N_\eta,\kappa,\Bbb Q)$-semi-generic for every 
$\underset\tilde {}\to \eta \in T$ 
\sn
\itemitem{ $(\gamma)$ } $q \Vdash ``\underset\tilde {}\to T \subseteq
T$ is a subtree and every $\eta \in \underset\tilde {}\to T$ for some $\nu,
\eta \trianglelefteq \nu \in { \text{\rm split\/}}
(\underset\tilde {}\to T,\bold I),
\bold I_\eta \le_{RK} \bold I_\nu$ and $\eta \in 
\underset\tilde {}\to T \Rightarrow 
N_\eta[{\underset\tilde {}\to G_{\Bbb Q}}] \cap \omega_1 = 
N_\eta \cap \omega_1"$.
\endroster}
\ermn
\ub{Then} 
$q \models ``N[{\underset\tilde {}\to G_{\Bbb Q}}]$ is $\lambda$-strictly
$(\Bbb I,\bold W)$-suitable for $(\chi,\langle x,
{\underset\tilde {}\to G_{\Bbb Q}} \rangle)"$ (in fact, if \nl
$\langle N_\eta:\eta \in (T,\Bbb I) \rangle$ was a witness for 
$N$ \ub{then} $(\langle N_\eta
[\underset\tilde {}\to G]:\eta \in (T,\bold I) \rangle$ is a witness for
$N[\underset\tilde {}\to G]$ being $(\Bbb I,\bold W)$-suitable for
$(\chi,x)$).
\endproclaim
\bigskip

\demo{Proof}  Let $G \subseteq \Bbb Q$ be generic over $\bold V$ such that $q \in G$.

In $\bold V^{\Bbb Q}$, i.e., in $\bold V[G]$,  
clearly $N_\eta[G] \prec ({\Cal H}(\chi)^{\bold V^{\Bbb Q}},
\in,<^*_\chi)$ (see
\cite[III,2.11,p.104]{Sh:f}, $N_\eta[G]$ is countable (trivially) and
$N_{\eta \restriction k}[G] \prec N_\eta[G]$ for $k \le \ell g(\eta)$.

As $q$ is $(N,\kappa,{\Bbb Q})$-semi-generic and $N
= N_{\langle \rangle}$, clearly
$q$ is $(N_{\langle \rangle},\Bbb Q)$-semi-generic.  First assume
$|\text{MAC}(\Bbb Q)| < \lambda$: now $\Bbb Q \in N$ hence MAC$(\Bbb
Q) \in N$ hence $|\text{MAC}(\Bbb Q)| \in N$, and as $N_{\langle
\rangle} cap \lambda <_\lambda
N_\eta \cap \lambda$ clearly $N_\eta \cap |\text{MAC}(\Bbb Q)| =
N_{\langle \rangle} \cap |\text{MAC}(\Bbb Q)|$ hence also $N_\eta \cap
\text{ MAC}(\Bbb Q) = N_{\langle \rangle} \cap \text{ MAC}(\Bbb Q)$ hence
every $\Bbb Q$-name $\underset\tilde {}\to \tau \in N_\eta$ of an ordinal belongs
to $N_{\langle \rangle}$ (essentially).  Hence $q$ is
$(N_\eta,\kappa,\Bbb Q)$-semi-generic, so $q \Vdash_{\Bbb Q} ``N_\eta[G] \cap \omega_1 =
N_\eta \cap \omega_1 = N \cap \omega_1"$, (even $N_{\eta \restriction \ell}
[G] <_\lambda N_\eta[G]$).  So $q$ is
$(N_\eta,\kappa,\Bbb Q)$-semi-generic for any $\eta\in T$ and even $\eta 
\in (\text{lim }T)^{\bold V}$ or $\eta \in 
(\text{lim }T)^{\bold V^{\Bbb Q}}$.

This almost shows that $\langle N_\eta[G]:\eta \in (T,\bold I) \rangle$ is 
a witness to $N[G]$ being $\lambda$-strictly $(\Bbb I,\bold
W)$-suitable for $(\chi,x)$.

The missing point is clause $(e)^-$ of Definition \scite{4.1}, that is,
that there may be $\eta \in T$, and $J \in \Bbb I \cap
N_\eta[G]$ such that $J \notin N_\eta$.  So there is a $\Bbb Q$-name
$\underset\tilde {}\to \tau \in N_\eta$ satisfying 
$J = \underset\tilde {}\to \tau[G]$; choose
a $\underset\tilde {}\to \tau$ like that with min$\{|Y|:Y \in N_\eta$ and
$\underset\tilde {}\to \tau[G] \in Y\}$ minimal, and let the set $Y$
be $\{J_i:i < \alpha\}$ (without loss of generality $Y \subseteq \Bbb I$), so
without loss of generality $\langle J_i:i < \alpha \rangle$ belongs to
$N_\eta$: by the minimality of $|Y| = |\alpha|$ and $q$ being 
$(N_\eta,\kappa,\Bbb Q)$-semi-generic we have $\alpha < \kappa$.

So as $\{J_i:i < \alpha\} \in N_\eta \cap \Bbb I$ and $\Bbb I$ is 
$\kappa$-closed there is $J \in N_\eta \cap \Bbb I$ such that 
$\dsize \bigwedge_{i < \alpha} J_i \le_{\text{RK}} J$
hence the set $\{\nu:\eta \triangleleft \nu \in \text{ split}(T)$ and
$J \le_{\text{RK}} \bold I_\nu\}$ is a front of $T^{[\eta]}$.  So 
$\langle N_\eta[G]:\eta \in (T,\Bbb I) \rangle$ is a
witness for ``$N[G]$ is $\lambda$-strictly $\Bbb I$-suitable" for 
$(\lambda,\chi,\langle x,G \rangle)$ (see Definition \scite{4.1}) so by 
\scite{4.4} we know that $N[G]$ is $\Bbb I$-suitable. \nl
Now if the second phrase of (E) holds, the proof is similar and if the third
holds, we can prove that \wilog \, 
the second holds.  \hfill$\square_{\scite{si.2}}$\margincite{si.2}
\enddemo
\bigskip

\noindent
The proof suggests some definitions, but we first consider:
\definition{\stag{si.3} Definition}  1) For a forcing notion $\Bbb Q$, family 
$\Bbb I$ of ideals and cardinal and $\Bbb Q$-names
$\underset\tilde {}\to \lambda$ of a cardinal $\kappa$ and stationary 
$\bold W \subseteq \omega_1$, we say $\Bbb Q$ satisfies 
UP$^6(\Bbb I,\kappa,\underset\tilde {}\to \lambda,\bold W)$ or
UP$_{\kappa,\underset\tilde {}\to \lambda}(\Bbb I,\bold W)$ \ub{if}:
\mr
\item "{$(*)$}"  for every $\chi$ regular large enough, $p \in \Bbb Q$ 
and $N$ a strictly $(\Bbb I^{[\kappa]},\bold W)$-suitable model for 
$\chi$ satisfying $\{p,\Bbb Q,\kappa,\underset\tilde {}\to \lambda\} \in N$, 
\ub{there is} $q$ satisfying $p \le_{\text{pr}} q \in \Bbb Q$ such that 
$q$ is $(N,\Bbb Q)$-semi-generic
and $q \Vdash_{\Bbb Q} ``N[{\underset\tilde {}\to G_{\Bbb Q}}]$ is strictly
$\Bbb I^{[{\underset\tilde {}\to \lambda}]}$-suitable model for 
$\chi,\lambda$.
\ermn
1A)  We say ``$q$ is $(N,{\underset\tilde {}\to {\Bbb I}},\Bbb Q)$-semi$_6$ generic" if
$q$ is $(N,\Bbb Q)$-semi genric and $q \Vdash ``N[{\underset\tilde
{}\to G_{\Bbb Q}}]$
is strictly ${\underset\tilde {}\to {\Bbb I}}$-suitable". \nl
2) In part (1), if we omit $\lambda$, we mean $\lambda = \text{ Max}\{\lambda:\Bbb I
\text{ is } \lambda$-complete$\}$, so $\lambda$ is regular $> \aleph_1$. \nl
3) For 
$\Bbb Q,\Bbb I,\bold W,\kappa,\underset\tilde {}\to \lambda$ as in part
(1) and $\underset\tilde {}\to \theta$ a $\Bbb Q$-name of a cardinal; 
we say that
$\Bbb Q$ satisfies UP$^6_{\kappa,\underset\tilde {}\to
\lambda,\underset\tilde {}\to \theta}
(\Bbb I,\bold W)$ \ub{if}:
\mr
\item "{$(**)$}"  if $\Bbb I^+$ is a set of partial orders ideal extending
$\Bbb I,\chi$ large enough, $p \in \Bbb Q \cap N$ and $N$ a strictly
$(\Bbb I,\bold W)$-suitable model for $\chi,\{p,\Bbb Q,\kappa,
\underset\tilde {}\to \lambda,\underset\tilde {}\to \theta\} \in N$, \ub{then}
there is $q$ satisfying $p \le_{pr} q \in \Bbb Q$ such that $q$ is $(N,
\Bbb I^{[\underset\tilde {}\to \lambda]} \cup (\Bbb I^+ \backslash \Bbb I)
^{[\underset\tilde {}\to \theta]},\Bbb Q)$-semi$_6$-generic 
(see (1A) above).
\ermn
4) We say $\Bbb Q$ satisfies UP$^5_{\kappa,{\underset\tilde {}\to \lambda}}
(\Bbb I,\bold W)$ \ub{if}: 
for any $(\Bbb I^{[\kappa]},\bold W)$-suitable tree
$\langle N_\eta:\eta \in (T,\Bbb I) \rangle$ of models and $p \in N_{<>}$
there are $q,\underset\tilde {}\to T$ such that $p \le_{\text{pr}} q
\in \Bbb Q,q \Vdash ``\underset\tilde {}\to T \subseteq T$ is a subtree and
$\langle N_\eta[G_P]:\eta \in \underset\tilde {}\to T \rangle$ is a
$(\Bbb I^{[\lambda]},\bold W)$-suitable tree of  models and
$N_\eta[G_{\Bbb Q}] \cap \omega_1 = N_{<>} \cap \omega_1$.  [Saharon
but \sciteu{4.x}!]
\enddefinition
\bn
Some variants of this Definition are equivalent by the following claim.
\proclaim{\stag{si.4} Claim}  1) If $\kappa$ is regular uncountable and
$\Bbb Q$ satisfies the $\kappa$-c.c. \ub{then}:
$q \in \Bbb Q$ is $(N,\kappa,\Bbb Q)$-semi-generic \ub{iff} $q \in
\Bbb Q$ is $(N,\Bbb Q)$-semi-generic. \nl
2) If $N$ is 
$(\Bbb I,\bold W)$-suitable for $(\chi,x)$ (see Definition \scite{4.2}),
$\Bbb I$ is $\lambda$-complete, $\lambda$ is regular and 
$(\forall \alpha < \lambda)(|\alpha|^{\aleph_0} < \lambda)$,
\ub{then} there is a $\lambda$-strictly $(\Bbb I,\bold W)$-suitable $N'$
for $(\chi,x)$ such that $N \prec N'$ and $N' \cap \omega_1 = N \cap
\omega_1$.  \nl
3)  If $q$ is $(N,\kappa,\Bbb Q)$-semi generic, 
$|{\text{\rm MAC\/}}(\Bbb Q)| < \lambda,
\Bbb I \in N$ is $\kappa$-closed $\lambda$-complete and $N$ is
$\Bbb I$-suitable \ub{then} $q$ is $(N,\Bbb I^{[\lambda]},\Bbb Q)$-semi$_6$ 
generic. \nl
4) If $|{\text{\rm MAC\/}}(\Bbb Q)| < \lambda$ 
and $\Bbb Q$ is ${\text{\rm UP\/}}^0(\Bbb I^{[\kappa]},
\bold W)$ \ub{then} $\Bbb Q$ is 
${\text{\rm UP\/}}^6_{\kappa,\lambda}(\Bbb I,\bold W)$.
Similarly for $\lambda \, 
\Bbb Q$-name such that for a dense set of $p$ we have
$p \Vdash \underset\tilde {}\to \lambda = \lambda$ and 
${\text{\rm MAC\/}}(\Bbb Q_{\ge p}) <
\lambda$.
\endproclaim
\bigskip

\demo{Proof}  Straight. \nl
1) Reflect. \nl
2) Use the partition theorem \scite{2.12}. \nl
3) By \scite{si.2} using possibility (A). \nl
4) By \scite{si.2}, too.  \hfill$\square_{\scite{si.4}}$\margincite{si.4}
\enddemo
\bigskip

\proclaim{\stag{si.4a} Claim}  1) If $\Bbb Q$ satisfies 
${\text{\rm UP\/}}^6_{\kappa,{\underset\tilde {}\to \lambda}}
(\Bbb I,\bold W)$ and $\bold W_1 \subseteq
\bold W,\kappa_1 \le \kappa$ and 
$\Vdash_{\Bbb Q} ``\underset\tilde {}\to \lambda \le
{\underset\tilde {}\to \lambda_1}"$, \ub{then} $\Bbb Q$ satisfies
${\text{\rm UP\/}}^6_{\kappa_1,
{\underset\tilde {}\to \lambda_1}}(\Bbb I,\bold W_1)$. \nl
2) If $\Bbb Q_0$ is a forcing notion satisfying 
${\text{\rm UP\/}}^6_{\kappa_0,
{\underset\tilde {}\to \kappa_1}}(\Bbb I,\bold W)$ and
${\underset\tilde {}\to {\Bbb Q}_1}$ is a $\Bbb Q$-name 
of a forcing notion satisfying
${\text{\rm UP\/}}^6_{{\underset\tilde {}\to \kappa_1},
{\underset\tilde {}\to \kappa_2}}
(\Bbb I,\bold W_1)$, \ub{then} $\Bbb Q_0 * {\underset\tilde {}\to
{\Bbb Q}_1}$ is a
forcing notion satisfying 
${\text{\rm UP\/}}^6_{\kappa_0,{\underset\tilde {}\to \kappa_2}}
(\Bbb I,\bold W_1)$. \nl
3) If $\Bbb Q$ satisfies 
${\text{\rm UP\/}}^6_{\kappa,\underset\tilde {}\to \lambda,
\underset\tilde {}\to \theta}(\Bbb I,\bold W_1)$ and $\Vdash_{\Bbb Q}
\underset\tilde {}\to \theta \le \theta$ and $\Bbb I \subseteq \Bbb I^+$
an $\Bbb I^+ \backslash \Bbb I$ is $\theta$-complete \ub{then} 
$\Bbb Q$ satisfies
${\text{\rm UP\/}}^6_{\kappa,\underset\tilde {}\to \lambda}
(\Bbb I',\bold W_1)$.
\endproclaim
\bigskip

\definition{\stag{si.5} Definition}  1) We say that $\bar{\Bbb Q} = \langle
\Bbb P_j,{\underset\tilde {}\to {\Bbb Q}_i},
{\underset\tilde {}\to \kappa_j}:j \le \alpha$ 
and $i < \alpha \rangle$ 
is a UP$^{6,e}(\Bbb I,\bold W,W)$-suitable iteration (with $e=6$ if not
mentioned explicitly) \ub{if}:
\medskip
\roster
\item "{$(a)$}"  $\langle 
\Bbb P_j,{\underset\tilde {}\to {\Bbb Q}_i}:j < \alpha,i <
\alpha \rangle$ is an $\aleph_1-\text{Sp}_e(W)$-iteration
\sn
\item "{$(b)$}"  $\Bbb I$ is a set of partial 
order ideals such that $\Bbb I^{[\kappa]}$ is $\kappa$-closed for any
regular $\kappa,\aleph_2 \le \kappa \le |\Bbb P_\alpha|^+$
\sn
\item "{$(c)$}"  $\bold W \subseteq \omega_1$ is stationary
\sn
\item "{$(d)$}"  for each $i < \alpha$ we have: 
${\underset\tilde {}\to \kappa_i}$ is a $\Bbb P_i$-name of a regular cardinal
$> \aleph_1$ in $\bold V$ [purity?],
\sn
\item "{$(e)$}"  $(i) \quad$  for $i < j$ we have $\Vdash_{{\Bbb P}_j} 
``{\underset\tilde {}\to \kappa_i} \le
{\underset\tilde {}\to \kappa_j}"$
\sn
\item "{${{}}$}"  $(ii) \quad$ if $\delta \le \alpha$ is a limit ordinal,
then $\Bbb P_\delta$ satisfies the local 
${\underset\tilde {}\to \kappa_\delta}$-c.c. purely \nl

$\qquad$ [and (?) if
$\nVdash {\underset\tilde {}\to \kappa_\delta} \ne \kappa$ then some
$\kappa$-complete ideal on \nl

$\qquad \kappa$ belongs to $\Bbb I]$
\sn
\item "{$(f)$}"  $\Vdash_{{\Bbb P}_i} 
``{\underset\tilde {}\to {\Bbb Q}_i}$ satisfies
UP$^6_{{\underset\tilde {}\to \kappa_i},{\underset\tilde {}\to \kappa_{i+1}}}
({\underset\tilde {}\to {\Bbb I}_i},\bold W)$ and
$({\underset\tilde {}\to {\Bbb Q}_i},\le_{\text{vpr}})$ 
is $\aleph_1$-complete".
\ermn
2) We say $\langle \Bbb P_j,{\underset\tilde {}\to {\Bbb Q}_i},
{\underset\tilde {}\to \kappa_j}:j \le \alpha,i < \alpha \rangle$ is a
UP$^4(\Bbb I,\bold W)$-iteration \ub{if}:
\mr
\item "{$(a)-(e)$}"  as above
\sn
\item "{$(f)$}"  $\Bbb Q_i$ 
satisfies UP$^5_{{\underset\tilde {}\to \kappa_i},
{\underset\tilde {}\to \kappa_{i+1}}}(\Bbb I,\bold W)$.
\endroster
\enddefinition
\bigskip

\definition{\stag{si.6} Definition}  1) We 
say that $\langle \Bbb P_j,{\underset\tilde {}\to {\Bbb Q}_i},
{\underset\tilde {}\to \kappa_j}:j \le \alpha$ and $i < \alpha \rangle$ is a 
weak UP$^{6,e}(\Bbb I,\bold W)$-iteration \ub{if} 
(when $e=6$ we may omit it):
\mr
\item "{$(a)$}"  $\langle 
\Bbb P_i,{\underset\tilde {}\to {\Bbb Q}_i}:i < \alpha \rangle$ is 
an $\aleph_1$-SP$_e(W)$-iteration (and, of course, $\Bbb P_\alpha =
\text{ Sp}_e(W)$-Lim$_\kappa(\bar{\Bbb Q}))$
\sn
\item "{$(b)$}"  $\Bbb I$ is a set of partial order ideals
\sn
\item "{$(c)$}"  $\bold W \subseteq \omega_1$ stationary
\sn
\item "{$(d)$}"  ${\underset\tilde {}\to \kappa_j}$ is a $\Bbb P_j$-name of a
member of RCar$^V \backslash \omega_2$, increasing with $j$
\sn
\item "{$(e)$}"  for $i < j \le \alpha,i$ nonlimit we have
$$
\Vdash_{P_i} ``\Bbb P_j/\Bbb P_i 
\text{ satisfies the \, UP}_{{\underset\tilde {}\to \kappa_i},
{\underset\tilde {}\to \kappa_j}}(\Bbb I,\bold W)"
$$
\item "{$(f)$}"  $\Vdash_{{\Bbb P}_i} ``({\underset\tilde {}\to {\Bbb Q}_i},
\le_{\text{vpr}})$ is $\aleph_1$-complete".
\ermn
2) We say $\bar{\Bbb Q} = \langle \Bbb P_j,{\underset\tilde {}\to
{\Bbb Q}_i},
{\underset\tilde {}\to \kappa_j}:j \le \alpha,i < \alpha \rangle$ is a
UP$^5(\Bbb I,\bold W)$-iteration \ub{if}:
\mr
\item "{$(a)-(d),(f)$}"  as above
\sn
\item "{$(g)$}"  if $i < j \le \alpha,i$ non limit we have \nl
$\Vdash_{{\Bbb P}_i} 
``\Bbb P_j/{\underset\tilde {}\to G_{{\Bbb P}_i}}$ satisfies
UP$^5_{{\underset\tilde {}\to \kappa_i},{\underset\tilde {}\to \kappa_j}}
(\Bbb I,\bold W)$.
\endroster
\enddefinition
\bigskip

\proclaim{\stag{si.7} Claim}  1) If $\langle \Bbb P_j,
{\underset\tilde {}\to {\Bbb Q}_i},
{\underset\tilde {}\to \kappa_j}:j \le \alpha$ and $i < \alpha \rangle$ is
a ${\text{\rm UP\/}}^6(\Bbb I,\bold W,W)$-iteration, \ub{then} it is a weak
${\text{\rm UP\/}}^6(\Bbb I,\bold W,W)$-iteration, moreover, in clause
(e) also limit $i$ is O.K. \nl
2) Assume $\bar{\Bbb Q} = \langle \Bbb P_j,{\underset\tilde {}\to
{\Bbb Q}_i},
{\underset\tilde {}\to \kappa_j}:j \le \alpha,i < \alpha \rangle$ is
an $\aleph_1-{\text{\rm Sp\/}}_e(W)$-iteration
\mr
\item "{$(a)$}"  if $\alpha$ is a limit ordinal and $\beta < \alpha
\Rightarrow \bar{\Bbb Q} \restriction \beta = \langle \Bbb P_j,
{\underset\tilde {}\to {\Bbb Q}_i},{\underset\tilde {}\to \kappa_j}:
j \le \beta,i < \beta \rangle$ is a weak 
${\text{\rm UP\/}}^6(\Bbb I,\bold W)$-iteration and
${\underset\tilde {}\to \kappa_\alpha} = \text{ sup}\{\kappa_j:j < \alpha\}$,
\ub{then} $\bar{\Bbb Q}$ is a weak 
${\text{\rm UP\/}}^6(\Bbb I,\bold W)$-iteration
\sn
\item "{$(b)$}"  if $\alpha = \beta +1,\bar{\Bbb Q} \restriction \beta$ is a
weak ${\text{\rm UP\/}}^6(\Bbb I,\bold W)$-iteration and in $\bold V^{{\Bbb
P}_\beta},\Bbb Q_\beta$ satisfies
${\text{\rm UP\/}}^6_{{\underset\tilde {}\to \kappa_\beta},
{\underset\tilde {}\to \kappa_{\beta +1}}}(\Bbb I,\bold W)$ \ub{then}
$\bar{\Bbb Q}$ is a weak ${\text{\rm UP\/}}^6(\Bbb I,\bold W)$-iteration. \nl
[Saharon - compare with \scite{si.8} and \scite{1.27A}(2)]
\endroster
\endproclaim
\bigskip

\demo{Proof}  Let $\gamma \le \beta \le \alpha$ and we need
\mr
\item "{$\boxtimes_{\gamma,\beta}$}"  Assume $G_\gamma \subseteq \Bbb P_\gamma$
is generic over $\bold V,\kappa_i = {\underset\tilde {}\to \kappa_i}[G_i]$ and
let $N \in \bold V[G_\gamma],N$ is strictly 
$\Bbb I^{[\kappa]}$-suitable, $N \cap \omega_1 \in \bold W$ (so
$N \prec ({\Cal H}(\chi),\in)$ is countable) and $p \in
\Bbb P_\beta/G_\gamma$ and $p \in \Bbb P_\beta \cap N$ and 
$\bar{\Bbb Q},\beta,\gamma \in N$. 
\ub{Then} we can find $q$ satisfying 
$p \le_{\text{pr}} q \in \Bbb P_\beta/G_\gamma$ and
$q \Vdash_{{\Bbb P}_\beta/G_\gamma} ``N[{\underset\tilde {}\to G_\beta}] \cap
\omega_1 = N \cap \omega_1$ and $N[{\underset\tilde {}\to G_\beta}]$ is
$\Bbb I^{[{\underset\tilde {}\to \kappa_\beta}]}$-suitable".\nl
Without loss of generality [??] $p$ forces a value to 
${\underset\tilde {}\to \kappa_\beta}$ moreover
$\kappa_\beta = \kappa_{\underset\tilde {}\to \kappa_\beta}(p,\Bbb P_\beta)$].
\ermn
Naturally, we prove this by induction on $\beta$ (for all $\gamma$).  The
case $\gamma = \beta$ holds trivially so assume $\gamma < \beta$.
If $\beta = 0$, we have nothing to prove.  If $\beta$ is a successor 
ordinal say $\gamma_1 +1$  so $\gamma \le \gamma_1$, now we use first 
$\boxtimes_{\gamma,\gamma_1}$ and then
the demand on ${\underset\tilde {}\to Q_{\gamma_1}}$ in definition
\scite{si.6}, in clause (f).
\enddemo
\bn
So from now on we shall assume that $\beta$ is a limit ordinal.  As in the
proof of \scite{5.5} we can note \nl
\ub{Fact A}:  If ${\underset\tilde {}\to \gamma_1} \le
{\underset\tilde {}\to \gamma_2}$ are simple $\bar{\Bbb Q}$-named 
$[0,\beta)$-ordinals then $\boxtimes_{{\underset\tilde {}\to \gamma_1},
{\underset\tilde {}\to \gamma_2}}$ holds.
\bn
\ub{Proof of the fact}:  Here we use ``$e = 6$" rather than ``$e=4$".  On
$\Bbb P_{\underset\tilde {}\to \zeta}$ see Definition \scite{1.13}(F)(g).  We prove it
by induction on the depth of ${\underset\tilde {}\to \zeta_2}$, see
Definition \scite{1.6}(5).  So we are given
$G_{\underset\tilde {}\to \zeta_2} \subseteq 
\Bbb P_{\underset\tilde {}\to \zeta_1}$ generic over $\bold V$ and in particular let
$\zeta_1 = {\underset\tilde {}\to \zeta_2}
[G_{\underset\tilde {}\to \zeta_1}]$, (so it is simpler to say that
$G_{\zeta_1} \subseteq \Bbb P_{\zeta_1}$ is generic over $\bold V,
{\underset\tilde {}\to \zeta_1}[G_{\zeta_1}] = \zeta_1$).  Let
$\kappa_1 = {\underset\tilde {}\to \kappa_{\zeta_1}}
[G_{\underset\tilde {}\to \zeta_1}]$ and we are also given $N$ which is
strictly $\Bbb I^{[\kappa]}$-suitable, $p \in
\Bbb P_{\underset\tilde {}\to \zeta_2}/G_{\underset\tilde {}\to \zeta_1},p \in N,
\{\bar{\Bbb Q},{\underset\tilde {}\to \zeta_1},G_{\underset\tilde {}\to \zeta_1}\}
\in N,{\underset\tilde {}\to \zeta_2} \in N$.  We have to find $q \in
\Bbb P_{\underset\tilde {}\to \zeta_2}/G_{\underset\tilde {}\to \zeta_1}$ such that
$p \le_{\text{pr}} q,q$ is $(N,\Bbb P_{\underset\tilde {}\to \zeta_2}/
\Bbb P_{\underset\tilde {}\to \zeta_1})$-generic and $q \Vdash N
[G_{{\Bbb P}_{\underset\tilde {}\to \zeta_2}/G_{\underset\tilde {}\to \zeta_1}}]$
is strictly $\Bbb I^{[\kappa_2]}$-suitable.

If the depth of ${\underset\tilde {}\to \zeta_2}$ is 0, then
$\Vdash ``{\underset\tilde {}\to \zeta_2} = \zeta_2"$
and we can use $\boxtimes_{\zeta_1,\zeta_2}$.  So assume the depth of
${\underset\tilde {}\to \zeta_2}$ is $> 0$, and so for some $\gamma^*$ and
a sequence $\langle {\underset\tilde {}\to \zeta_{2,\varepsilon}}:
\varepsilon < \varepsilon^* \rangle$ of simple $\bar{\Bbb Q}$-named
$[\text{Max}\{\gamma^*,\gamma\},\beta)$-ordinals and $\Bbb P_{\gamma^*}$-name
$\underset\tilde {}\to \varepsilon$ we have
$\Vdash_{\bar{\Bbb Q}} ``{\underset\tilde {}\to \zeta_2} =
{\underset\tilde {}\to \zeta_{2,{\underset\tilde {}\to \varepsilon}}}"$.
So \wilog \, $\{\gamma^*,\langle
{\underset\tilde {}\to \zeta_{2,\varepsilon}}:
\varepsilon < \varepsilon^* \rangle,\underset\tilde {}\to \varepsilon\} \in
N$.  Let ${\underset\tilde {}\to \zeta_1} = \text{ Max}\{\gamma,\gamma^*,
{\underset\tilde {}\to \zeta_1}\},
\Vdash_{\Bbb Q} ``{\underset\tilde {}\to \zeta_1} \le
{\underset\tilde {}\to \zeta'_1} \le \zeta_2"$.

Now clearly ${\underset\tilde {}\to \zeta'_1}$ is a simple $\bar{\Bbb Q}$-named
$[0,\beta)$-ordinal, $\Vdash_{\bar{\Bbb Q}} ``{\underset\tilde {}\to \zeta_1} \le
{\underset\tilde {}\to \zeta'_1} \le {\underset\tilde {}\to \zeta'_2}"$
and $\boxtimes_{{\underset\tilde {}\to \zeta_1},
{\underset\tilde {}\to \zeta'_1}} \and
\boxtimes_{{\underset\tilde {}\to \zeta'_1},
{\underset\tilde {}\to \zeta_2}} \Rightarrow
\boxtimes_{{\underset\tilde {}\to \zeta_1},{\underset\tilde {}\to \zeta_2}}$
and $\boxtimes_{{\underset\tilde {}\to \zeta_1},
{\underset\tilde {}\to \zeta_2}}$ easily holds (by the cases proved above)
so it is enough to prove
$\boxtimes_{{\underset\tilde {}\to \zeta'_1},
{\underset\tilde {}\to \zeta_2}}$.  This just means that \wilog \,
$\zeta_1 \le \gamma^*$ and even ${\underset\tilde {}\to \zeta_1} =
\gamma^*$.  Now $\underset\tilde {}\to \varepsilon
[G_{\underset\tilde {}\to \zeta_1}] \in N$ so we use the induction hypothesis
to get the desired $q$.  \hfill$\square_{\text{fact}}$
\bn
\ub{Fact B}:  If $\underset\tilde {}\to \xi$ is a simple $\bar{\Bbb Q}$-named
$[0,\beta)$-ordinal, $p \in \Bbb P_\beta$ and $\underset\tilde {}\to \tau$ is
a $\Bbb P_\beta$-name of a countable ordinal, \ub{then} there are
$\underset\tilde {}\to \varepsilon$ and $q$ such that:
\mr
\item "{$(*)(i)$}"  $\Bbb P_\beta \models p_{\text{pr}} q$
\sn
\item "{$(ii)$}"  $q \restriction \underset\tilde {}\to \xi = p
\restriction \underset\tilde {}\to \xi$
\sn
\item "{$(iii)$}"  $\underset\tilde {}\to \varepsilon$ is a
$\Bbb P_{\underset\tilde {}\to \xi}$-name of a countable ordinal
\sn
\item "{$(iv)$}"  $q \Vdash_{{\Bbb P}_\beta} ``\underset\tilde {}\to \tau <
\underset\tilde {}\to \varepsilon"$.
\endroster
\bigskip

\demo{Proof}  Let $\langle {\underset\tilde {}\to \zeta_n}:n < \omega
\rangle$ be a witness for $p$, so each ${\underset\tilde {}\to \zeta_n}$
is a simple $\bar{\Bbb Q}$-named $[0,\beta)$-ordinal.  For each $n$ we define a
$P_{\underset\tilde {}\to \zeta_n}$-name of $\bold t_n$ of a truth value
and ${\underset\tilde {}\to r_n}$ of a member of $\Bbb P_\beta,
{\underset\tilde {}\to r_n} = {\underset\tilde {}\to r_n} \restriction
[{\underset\tilde {}\to \zeta_n},\beta)$, as follows: if $G^n \subseteq
\Bbb P_{\underset\tilde {}\to \zeta_n}$ is generic over $\bold V$ and $\neg
({\underset\tilde {}\to \zeta_n}[G^n] \ge \underset\tilde {}\to \xi
[G^n])$, and there are $q \in \Bbb P_\beta/G^n$ and $\varepsilon < \omega_1,
\Bbb P_\beta \Vdash ``p \le_{\text{pr}} q"$ and $q \Vdash_{{\Bbb P}_\beta/G^n}
``\underset\tilde {}\to \tau < \varepsilon"$, \ub{then}
${\underset\tilde {}\to {\bold t}_n}[G^n] =$ truth and then
${\underset\tilde {}\to r_n}[G^n]$ is $q \restriction
[{\underset\tilde {}\to \zeta_n},\beta)$ for some such $q$ otherwise
${\underset\tilde {}\to {\bold t}_n}[G^n] =$ false,
$\underset\tilde {}\to r[G^n] = \emptyset$.  Let $r'_n$ be the following
$\Bbb P_{\underset\tilde {}\to \zeta_n}$-name: \nl
for $G^n \subseteq \Bbb P_{\underset\tilde {}\to \zeta_n}$ \ub{if}:
\mr
\item "{$(a)$}"  ${\underset\tilde {}\to {\bold t}_n}[G^n] =$ truth and
\sn
\item "{$(b)$}"  for no $m < \omega$, for some $r \in G^n$ forces a value
to ${\underset\tilde {}\to \zeta_m}$, say $\zeta_m,\zeta_m < \zeta_n[G^n]
\vee (\zeta_m = {\underset\tilde {}\to \zeta_n}[G^n] \and m < n)$
\ermn
\ub{then}
${\underset\tilde {}\to \zeta'_n}[G^n]$ is ${\underset\tilde {}\to r_n}$,
otherwise ${\underset\tilde {}\to \zeta'_n}[G^n] = \emptyset$.  Let
$p^* = p \cup \dbcu_n {\underset\tilde {}\to r'_n}$, easily $q \in
\Bbb P_\beta$ and $p \le_{\text{pr}} p^*$.
\nl
Let ${\underset\tilde {}\to \xi_n} = \text{ Min}\{\underset\tilde {}\to \xi
+1,{\underset\tilde {}\to \zeta_0} +1,\dotsc,
{\underset\tilde {}\to \zeta_{n-1}} +1\}$, ${\underset\tilde {}\to \xi_n}$
is a simple $\bar{\Bbb Q}$-named $[0,\beta)$-ordinal, 
${\underset\tilde {}\to \xi_0} = \underset\tilde {}\to \xi,
{\underset\tilde {}\to \zeta_n} < {\underset\tilde {}\to \xi_{n+1}},
{\underset\tilde {}\to \xi_n} \le {\underset\tilde {}\to \xi_{n+1}}$.  Let
$N$ be a strictly $(\Bbb I,\bold W)$-suitable model for $\chi$ such that
$\{\bar{\Bbb Q},p,\langle {\underset\tilde {}\to \zeta_n},
{\underset\tilde {}\to \xi_n}:n < \omega \rangle,\langle
{\underset\tilde {}\to r_n},{\underset\tilde {}\to r'_n}:n < \omega\}$
belongs to $N$, it exists by \scite{4.8}.  Now we choose $q_n$ by induction
on $n < \omega$ such that

$$
q_n \in \Bbb P_{\underset\tilde {}\to \zeta_n} \text{ is }
(N,{\Bbb I}^{[{\underset\tilde {}\to \kappa_{\underset\tilde {}\to \zeta_n}}]},
P_{\underset\tilde {}\to \zeta_n})-\text{semi}_6 \text{ generic}
$$

$$
p^* \restriction \zeta_n,\beta \le_{\text{pr}} q_n
$$

$$
q_{n+1} \restriction {\underset\tilde {}\to \zeta_n} = q_n.
$$
\mn
This is possible by Fact A.  Now (see \scite{1.27}) $q^* =: \dbcu_{n < \omega}
q_n \cup p^*$ belongs to $\Bbb P_\beta,p^* \le_{\text{pr}} q^*$.  It is enough
to show that $q^* \Vdash \underset\tilde {}\to \tau \in N \cap \omega_1$,
assume not so there is $r,q^* \le r \in \Bbb P_\beta$ such that $r$ forces a value
$\varepsilon^* \in \omega_1 \backslash (N \cap \omega_1)$ to
$\underset\tilde {}\to \tau$.

By \scite{1.15}(1) \wilog \, $q^* \le r$ above $\{\Upsilon_1,\dotsc,
\Upsilon_m\}$ for some $m < \omega$ and $\Upsilon_\ell < \beta$.  There are
$r',k,\xi_k$ such that: $r \le r',\xi_k < \omega,r' \restriction [\xi_k,
\beta) = r \restriction [\xi_k,\beta),r'$ forces
${\underset\tilde {}\to \xi_k}$ is $\xi_k$ and for each $\ell = 1,2,\dotsc,m$
we have $\Upsilon_\ell < \xi_k$ \ub{or} 

$$
\neg(\exists r')(\exists \Upsilon < \beta)(\exists k')[r \restriction
\Upsilon \le r' \in \Bbb P_\Upsilon \and r' \Vdash {\underset\tilde {}\to \xi_{k'}}
= \Upsilon \and \Upsilon_\ell < \Upsilon].
$$
\mn
Clearly $r'$ forces that $(\forall n)[{\underset\tilde {}\to \zeta_n} \ge
\xi_k \rightarrow {\underset\tilde {}\to {\bold t}_n}$ is truth] and we
easily finish.  \hfill$\square_{\text{Fact B}}$.

Now we do not just have to find $q$ satisfying $\Bbb P_\beta/G_\gamma \models 
p \le_{\text{pr}} q$ and $q$ is $(N,\Bbb P_\beta/G_\gamma)$-semi-generic, but 
we need more in the $(N,\Bbb I^{[{\underset\tilde {}\to \kappa_\beta}]}
\restriction \Bbb P_\beta/G_\gamma)$-semi$_6$ generic.  
Now for $G_\beta \subseteq \Bbb P_\beta$ generic over $\bold V$, in
$\bold V[G_\beta]$
for every countable elementary submodel $M$ of $({\Cal H}(\chi)^{\bold V},G_\beta
\in,<^*),\langle \bar{\Bbb Q},\gamma,\beta,\Bbb I \rangle \in M$, we 
have (in \scite{4.5B}) defined Dp$(M)$, an ordinal or $\infty$ and $I_M$ 
such that
\mr
\item "{$(A)$}"  Dp$(M) = \infty$ iff $M[G_\beta] = \{
\underset\tilde {}\to \tau[G_\beta]:\underset\tilde {}\to \tau \in M
\text{ a } \Bbb P_\beta \text{-name}\}$ includes $M$, 
has the same countable ordinals, is $\prec({\Cal H}(\chi)^{\bold V[G_\beta]},\in)$ 
and $M[G_\beta]$ is strictly
$\Bbb I^{[{\underset\tilde {}\to \kappa_\beta}[G_\beta]]}$-suitable
\sn
\item "{$(B)$}"  if Dp$(M) = \alpha < \infty$ then $I_M$ is a member of 
$\Bbb I \cap M$ which
is ${\underset\tilde {}\to \kappa_\beta}[G_\beta]$-complete, and
$$
\align
Y_M = \{t \in \text{ Dom}(I):&\text{there is } N \text{ as above},
M \subseteq N,M \cap \omega_1 = N \cap \omega_1 \\
  &\text{and } t \in N \text{ and Dp}(N) \ge \alpha\} = 
\emptyset \text{ mod } I.
\endalign
$$
\ermn
So we have $\Bbb P_\beta$-names ${\underset\tilde {}\to {\text{Dp}}},
{\underset\tilde {}\to I_M}$. \nl
Consider 

$$
\align
{\frak K} = \bigl\{ (\zeta,G_\zeta,N):&\gamma \le \zeta < \beta,\zeta
\text{ nonlimit}, G_\zeta \subseteq \Bbb P_\zeta \\
  &\text{generic over } \bold V, \text{ in } \bold V[G_\zeta],N[G_\zeta] \text{ is } \\
  &\Bbb I^{[{\underset\tilde {}\to \kappa_\zeta}[G_\zeta]]} \text{-suitable}
\bigr\}.
\endalign
$$
\mn
We now define by induction $({\underset\tilde {}\to \zeta_n},q_n,
p_n,{\underset\tilde {}\to N_n})$ such that
\mr
\item "{$(a)$}"  ${\underset\tilde {}\to \zeta_n}$ is a simple $\bar{\Bbb Q}$-named
$[\gamma,\beta)$-ordinal, (as $e=6$, it is full)
\sn
\item "{$(b)$}"  ${\underset\tilde {}\to N_n}$ is a 
$\Bbb P_{\underset\tilde {}\to \zeta_n}$-name,
$q_n \in \Bbb P_{\underset\tilde {}\to \zeta_n},p_n$ is a
$\Bbb P_{\underset\tilde {}\to \zeta_m}$-name of a member of \nl
${\underset\tilde {}\to N_n}[G_{\underset\tilde {}\to \zeta_n}] \cap
(\Bbb P_\beta/{\underset\tilde {}\to G_{\zeta_n}})$
\sn
\item "{$(c)$}"  ${\underset\tilde {}\to \zeta_0} = \gamma,N_0 = N$
\sn
\item "{$(d)$}"  if $G^n \subseteq \Bbb P_{\underset\tilde {}\to \zeta_n}$ is
generic over $\bold V,q_n \in G^n,G_\gamma \subseteq G^n,
{\underset\tilde {}\to \zeta_n}[G^n] = \zeta_n$ (so essentially
$G_n$ is just a generic
subset of $G_{\zeta_n}$ over $\Bbb V$ such that ${\underset\tilde {}\to \zeta_n}
[G^n] = \zeta_n$), \ub{then} ${\underset\tilde {}\to N_n}[G^n]$ is a
countable elementary submodel of $({\Cal H}(\chi)^{\bold V[G^n]},\in)$ to which
$\bar{\Bbb Q},\gamma,\beta,\Bbb I$ belongs and is strictly
$\Bbb I^{[{\underset\tilde {}\to \kappa_{\zeta_n}}[G^n]]}$-suitable
\sn
\item "{$(e)$}"  ${\underset\tilde {}\to N_n} \subseteq N_{n+1},
{\underset\tilde {}\to N_n} \cap \omega_1 = 
{\underset\tilde {}\to N_{n+1}} \cap \omega_1,q_{n+1} \restriction
{\underset\tilde {}\to \zeta_n} = q_n,p_n \restriction
{\underset\tilde {}\to \zeta_n} \le q_n,p_n \le_{\text{pr}} p_{n+1}$
\sn
\item "{$(f)$}"  for $G^n,\zeta_n,N_n$ as in (d) and $I \in \Bbb I \cap
N_n$ there is $k > n$ such that: if $G^k,\zeta_k,N_k$ are as in (d),
$G^n \subseteq G^k$, \ub{then} there is $t \in \text{ Dom}(I) \cap N_k
\backslash {\underset\tilde {}\to Y_{({\underset\tilde {}\to N_n}[G^n])
[{\underset\tilde {}\to G_\beta}]}}$ (i.e. forced to be there) \nl
(recall $p$ forces a value to ${\underset\tilde {}\to \kappa_\beta})$
sn
\item "{$(g)$}"  for $G^n_n,\zeta_n,N_n$ as in clause (d) above and
$\underset\tilde {}\to \gamma \in N_n$ a $\Bbb P_\beta$-name of a countable
ordinal there is $k > n$ such that: if $G^k,\zeta_k,N_k$ are as in clause
(d), $G^n \subseteq G^k$ then $p_k$ forces $(\Vdash_{{\Bbb P}_\beta/G^k})$ to
be $< N \cap \omega_1$.
\ermn
No problem to carry the definition.  Now
\mr
\item "{$(*)_a$}"  $q = \dbcu_{n < \omega} q_n \in \Bbb P_\beta$ \nl
Here we use the $q_{n+1} \restriction {\underset\tilde {}\to \zeta_n} = q_n$
below $\dbcu_{n < \omega} {\underset\tilde {}\to \zeta_n}$ and
``$({\underset\tilde {}\to {\Bbb Q}_i},\le_{\text{vpr}})$ is $\aleph_1$-complete"
above $\dbcu_{n < \omega} {\underset\tilde {}\to \zeta_n}$
\sn
\item "{$(*)_b$}"  for $G^n,N_n,\zeta_n$ as in clause (d), $q$ is
$(N_n,\Bbb P_\beta/G^n)$-semi generic and above $p \restriction \zeta_n$ \nl
[why?  note clause (g)]
\sn
\item "{$(*)_c$}"  if $q \in G_\beta \subseteq \Bbb P_\beta,G_\beta$ is generic
over $\bold V$ extending $G_\gamma$ \ub{then} in $\bold V[G_\beta]$, 
Dp$({\underset\tilde {}\to N_n}[G_\beta])$ is
well defined
\sn
\item "{$(*)_d$}"  Dp$((N_0[G_\beta]))$ is $\infty$ \nl
(use clause (B) in the demands Dp and clause (f) above).  We use
$I^{[\kappa]}$ is $\kappa$-closed for the relevant $\alpha$'s.
\ermn
So we are done.  \hfill$\square_{\scite{si.7}}$\margincite{si.7}
\enddemo
\bigskip

\proclaim{\stag{si.8} Claim}  1) Assume 
that $\bar{\Bbb Q} = \langle \Bbb P_j,
{\underset\tilde {}\to {\Bbb Q}_i},
{\underset\tilde {}\to \kappa_j}:j \le \alpha,i < \alpha \rangle$ satisfies:
\mr
\item "{$(\alpha)$}"  $\alpha$ is a limit ordinal
\sn
\item "{$(\beta)$}"  if $\beta < \alpha$ then $\langle \Bbb P_j,
{\underset\tilde {}\to Q_i},{\underset\tilde {}\to \kappa_j}:j \le \beta,
i < \alpha \rangle$ is a weak 
${\text{\rm UP\/}}^6(\Bbb I,\bold W,W)$-iteration
\sn
\item "{$(\gamma)$}"  $\Bbb P_\alpha$ is the 
$\aleph_1-{\text{\rm Sp\/}}_6(W)$-limit of
$\langle \Bbb P_i,{\underset\tilde {}\to {\Bbb Q}_i}:i < \alpha \rangle$
\sn
\item "{$(\delta)$}"  ${\underset\tilde {}\to \kappa_\alpha}$, a
$\Bbb P_\alpha$-name is 
$\sup\{{\underset\tilde {}\to \kappa_i}:i < \alpha\}$.
\ermn
\ub{Then} $\bar{\Bbb Q}$ is a weak 
${\text{\rm UP\/}}^6(\Bbb I,\bold W,W)$-iteration. \nl
2) Assume that
\mr
\item "{$(\alpha)$}"  $\langle 
\Bbb P_j,\Bbb Q_i,{\underset\tilde {}\to \kappa_i}:
j \le \alpha,i < \alpha \rangle$ is a weak ${\text{\rm UP\/}}^6
(\Bbb I,\bold W,W)$-iteration
\sn
\item "{$(\beta)$}"  in 
$\bold V^{\Bbb P},\Bbb Q$ is a forcing notion satisfying
${\text{\rm UP\/}}^6\Bbb I,\kappa_\alpha,\underset\tilde {}\to \kappa,\bold W)$, where
$\kappa_\alpha$ is the interpretation of ${\underset\tilde {}\to
\kappa_\alpha}$) and let
${\underset\tilde {}\to {\Bbb Q}},\underset\tilde {}\to \kappa$ be
$\Bbb P_\alpha$-names
of those objects.
\ermn
\ub{Then} there is a ${\text{\rm UP\/}}^6(\Bbb I,\bold W,W)$-iteration
$\langle P_i,{\underset\tilde {}\to {\Bbb Q}_j},
{\underset\tilde {}\to \kappa_i}:
i \le \alpha +1,j < \alpha \rangle$ 
with ${\underset\tilde {}\to {\Bbb Q}_\alpha}
= {\underset\tilde {}\to {\Bbb Q}},{\underset\tilde {}\to \kappa_{\alpha +1}} =
\underset\tilde {}\to \kappa$.
\endproclaim
\bigskip

\demo{Proof}  1) By the proof of \scite{si.7}. \nl
2) Straightforward.
\enddemo
\bigskip

\proclaim{\stag{si6.9} Claim}  Assume
\mr
\item "{$(a)$}"  $\bar{\Bbb Q} = \langle \Bbb P_j,
{\underset\tilde {}\to {\Bbb Q}_i},
{\underset\tilde {}\to \kappa_j}:j \le \alpha^*,i < \kappa \rangle$ is a
weak ${\text{\rm UP\/}}^6(\Bbb I,\bold W,W)$-iteration
\sn
\item "{$(b)$}"  $\gamma < \beta \le \alpha^*$
\sn
\item "{$(c)$}"  $G_\gamma \subseteq \Bbb P_\gamma$ is generic over $V$
\sn
\item "{$(d)$}"  $N$ is a strictly $(\Bbb I,\bold W)$-suitable model $N$
for $(\chi,\langle \bar{\Bbb Q},\gamma,\beta \rangle)$ in 
$\bold V[G_\gamma]$
\sn
\item "{$(e)$}" $p \in N \cap (\Bbb P_\beta/G_\gamma)$.
\ermn
\ub{Then} there is $q$ such that:
\mr
\item "{$(\alpha)$}"  $p \le_{\text{pr}} q \in \Bbb P_\beta/G_\gamma$
\sn
\item "{$(\beta)$}"  $p \restriction \gamma = q \restriction \gamma$
\sn
\item "{$(\gamma)$}"  $q$ is $(N,\Bbb P_\beta/G_\gamma)$-semi generic
\sn
\item "{$(\delta)$}"  $q$ has 
a witness listing $\{\underset\tilde {}\to \zeta
\in N:\underset\tilde {}\to \zeta$ a simple $\bar{\Bbb Q}$-named
$[\gamma,\beta)$-ordinal$\}$.
\endroster
\endproclaim
\bigskip

\demo{Proof}  Like \scite{si.7}, forgetting to take care of the
$I \in \Bbb I \cap {\underset\tilde {}\to N_n}$ so we can use
$\underset\tilde {}\to N = 
N[{\underset\tilde {}\to G_{{\Bbb P}_{\underset\tilde {}\to
\zeta_n}}}]$. \nl
${{}}$  \hfill$\square_{\scite{si6.9}}$\margincite{si6.9}
\enddemo
\bigskip

\proclaim{\stag{si6.10} Claim}  Suppose $F_f,F_c$ are functions (possibly
classes), $\bold W \subseteq \omega_1$ is stationary, $\Bbb I$ is a class
of ($\aleph_2$-complete) quasi order ideals, $W$ a class of strongly
inaccessible cardinals.

\ub{Then} for every ordinal $\alpha$ there is a unique $\bar{\Bbb Q} =
\langle \Bbb P_j,
{\underset\tilde {}\to {\Bbb Q}_i},{\underset\tilde {}\to \kappa_j}:j \le 
\alpha^\dag,i < \alpha \rangle$ such that:
\mr
\item "{$(a)$}"  $\bar{\Bbb Q}$ is a ${\text{\rm UP\/}}^6
(\Bbb I,\bold W,W \cap (\alpha^\dag
+1))$-iteration
\sn
\item "{$(b)$}"  for every $i < \alpha^\dag$ we have 
${\underset\tilde {}\to {\Bbb Q}_i} = F_f(\bar{\Bbb Q} \restriction i),
{\underset\tilde {}\to \kappa_i} = F_c(\bar{\Bbb Q} \restriction i)$
\sn
\item "{$(c)$}"   $\alpha^\dag \le c$
\sn
\item "{$(d)$}"  for limit $\beta \le \alpha^\dag$ we have
${\underset\tilde {}\to \kappa_\beta} = \sup
\{{\underset\tilde {}\to \kappa_\gamma}:\gamma < \beta\}$
\sn
\item "{$(e)$}"  if $\alpha^\dag < \alpha$ then the following is impossible
{\roster
\itemitem{ $(\alpha)$ }   $F_f(\bar{\Bbb Q} \restriction i)$ is a
$\Bbb P_{\alpha^\dag}$-name of a forcing notion
\sn
\itemitem{ $(\beta)$ }  $F_c(\bar{\Bbb Q} \restriction i)$ is a
$\Bbb P_{\alpha^\dag}$-name of a $F_f(\bar{\Bbb Q} \restriction i)$-name of a
$\bold V$-cardinal $\ge \aleph_1$
\sn
\itemitem{ $(\gamma)$ }  $\Vdash_{{\Bbb P}_{\alpha^\dag}} 
``F_f(\bar{\Bbb Q} \restriction i)$ is 
${\text{\rm UP\/}}^6_{{\underset\tilde {}\to \kappa_{\alpha^\dag}},
{\underset\tilde {}\to \kappa}}(\Bbb I,\bold W)"$.
\endroster}
\endroster
\endproclaim
\bigskip

\demo{Proof}  Straight.
\enddemo
\bn
Next we give sufficient conditions for $\dbcu_{i < \delta} \Bbb P_i$ being a
dense subset of $\Bbb P_j$
\proclaim{\stag{si6.11} Claim}  Assume that $\langle \Bbb P_j,
{\underset\tilde {}\to {\Bbb Q}_i},{\underset\tilde {}\to \kappa_j}:
j \le \alpha^*,i < \alpha^* \rangle$ is a 
${\text{\rm UP\/}}^6(\Bbb I,\bold W,W)$-iteration. \nl
1) Assume
\mr
\item "{$(*)_1$}"  $(\forall i < \alpha^*)[\Vdash_{{\Bbb P}_i}
``\le^{{\Bbb Q}_i}_{\text{vpr}}$ is equality"].
\ermn
If ${\text{\rm cf\/}}(\delta) = \aleph_1 \and \delta \le \alpha^*$ \ub{then}
$\dbcu_{i < \delta} \Bbb P_i$ is a dense subset of $\Bbb P_\delta$ even under
$\le_{\text{pr}}$. \nl
2) Assume
\mr
\item "{$(*)_2$}"  $(\forall i < \alpha^*)(\Vdash_{{\Bbb P}_i}
``(\{r \in {\underset\tilde {}\to {\Bbb
Q}_i}:\emptyset_{\underset\tilde {}\to {\Bbb Q}_i}
\le_{\text{vpr}} r\},\le^{\underset\tilde {}\to {\Bbb Q}_i}_{\text{vpr}})$ is
directed).
\ermn
If $\delta \le \alpha^*,{ \text{\rm cf\/}}(\delta) = \aleph_1$ \ub{then}
\mr
\item "{$(a)$}"  $\dbcu_{i < \delta} {\Cal D}_{\delta,i}$ is a dense 
subset of $\Bbb P_\delta$ even under $\le_{\text{pr}}$ where
${\Cal D}_{\delta,i} = \{p \in \Bbb P_\delta:\Bbb P_\delta \models ``p \restriction
i \le_{\text{vpr}} p"\}$
\sn
\item "{$(b)$}"  if $i < \delta,\{p_0,p_1\} \subseteq {\Cal D}_{\delta,i}$
and $p_0 \restriction i = p_j \restriction i$ then $p_0,p_1$ has an upper
bound $p$ even in ${\Cal D}_{\delta,i},p \restriction i = p_0 \restriction
i = p_1 \restriction i$.
\ermn
3) If $\delta \in W \and \delta \le \alpha^* \and (\forall i < \delta)$
[${\text{\rm density\/}}(\Bbb P_i) 
< \delta$] \ub{then} $\Bbb P_\delta = \dbcu_{i < \delta} P_i$.
\endproclaim
\bigskip

\demo{Proof}  1), 2)  Let $p \in \Bbb P_\delta$, choose $\chi$ large enough.
There is a strictly $(\Bbb I,\bold W)$-suitable countable model
$N \prec ({\Cal H}(\chi),\in)$ to which $\{\bar{\Bbb Q},\delta,p\}$ belongs.
Applying \scite{si.7} for $\gamma = 0,\beta = \delta$
(i.e. $\boxtimes_{\gamma,\beta}$ from the proof) we can find $q \in
\Bbb P_\delta,p \le_{\text{pr}} q$ which is 
$(N,\Bbb I^{[{\underset\tilde {}\to \kappa_\delta}]},\Bbb P_\delta)$-semi$_6$ 
generic
and $q$ has a witness $\subseteq \{\underset\tilde {}\to \zeta \in N:
\underset\tilde {}\to \zeta$ a simple $\bar{\Bbb Q}$-named $[0,\delta)$-ordinal$\}$.
As cf$(\delta) = \aleph_1$ there is an increasing continuous $\bar \beta =
\langle \beta_\varepsilon:\varepsilon < \omega_1 \rangle$ with limit
$\delta$, \wilog \, $\bar \beta \in N$ so $q \Vdash ``
N[{\underset\tilde {}\to G_{{\Bbb P}_\delta}}] \cap \delta = \cup
\{\beta_\varepsilon:\varepsilon \in N
[{\underset\tilde {}\to G_{{\Bbb P}_\delta}}] \cap \omega_1\}"$ but
$q \Vdash \cup \{\beta_\varepsilon:\varepsilon \in 
N[{\underset\tilde {}\to G_{{\Bbb P}_\delta}}] \cap \omega_1\} = \cup
\{\beta_\varepsilon:\varepsilon \in N \cap \omega_1\} =
\beta_{N \cap \omega_1}$, so clearly $\Bbb P_\delta \models q \restriction
\beta_{N \cap \omega_1} \le_{\text{vpr}}$.  For (1) it follows that
$q \in P_{N \cap \delta}$ and we are done. \nl
For (2) just reflect. \nl
3) Straight.
\enddemo
\bigskip

\demo{\stag{si6.12} Conclusion}  Let $\bar{\Bbb Q} = \langle \Bbb P_j,
{\underset\tilde {}\to {\Bbb Q}_i},{\underset\tilde {}\to \kappa_j}:j \le \alpha^*,
i < \alpha^* \rangle$ be a UP$^6(\Bbb I,\bold W,W)$-iteration and
$\kappa = \text{ cf}(\kappa) \le \alpha^*$ and $(\forall i < \delta)$
(density$(\Bbb P_i) < \kappa)$. \nl
1) If $\{\theta < \kappa:\theta = \text{ cf}(\theta) \in W\}$ is
stationary.  
\ub{Then} $\Bbb P_\kappa$ satisfies the $\kappa$-c.c. (in fact a strong
version and even under $\le^{{\Bbb P}_\kappa}_{\text{vpr}}$. \nl
2) If $(*)_1$ from \scite{si6.11}(1), i.e. $\Vdash_{{\Bbb P}_i} 
``\le^{\underset\tilde {}\to {\Bbb Q}_i}_{\text{pr}}$ is equality" 
for $i < \kappa$ \ub{then} 
$\Bbb P_\kappa$ satisfies the $\kappa$-c.c. (in fact a strong version)
even for $\le^{{\Bbb P}_\kappa}_{\text{vpr}})(\theta \le_{\text{vpr}})$ is
$(2^{|\delta|})^+$-complete. \nl
3) Assume $\kappa \notin W$ and $S \subseteq \{\delta < \kappa:
\text{cf}(\delta) = \aleph_1\}$ is stationary and $i \ge \delta \and \delta
\in S \Rightarrow \Vdash_{{\Bbb P}_i} 
``(\{r:\emptyset_{\underset\tilde {}\to {\Bbb Q}_i}
\le_{\text{vpr}} r\}$.  \ub{Then} forcing with $\Bbb P_\kappa$ does not add a
function in ${}^\omega\text{Ord}$ not in $\dbcu_{\beta < \kappa}
\bold V^{{\Bbb P}_\beta}$, even any function in ${}^{\alpha(*)}$Ord $\backslash
\dbcu_{\beta < \kappa} \bold V^{{\Bbb P}_\beta},\alpha(*) < \kappa$.
\enddemo
\bigskip

\demo{Proof}  1) Let $S = \{\theta < \kappa:\theta = \text{ cf}(\theta) 
\in W\}$ 
and let $\langle \Bbb P_\theta:\theta \in S \rangle$ be a sequence of
members of 
$\Bbb P_\kappa$.  So for each $\theta \in S$ for some $i(\theta) < 
\theta$ 
we have $p_\theta \restriction \theta \in \Bbb P_{i(\theta)}$ and we can
find a pressing down function $h$ on $S$ such that $h(\theta_1) =
h(\theta_2) \Rightarrow p_{\theta_1} \restriction \theta_1 = p_{\theta_2}
\restriction \theta_2$.  Clearly there is a club $E$ or $\kappa$ such that
$\theta_i \in \theta_2 \cap s \and \theta_2 \in S \cap S \cap E \Rightarrow
p_{\theta_1} \in P_{\theta_2}$.

Lastly, if $\theta_1,\theta_2 \in E \cap S$ and $f(\theta_1) = f(\theta_2)$
then $p_{\theta_1} \cup p_{\theta_2}$ is a common upper bound of 
$p_{\theta_1},p_{\theta_2}$ (even a $\le_{\text{vpr}}$ one). \nl
2) Similar using $S = \{\delta < \kappa:\text{cf}(\delta) = \aleph_1$ and
$\alpha < \delta \Rightarrow$ density$(\Bbb P_\alpha) < \delta\}$. \nl
3) If $W \cap \kappa$ is stationary in $\kappa$ use part (2), so let $E$ be
a club of $\kappa$ disjoint to $W$.  Assume toward contradiction that
$p \Vdash_{{\Bbb P}_\kappa} ``\text{the function } \underset\tilde {}\to \tau:
\omega_1 \rightarrow {}^\kappa\text{Ord}$ is not in $\dbcu_{\beta < \kappa}
\bold V^{{\Bbb P}_\beta}"$, let $S$ be as in part (2).  We choose by induction on
$j < \kappa,(p_j,\alpha_j)$ and if $\alpha _j \in \text{ Salso}(q_j,
\varepsilon_j,\gamma_j,\beta_j)$ such that:
\mr
\widestnumber\item{$(iii)$}
\item "{$(i)$}"  $\alpha_j \in E$ is increasing continuous
\sn
\item "{$(ii)$}"  $p_j \in \Bbb P_\kappa,p_0 = p$
\sn
\item "{$(iii)$}"  $i < j \Rightarrow p_j  \restriction \alpha_i = p_i$
\sn
\item "{$(iv)$}"  $i < j \and \beta \in [\alpha_i,\alpha_j) \Rightarrow
\Vdash_{{\Bbb P}_\beta} ``\emptyset \le_{\text{vpr}} p_j \restriction \{\beta\}$
in ${\underset\tilde {}\to {\hat{\Bbb Q}}_\beta}"$
\sn
\item "{$(v)$}"  if $\alpha_i \in S,j < \alpha_i$ \ub{then} for every
$p_i \le q \in {\Cal D}_{\alpha_i,j}$ there is $\varepsilon_i(q) < \omega_i$
such that: \ub{if} there are $\varepsilon,r$ such that $p_{i+1} \cup 
\le_{\text{vpr}} r,r \restriction \alpha_i = q,\varepsilon < \omega_1,r$
forcing a value to $\underset\tilde {}\to \tau(\varepsilon)$ but for no
$r,q \le r' \in \Bbb P_{\alpha_i}$, does $P_{i+1} \cup r'$ forces a value to
$\underset\tilde {}\to \tau(\varepsilon)$ then $\varepsilon = \varepsilon_i
(q)$ satisfies this
\sn
\item "{$(vi)$}"  $\alpha_i \in S$ then $p_{i+1} \le q_i,q_i \Vdash
``\underset\tilde {}\to \tau(\varepsilon_i) = \gamma_i",\varepsilon_i <
\alpha(*),\beta_i < \alpha_i,q_i \restriction \beta_i \le_{\text{vpr}} q_i
\restriction \alpha_i$
\sn
\item "{$(vii)$}"  there is no $q,q_i \restriction \alpha_i \le q \in
\Bbb P_{\alpha_i}$ such that $p_{i+1} \cup (q_i \restriction \alpha_i)$ forces a
value to $\underset\tilde {}\to \tau(\varepsilon_i)$.
\ermn
For any $j$ we choose $(p_j,\alpha_j)$. \nl
For $j=0$ let $p_0=p$ and by \scite{si.8} for some $\alpha_0 < \kappa,p_0
\in \Bbb P_{\alpha_0}$.  For $j=i+1$, first choose $p_j$ to satisfy clause (v)
and then $\alpha_j$ such that $p_j,q \in \Bbb P_{\alpha_j}$.  Lastly for $j$
limit let $p_j = \dbcu_{i < j} p_i,\alpha_j = \dbcu_{i < j} \alpha_i$ and
check.
The contradiction is easy.

Let $G_{\alpha_i} \subseteq \Bbb P_{\alpha_i}$ be generic over $\bold V$ such that
$p_i \in G_{\alpha_i}$.  Clearly for some $\varepsilon < \alpha(*)$ for no
$q \in \Bbb P_{\alpha_i}$ do we have $q \cup p_{i+1}$ forces a value to
$\underset\tilde {}\to \tau(\varepsilon)$ as otherwise
$p_{\alpha_{i+1}} \Vdash_{{\Bbb P}_\kappa/G_{\alpha_i}} \underset\tilde {}\to \tau
\in \bold V[G_{\alpha_i}]$.  Choose $\varepsilon_i < \alpha(*)$ as above, choose
$q'_i \in G_{\alpha_i}$ which forces this choose $q_i \in \Bbb P_\kappa$ above
$q'_i \cup p_{i+1}$ which forces a value to $\underset\tilde {}\to \tau
(\varepsilon_i)$ and \wilog \, there is $\gamma _i < \alpha_i$ such that
$q_i \restriction \gamma_i  \le_{\text{vpr}} q_i \restriction \alpha_i$.
Lastly let $\alpha_{i+1}$ be such that $\beta \in [\alpha_{i+1},\kappa)
\Rightarrow \Vdash_{P_\beta} ``\emptyset_{\underset\tilde {}\to {\Bbb Q}_\beta}
\le_{\text{vpr}} p_{i+1} \restriction \{\beta\}$. \nl
[Saharon: the role of $W$].
\enddemo
\newpage

\head {\S6B On UP$^2$-iteration} \endhead  \resetall \sectno=6
\bigskip

\proclaim{\stag{6.7}(?)Lemma}  Assume that $\bold W \subseteq \omega_1$ is 
stationary and $\bar{\Bbb Q} = \langle \Bbb P_i,{\underset\tilde
{}\to {\Bbb Q}_i},
{\underset\tilde {}\to {\Bbb I}_i},{\underset\tilde {}\to \kappa_i}:
i < \alpha \rangle$ is a ${\text{\rm UP\/}}^{2,e}(\bold W,W)$-suitable 
iteration, and $\Bbb P_\alpha 
= { \text{\rm Sp\/}}_e(W)-{\text{\rm Lim\/}}_\kappa
(\bar{\Bbb Q})$ be the limit and 
$\underset\tilde {}\to \kappa(\beta) = \sup
\{{\underset\tilde {}\to \kappa_\gamma}:\gamma < \beta\}$ (this is a
$\Bbb P_\beta$-name of a $V$-cardinal) and

$$
\kappa^-(\beta) = {\text{ \rm Sup\/}}
\{\kappa:\,\,\nVdash_{{\Bbb P}_\gamma} ``\kappa \ne
\underset\tilde {}\to \kappa(\beta)" \text{ for some } \gamma < \beta\}.
$$
\mn
1) For each $\beta \le \alpha,\Bbb P_\beta$ satisfies 
${\text{\rm UP\/}}^{0,e}_{\underset\tilde {}\to \kappa
(\beta)}(\Bbb I'_\beta,\bold W)$ for 
some $\aleph_2$-complete $\Bbb I'_\beta \in \bold V$. \newline
2) In fact, $\Bbb I_\beta$ is $\kappa^\beta$-complete where
$\kappa^\beta = \text{ min}\{\kappa:\text{for some } \gamma < \beta$ we
have \nl
$\nVdash_{{\Bbb P}_{\gamma +1}} ``{\underset\tilde {}\to \kappa_\gamma} \ne
\kappa"\}$, and each $I \in \Bbb I_\beta$ has domain of cardinality \nl
$\le (\underset{\gamma < \beta} \to \sup\{ \lambda < \kappa_{\delta +1}:\,
\nVdash_{{\Bbb P}_\gamma} ``\neg(\exists I \in 
{\underset\tilde {}\to {\Bbb I}_\gamma})(\lambda = |\text{Dom}(I)|)"\})$ and
\nl
$|\Bbb I'_\beta| 
\le \dsize \sum_{\gamma < \beta} (\aleph_0 + |\Bbb P_\gamma| +
\text{ min}\{\lambda:\Vdash_{{\Bbb P}_\gamma} 
``|{\underset\tilde {}\to {\Bbb I}_\gamma}| \le \lambda)^{< \kappa}\}$. \nl
3) Similarly for ${\text{\rm UP\/}}^{0,e}_a$ 
and weak ${\text{\rm UP\/}}^{0,e}_y(\bold W)$-iterations.
\endproclaim
\bigskip

\remark{\stag{6.7A} Remark}  We can also get the preservation version of this
Lemma.
\endremark
\bigskip

\demo{Proof}  For each $\gamma < \alpha$ let ${\Cal J}_\gamma =:
\{q \in \Bbb P_{\gamma +1}:q \text{ forces a value to } 
{\underset\tilde {}\to \kappa_\gamma}$, called 
$\kappa_{\gamma,q}$ and $q$ forces
${\underset\tilde {}\to {\Bbb I}_\gamma}$ to be equal to a $\Bbb P_\gamma$-name
${\underset\tilde {}\to {\Bbb I}_{\gamma,q}}$ and
$q \restriction \gamma$ forces a value to
$|\Bbb I_\gamma|$ says $\mu_{\gamma,q}$ is this is purely decidable, if not,
just an upper bound to it$\}$; let ${\Cal J}'_\gamma \subseteq
{\Cal J}_\gamma$ a maximal antichain.  Let $\mu_\gamma =
\underset{q \in {\Cal J}'_\gamma} \to \sup \mu_{\gamma,q}$. \newline
Let $q \Vdash_{{\Bbb P}_\gamma} ``{\underset\tilde {}\to {\Bbb I}_\gamma} = \{
{\underset\tilde {}\to I_{\gamma,\zeta}}:\zeta < \mu_{\gamma,q}\}"$ for
$q \in {\Cal J}_\gamma$ and let
${\Cal J}_{\gamma,\zeta} = \{q \in {\Cal J}_\gamma:\mu_{\gamma,q} > \zeta
\text{ and } q \Vdash \text{ Dom}({\underset\tilde {}\to I_{\gamma,\zeta}})
\text{ is } \lambda_{\gamma,q,\zeta}\}$ and let
${\underset\tilde {}\to I_{\gamma,\zeta}} \text{ be id}_{L_{\gamma,\zeta}}$, 
so $L_{\gamma,q,\zeta}$ is a
$\Bbb P_\gamma$-name of a ${\underset\tilde {}\to \kappa_{\gamma,q}}$-directed 
partial order on $\lambda_{\gamma,q,\zeta}$ (but $\Vdash_{{\Bbb P}_\gamma}$ ``if
$|{\underset\tilde {}\to {\Bbb I}_\gamma}| \le \zeta < \mu_\gamma$ then let
$L_{\gamma,\zeta}$ be trivial").

For $q \in {\Cal J}_\gamma$ let 
$L^*_{\gamma,q,\zeta}$ be ap$_{\kappa_{\gamma,q}}(
{\underset\tilde {}\to L_{\gamma,\zeta}})$ for the forcing notions \nl
$\Bbb P_\gamma^{[q]} = \{p \in \Bbb P_\gamma:q \restriction \gamma 
\le^{{\Bbb P}_\gamma}_{pr} p\}$ from Definition \scite{3.9}.  So by Claim \scite{3.10}
\medskip
\roster
\widestnumber\item{$(iii)$}
\item "{$(i)$}"   $L^*_{\gamma,q,\zeta}$ is $\kappa_{\gamma,q}$-directed
partial order on $[\lambda_{\gamma,q,\zeta}]^{< \kappa_{\gamma,q}}$
\sn
\item "{$(ii)$}"  $|L^*_{\gamma,q,\zeta}| \le (\lambda_{\gamma,q,\zeta})
^{< \kappa_{\gamma,q}}$
\sn
\item "{$(iii)$}"  $q \restriction \gamma \Vdash_{{\Bbb P}_\gamma} 
``{\underset\tilde {}\to I_{\gamma,\zeta}} = 
\text{ id}_{\underset\tilde {}\to L_{\gamma,\zeta}} \le_{\text{RK}} \text{ id}
_{L^*_{\gamma,q,\zeta}}"$.
\endroster
\medskip

\noindent
Let $\kappa_\beta = \sup\{\kappa_{\gamma,q}:\gamma < \beta,q \in
{\Cal J}_\gamma\}$. 
\sn
Let $\Bbb I^*_\beta$ be the $(< \kappa_\beta)$-closure of $\{
\text{id}_{L^*_{\gamma,q,\zeta}}:\gamma < \beta,q \in {\Cal J}_\gamma,
\zeta < \mu_{\gamma,q}\}$ (see Definition \scite{3.13}(1)).
\smallskip

Let $N$ be $(\Bbb I^*_\alpha,\bold W)$-suitable model for $(\chi,\lambda),x$
code enough informtion so for some $\bar N,N = N_{\langle \rangle}$ and 
$\bar N = \langle N_\eta:\eta \in (T,\bold I) \rangle$ be a strict 
$(\Bbb I^*_\alpha,\bold W)$-suitable tree of models for $(\chi,x),x$ 
coding enough information
(so $\bar{\Bbb Q},\Bbb I^*_\alpha,{\underset\tilde {}\to {\bold S}},\bold W \in
N_{\langle \rangle}$).

Let ${\Cal T}_{\bar N}$ be the set of pairs $(\gamma,p)$ such that:
\medskip
\roster
\item "{$\bigotimes$}"  $\gamma \le \alpha,p \in \Bbb P_\gamma$, and for some
$\kappa$, \newline
$p \Vdash_{P_\gamma} ``N[{\underset\tilde {}\to G_{P_\gamma}}] \cap \omega_1 
= N_{\langle \rangle} \cap \omega_1$ and $\gamma \in N[G_{P_\gamma}]$".
\endroster
\medskip

\noindent
${\Cal T}'_{\bar N}$ is defined similarly as the set of pairs
$(\underset\tilde {}\to \gamma,p)$ such that:
$\underset\tilde {}\to \gamma$ is a simple $(\bar{\Bbb Q},W)$-named ordinal, $p \in
\Bbb P_{\underset\tilde {}\to \gamma}$.  (I.e. if $\zeta <
\beta,G_{{\Bbb P}_\gamma}
\subseteq \Bbb P_\zeta$ is generic over $\bold V$ and $\zeta = {\underset\tilde {}\to 
\gamma_n}[G_{{\Bbb P}_\zeta}]$ then $r \in q_n \Rightarrow
{\underset\tilde {}\to \zeta_n}[G_\zeta] < \zeta$, i.e. is well defined
$< \zeta$ \ub{or} is forced ($\Vdash_{{\Bbb P}_\alpha/G_\zeta}$) to be not 
well defined, and $p \Vdash_{{\Bbb P}_{\underset\tilde {}\to \gamma}} 
``\eta \in \text{ lim } T"$).
\medskip

We consider the statements, for $\gamma \le \beta < \alpha$ (or restrict
ourselves to $\gamma$ non-limit)
\medskip
\roster
\item "{$\boxtimes_{\gamma,\beta}$}"  for any 
$(\gamma,p) \in {\Cal T}_{\bar N}$ and $\rho$ such that \newline
$p \Vdash_{{\Bbb P}_\gamma} ``\rho \triangleleft \underset\tilde {}\to \eta"$
\newline
and ${\underset\tilde {}\to p'}$ a $\Bbb P_\gamma$-name such that 
$p \Vdash_{{\Bbb P}_\gamma} ``{\underset\tilde {}\to p'}[G_{{\Bbb P}_\gamma}] \in
N[{\underset\tilde {}\to G_{{\Bbb P}_\gamma}}] \cap P_\beta/G_{{\Bbb P}_\gamma}$ and
$({\underset\tilde {}\to p'}[G_{{\Bbb P}_\gamma}]) \restriction \gamma \le_{pr} p"$ 
\ub{there is} $(\beta,q) \in {\Cal T}$ such that 
${\underset\tilde {}\to p'} \le_{pr} q$ \nl
(i.e. $p \Vdash_{{\Bbb P}_\gamma} ``{\underset\tilde {}\to p'}
[{\underset\tilde {}\to G_{{\Bbb P}_\gamma}}] \le_{pr} q"$) and 
$q \restriction \gamma = p"$.
\endroster
\medskip

We prove by induction on $\beta \le \alpha$ that $(\forall \gamma \le \beta)
\, \boxtimes_{\gamma,\beta}$ (but for \scite{6.7}(3), we use
($\forall$ non-limit $\gamma \le \beta$) $\boxtimes_{\gamma,\beta}$), note 
that for $\gamma = \beta$ the statement is trivial hence we shall consider 
only $\gamma < \beta$.
\bigskip

\noindent
\underbar{Case 1}:  $\beta = 0$.

Trivial.
\bigskip

\noindent
\underbar{Case 2}:  $\beta$ a successor ordinal (for part (3), $\beta$
successor of non-limit ordinal).

As trivially $\boxtimes_{\gamma_0,\gamma_1} \and \boxtimes_{\gamma_1,\gamma_2}
\Rightarrow \boxtimes_{\gamma_0,\gamma_2}$, clearly without loss of
generality \nl
$\beta = \gamma + 1$.

Let $G_{{\Bbb P}_\gamma}$ be such that $p \in G_{{\Bbb P}_\gamma} 
\subseteq \Bbb P_\gamma$ and
$G_{{\Bbb P}_\gamma}$ generic over $V$.

Let $\bar N' = \langle N_\nu[G_{{\Bbb P}_\gamma}]:\eta \in (T',\bold I) \rangle$. \nl
In $\bold V[G_{{\Bbb P}_\gamma}]$ we apply \sciteu{6.2} for $\lambda = \aleph_2$ to
$\bar N'$ and find $T' \subseteq T$ such that $\bold V[G_{{\Bbb P}_\gamma}] \models 
``\langle N'_\eta[G_{{\Bbb P}_\gamma}]: \eta \in (T'',\bold I') \rangle$ is strict 
$(\Bbb I^*_\gamma)^{[\kappa(\gamma)]},\bold W)$-suitable".  So we can apply
clause (f) of \sciteu{6.4A}.
\bn
\ub{Discussion}:  This is a question whether there is an $\Bbb I$-tree of
model $\langle N_\eta:\eta \in (T,\Bbb I) \rangle$ such that:
\mr
\item "{{}}"  if $\Bbb I_\eta$ is $\lambda(\Bbb I_\eta)$-complete,
$\lambda(\bold I_\eta)$ regular, $\alpha < \lambda(\bold I_\eta)
\Rightarrow |\alpha|^{\aleph_0} < \lambda(\bold I_\eta)$, \ub{then}
$\nu \in \text{ Suc}_T(\eta) \Rightarrow N_\eta <_\lambda N_\nu$.
\ermn
This would make \scite{si.2}? more effective.
\bigskip

\noindent
\underbar{Case 3}:   $\beta$ is a limit ordinal.

We shall choose by induction on $n < \omega,{\underset\tilde {}\to \gamma_n},
q_n,{\underset\tilde {}\to p_n}$ such that:
\medskip
\roster
\item "{$(a)$}"  $({\underset\tilde {}\to \gamma_n},q_n) \in 
{\Cal T}'_{\bar N}$ \newline
(so ${\underset\tilde {}\to \gamma_n}$ is a simple $(\bar{\Bbb Q},W)$-named ordinal)
\sn
\item "{$(b)$}"  ${\underset\tilde {}\to \gamma_0} = \gamma$ and
$\Vdash_{\bar{\Bbb Q}}``{\underset\tilde {}\to \gamma_n} < 
{\underset\tilde {}\to \gamma_{n+1}} < \beta"$ \newline
i.e. if $\zeta < \beta$ and $G_{{\Bbb P}_\gamma}
\subseteq \Bbb P_\zeta$ is generic over $V$ and $\zeta = {\underset\tilde {}\to 
\gamma_n}[G_{{\Bbb P}_\zeta}]$ then \nl
$r \in q_n \Rightarrow {\underset\tilde {}\to \zeta_n}[G_\zeta] < \zeta$ 
(i.e. is well defined $< \zeta$ or is forced to be not well defined),
\sn
\item "{$(c)$}"  $q_{n+1} \restriction {\underset\tilde {}\to \gamma_n} =
q_n$
\sn
\item "{$(d)$}"   ${\underset\tilde {}\to p_n}$ is a
$\Bbb P_{\underset\tilde {}\to \gamma_n}$-name, $p_0 = p,
{\underset\tilde {}\to p_n} \restriction {\underset\tilde {}\to \gamma_n}
\le_{pr} q_n$ and \newline
$q_n \Vdash_{{\Bbb P}_{\underset\tilde {}\to \gamma_n}} ``
{\underset\tilde {}\to p_n} \in N_{\underset\tilde {}\to \rho_n}
[G_{{\Bbb P}_{\underset\tilde {}\to \gamma_n}}] \cap \Bbb P_\beta$ and 
$p_n \restriction {\underset\tilde {}\to \gamma_n} \in G_{{\Bbb P}_
{\underset\tilde {}\to \gamma_n}}"$
\sn
\item "{$(e)$}"  $q_n \Vdash_{{\Bbb P}_{\underset\tilde {}\to \gamma_n}}
``{\underset\tilde {}\to p_n} \le^{{\Bbb P}_\beta}_{pr}
{\underset\tilde {}\to p_{n+1}} \in N[G_{{\Bbb P}_{\gamma_n}}] \cap P_\beta"$
\sn
\item "{$(f)$}"  letting $\langle {\underset\tilde {}\to \tau_\ell}:
\ell < \omega \rangle$ list the $P_\beta$-names of ordinals from $N$: 
for $m,\ell \le n$ we have:
\endroster
\medskip
$$
\align
q_n \Vdash_{{\Bbb P}_{\underset\tilde {}\to \gamma_n}} ``{\underset\tilde {}\to 
p_{n+1}} \text{ force } (\Vdash_{{\Bbb P}_{\underset\tilde {}\to \gamma_{n+1}}})
\text{ that:} \, &\text{ if } 
{\underset\tilde {}\to \tau_\ell} 
\text{ is a countable ordinal, then it is smaller} \\
  &\quad \text{than some }  
{\underset\tilde {}\to \tau'_\ell} \in N
[G_{{\Bbb P}_{\underset\tilde {}\to \gamma_{n+1}}}], \\
  &\quad \text{a } 
\Bbb P_{\underset\tilde {}\to \gamma_n} \text{-name of a countable ordinal}".
\endalign
$$
\medskip

\noindent
The induction is straight and $\dsize \bigcup_{n < \omega} q_n$
are as required noting we need and have $(*)_1$ or $(*)_2$ below:
\mr
\item "{$(*)_1$}"  Assume $\le_{\text{pr}},\le_{\text{vpr}}$ 
are equal to $\le$ \nl
(i.e., $\Vdash_{{\Bbb P}_\beta} 
``\le^{\underset\tilde {}\to {\Bbb Q}_\beta}_{\text{vpr}}$ is
$\le^{\underset\tilde {}\to {\Bbb Q}_\beta}"$ for each $\beta < \alpha$),
if $p \in \Bbb P_\beta,\gamma < \beta,\underset\tilde {}\to \tau$
a $\Bbb P_\beta$-name of an ordinal \ub{then} there are $p',
{\underset\tilde {}\to \tau'}$ such that:
{\roster
\itemitem{ (i) }  ${\underset\tilde {}\to \tau'}$ 
is a $\Bbb P_\gamma$-name of an ordinal
\sn
\itemitem{ (ii) }  $p \le_{\text{pr}} 
p' \in \Bbb P_\beta$ and $p \restriction \gamma = p' \restriction \gamma$
\sn
\itemitem{ (iii) }  $p' \Vdash_{{\Bbb P}_\beta} ``
\underset\tilde {}\to \tau = {\underset\tilde {}\to \tau'}"$. \nl
[why?  straight by \scitet{1.16}].
\endroster}
\sn
\item "{$(*)_2$}"  if $p \in \Bbb P_\beta,\gamma < \beta,
\underset\tilde {}\to \tau$ is a $\Bbb P_\beta$-name of a countable ordinal,
\ub{then} there are $p',\tau'$ such that
{\roster
\itemitem{ (i) }  $\tau'$ is a $\Bbb P_\gamma$-name of a countable ordinal
\sn
\itemitem{ (ii) }  $p \le_{\text{pr}} 
p' \in \Bbb P_\beta$ and $p' \restriction \gamma = p \restriction \gamma$
\sn
\itemitem{ (iii) }  $p' 
\Vdash_{{\Bbb P}_\beta} ``\underset\tilde {}\to \tau \le \tau'\,"$ 
\endroster}
\sn
[why $(*)_2$?  let $\underset\tilde {}\to \zeta$ be the following simple
$(\bar{\Bbb Q},W)$-named $[\gamma,\beta)$-ordinal: \nl
$G_\zeta \subseteq \Bbb P_\zeta$ is generic over $\bold V$ for 
$\zeta \in [\gamma,\beta)$ we let
$\underset\tilde {}\to \zeta[G_\zeta] = \zeta$ if
{\roster
\itemitem{ $(a)$ }  $p \restriction \zeta \notin G_\zeta$ \ub{or}: \nl
for some $p' \in \Bbb P_\beta$ we have $p' \restriction \zeta
= p$ and $\Bbb P_\beta \models p \le_{\text{pr}} p'$ \nl
and $p' \Vdash_{{\Bbb P}_\beta/G_\zeta} ``\underset\tilde {}\to \tau < \tau^*"$ for 
some countable ordinal $\underset\tilde {}\to \tau$ 
\sn
\itemitem{ $(b)$ }  for no $\xi \in [\gamma,\zeta)$ does clause (a) hold for
$\xi,G_\zeta \cap \Bbb P_\xi$. 
\endroster}
\sn
Now if $p \Vdash_{{\Bbb P}_\alpha} ``\underset\tilde {}\to \zeta = \gamma"$ we are
done.  Also $\Vdash_{{\Bbb P}_\alpha} ``\underset\tilde {}\to \zeta
[{\underset\tilde {}\to G_{{\Bbb P}_\alpha}}]$ is well defined" as if $p \in
G_\alpha \subseteq \Bbb P_\alpha$ and $G_\alpha$ is generic over
$\bold V$, then for
some $q \in G_\alpha$ and countable ordinal $\tau^*$ we have
$q \Vdash ``\underset\tilde {}\to \tau = \tau^*$.  By the definition of
$\aleph_1-Sp_e(W)$-iteration for some $\zeta \in [\gamma,\beta)$ we have
$\xi \in [\zeta,\beta) \Rightarrow [p \restriction \{\xi\}
\le^{\underset\tilde {}\to {\Bbb Q}_\xi}_{\text{pr}} 
q \restriction \{\xi\}$ or
$e = 4 \and p \restriction \{\xi\}$ not defined]. \nl
Define $p'$ by: $p' \restriction \underset\tilde {}\to \zeta = p 
\restriction \underset\tilde {}\to \zeta$, and for
$\xi \in [\zeta,\beta)$ we let $p' \restriction \{\xi\}$ be $q \restriction
\{\xi\}$ if: $p \restriction \{\xi\} 
\le^{\underset\tilde {}\to {\Bbb Q}_\xi}_{\text{pr}}
q \restriction \{\xi\}$ or $e=4 \and p \restriction \{\xi\}$ not defined.
Now $p'$ is as required.
\sn
So there is a $\Bbb P_{\underset\tilde {}\to \zeta}$-name of
${\underset\tilde {}\to p'}$ as appearing in the definition of
$\underset\tilde {}\to \zeta$ and it is, essentially, a member of
$\Bbb P_\beta$.
Now as we have finite apure support, the proof of
``$\underset\tilde {}\to \zeta[{\underset\tilde {}\to G_{{\Bbb P}_\alpha}}]$ is well
defined" gives $\Vdash_{{\Bbb P}_\alpha} ``\underset\tilde {}\to \zeta$ is not a
limit ordinal $> \gamma$".  Lastly $\Vdash_{{\Bbb P}_\alpha} ``\underset\tilde {}\to 
\zeta$ is not a successor ordinal $> \gamma"$ is proved by the property of
each ${\underset\tilde {}\to {\Bbb Q}_\xi}$.]
\ermn
Finishing the induction we define $q_\omega \restriction
{\underset\tilde {}\to \gamma_n} = q_n,
q_\omega \restriction [\dsize \bigcup_{n < \omega} 
{\underset\tilde {}\to \gamma_n},\beta)$ is defined as $\le_{vpr}$-upper 
bound of $\langle {\underset\tilde {}\to p_m} \restriction
[\dsize \bigcup_{n < \omega} \gamma_n,\beta):m < \omega \rangle$.
\medskip

More formally, let $\gamma^*,\beta$ and $G_{\gamma^*} \subseteq
\Bbb P_{\gamma^*}$ be such that: $G_{\gamma^*}$ is generic over $\bold V,\gamma^* =
\dsize \bigcup_{n < \omega} \gamma^*_n,\gamma^*_n = 
{\underset\tilde {}\to \gamma_n}[G_{\gamma^*}]$, let $p'_n =
{\underset\tilde {}\to p_n}[G_{\gamma^*} \cap \Bbb P_{\gamma^*_n}]$, let
$p'_n \restriction [\gamma^*,\alpha) = \{{\underset\tilde {}\to r^n_\zeta}:
\zeta < \zeta^*_n\}$ where ${\underset\tilde {}\to r^n_\zeta}$ is a simple
$[\gamma^*,\alpha)$-named atomic condition. \newline
Now we define ${\underset\tilde {}\to s^n_\zeta}$, a simple
$(\bar{\Bbb Q},W)$-named $[\gamma^*,\alpha)$-ordinal atomic condition as follows:
\medskip
\roster
\item "{$(a)$}"  $\zeta_{\underset\tilde {}\to s^n_\zeta} =
\zeta_{\underset\tilde {}\to r^n_\zeta}$
\sn
\item "{$(b)$}"  if $\zeta \in [\gamma^*,\alpha),G_{\gamma^*} \subseteq
G_\gamma \subseteq \Bbb P_\gamma,G_\gamma$ generic over $\bold V,
\zeta_{\underset\tilde {}\to r^n_\zeta}[G_\zeta] = \zeta$ then
${\underset\tilde {}\to s^n_\zeta}[G_\zeta]$ is the $<^{*{\bold V}[G_\zeta]}$-first
elements of ${\underset\tilde {}\to {\Bbb Q}_\zeta}[G_\zeta]$ which satisfies the
following:
{\roster
\itemitem{ $(*)(a)$ }  $({\underset\tilde {}\to p_n} \restriction \{\zeta\})$
\sn
\itemitem{ $(b)$ }  if 
$\langle \emptyset_{{\underset\tilde {}\to {\Bbb Q}_\zeta}[G_\zeta]} \rangle
\char 94 \langle ({\underset\tilde {}\to p_m} \restriction \{\xi\})
[G_\zeta]:m \in (n,\omega) \rangle$ has a $\le_{vpr}$-upper bound then
${\underset\tilde {}\to s^n_\zeta}[G_\zeta]$ is such upper bound.
\endroster}
\ermn
Now actually such $\le_{vpr}$-upper bound actually exists, and $q_\omega$
is as required.  \hfill$\square_{\scite{6.7}}$\margincite{6.7}
\enddemo
\bn
Now we can refine \scite{1.16A} to our iteration theorem. \nl
\ub{Saharon} revise.
\bigskip

\proclaim{\stag{6.8} Claim}  1) Suppose $F$ is a function, 
$\bold W \subseteq \omega_1$ stationary, $\Bbb I$ a class of quasi order
ideals, \ub{then} for every ordinal $\alpha$ there is 
UP$^6(\Bbb I,\bold W)$-iteration $\bar{\Bbb Q} = \langle \Bbb P_j,
{\underset\tilde {}\to {\Bbb Q}_i},{\underset\tilde {}\to \kappa_j}:
j \le \alpha^\dag,i < \alpha^\dag \rangle$, such that:
\mr
\item "{$(a)$}"  for every $i < \alpha^\dag$ we have 
${\underset\tilde {}\to {\Bbb Q}_i} = F(\bar Q \restriction i)$,
\sn
\item "{$(b)$}"  $\alpha^\dag \le \alpha$
\sn
\item "{$(c)$}"  for $\delta \le \alpha^\dag,{\underset\tilde {}\to \kappa
_\delta}$ is as in clause (d) of Definition \sciteu{6.3}
\sn
\item "{$(d)$}"  either $\alpha^\dag = \alpha$ or $F(\bar{\Bbb Q})$ is not a
pair 
$({\underset\tilde {}\to {\Bbb Q}},\underset\tilde {}\to \kappa)$ such that:
$\underset\tilde {}\to \kappa$ is a $\Bbb P_{\alpha^\dag} *
{\underset\tilde {}\to {\Bbb Q}}$-name of a cardinal from $\bold V$ and
$\Vdash_{{\Bbb P}_{\alpha^\dag}} ``{\underset\tilde {}\to {\Bbb Q}}$ 
satisfies
UP$^6_{{\underset\tilde {}\to \kappa_{\alpha^\dag}},
\underset\tilde {}\to \kappa}(\Bbb I,\bold W)$.
\ermn
2) Suppose $\beta < \alpha,G_\beta \subseteq \Bbb P_\beta$ is generic over
$\bold V$, then in $\bold V[G_\beta],\bar{\Bbb Q}/G_\beta = \langle
\Bbb P_i/G_\beta:
{\underset\tilde {}\to {\Bbb Q}_i},{\underset\tilde {}\to \kappa_i}:
\beta \le i < \alpha \rangle$ is an
${\text{\rm UP\/}}^6
(\Bbb I^{[{\underset\tilde {}\to \kappa_\beta}[G_\beta]]},
\bold W)$-iteration. \nl
3) If $\bar{\Bbb Q}$ is an ${\text{\rm UP\/}}^6(\Bbb I,\bold W)$-iteration, 
$p \in { \text{\rm Sp\/}}_e(W)-{\text{\rm Lim\/}}(\bar{\Bbb Q}),
\Bbb P'_i = \{q \in \Bbb P_i:q \ge p \restriction i\},{\underset\tilde
{}\to {\Bbb Q}'_i} =
\{p \in {\underset\tilde {}\to {\Bbb Q}_i}:
p \ge p \restriction \{i\}\}$, \ub{then}
$\bar{\Bbb Q} = \langle \Bbb P'_i,\Bbb Q'_i:i < \ell g(\bar{\Bbb Q}) 
\rangle$ is (essentially) an
${\text{\rm UP\/}}^6(\Bbb I,\bold W)$-iteration.
\endproclaim
\newpage

\head {\S7 No New Reals, Replacements for Completeness} \endhead  \resetall \sectno=7
\bn
Now we turn to ``No New Reals", there are versions corresponding to
\cite[Ch.V,\S1-\S3]{Sh:f} ($\bold W$-complete), \cite[Ch.V,\S5-\S7]{Sh:f}
($\dsize \bigwedge_{\alpha < \omega_1} \alpha$-proper 
$+ {\Cal D}$-completeness) and better \cite[Ch.VIII, \S4]{Sh:f} 
(making the previous preserved) and in different directions
\cite[Ch.XVIII,\S2]{Sh:f} and \cite{Sh:656}.
\bn
We deal here with the first (here we are interested in the cases
$\le_{pr} = \le$)

\definition{\stag{stc.0} Definition}  For $p \in \Bbb Q$ let
$\Bbb Q^{\text{pr}}_p = 
\{q \in \Bbb Q:p \le_{\text{pr}} q\}$.  A point which may confuse is that the
pure extension notion used in 
Definition \scite{stc.1}, is not necessarily the
one used seriously in the iteration.  This is the reason for the case
$e=5$ in \S1.  [Saharon check: main question: do we really need the purity
in the iteration for Nm$'$.  For Nm it is not needed (as in
\cite[Ch.XI]{Sh:f}.]
\enddefinition
\bigskip

\definition{\stag{stc.1} Definition}  1) UP$^4_{\text{com}}(\Bbb I,\bold W)$ is
satisfied by the forcing notion $\Bbb Q$, \ub{iff}: 
for any $\langle N_\eta:\eta \in
(T,\bold I) \rangle$ a strict $(\Bbb I,\bold W)$-suitable tree of models 
for $(\chi,x),x$ coding enough information, we have $(*)_1 \Rightarrow (*)_2$
where:
\medskip
\roster
\item "{$(*)_1$}"  for every $\eta,\nu \in T$, of the same length we have
$(N_\eta,\Bbb Q) \cong (N_\nu,\Bbb Q)$ and letting $h_{\eta,\nu}$ be the isomorphism 
from $N_\eta$ onto $N_\nu$ we have $h_{\eta,\nu}(x) = x$ and
$\ell < \ell g(\eta) \Rightarrow h_{\eta \restriction \ell,\nu 
\restriction \ell} \subseteq h_{\eta,\nu}$; (for $\eta,\nu \in \text{ lim}
(T)$ let $h_{\nu,\eta} = \dbcu_{\ell < \omega}
h_{\nu \restriction \ell,\eta \restriction \ell})$
\sn
\item "{$(*)_2$}"  if $\eta^* \in \text{ lim}(T),p \in N_{\langle \rangle}
\cap G$ and $G_{\eta^*}$ is a $\le_{pr}$-directed subset of 
$N_{\eta^*} \cap \Bbb Q^{pr}_p = \dbcu_{\ell < \omega} N_{\eta^* \restriction
\ell} \cap \Bbb Q^{pr}$, not disjoint to any dense subset of $\Bbb Q^{pr}_p \cap
\dsize \bigcup_{m < \omega} N_{\eta^* \restriction m}$ definable in
$(\dbcu_{m < \omega} N_{\eta^* \restriction m},N_{\eta^* \restriction m},
\bold I_{\eta^* \restriction m})_{m < \omega}$, \underbar{then}
there is $q \in \Bbb Q$ such that $p \le_{pr} q$ and 
$q \Vdash_{\Bbb Q}$ ``there is $\nu \in \text{ lim}(T)$
(in $\bold V^{\Bbb Q}$) such that $\dsize \bigcup_{\ell < \omega}
h_{\eta^* \restriction \ell,\nu \restriction \ell} (G \cap N_{\nu \restriction
\ell})$ is a subset of ${\underset\tilde {}\to G_{\Bbb Q}}$".
\ermn
2) UP$^4_{\text{stc}}(\Bbb I,\bold W)$ is satisfied by the forcing
notion $\Bbb Q$
\ub{iff} for any $\bar N = \langle N_\eta:\eta \in (T,\bold I))$ a strict
$(\Bbb I,\bold W)$-suitable tree of models for $(\chi,x),x$ coding enough
information we have $(*)_1 \Rightarrow (*)_2$ where
\mr
\item "{$(*)_1$}"  for every $\eta,\nu \in T$, of the same length we have
$(N_\eta,\Bbb Q) \cong (N_\nu,\Bbb Q)$ and letting $h_{\eta,\nu}$ be the isomorphism 
from $N_\eta$ onto $N_\nu$ we have $h_{\eta,\nu}(x) = x$ and
$\ell < \ell g(\eta) \Rightarrow h_{\eta \restriction \ell,\nu 
\restriction \ell} \subseteq h_{\eta,\nu}$; (for $\eta,\nu \in \text{ lim}
(T)$ let $h_{\nu,\eta} = \dbcu_{\ell < \omega}
h_{\nu \restriction \ell,\eta \restriction \ell})$
\sn
\item "{$(*)_2$}"  if $\eta^* \in \text{ lim}(T);\rho \in N_{<>} \cap G$ then
in the following game $\Game = \Game_{{\Bbb Q},p,\bar N,\eta}$; player complete has
a winning strategy.  The plays last $\omega$-moves, in the $n$-th move (we
can incorporate the last into the game) where
$\bold t_{\eta \restriction n} = T * x(\text{Suc}(\eta^* \restriction n) \in
\bold I^+_\eta)$ a condition $p_n \in P \cap N_{\eta^*}$ is chosen such that
$p \le p_0,[n > 0 \Rightarrow p_{n-1} \le p_n]$.  In the $n$-th move, the
anti-completeness player chooses $q_n \in \Bbb Q \cap N_{\eta^*}$ such that
$n = 0 \Rightarrow p = q_n$ and $n > 0 \Rightarrow p_{n-1} 
\le_{\text{pr}} q_n$ and the completeness player chooses $p_n \in \Bbb Q \cap N$
such that $q_n \le p_n \in \Bbb Q \cap N_{\eta^*}$.  In the end the completeness
player wins iff one of the following occurs:
{\roster
\itemitem{ $(\alpha)$ }  for some $\Bbb Q$-name $\underset\tilde {}\to \tau \in
N_{\eta^*}$ of a countable ordinal, no $p_n$ forces a value
\sn
\itemitem{ $(\beta)$ }  not $(\alpha)$ but there is $q$ such that:
$p \le_{\text{pr}} q \in \Bbb Q$ and $q \Vdash_{\Bbb Q}$ ``there is $\eta \in
\text{ lim}(T)$ such that $n < \omega \Rightarrow
h_{{\underset\tilde {}\to \eta} \restriction n,\eta^* \restriction n}(p_n)
\in {\underset\tilde {}\to G_{\Bbb Q}}"$.
\endroster}
\ermn
3) We define ``$\Bbb Q$ satisfies UP$^1_{\text{com},\underset\tilde {}\to \kappa}
(\Bbb I,W)$ is defined as in (1) but we replace the conclusion of $(*)_2$
by: there is $q \in \Bbb Q$ such that $p \le_{\text{pr}} q$ and

$$
\align
q \Vdash_{\Bbb Q} ``&\text{there is } (T',\bold I) \text{ satisfying }
(T,\bold I) \le^{\underset\tilde {}\to \kappa} (T',\bold I) \\
  &\text{(see \scite{2.4}(d)) such that for every} \\
  &\nu \in \text{ lim}(T') \text{ we have } 
h_{\nu,\eta^*}(G^*_\eta) \subseteq {\underset\tilde {}\to G_{\Bbb Q}}".
\endalign
$$
\mn
4) Similarly UP$^1_{\text{stc},\underset\tilde {}\to \kappa}(\Bbb I,W)$.
\nl
5) We say UP$^4_{\text{stc}}(\Bbb I,\bold W)$ if letting ${\frak p} = 
\langle({\underset\tilde {}\to \tau_n},I_n,\eta^*(n):I_{\eta^* \restriction
n},(\bold t_{\eta^* \restriction n}),n < \omega \rangle$ be such that
$\{{\underset\tilde {}\to \tau_n}:n < \omega\}$ list the $\Bbb Q$-names of
ordinals in $N_{\eta^*},\{I_n:n < \omega\}$ lists $\Bbb I \cap
N_{\eta^*}$, the winning strategy on each stage depends just on
$\theta,N_{<>} \cap \omega_1$ and in ${\frak p}$ continuously.
\enddefinition
\bigskip

\remark{\stag{stc.2} Remark}  1) This property relates to the UP$(\Bbb I,
\bold W)$ just as $E$-complete relates to $E$-proper (see 
\cite[Ch.V,\S1]{Sh:f}). \newline
2) Who satisfies this condition?  See section 8, so $\bold W$-complete
forcing notions, Nm$'(D)$,Nm$(D)(D$ is $\aleph_2$-complete) 
Nm$^{(\prime)}(T,\bold I)$ (when $\bold I$ is $\aleph_2$-complete), 
and shooting a club through a stationary subset of
some $\lambda = \text{ cf}(\lambda) > \aleph_1$ consisting of ordinals of
cofinality $\aleph_0$ (and generally those satisfying the $\Bbb I$-condition
from \cite[Ch.XI]{Sh:f}).
\endremark
\bigskip

\proclaim{\stag{stc.3} Claim}  If $\Bbb Q$ satisfies 
${\text{\rm UP\/}}^1_{\text{com}}
(\Bbb I,\bold W)$ or ${\text{\rm UP\/}}^1_{\text{stc}}(\Bbb I,\bold W)$ and 
$\Bbb I$ is $(2^{\aleph_0})^+$-complete, and $\Bbb Q$ has
$(\omega_1,2)$-pure decidability, \ub{then} 
forcing by $\Bbb Q$ add no new real.
\endproclaim
\bigskip

\demo{Proof}  Immediate.
\enddemo
\bigskip

\proclaim{\stag{stc.4} Claim}  Suppose:
\mr
\item "{$(a)$}"  $\Bbb Q$ is a 
forcing notion satisfying the ${\text{\rm UP\/}}^1_{\text{com}}
(\Bbb I,\bold W)$ and the local $\underset\tilde {}\to \kappa$-c.c. where
$\underset\tilde {}\to \kappa$ is a purely decidable $\Bbb Q$-name 
\sn
\item "{$(b)$}"  $\bar N = \langle N_\eta:\eta \in (T,\bold I) \rangle$
is a strict $(\Bbb I,\bold W)$-suitable tree of models (for $\chi$ and 
$x = \langle \Bbb Q,\underset\tilde {}\to \kappa,\Bbb I,\bold W) \rangle)$ 
satisfying $(*)_1$ of Definition \scite{stc.1}
\sn
\item "{$(c)$}"  the family 
${\underset\tilde {}\to {\Bbb I}'} =: \{I \in \Bbb I:I
\text{ is } \underset\tilde {}\to \kappa \text{-complete}\}$ is 
$(< \underset\tilde {}\to \kappa)$-closed.
\ermn
\ub{Then} $\Bbb Q$ satisfies ${\text{\rm UP\/}}^1_{\text{stc},
\underset\tilde {}\to \kappa}
(\Bbb I,\bold W)$.
\endproclaim
\bigskip

\demo{Proof}  Let $(T^*,\bold I^*)$ and $\bar N = \langle N_\eta:\eta \in
(T^*,\bold I^*) \rangle,\langle h_{\eta,\nu}:\eta,\nu \in T^* \cap
{}^n\text{Ord for some } n \rangle$ be as in Definition \scite{stc.1}. \nl
Let $\eta^* \in \text{ lim}(T^*),p \in \Bbb Q \cap N_{<>}$ be given and we choose
as our strategy for proving UP$^1_{\text{stc}}(\Bbb I,\bold W)$ the same
strategy that exists as UP$^1_{\text{stc}}(\Bbb I,\bold W)$ and let
$\langle p_n:n < \omega \rangle$ be a play as in Definition \scite{stc.1} in
which the completeness player uses his winning strategy. \nl
Let ${\Cal T} = \{T:(T^*,\bold I^*) \le^* (T,\bold I^* \restriction T)\}$.
As we can replace $p$ by $p'$ if $p \le_{\text{pr}} p' \in \Bbb Q \cap N_{<>}$,
\wilog \, $p$ forces a value to $\kappa$.  So for every $T \in {\Cal T}$
there are $q$ and $\underset\tilde {}\to \eta$ such that

$$
\align
p \le_{\text{pr}} q \in \Bbb Q \text{ and } q \Vdash_Q 
``&\underset\tilde {}\to \eta
\in \text{ lim } \underset\tilde {}\to T \text{ and } h = \bigcup
h_{\underset\tilde {}\to \eta} \restriction n,\eta^* \restriction n \\
  &\text{satisfies } n < \omega \Rightarrow h(\beta) \in
{\underset\tilde {}\to G_{\Bbb Q}} \\
  &\text{(hence } N_{{\underset\tilde {}\to \eta} \restriction n} \cap
\omega_1 = N_{<>} \cap \omega_1)".
\endalign
$$
\mn
Remember that we can replace $\eta^*$ by any $\eta^{**} \in \text{ lim}(T)$.
Let 

$$
\align
{\underset\tilde {}\to T^*}[G_{\Bbb Q}] = \{\nu \in T^*:&G_{\Bbb Q} \cap N_\nu \supseteq
\{h_{\eta,\eta^* \restriction \ell g(\eta)}(r), \\
  &r \in \Bbb Q \cap N_{\eta^* \restriction (\ell g(\nu))} \text{ and }
r \le p_n \text{ for some } n\}
\endalign
$$
\mn
clearly it is a subtree.  We continue as in the proof of \scite{5.2}.
\hfill$\square_{\scite{stc.4}}$\margincite{stc.4}
\enddemo
\bigskip

\proclaim{\stag{stc.5} Claim}  Suppose:
\mr
\item "{$(a)$}"  $\bar{\Bbb Q} = \langle \Bbb P_i,{\underset\tilde
{}\to {\Bbb Q}_i},
{\underset\tilde {}\to {\Bbb I}_i},{\underset\tilde {}\to \kappa_i}:
i < \alpha \rangle$ is a ${\text{\rm UP\/}}^{1,e}(\bold
W,W)$-iteration with \nl 
${\text{\rm Sp\/}}_e(W)$-limit $\Bbb P_\alpha$
\sn
\item "{$(b)$}"  $\Vdash_{{\Bbb P}_i} ``{\underset\tilde {}\to {\Bbb Q}_i}$ 
satisfies ${\text{\rm UP\/}}^1_{\text{stc}}
({\underset\tilde {}\to {\Bbb I}_i},\bold W)$ 
\sn
\item "{$(c)$}"  ${\underset\tilde {}\to \kappa^-}(\beta) = {
\text{\rm Min\/}} 
\{{\underset\tilde {}\to \kappa_\gamma}:\gamma < \gamma\}$, \newline
$\kappa^+(\beta) = { \text{\rm Sup\/}}
\{ \kappa:\text{for some } \gamma < \beta,
\nVdash_{{\Bbb P}_\gamma} 
{\underset\tilde {}\to \kappa_\gamma} \ne \kappa\}$, \nl
$\kappa^-(\beta) = { \text{\rm Min\/}}
\{ \kappa:\text{for some } \gamma < \beta
\text{ we have } \nVdash_{{\Bbb P}_\gamma} 
{\underset\tilde {}\to \kappa_\gamma} \ne \kappa\}$.
\ermn
\underbar{Then} \,1) for each $\beta \le \alpha,\Bbb P_\beta$ satisfies
${\text{\rm UP\/}}^1_{\text{com}}({\Bbb I}'_\beta,\bold W)$ for some 
$\kappa^-(\beta)$-complete set ${\Bbb I}'_\beta$ of (partial order)
ideals. \newline 
2) In fact, ${\Bbb I}_\beta$ is $\kappa^\beta$-complete where \newline
$\kappa^\beta = \min\{\kappa:\text{for some } \gamma < \beta$ we
have $\nVdash_{{\Bbb P}_{\gamma +1}} ``{\underset\tilde {}\to \kappa_\gamma} \ne
\kappa"\}$, and each $I \in \Bbb I_\beta$ \newline
has domain of cardinality
$\le (\underset{\gamma < \beta} \to \sup\{ \lambda < \kappa_{\delta +1}:\,
\nVdash_{{\Bbb P}_\gamma} ``\neg(\exists I \in 
{\underset\tilde {}\to {\Bbb I}_\gamma})(\lambda = |{\text{\rm Dom\/}}
(I)|)"\})$
\newline
and $|{\Bbb I}'_\beta| \le \dsize \sum_{\gamma < \beta} 
(\aleph_0 + |\Bbb P_\gamma| + \min\{(\lambda:\,\Vdash_{{\Bbb P}_\gamma} 
``|{\underset\tilde {}\to {\Bbb I}_\gamma}| \le \lambda)^{< \kappa}\}$.
\endproclaim
\bigskip

\remark{\stag{stc.5A} Remark}  We can 
use this for iteration as in \scite{5.8},
the version with clauses (b), (d) or (d)$'$, $W \cap \alpha = \emptyset$.
To prove $\Bbb P_\alpha$ does not add reals, it is enough to prove
that for each
$\beta < \alpha$, forcing with $\Bbb P_\beta$ does not add reals.  By
$\{p \in \Bbb P_\beta:(\forall \gamma < \beta) \Vdash_{{\Bbb P}_\gamma}
``\emptyset_{\underset\tilde {}\to {\Bbb Q}_\gamma} < p \restriction 
\{\gamma\}"\}$ is $\le_{\text{vpr}}$-dense.  This should be useful in \cite{GoSh:511}.
\endremark
\bn
SAHARON:  1) Use less $\underset\tilde {}\to \kappa$. \nl
2) What requirements will resurrect $\le_{\text{vpr}}$?
\bigskip

\demo{Proof}  Similar to the one of \scite{5.5}. \newline

For each $\gamma < \alpha$ let ${\Cal J}_\gamma =:
\{q \in \Bbb P_{\gamma +1}:q \text{ forces a value to } 
{\underset\tilde {}\to \kappa_\gamma}$, called 
$\kappa_{\gamma,q}$ and $q$ forces
${\underset\tilde {}\to {\Bbb I}_\gamma}$ to be equal to a $\Bbb P_\gamma$-name 
${\underset\tilde {}\to {\Bbb I}_{\gamma,q}}$ and
$q \restriction \gamma$ forces a value to
$|\Bbb I_\gamma|$ says $\mu_{\gamma,q}\}$; let ${\Cal J}'_\gamma \subseteq
{\Cal J}_\gamma$ be a maximal antichain.  Let $\mu_\gamma =
\underset{q \in {\Cal J}'_\gamma} \to \sup \mu_{\gamma,q}$. \newline
Let $q \Vdash_{{\Bbb P}_\gamma} ``{\underset\tilde {}\to {\Bbb I}_\gamma} = \{
{\underset\tilde {}\to I_{\gamma,\zeta}}:\zeta < \mu_{\gamma,q}\}"$ for
$q \in {\Cal J}_\gamma$ and let
${\Cal J}_{\gamma,\zeta} = \{q \in {\Cal J}_\gamma:\mu_{\gamma,q} > \zeta
\text{ and } q \Vdash ``\text{ Dom}({\underset\tilde {}\to I_{\gamma,\zeta}})
\text{ is } \lambda_{\gamma,q,\zeta}"$ if this is purely decidable$\}$ and let
${\underset\tilde {}\to I_{\gamma,\zeta}} \text{ be id}_{L_{\gamma,\zeta}}$, 
so $L_{\gamma,q,\zeta}$ is a
$\Bbb P_\gamma$-name of a ${\underset\tilde {}\to \kappa_{\gamma,q}}$-directed 
partial order on $\lambda_{\gamma,q,\zeta}$ (but $\Vdash_{{\Bbb P}_\gamma}$ ``if
$|{\underset\tilde {}\to {\Bbb I}_\gamma}| \le \zeta < \mu_\gamma$ then let
$L_{\gamma,\zeta}$ be trivial").

For $q \in {\Cal J}_\gamma$ let 
$L^*_{\gamma,q,\zeta}$ be ap$_{\kappa_{\gamma,q}}(
{\underset\tilde {}\to L_{\gamma,\zeta}})$ for the forcing notions \nl
$\Bbb P_\gamma^{[q]} = \{p \in P_\gamma:q \restriction \gamma
\le^{{\Bbb P}_\gamma}_{pr} p\}$ from Definition \scite{3.9}.  So by Claim \scite{3.10}
\medskip
\roster
\widestnumber\item{$(iii)$}
\item "{$(i)$}"   $L^*_{\gamma,q,\zeta}$ is $\kappa_{\gamma,q}$-directed
partial order on $[\lambda_{\gamma,q,\zeta}]^{< \kappa_{\gamma,q}}$
\sn
\item "{$(ii)$}"  $|L^*_{\gamma,q,\zeta}| \le (\lambda_{\gamma,q,\zeta})
^{< \kappa_{\gamma,q}}$
\sn
\item "{$(iii)$}"  
$q \restriction \gamma \Vdash_{{\Bbb P}_\gamma} 
``{\underset\tilde {}\to I_{\gamma,\zeta}} = 
\text{ id}_{\underset\tilde {}\to L_{\gamma,\zeta}} \le_{\text{RK}} \text{ id}
_{L^*_{\gamma,q,\zeta}}"$.
\endroster
\medskip

\noindent
Note: $\kappa^-(\beta) = \text{ Min}\{\kappa_{\gamma,q}:\gamma < \beta,q \in
{\Cal J}_\gamma\}$. 
\sn
Let $\Bbb I^*_\beta$ be the $(< \kappa^-(\beta))$-closure of $\{
\text{id}_{L^*_{\gamma,q,\zeta}}:\gamma < \beta,q \in {\Cal J}_\gamma,
\zeta < \mu_{\gamma,q}\}$ (see Definition \scite{3.13}(1)).
\smallskip

Let $\bar N = \langle N_\eta:\eta \in (T^*,\bold I) \rangle$ be a strict 
$(\Bbb I^*_\alpha,\bold W)$-suitable tree of models for $(\chi,x),x$ coding 
enough information (so $\bar{\Bbb Q},\Bbb I^*_\alpha,\bold W,W \in 
N_{\langle \rangle}$).  
For any $\gamma < \alpha$ and $G_\gamma \subseteq \Bbb P_\gamma$
generic over $\bold V,
T \subseteq T^*$ and $\bold I^{[\kappa^-(\gamma)]}$-tree and $\nu \in T$ and
$\eta^* \in \text{ lim}(T)$ and $p \in N_\nu[G_\gamma] \cap (\Bbb P_\alpha/
G_\gamma)$ then let $\bar N_{\nu,T}[G] = \langle N_{\nu \char 94 \rho}[G]:
\nu \char 94 \rho \in T \rangle$ then we can find a winning strategy
{\bf St} for the completeness player in the game $\Game = 
\Game_{\bar N_{\nu,T}}[G],p,\eta^*,\Bbb P_\alpha/G_\gamma$ of \scite{stc.1}(2).
Without loss of generality if $\eta^*_1,\eta^*_2 \in \text{ lim}(T)$ the
isomorphism from $N_{\eta^*_1}[G]$ onto $N_{\eta^*_2}[G]$ commutes with the
winning strategies; so the choice of $\eta^*$ is not important.  Of course,
we have a name ${\underset\tilde {}\to {\bold{St}}} = 
\bold{St}_{\underset\tilde {}\to \nu,\underset\tilde {}\to T,
{\underset\tilde {}\to \eta^*}}$.  
\sn
Now fix $\eta^* \in \text{ lim}(T^*)$ and we define a strategy {\bf St} for
the game.  For each simple $(\bar{\Bbb Q},W)$-name of an $[0,\alpha)$-ordinal
$\underset\tilde {}\to \gamma \in N_{\eta^*}$, let $\langle
(\tau^{\underset\tilde {}\to \gamma}_n,{\underset\tilde {}\to I^\gamma_n},
\eta^*(n),\bold I_{\eta^* \restriction n},\bold t_{\eta^* \restriction m}):
n < \omega \rangle$ be as in \scite{stc.1}(2) for
${\underset\tilde {}\to {\Bbb Q}_{\underset\tilde {}\to \gamma}}$.
\sn
We define {\bf St} such that if $\bar p = \langle p_n:n < \omega \rangle$ is
a play in which the completeness player uses his winning strategy then this
holds for $\langle p_n(\underset\tilde {}\to \gamma):n < \omega \rangle$
for each $\gamma$, i.e.,
\mr
\item "{$(*)$}"  if $\gamma < \alpha,G_\gamma \subseteq \Bbb P_\gamma$ is
generic over $\bold V$ and $p \restriction \gamma \in G$ and
$\underset\tilde {}\to \gamma[G_\gamma] = \gamma$ and $T \in \bold V[G_\gamma]$ is
a subtree of $T^*,\Bbb I^{[\kappa^-(\gamma)]}$-large and $\eta^{**} \in
\text{ lim}(T)^{V[G_\gamma]}$ and $p'_\gamma = (h^*_{\eta^{**},\eta^*}(p_n))
(\gamma)[G_\gamma] \in {\underset\tilde {}\to {\Bbb Q}_\gamma}[G_\gamma]$, \ub{then}
in $\langle p'_n:n < \omega \rangle$ the completeness player uses his
winning strategy from above.
\ermn
So fix such $\bar p^* = \langle p^*_n:n < \omega \rangle$, we would like to
find $q$ as in Definition \scite{stc.1}(2).
\sn
Let ${\Cal T}_{\bar N}$ be the set of quadruples 
$(\gamma,q,\underset\tilde {}\to \nu,\underset\tilde {}\to T)$ such that:
\medskip
\roster
\item "{$\bigotimes_1$}"  $\gamma \le \alpha,q \in \Bbb P_\gamma,
\Bbb P_\gamma \models p \restriction \gamma \le_{\text{pr}} q$ and \newline
$q \Vdash_{{\Bbb P}_\gamma} ``(\alpha) \quad 
\underset\tilde {}\to \nu \in \underset\tilde {}\to T
\subseteq T^*$, where $\underset\tilde {}\to \nu,\underset\tilde {}\to T$
are $\Bbb P_\gamma$-names \nl

$\qquad (\beta) \quad 
N_{\langle \rangle}[{\underset\tilde {}\to G_{{\Bbb P}_\gamma}}]
\cap \omega_1 = N_{\langle \rangle} \cap \omega_1$ 

$\qquad (\gamma) \quad \gamma \in \dbcu_{\ell < \omega} 
N_{\underset\tilde {}\to \nu}[G_{{\Bbb P}_\gamma}]$ \nl

$\qquad (\delta) \quad \langle N_{\eta,\ell}
[{\underset\tilde {}\to G_{{\Bbb P}_\gamma}}]:\eta \in \underset\tilde {}\to T$ is a
strictly 
$(\Bbb I^*_\gamma)^{[{\underset\tilde {}\to \kappa}(\gamma)]}$-suitable 
tree, \nl

$\qquad (\varepsilon) \quad$ for every
$\underset\tilde {}\to \eta \in \text{ lim}(\underset\tilde {}\to T)$ we have
\nl

$\qquad \qquad 
\{h_{{\underset\tilde {}\to \eta},\eta^*}(p^*_n) \restriction \gamma:n <
\omega\}$ is a subset of ${\underset\tilde {}\to G_{{\Bbb P}_\gamma}}"$.
\endroster
\medskip

\noindent
Now ${\Cal T}'_{\bar N}$ is defined similarly as the set of quadruples
$(\underset\tilde {}\to \gamma,q,\underset\tilde {}\to \nu,
\underset\tilde {}\to T)$ such that: as in $\otimes_1$ but we have
$\underset\tilde {}\to \gamma$ is a simple $(\bar{\Bbb Q},W)$-named ordinal, $q \in
\Bbb P_{\underset\tilde {}\to \gamma}$ and in clause $(\gamma)$ \,
$\Vdash \underset\tilde {}\to \gamma
\in 
N_{\underset\tilde {}\to \nu}[G_{{\Bbb P}_{\underset\tilde {}\to \gamma}}]$.  
(I.e., if $\zeta < \beta,G_{{\Bbb P}_\gamma}
\subseteq \Bbb P_\zeta$ is generic over $\bold V$ and $\zeta = {\underset\tilde {}\to 
\gamma_n}[G_{{\Bbb P}_\zeta}]$ then $r \in q_n \Rightarrow
{\underset\tilde {}\to \zeta_n}[G_\zeta] < \zeta$, i.e., is well defined
$< \zeta$ \ub{or} is forced ($\Vdash_{{\Bbb P}_\alpha/G_\zeta}$) to be not 
well defined, and $p \Vdash_{{\Bbb P}_{\underset\tilde {}\to \gamma}} 
``\eta \in \text{ lim}(T)"$).
\medskip

We consider the statements, for $\gamma \le \beta < \alpha$
\medskip
\roster
\item "{$\boxtimes_{\gamma,\beta}$}"  \ub{for any} 
$(\gamma,p,\underset\tilde {}\to \eta,\underset\tilde {}\to T) \in 
{\Cal T}_{\bar N}$ and $\underset\tilde {}\to \rho$ such that \newline
$p \Vdash_{{\Bbb P}_\gamma} ``\underset\tilde {}\to \eta \triangleleft 
\underset\tilde {}\to \rho \in \underset\tilde {}\to T"$
\newline
and ${\underset\tilde {}\to p'}$ a $\Bbb P_\gamma$-name such that \nl
$p \Vdash_{{\Bbb P}_\gamma} ``{\underset\tilde {}\to p'}
[{\underset\tilde {}\to G_{{\Bbb P}_\gamma}}] \in
N_{\underset\tilde {}\to \rho}
[{\underset\tilde {}\to G_{{\Bbb P}_\gamma}}] \cap 
P_\beta/{\underset\tilde {}\to G_{{\Bbb P}_\gamma}}$ and \nl
$({\underset\tilde {}\to p'}[G_{P_\gamma}]) \restriction \gamma \le p"$ 
\ub{there is} $(\beta,q,\underset\tilde {}\to \nu,{\underset\tilde {}\to T'}) 
\in {\Cal T}$ such that \nl
${\underset\tilde {}\to p'} \le q$ 
(i.e., $p \Vdash_{{\Bbb P}_\gamma} ``{\underset\tilde {}\to p'}
[{\underset\tilde {}\to G_{{\Bbb P}_\gamma}}] \le q"$) and 
$q \restriction \gamma = p$ and ${\Cal T}' \subseteq {\Cal T}"$.
\endroster
\medskip

We prove by induction on $\beta \le \alpha$ that $(\forall \gamma \le \beta)
\, \boxtimes_{\gamma,\beta}$ (or, for strong preservation), that 
($\forall$ non-limit
$\gamma \le \beta) \boxtimes_{\gamma,\beta})$, note that for $\gamma = \beta$ 
the statement is trivial hence we shall consider only $\gamma < \beta$.
\bigskip

\noindent
\underbar{Case 1}:  $\beta = 0$.

Trivial.
\bigskip

\noindent
\underbar{Case 2}:  $\beta$ a successor ordinal.

As trivially $\boxtimes_{\gamma_0,\gamma_1} \and \boxtimes_{\gamma_1,\gamma_2}
\Rightarrow \boxtimes_{\gamma_0,\gamma_2}$, clearly without loss of
generality \nl
$\beta = \gamma + 1$.

Let $G_{{\Bbb P}_\gamma}$ be such that $p \in G_{{\Bbb P}_\gamma} 
\subseteq \Bbb P_\gamma$ and
$G_{{\Bbb P}_\gamma}$ generic over $V$.

Let $T' = \{\nu:\rho \char 94 \nu \in \underset\tilde {}\to T
[G_{{\Bbb P}_\gamma}]\},\bar N' = \langle N'_\nu:\nu \in
(T',\bold I') \rangle$ where $N'_\nu = N_{\rho \char 94 \nu}[G_{{\Bbb P}_\gamma}]$,
\nl
$\bold I'_\nu = \bold I^*_{\rho \char 94 \nu}$. \newline
By \scite{stc.4} applied to $\bar N'$ we can find 
$\underset\tilde {}\to p',{\underset\tilde {}\to T''}$ as required.
\bigskip

\noindent
\underbar{Case 3}:   $\beta$ is a limit ordinal.

By \scite{stc.4} it suffices to prove $\bigotimes_2$ there are $q$ and
$\underset\tilde {}\to \eta$ such that: $\underset\tilde {}\to \eta$ is a
$\Bbb P_\beta$-name, $q \in P_\beta,q \restriction \gamma = p \restriction \gamma,
p \le q$ and $q \Vdash ``\underset\tilde {}\to \eta \in \text{ lim}
(\underset\tilde {}\to T)"$ and $\dbcu_{\ell < \omega}
N_{{\underset\tilde {}\to \eta} \restriction \ell}
[{\underset\tilde {}\to G_{{\Bbb P}_\beta}}] \cap \omega_1 = N_{\langle \rangle}
\cap \omega_1$ and $\{ h_{{\underset\tilde {}\to \eta},\eta^*}(r)
\restriction \beta:r \in G_{\eta^*}\} \subseteq 
{\underset\tilde {}\to G_{{\Bbb P}_\beta}}$.

We should choose by induction on $n < \omega,{\underset\tilde {}\to \gamma_n},
q_n,{\underset\tilde {}\to \rho_n},{\underset\tilde {}\to \eta_n},k_n$ 
such that:
\medskip
\roster
\item "{$(a)$}"  $({\underset\tilde {}\to \gamma_n},q_n,
{\underset\tilde {}\to \eta_n}) \in {\Cal T}'_{\bar N}$ \newline
(so ${\underset\tilde {}\to \gamma_n}$ is a $\bar{\Bbb Q}$-named ordinal)
\sn
\item "{$(b)$}"  ${\underset\tilde {}\to k_n}$ is a
$P_{\underset\tilde {}\to \gamma_n}$-name of a natural number
\sn
\item "{$(c)$}"  ${\underset\tilde {}\to \rho_n}$ is a $\Bbb P_{\gamma_n}$-name
(of a member of $T$)
\sn
\item "{$(d)$}"  $q_n \Vdash_{{\Bbb P}_{\underset\tilde {}\to \gamma_n}}
``{\underset\tilde {}\to \eta_n} \restriction {\underset\tilde {}\to k_n} =
{\underset\tilde {}\to \rho_n}"$
\sn
\item "{$(e)$}"  ${\underset\tilde {}\to \gamma_0} = \gamma$ and
$\Vdash_{\bar{\Bbb Q}}``{\underset\tilde {}\to \gamma_n} < 
{\underset\tilde {}\to \gamma_{n+1}} < \beta$ and 
${\underset\tilde {}\to \gamma_{n+1}}$ non-limit" \newline
i.e., if $\zeta < \beta$ and $G_{{\Bbb P}_\gamma}
\subseteq \Bbb P_\zeta$ is generic over $V$ and $\zeta = {\underset\tilde {}\to 
\gamma_n}[G_{{\Bbb P}_\zeta}]$ then \nl
$r \in q_n \Rightarrow {\underset\tilde {}\to \zeta_n}[G_\zeta] < \zeta$ 
(i.e. is well defined $< \zeta$ or is forced to be not well defined),
\sn
\item "{$(f)$}"  $q_{n+1} \restriction {\underset\tilde {}\to \gamma_n} =
q_n$
\sn
\item "{$(g)$}"  $q_{n+1} 
\Vdash_{{\Bbb P}_{\underset\tilde {}\to \gamma_{n+1}}} 
``{\underset\tilde {}\to \rho_n} \triangleleft
{\underset\tilde {}\to \rho_{n+1}}$, so 
${\underset\tilde {}\to k_n} < {\underset\tilde {}\to k_{n+1}}"$.
\ermn
Finishing the induction we let
$\underset\tilde {}\to \eta = \dsize \bigcup_{n < \omega}
{\underset\tilde {}\to \rho_n}$ and we define $q_\omega \restriction
{\underset\tilde {}\to \gamma_n} = q_n$ and \nl
$q_\omega \in \Bbb P
\dbcu_{n < \omega} {\underset\tilde {}\to \gamma_n}$. \nl
We shall check that $\bigotimes_2$ holds which is straight.
\enddemo
\bigskip

\demo{\stag{7.6} Discussion}  1) As in \S6 (not \S5)? \nl
2) The other $NNR$.

Like V and like XVIII. \nl
A. \,\, Like XVIII,\S2 - seem straight but check. \nl
B. \,\, Like V,\S6 - think. \nl
3) Explain the specific choice for \scite{stc.1}. \nl
4) Think Ch.VI,\S1, $\le = \le_{pr}$. \nl
\S3 not necessarily.
\enddemo 
\newpage

\head {\S8 Examples} \endhead  \resetall \sectno=8
\bigskip

\noindent
Namba [Nm] defines Nm$(J^{bd}_\lambda)$ (and also with $\omega$ ideals)
as examples of forcing notion preserving $\aleph_1$ but changing the 
cofinality of some $\lambda = \text{ cf}(\lambda)$ to $\aleph_0$. \newline
More \cite{RuSh:117}, \cite[X,XI,XV,XIV,\S5]{Sh:f}.
\bigskip

\definition{\stag{8.1} Definition}  1) For an ideal $I$ on a cardinal
$\lambda$, let the forcing notion Nm$(I)$ be

$$
\align
\text{Nm}(I) = \biggl\{ T:&T \subseteq {}^{\omega >}\lambda 
\text{ is non-empty, closed under initial segments and} \\
  &(\forall \eta \in T)(\exists \nu)[\eta \triangleleft \nu \in T \and
   (\exists^{I^+} \alpha < \lambda)(\eta \char 94 \langle \alpha \rangle \in
    T)] \biggr\}
\endalign
$$
\mn
where $(\exists^{I^+} \alpha < \lambda)\text{Pr}(\alpha) \text{ means }
\{\alpha < \lambda:\text{Pr}(\alpha)\} \in I^+ 
\text{ and } I^+ = \{A \subseteq \lambda:A \notin I\}$ ordered by 
inverse inclusion and let $<_{\text{pr}} = \le$ 
and $\le_{\text{vpr}}$ be the equality
$p \in \text{ Nm}(I)$ 
is normal if $\forall \eta \in p \Rightarrow |\text{Suc}
(\eta)| = 1 \vee \text{ Suc}_T(\eta) \ne \emptyset \text{ mod } I$. \nl
2) For an ideal $I$ and a cardinal $\lambda$, let the forcing notions
Nm$'(I)$ be

$$
\align
\text{Nm}'(I) = \biggl\{ T:\,&T \subseteq {}^{\omega >}\lambda 
\text{ is non-empty, closed under initial segments and for some} \\
  &\,n = n(T) < \omega \text{ we have}: \\
  &\,(i) \quad \ell \le n \Rightarrow |T \cap {}^\ell \lambda| = 1 \\
  &\,(ii) \quad \eta \in T \and \ell g(\eta) \ge n \Rightarrow
(\exists^{I^+} \alpha < \lambda)[\eta \char 94 \langle \alpha \rangle \in
T] \biggr\}
\endalign
$$
\mn
we call the $\eta \in T \cap {}^{n(T)}\lambda$ 
the trunk of $T$ and denote it by tr$(T)$) 
\sn
\underbar{ordered} by inverse inclusion 
and let $\le_{\text{pr}} = \le^*$ (see \S2)
and $\le_{\text{vpr}}$ be the equality.
\sn
3) Writing a filter $D$ means the dual ideal.
\enddefinition
\bigskip

\proclaim{\stag{8.2} Claim}  Let $I$ be a $\kappa$-complete ideal on
$\lambda,\lambda \ge \kappa \ge \aleph_2,I \in \Bbb I,\Bbb I$ is 
(restriction closed and) $\kappa$-complete. \newline
1) \text{\rm Nm}$(I)$ and 
\text{\rm Nm}$'(I)$ satisfies ${\text{\rm UP\/}}^1(\Bbb I)$ 
and ${\text{\rm UP\/}}^4_{\lambda^+}(\Bbb I)$
and ${\text{\rm UP\/}}^6_{\lambda^+}(\Bbb I)$ 
so does not collapse $\aleph_1$. \newline
) If $I$ is uniform, then $\Vdash_{\text{Nm}(I)} ``{\text{\rm cf\/}}
(\lambda) = 
\aleph_0"$ and 
$\Vdash_{\text{Nm}'(I)} ``{\text{\rm cf\/}}(\lambda) = \aleph_0"$, 
in fact if $[A \in I^+ \Rightarrow I \restriction A$ is 
$\lambda'$-decomposable] and $\lambda'$ is regular, 
\ub{then} the same holds for $\lambda'$. \newline
3) $|{\text{\rm Nm\/}}(I)|,|{\text{\rm Nm\/}}'(I)| \le 2^\lambda$. \newline
4) If in addition $2^{\aleph_0} < \kappa$, \underbar{then} forcing with
${\text{\rm Nm\/}}(I)$ does not add reals, moreover it satisfies the 
condition from \scite{stc.1}, ${\text{\rm UP\/}}^4_{\text{com}}(\Bbb I)$. \nl
5) If in addition $2^{\aleph_0} < \kappa$ \ub{then} forcing with
${\text{\rm Nm\/}}(I),{\text{\rm Nm\/}}'(I)$ 
does not add reals; moreover, they satisfy the condition
${\text{\rm UP\/}}^{4,+}_{\text{stc}}(\Bbb I)$.
\endproclaim
\bigskip

\demo{Proof}  1)  We will use the following fact about $\Bbb Q = \text{ Nm}(I)$ 
and $\Bbb Q = \text{ Nm}'(I)$:
\medskip
\roster
\item "{$(*)$}"  If $p \in \Bbb Q,\underset\tilde {}\to \alpha$ is a
$\Bbb Q$-name of 
an ordinal, \newline
then there is $q,p \le_{pr} q$ such that the set \newline
$\{\eta \in q:\text{for some } \beta \text{ we have } q^{[\eta]} \Vdash
``\underset\tilde {}\to \alpha = \beta"\}$ contains a front.
\endroster
\medskip

This fact follows easily from \scite{2.14} (let $H:p \rightarrow \{0,1\}$
(i.e., Dom$(H) = T^p$) be defined by $H(\eta) = 1$ iff $(\exists q)[p
^{[\eta]} \le_{pr} q \and q$ decides $\underset\tilde {}\to \alpha]$, define
$H(\eta) =$ \nl
$\text{lim}_{n < \omega}(H(\eta \restriction n))$ for
$\eta \in \text{ lim}(p)$, and find $q$ such that $H$ is constant on
lim$(q)$).  Let $Y = \{\eta \in T^q:H(\eta) = 1 \and (\forall \nu)(\nu
\triangleleft \eta \rightarrow H(\eta) = 0]\}$, so $Y$ is a front of $q$.
For $\eta \in Y$ let $q_\eta$ be such that $p^{[\eta]} \le_{\text{pr}} 
q_\eta$ and $q_\eta$ forces a value 
to $\underset\tilde {}\to \alpha$ let $r$ be such that
$T^r = \dbcu_{\eta \in Y} T^{q_\eta}$.  So clearly $r$ is as required, $Y$
such a front.

Now let $\langle N_\eta:\eta \in (T,\bold I) \rangle$ be a strictly 
$\Bbb I$-suitable tree of models for $\chi,x$ satisfying 
$\{p,I,\Bbb I\} \in N_{\langle \rangle}$ where $p \in \Bbb Q \cap N_{<>}$ is a
condition.  We can now find a condition $q,p \le_{pr} q$, a family $\langle
p_\eta:\eta \in p \rangle$ of conditions and a function $f:q \rightarrow T$
satisfying the following:
\medskip
\roster
\item  If $\eta \triangleleft \nu$ in $q$, then $f(\eta) \triangleleft
f(\nu)$.
\sn
\item  For all $\eta$ in $q$, Suc$_T(f(\eta)) \ne 0 \text{ mod } I$ and
$\bold I_\eta = I$.
\sn
\item  For all $\eta$ in $q$, Suc$_q(\eta) \subseteq \text{ Suc}_T(f(\eta))$.
\sn
\item  For all $\eta$ in $q,p_\eta \in N_{f(\eta)},\text{tr}(p_\eta) =
\eta,p^{[\eta]} \le_{pr} p_\eta$.
\sn
\item  For all $\eta$ in $q,p_\eta \le_{pr} q^{[\eta]}$.
\sn
\item  For all $\eta$ in $q$, all names $\underset\tilde {}\to \alpha$ in
$N_{f(\eta)}$, the set \newline
$\{\nu \in q:p_\nu \text{ decides } \alpha\}$ contains a front of
$p^{[\eta]}$.
\endroster
\medskip

We can do this as follows: by induction on $\eta \in p$ we define
$f(\eta),p_\eta$ and Suc$_q(\eta)$.  We can find $f(\eta)$ satisfying (2)
+ (3) because $T$ is $\Bbb I$-suitable and $I \in \Bbb I$ and $\Bbb I$ is
restriction closed.  We choose $p_\eta$ using a bookkeeping
argument to take care of a case of (6), using $(*)$.  Then we 
choose Suc$_q(\eta)$ such that (3) are satisfied.  

Lastly, let $q = \{\nu$: for some $\eta,p_\eta$ is well defined and $\nu
\trianglelefteq \eta\}$.  Clearly $p \le_{\text{pr}} q \in \Bbb Q$.

Now let $G$ be $\Bbb Q$-generic, $q \in G$.  Now $G$ defines a generic branch
$\eta$ through $q$.  This induces a branch $\nu$ through $T:\nu =
\dsize \bigcup_{n < \omega} f(\eta \restriction n)$.  Let
$\underset\tilde {}\to \alpha \in N_{\nu \restriction k}$, then there is
$\ell$ such that $p_{\eta \restriction \ell} \Vdash ``
\underset\tilde {}\to \alpha = \beta$ and $\beta \in N_{f(\eta \restriction
\ell)} \subseteq N_\nu"$. \nl
2) It is enough to prove the second version for any condition $p$, let
$\bold I^p_\eta$ be $I$ ``mapped" to Suc$_p(\eta)$. \nl
For any condition $p$, for each $\eta \in p$ such that Suc$_p(\eta) \ne
\emptyset$ mod $I$ let $h_\eta:\text{Suc}_p(\eta) \rightarrow \lambda'$ be
such that $(\forall \alpha < \lambda')[\{\nu \in \text{ Suc}_p(\eta):h_\eta
(\nu) < \alpha\} \in \bold I^p_\eta]$.  Now letting 
$\underset\tilde {}\to \eta \in {}^\omega \text{Ord},
\underset\tilde {}\to \eta \restriction \ell =
{}^n \text{Ord} \cap r$ for any $r \in {\underset\tilde {}\to G_{\Bbb Q}}$ large
enough, we let $\underset\tilde {}\to A = \{
h_{{\underset\tilde {}\to \eta} \restriction \ell}(\underset\tilde {}\to \eta
\restriction (\ell + 1)):\ell < \omega\}$. \nl
So easily $\Vdash_{\Bbb Q} ``\underset\tilde {}\to A \subseteq \lambda'$ is
unbounded. \nl
3) Trivial. \nl
4) Without loss of generality $\Bbb I$ is $\kappa$-complete (as we can
decrease it).  So by \scite{stc.3} it suffices to prove
UP$^1_{\text{com}}(\Bbb I,\bold W)$.  So assume $\langle N_\eta:\eta \in
(T,\bold I) \rangle,h_{\eta,\nu}$ (for $\eta,\nu \in T \cup \text{ lim}(T),
\ell g(\eta) = \ell g(\nu))$ are as in Definition \scite{stc.1}, $(*)_1$ and
$\eta^* \in \text{ lim}(T),G_{\eta^*}$ as in the assumption of $(*)_2$ there.
Now choose inductively on $n < \omega,p_n$ and $k_n$ such that:
$p_0 = p,k_0 = 0,p_n \in G_{\eta^*},p_n \in N_{\eta^* \restriction
k_{n+1}}$, and $p_n < p_{n+1},k_n < k_{n+1},\eta_n$ is the trunk of
$p_n,\eta_n \triangleleft \eta_{n+1}$, Suc$_{p_n}(\eta_n) \ne \emptyset
\text{ mod } I$ (as in proof of part (1)) and Suc$_T(\eta \restriction
k_{n+1}) = \text{ Suc}_{p_n}(\eta_n)$ and 
$p^{[\eta_n \char 94 \langle \eta^*(k_{n+1})]}_n \le_{pr} p_{n+1}$ and if
$\underset\tilde {}\to \tau \in N_{\eta^*}$ is a Nm$(I)$-name of a countable
ordinal then for some $n,p_n$ decides its value. \nl
5) The winning strategy of the completeness player is, given $q_n$, let
$\nu = tr(T)$ and let $n$ be minimal such that $q_n \in N_{\eta^* \restriction
n}$ and $I \restriction \text{ Suc}_{q_n}(\nu) = I_{\eta^* \restriction n}$
and let $p_n = (q_n)^{[\nu \restriction \eta^*(n)]}$. 
\hfill$\square_{\scite{8.2}}$\margincite{8.2}
\enddemo
\bigskip

\definition{\stag{8.3} Definition}  1) We can consider an $\Bbb I$-suitable 
tree of models \newline
$\bar N = \langle N_\eta:\eta \in (T^*,{\Bbb I}^*) \rangle$, and let

$$
\align
a) \qquad \Bbb Q_{\bar N} = \biggl\{ T \subseteq T^*:\,&T \text{ non-empty, closed
under initial segments} \\
  &\,\text{such that } \langle N_\eta:\eta \in (T,\bold I \restriction T)
\rangle \text{ is an} \\
  &\,{\Bbb I} \text{-suitable tree of models} \biggr\}
\endalign
$$

\noindent
ordered by inverse inclusion. \newline
2) We can consider for any tagged tree $(T^*,\bold I^*)$

$$
\align
\Bbb Q^0_{(T^*,\bold I^*)} = \biggl\{T \subseteq T^*:\,&T \text{ non-empty, closed
under initial segments} \\
  &\,\text{such that for some } n = n(T), \\
  &\,(i) \qquad \ell \le n \Rightarrow |T \cap {}^n\text{Ord}| = 1 \\
  &\,(ii) \qquad \text{if } \eta \in T \and \ell g(\eta) 
\ge n \and \text{ Suc}_{T^*}
(\eta) \ne \emptyset \text{ mod } \bold I_\eta \\
  &\qquad \qquad \,\,\text{ then Suc}_T(\eta) \ne 
\emptyset \text{ mod } \bold I_\eta \biggr\}
\endalign
$$
\mn
is ordered by inverse inclusion.

$$
\align
\Bbb Q^1_{(T,{\bold I}^*)} = \biggl\{ (T,\bold I):&(T^*,\bold I^*) \le
(T,\bold I), \text{ and for every } \eta \in \text{ lim}(T) \\
  &\text{ we have } (\forall k)(\exists^\infty n)[\eta \restriction n
\text{ is a splitting point of } (T,\bold I) \\
 &\text{ and } \bold I_{\eta \restriction k} \le_{\text{RK}}
\bold I_{\eta \restriction n} \biggr\}
\endalign
$$

\noindent
ordered by inverse inclusion.  [Saharon" [$\Bbb Q^0 \ne \Bbb Q^2$?]

$$
\Bbb Q^2_{(T,{\bold I}^*)} = \biggl\{ (T,\bold I):(T^*,\bold I^*)^{[\eta]} \le^*
(T,\bold I), \text{ for some } \eta \in T^* \biggr\}
$$

\noindent
ordered by inverse inclusion. 
\enddefinition
\bigskip

\proclaim{\stag{8.3A} Claim}  For the forcing notions defined in Definition
\scite{8.3} for $\Bbb I$ being $\aleph_2$-complete, of course, we have: if
$P \in \{\Bbb Q_{\bar N},\Bbb Q^0_{(T^*,\bold I^*)},\Bbb Q^1_{(T,\bold I^*)}\}$, \ub{then}
\medskip
\roster
\item "{$(a)$}"  $P$ satisfies ${\text{\rm UP\/}}^1(\Bbb I)$
\sn
\item "{$(b)$}"  if $I \in \Bbb I \Rightarrow |{\text{\rm Dom\/}}
(I)| < \lambda =
{ \text{\rm cf\/}}(\lambda)$, \ub{then} $|\Bbb P| \le 2^{< \lambda}$ and even
$\le 2^\mu$ for some $\mu < \lambda$
\sn
\item "{$(c)$}"  if for $\lambda$ regular \newline
$(\forall I \in \Bbb I)(\forall A \in (I)^+[I \restriction A$ is not
$\lambda$-indecomposable] \newline
\underbar{then} $\Vdash_{\Bbb P} ``{\text{\rm cf\/}}
(\lambda) = \aleph_0"$
\sn
\item "{$(d)$}"  if $\Bbb P = \Bbb Q^0_{(T^*,\bold I^*)}$ then 
$(\forall \eta \in \text{ lim}(T^*)) \, 
\exists^\infty n \dsize \bigwedge_{m \ge n} \forall A \in
(\bold I^*_{\eta \restriction n})^+[\bold I_\eta \restriction A$ is not
$\lambda$-indecomposable] 
\underbar{then} $\Vdash_{\Bbb P} ``{\text{\rm cf\/}}
(\lambda) = \aleph_0"$
\sn
\item "{$(e)$}"  if $\Bbb I$ is $(2^{\aleph_0})^+$-complete then forcing with
$\Bbb P$ add no new reals, moreover it satisfies 
${\text{\rm UP\/}}^{4,+}_{\text{stc}}(\Bbb I)$
and if $p \in \{{\Bbb Q}_{\bar N},\Bbb Q^0_{(T^*,\bold I^*)}\}$ \ub{then} it 
satisfies ${\text{\rm UP\/}}^4_{\text{com}}(\Bbb I)$.
\endroster
\endproclaim
\bigskip

\demo{Proof}  Left to the reader. \hfill$\square_{\scite{8.3A}}$\margincite{8.3A}
\enddemo
\bigskip

\definition{\stag{8.4} Definition}  Let $\lambda = \text{ cf}(\lambda) >
\aleph_1,S \subseteq \{\delta < \lambda:\text{cf}(\delta) = \aleph_0\}$ be
stationary and

$$
\align
\text{club}_S(S) = \biggl\{ h:\,&\text{for some non-limit } \alpha 
< \omega_1,\\
  &\,h \text{ is an increasing function from } \alpha \text{ to } S \biggl\}
\endalign
$$

\noindent
ordered by inverse inclusion, $\le_{\text{pr}} = \le,\le_{vpr}$ is
equality.
\enddefinition
\bigskip

\proclaim{\stag{8.5} Claim}  For $\lambda,S$ as in Definition \scite{8.4}
we have (for any $I,I$ is an $\aleph_2$-complete ideal on $\lambda$ extending
$J^{bd}_\lambda$). \newline
1) Club$(S)$ satisfies ${\text{\rm UP\/}}^1(\{I\})$ 
of cardinality $\le \lambda^{\aleph_0}$.
\newline
2) If $I$ is $(2^{\aleph_0})^+$-complete and $I \in \Bbb I$, \ub{then} 
Club$(S)$ satisfies ${\text{\rm UP\/}}^4_{\text{com}}(\Bbb I)$, hence 
${\text{\rm UP\/}}^{4,+}_{\text{stc}}(\Bbb I)$.
\endproclaim
\bigskip

\demo{Proof}  Left to the reader \ub{or} follows from \cite[XI,4.6]{Sh:f} by
\scite{8.7} below. \hfill$\square_{\scite{8.5}}$\margincite{8.5}
\enddemo
\bigskip

\proclaim{\stag{8.5A} Lemma}  Let $\bar W = \langle W_i:i < \omega_1 \rangle$ 
be a sequence of stationary subsets of $\{\alpha < \lambda:
{\text{\rm cf\/}}(\alpha)
= \omega\}$ where $\lambda = { \text{\rm cf\/}}(\lambda) > 
\aleph_0$ and let the
forcing notion $\Bbb P[\bar W]$ be defined by

$$
\align
\Bbb P[\bar W] =: \biggl\{ f:\,&f \text{ is an increasing and continuous function
from} \\
  &\,\alpha + 1 \text{ into } W_0 \text{ for some } \alpha < \omega_1, \\
  &\,\text{ such that for every } i \le \alpha \text{ we have } f(i) \in
W_i \biggr\}
\endalign
$$
\mn
(ordered by inclusion).  If $I \supseteq J^{bd}_\lambda$ be
$\aleph_2$-complete, \ub{then} $\Bbb P[\bar W]$ satisfies 
${\text{\rm UP\/}}^{4,+}(\Bbb I)$ for
any $\Bbb I$ such that $I \in \Bbb I$ and if $\Bbb I$ is
$(2^{\aleph_0})^+$-complete it also satisfies 
${\text{\rm UP\/}}^4_{\text{com}}(\Bbb I)$ 
hence ${\text{\rm UP\/}}^{4,+}_{\text{stc}}(\Bbb I)$.
\endproclaim
\bigskip

\demo{Proof}  Left to the reader or follows from \cite[Ch.XI,4.6A]{Sh:f} by
\scite{8.7} above.  \hfill$\square_{\scite{8.5A}}$\margincite{8.5A}
\enddemo
\bn
Concerning $\bold W$-completeness (see \cite[Ch.V]{Sh:f}):
\proclaim{\stag{8.6} Claim}  Assume $\bold W \subseteq \omega_1$ is 
stationary and $\Bbb Q$ is $\bold W$-complete forcing notion (i.e., 
if $\chi$ is large enough, $\Bbb Q \in N \prec ({\Cal H}(\chi),\in,
<^*_\chi),N$ countable, 
$p_n \in \Bbb Q \cap N$ 
is $\le_{\Bbb Q}$-increasing and $(\forall {\Cal I} \in N)
({\Cal I} \subseteq \Bbb Q$ is dense $\rightarrow \dsize \bigvee_n p_n \in
{\Cal I})$ \ub{then} $\{p_n:n < \omega\}$ has an upper bound in $\Bbb Q$).

\underbar{Then} $\Bbb Q$ satisfies 
${\text{\rm UP\/}}^{4,+}_{\text{com}}(\Bbb I,\bold W)$ (i.e.,
for any $\Bbb I$).
\endproclaim  
\bigskip

\demo{Proof}  Trivial (see Definition \sciteu{7.1}), for the $\omega$-branch
$\eta^*$ of $T$ there is $q,p \le q,q$ is an upper bound of $\Bbb Q^*$ hence is
$(\dsize \bigcup_{\ell < \omega} N_{\eta^* \restriction \ell},\Bbb Q)$-generic.
\hfill$\square_{\scite{8.6}}$\margincite{8.6}
\enddemo
\bn
Comparing to \cite{Sh:f} we have
\proclaim{\stag{8.7} Claim}  1) If $\Bbb Q$ satisfies the $\Bbb I$-condition
of \cite[Ch.XI,Def.2.6,2.7]{Sh:f} $\Bbb I$ is $(2^{\aleph_0})^+$-complete,
\ub{then} $\Bbb Q$ satisfies ${\text{\rm UP\/}}^4_{\text{com}}(\Bbb I)$. \nl
2) If $\Bbb Q$ satisfies condition 
${\text{\rm UP\/}}(\Bbb I,\bold W)$ of \cite[Ch.XV, Definition 2.7A]{Sh:f}, 
\ub{then} it satisfies ${\text{\rm UP\/}}^4
(\Bbb I,\bold W)$ here. \nl
3) If $\Bbb Q$ is a proper or just semi-proper forcing notion,
\ub{then} $\Bbb Q$
satisfies ${\text{\rm UP\/}}^6(\Bbb I)$ and all the UP$^\ell(\Bbb I)$.
\endproclaim
\bigskip

\demo{Proof}  Part (2) holds by \cite[XV,2.11]{Sh:f}. \nl
Part (1) follows by part (2) and \cite[Ch.XV,2.11]{Sh:f}. \nl
Part (3) is immediate by the definition (can use any fix branch).
\hfill$\square_{\scite{8.7}}$\margincite{8.7}
\enddemo
\newpage

\head {\S9 Reflection in $[\omega_2]^{\aleph_0}$} \endhead  \resetall \sectno=9
\bn
As an exercise 
we answer a question (and variants) of Jech \footnote{Done 10/97}.
\proclaim{\stag{11.1} Theorem}  Assume
\mr
\item "{$(A)$}"  $\kappa$ is large enough (supercompact or just Woodin).
\ermn
\ub{Then} \nl
1) there is a $\kappa$-c.c. semi-proper forcing notion $\Bbb P$ of cardinality
$\kappa$ such that in $\bold V^{\Bbb P}$
\mr
\item "{$(\alpha)$}"  $\aleph_2 = \aleph^V_1,\aleph_2 = \kappa$ and cardinal
$\ge \kappa$ are the same as in $\bold V$ and $2^{\aleph_0} = \aleph_2$
\sn
\item "{$(\beta)$}"  every stationary subset of $[\omega_2]^{\le \aleph_0}$
reflect (in a set of cardinality $\aleph_1$)
\sn
\item "{$(\gamma)$}"  $D_{\omega_1}$ is $\aleph_2$-saturated
\sn
\item "{$(\delta)$}"  there is a projectively stationary $S \subseteq
[\omega_2]^{\le \aleph_0}$, see below 
such that there is no sequence $\langle a_i:i < \omega_1 \rangle$ increasing
continuous, $a_i \in S,a_i \ne a_{i+1}$.
\ermn
2) Assume in addition
\mr
\item "{$\boxtimes$}"   $\{\lambda:\lambda \text{ measurable}\}$ is
not in the weakly compact ideal of $\kappa$. 
\ermn
We can add to (1) the statement (on $\bold V^{\Bbb P}$)
\mr
\item "{$(*)$}"  every stationary $S \subseteq S^2_0 = \{\delta < \omega_2:
\text{cf}(\delta) = \aleph_2\}$ contains a closed copy of $\omega_1$.
\ermn
3) We may strengthen clause $(\delta)$ of (1) to ``$S$ is 
$S^2_1$-projectively stationary".
\endproclaim 
\bigskip

\definition{\stag{11.1A} Definition}  1) We call 
$S \subseteq [\omega_2]^{\aleph_0}$ projectively stationary if: \nl
for every club $E$ of $[\omega_2]^{\aleph_0}$ and stationary
co-stationary $W \subseteq \omega_1$ we can find a sequence
$\langle a_i:i < \omega_1 \rangle$ increasing continuous, 
$a_i \in [\omega_2]^{\aleph_0},a_i \in E$
and $\{i \in W:a_i \notin S\}$ is not a stationary subset of $\omega_1$. \nl
2) We say $S \subseteq [\omega_2]^{\aleph_0}$ is $S^2_1$-projectively 
stationary for $W \subseteq \omega_1$ stationary co-stationary, for 
stationarily many $\delta \in S^2_1$ if we let 
$\delta = \dbcu_{i < \omega_1} a^\delta_i,a^\delta_i$
countable increasing continuous we have 
$\{i \in W:a^\delta_i \notin S\}$ non-stationary.
\enddefinition
\bigskip

\demo{Proof}  Like \cite[Ch.XVI,2.4]{Sh:f}'s proof.
\mn
1), 2)  We 
define a RCS iteration $\langle \Bbb P_i,{\underset\tilde {}\to {\Bbb Q}_j}:i \le
\kappa,j < \kappa \rangle$ such that:

$$
|\Bbb Q_j| \le (2^{\aleph_2})^{\bold V^{{\Bbb P}_j}}
$$
\mn
in $\bold V^{{\Bbb P}_j},\Bbb Q_j$ is the disjoint union of the following (so the choice by
which of them we force is generic), \ub{but} if $j$ is non-limit only (a)
is allowed
\mr
\item "{$(a)$}"  $\Bbb Q^0_j = \text{ Levy}(\aleph_1,2^{\aleph_1}) *$ Cohen
\sn
\item "{$(b)$}"  $\Bbb Q^1_j =$ sealing all semi-proper maximal antichanges of
$D_{\omega_1}$ provided that strong Chang conjecture holds in
$\bold V^{{\Bbb P}_j}$ (true if $j$ is measurable $> \aleph_0$)
\sn
\item "{$(c)$}"  if we like to have $(*)$ and in $\bold V^{{\Bbb P}_j}$, strong
Chang conjecture holds then allow: $\Bbb Q^2_{j,S}$ where 
$S \subseteq S^2_0$ stationary not containing a closed copy of $\omega_1$
and $\Bbb Q^2_{j,S}$ semi-proper where $\Bbb Q^2_{j,S}$ shoot an $\omega_1$-increasingly
continuous chain i.e.
$$
\align
\Bbb Q^2_{j,S} = \{(S,f):&S \subseteq S^2_0 \text{ stationary, Dom}(f)
\text{ is a successor} \\
  &\text{countable ordinal, } f \text{ is increasingly continuous into } S\}
\endalign
$$

$$
(S,f_1) \le (S,f_2) \Leftrightarrow S_1 = S_2 \wedge f_1 \subseteq f_2.
$$
\ermn
So $\Vdash_{{\Bbb P}_j} ``{\underset\tilde {}\to {\Bbb Q}_j}$ is semi-proper of cardinality
$\le (2^{\aleph_1})^{V^{P_j}}"$. \nl
So by \cite[Ch.XVI,\S2,2.4,2.5]{Sh:f}
\mr
\item "{$\bigotimes_1$}"  for $i < j \le \kappa,\Bbb P_j/\Bbb P_i$ is semi-proper,
so $\aleph_1$ is not collapsed
\sn
\item "{$\bigotimes_2$}"  $\Bbb P_\kappa$ collapses every 
$\theta \in (\aleph_1,\kappa)$, satisfies the $\kappa$-c.c. and has
cardinality $\kappa$
\sn
\item "{$\bigotimes_3$}"  $\Vdash_{{\Bbb P}_\kappa} 
``{\underset\tilde {}\to D_{\omega_1}}$ is $\aleph_2$-saturated".
\ermn
By preliminary forcing, \wilog \, there is 
$S_0 = \{\delta < \kappa:\delta \text{ strong limit, cf}(\delta) =
\aleph_0\}$, stationary in $\kappa$, reflecting only in inaccessibles.
Let $S_1 = \{\lambda < \kappa:\lambda
\text{ is measurable}\}$ so we know $S_1$ is stationary.  If we are proving
the version with $(*)$, note that $\lambda \in S_1 \Rightarrow$ in
$\bold V^{{\Bbb P}_\lambda}$, the strong Chang conjecture holds (\cite[Ch.XIII,1.9]{Sh:f})
hence $\Bbb Q^2_{\lambda,S}$ is semi-proper for every stationary $S \subseteq
\{\delta < \lambda:\bold V^{{\Bbb P}_\lambda} \models \text{ cf}(\delta) = \aleph_0\}$.
Also if $\boxtimes$ holds then

$$
\align
\Vdash_{{\Bbb P}_\kappa} ``&\underset\tilde {}\to S \subseteq \{\delta < \lambda:
\text{in } \bold V^{{\Bbb P}_\kappa},\text{cf}(\delta) = \aleph_0\} \text{ is
stationary } \Rightarrow \\
  &S'_1 = \{\lambda \in S_1:\underset\tilde {}\to S \restriction \lambda
\text{ is a } \Bbb P_\lambda \text{-name of a stationary} \\  
  &\text{subset of } \{\delta < \lambda:\text{cf}(\delta) = \aleph_0
\text{ in } \bold V^{{\Bbb P}_\lambda}\}\} \text{ is stationary}".
\endalign
$$
\mn
So
\mr
\item "{$\bigotimes_4$}"  if we are proving $(*)$, then in $\bold
V^{{\Bbb P}_\kappa}$ 
every stationary $S \subseteq \{\delta < \kappa:\text{cf}(\delta) = 
\aleph_0\}$ contains a close copy of $\omega_1$.
\ermn
Now we have to deal with the ``projectively stationary".  We can find
function $h$, Dom$(h) = S_0,h(\delta)$ is a $\Bbb P_\delta$-name of a stationary
co-stationary subset of $\omega_1$ (even $\Bbb P_\alpha$-name for some
$\alpha < \delta$) such that: every such name appears
stationarily often.  Let $\langle {\underset\tilde {}\to a^\delta_i}:i <
\omega_1 \rangle$ be a $\Bbb P_\delta$-name such that

$$
\align
\Vdash_{{\Bbb P}_\delta} ``&{\underset\tilde {}\to a^\delta_i} 
\subseteq \delta \text{ is countable unbounded
in } \delta \text{ increasingly continuous in } i \\
  &\text{and } \delta = \dbcu_{i < \omega_1}
{\underset\tilde {}\to a^\delta_i}".
\endalign
$$
\mn
Let

$$
{\underset\tilde {}\to {\Cal W}_\delta} = \{a^\delta_i:i \in h(\delta)\}
$$

$$
{\underset\tilde {}\to {\Cal W}_{< \alpha}} = \cup\{{\Cal W}_\delta:
\delta \in \alpha \cap S_0\}
$$

$$
\underset\tilde {}\to {\Cal W} = {\Cal W}_{< \kappa}.
$$
\mn
So 
\mr
\item "{$\bigotimes_5$}"  $\underset\tilde {}\to {\Cal W}$ is a
$\Bbb P_\kappa$-name of a subset of $[\kappa]^{\aleph_0}$ and $S_0$ is
stationary in $\bold V^{{\Bbb P}_\kappa}$ (as $\Bbb P_\kappa \models \kappa$-c.c.)
\sn
\item "{$\bigotimes_6$}"  $\Vdash_{{\Bbb P}_\kappa} ``\underset\tilde {}\to {\Cal W}$
is stationary.
\ermn
[why?  If 
$\Vdash_{{\Bbb P}_\kappa} ``\underset\tilde {}\to M$ is a model with countable
vocabulary and universe $\kappa$" \ub{then} $\underset\tilde {}\to E = \{
\delta < \lambda:\underset\tilde {}\to M \restriction \delta \text{ is a }
\Bbb P_\delta$-name and is an elementary submodel of $M\}$ is a $\Bbb P_\kappa$-name
of a club of $\kappa$ hence contains a club $E^*$ of $\kappa$ from
$\bold V$.  
So for a club $i < \omega_1,M \restriction a^\delta_i$ is an 
elementary submodel of $M$.  But for stationarily many $i < \omega_1,
a^\delta_i \in {\Cal W}_\delta \subseteq {\Cal W}$, so really ${\Cal W}$ is
stationary. If
$\underset\tilde {}\to W$ is a $\Bbb P_\kappa$-name of a stationary co-stationary
subset of $\omega_1$ then for some, even for stationarily many 
$\delta \in E^* \cap S_0$ we have $h(\delta) =
\underset\tilde {}\to W$ and so easily
\mr
\item "{$\bigotimes_7$}"  $\Vdash_{{\Bbb P}_\kappa} 
``\underset\tilde {}\to {\Cal W}$ is projectively stationary".
\ermn
Lastly, why would ${\Cal W}$ contain no increasing $\omega_1$-chains?  Assume
$p^* \Vdash ``\langle {\underset\tilde {}\to a_i}:i < \omega_1 \rangle$ is
increasing continously and ${\underset\tilde {}\to a_i} \in W"$.  
So \wilog \, for some $\delta^*$ either
\mr
\item "{$(\alpha)$}"  $p^* \Vdash ``\sup(a_i)$ is strictly increasing with
limit $\delta^*"$ or
\sn
\item "{$(\beta)$}"  $p^* \Vdash ``\sup {\underset\tilde {}\to a_i}$ is
constantly $\delta^*$ for $i \ge i^*,i^* < \omega_1$ so \nl
without loss of generality $i^*=0"$.
\ermn
\ub{Case A}:  The possibility $(\beta)$ holds. \nl
Necessarily $\delta^* \in S_0,p^* \Vdash_{\Bbb P} ``\{{\underset\tilde {}\to a_i}:
i < \omega_1\} \subseteq {\Cal W}_{\delta^*}"$ and as $\Vdash_{P_{\delta^*}}
``h(\delta^*)$ is co-stationary subset of $\omega_1"$ and $\Bbb P_\kappa/
\Bbb P_{\delta^*}$ is semi-proper hence preserve stationarity of subsets of
$\omega_1$ we are done. 
\mn
\ub{Case B}:  Possibility $(\alpha)$ holds and $\delta^*$ not 
strongly inaccessible.  So $S_0 \cap \delta^*$ 
is not stationary in $\delta^*$ hence ${\underset\tilde {}\to {\Cal W}} \cap
[\delta^*]^{\aleph_0} = \dbcu_{\delta \in \delta^* \cap S_0}{\Cal W}_\delta$
is not even stationary.
\mn
\ub{Case C}:  Possibility $(\alpha)$ holds and not case B, in 
$V^{P_{\delta^*}}$ strong Chang conjecture fails.

Then ${\underset\tilde {}\to {\Bbb Q}_\delta}$ is Levy$(\aleph_1,2^{\aleph_1}) *$
Cohen (as in clauses (b) and (c) in $\bold V^{{\Bbb P}_\delta}$ strong Chang conjecture 
holds), so as
clearly in $\bold V^{{\Bbb P}_\delta},2^{\aleph_0} = \aleph_2$ (by the Cohen in (a), i.e.
$\Bbb Q^0_j$), then in $\bold V^{{\Bbb P}_{\delta^*}}$ for every club $E'$ of $[\delta^*]
^{\aleph_0}$, we can find some $\delta < \delta^*$ and $2^{\aleph_0}$
members of $E' \cap [\delta]^{\aleph_0} \backslash
{\underset\tilde {}\to W_{< \delta^*}}$.  So in $\bold V^{{\Bbb P}_{\delta^*}},
[\delta^*]^{\aleph_0} \backslash 
{\underset\tilde {}\to {\Cal W}_{< \delta^*}}$ is stationary
and ${\underset\tilde {}\to {\Bbb Q}_\delta}$ is proper so this
holds in $\bold V^{{\Bbb P}_{\delta^*+1}}$.  But 
$\Bbb P_\kappa/\Bbb P_{\delta^*+1}$ preserves
stationarity of subsets of $\omega_1$ hence in $\bold V^{{\Bbb P}_\kappa}
[\delta^*]^{\aleph_0} < {\Cal W} < \delta$ is stationary, so we are done.
\mn
\ub{Case D}:  Possibility $(\alpha)$ holds, not case B and in 
$\bold V^{{\Bbb P}_{\delta^*}}$ strong Chang conjecture holds.
\mn
Just note: in $\bold V^{{\Bbb P}_{\delta^*}}$, let $p \in \Bbb Q_{\delta^*}$, let 
$\chi$ large enough $N \prec ({\Cal H}(\chi),\in)$ is countable
to which $\bar{\Bbb Q},\delta^*,G_{{\Bbb P}_{\delta^*}},p,\langle {\Cal W}_\delta:
\delta \in S_0 \cap \delta^* \rangle$ belong, \ub{then} we can find 
(see \cite[Ch.XIII]{Sh:f}) $T \subseteq {}^{\omega >} \omega_2$ closed 
under initial segments $T \cap N = \emptyset$, 
satisfying $(\forall \eta \in T)(\exists^{\aleph_2}\alpha)
(\eta \char 94 \langle \alpha \rangle \in T)$ and $\langle N_\eta:\eta \in T
\rangle$ such that
\mr
\widestnumber\item{$(iii)$}
\item "{$(i)$}"  $N_{<>} = N$
\sn
\item "{$(ii)$}"  $N_\eta \prec ({\Cal H}(\chi),\in)$ is countable
\sn
\item "{$(iii)$}"  $\eta \in N_\eta,N_\eta \cap \omega_1 = N \cap \omega_1$
\sn
\item "{$(iv)$}"  $\nu \triangleleft \eta \Rightarrow N_\eta \subseteq N_\nu$
\sn
\item "{$(v)$}"  if ${\Cal I} = \{A_\zeta:\zeta < \zeta^*\}$ is a maximal
antichain of ${\Cal D}_{\omega_1}$ which is semi-proper 
and ${\Cal I} \in N_\eta$ \ub{then} for
some $k < \omega,\eta \triangleleft \nu \in T \and \ell g(\nu) \ge k
\Rightarrow N \cap \omega_1 \in \dbcu_{\zeta \in N_\nu} A_\zeta$.
\ermn
Let

$$
\align
E = \{\delta < \omega_2:&\text{ if } \eta \in {}^{\omega >}\delta
\text{ then } N_\eta \cap \omega_2 \\  
  &\text{is a bounded subset of } \delta\}.
\endalign
$$
\mn
Now if $p \in \Bbb Q^0_\delta$ we do as in Case C.  If $p \in \Bbb Q^1_{\delta^*}$, 
choose $\delta \in E$, cf$(\delta) = \aleph_0$, and
such that for every $\eta \in T \cap {}^{\omega >} \delta,\delta =
\text{ otp}\{\beta < \delta:\eta \char 94 \langle \beta \rangle \in T\}$ and
$\eta \char 94 \langle \alpha \rangle \in T \and \alpha < \delta \Rightarrow 
\text{ sup}(N_\alpha \cap \omega_2) < \delta\}$.  Now 
we can by cardinality  considerations
$(2^{\aleph_0} > \aleph_1)$ find $\eta \in \lim(T) \cap {}^\omega \delta$
such that letting $M = \dbcu_{\ell < \omega} N_{\eta \restriction \ell},
M \cap \omega_2 = M \cap \delta \notin {\Cal W}_{< \delta^*}$.  
So there is $q \in \Bbb Q_\delta$ which is $(M,\Bbb Q_\delta)$-generic, 
$p \le q$ (by the definition of $\Bbb Q^1_{\delta^*}$).  Now $q$  forces
${\underset\tilde {}\to a_{M \cap \omega_1}} = a_{M \cap \omega_1}$ to be
$M \cap \omega_2$ which is not in 
${\underset\tilde {}\to {\Cal W}_{< \delta^*}}$.
\sn
Lastly if $q \in \Bbb Q^2_{j,S}$ (in $\bold V^{{\Bbb P}_{\delta^*}}$) as $S$
does not reflect we can find $\delta \in E$ as above, $\delta \in S$, cf$(\delta) =
\aleph_0$ and choose $\eta,M$ as above. \nl
3) We may like to adapt the proof above.  \nl
We omit the choice of $\langle {\underset\tilde {}\to a^\delta_i}:i <
\omega_1 \rangle$, \ub{but} in ${\underset\tilde {}\to {\Bbb Q}_j}$ if $j \in S_0$
we also choose a $\Bbb P_j$-name of a countable unbounded subset of $\delta,
{\underset\tilde {}\to a_\delta}$ and let 
${\underset\tilde {}\to {\Cal W}_\delta} = \{{\underset\tilde {}\to a_\delta}
\}$ so $\Bbb Q_j$ is replaced by 
$\Bbb Q_j \times \{\underset\tilde {}\to a:\underset\tilde {}\to a$
a name as above$\}$.  Now $h_0$ has domain $S^* = \{\delta:\delta$
strongly inaccessible, in $\bold V^{{\Bbb P}_\delta}$, strong Chang conjecture holds$\},
h_0(\delta)$ a $\Bbb P_\delta$-name of a stationary co-stationary subset of
$\omega_1$ and we add to clauses (a), (b), (c) above also 
\smallskip

(d)   define in $\bold V^{{\Bbb P}_\delta}$:

$$
\align
\Bbb Q^3_\delta = \biggl\{ \langle M_i:i \le j \rangle:&\text{the ordinal }
j \text{ is countable and} \\
  &M_i \prec ({\Cal H}
((2^\delta)^+),\in) \text{ is countable increasing} \\
  &\text{continuous, and: \ub{if} } M_i \cap \omega_1 \in h_0(\delta)
\text{ \ub{then}} \\
  &M_i \in {\Cal W}_{< \delta} \text{ and if } 
M_i \cap \omega_1 \notin h(\delta) \text{ then } 
M_i \notin {\Cal W}_{< \delta} \biggr\}.
\endalign
$$
\mn
Now again we use $\langle N_\eta:\eta \in T \rangle$ and choosing $M$ it is
enough to show that
\mr
\item "{$\boxtimes_1$}"  for some $\eta \in \lim T,\dbcu_{i < \omega}
N_{\eta \restriction \ell} \cap \omega_2 \in {\Cal W}_{< \delta}$
\sn
\item "{$\boxtimes_2$}"  some $\eta \in \lim T,\dbcu_{\ell < \omega}
N_{\eta \restriction \ell} \cap \omega_2 \notin {\Cal W}_{< \delta}$.
\ermn
Now $\boxtimes_2$ is as before, $\boxtimes_1$ O.K. by the way
${\underset\tilde {}\to {\Cal W}_{< \delta}}$ is defined.
\hfill$\boxtimes_{\scite{11.1}}$
\enddemo
\newpage

\head {\S10 Mixing finitary norms and ideals} \endhead  \resetall \sectno=10
\bn
We may consider replacing families of ideals by families of creatures
see \cite{RoSh:470} on creatures: \nl
We hope it will gain something
\definition{\stag{0.A} Definition}:  1)  A $\lambda$-creature ${\frak c}$ 
consists of $(D^{\frak c},\le,\text{val}^{\frak c},\text{nor}^{\frak c},
\lambda^{\frak c})$, where: 
\mr
\item "{{}}"  $\lambda^{\frak c} = \lambda$
\sn
\item "{{}}"  $D^{\frak c}$ the domain,
\sn
\item "{{}}"  $\le$ a partial order on $D^{\frak c}$,
\sn
\item "{{}}"  val$^{\frak c}$ is a function from $D^{\frak c}$ to
${\Cal P}(\lambda) \backslash \{\emptyset\}$
\sn
\item "{{}}"  nor$^{\frak c}:D^{\frak c} \rightarrow \omega$ or to??
\ermn
2) It is called simple if nor$^{\frak c}$ is always $>0$ (without loss 
of generality constant, e.g. Rang(val$^{\frak c}$) = $I^+,I$ an ideal 
on $\lambda$).  A creature is a $\lambda$-creature for some $\lambda$. \nl
$\Bbb I$ will be a set of creatures. 
\enddefinition
\bigskip

\definition{\stag{0.B} Definition}  1) An $\Bbb I$-tree is 
$(T,\bold I,\bold d)$ such that: \newline
for some ordinal $\alpha,T \subseteq {}^{\omega >}\alpha$ closed 
under initial segments, $\ne \emptyset$ \newline
$\bold I$ is a partial function Dom$(\bold I) \subseteq T,\bold I_\eta \in
\Bbb I$, \newline
$\bold d$ has domain Dom$(\bold I),\bold d(\eta) \in D^{\bold I_\eta}$ \nl
val$^{\bold I_\eta}(\bold d(\eta)) = 
\{\alpha:\eta \char 94 \langle \alpha \rangle \in T\}$. \newline
2) Let $(T^*,\bold I^*,\bold d^*)$ be an $\Bbb I$-tree, such that
\medskip
\roster
\item "{$(*)$}"  $(\forall \eta \in \text{ lim } T^*)
[\text{lim } \sup {n < \omega} \text{ nor}^{\bold I^*_{\eta \restriction n}}
(\bold d^*_\eta) = \infty]$ \newline
and Dom$(\bold I_\eta) = T^*$.
\endroster
\medskip

\noindent
We define a forcing notion $\Bbb Q = \Bbb Q_{(T^*,\bold I^*,\bold d^*)}$:

$$
\align
\Bbb Q = \biggl\{(T,\bold I,\bold d):&T \subseteq T^*,\bold I = \bold I^*
\restriction T, \\
  &(T,\bold I,\bold d) \text{ an } \Bbb I \text{-tree}, \\
  &(\forall \eta \in T)(\bold d^*_\eta \le^{\bold I_\eta} \bold d_\eta) \\
  &\text{and } ((\forall \eta \in \text{ lim } T^*) \text{ lim sup nor}
^{{\bold I}^*_{\eta \restriction n}}(\bold d_\eta) = \infty) \biggr\}.
\endalign 
$$
\mn
Order: natural.
\sn
3) Let $(T^*,I^*,\bold d^*)$ be an $\Bbb I$-tree such that
\medskip
\roster
\item "{$(**)$}"  $(\forall \eta \in \text{ lim }T^*)(\forall n)
(\forall^* \ell)(\text{nor}^{I^*_{\eta \restriction n}}(\bold d^*_\eta) \ge
n)$ \newline
(i.e. lim inf $= \infty$).
\endroster
\medskip

\noindent
and define $\Bbb Q' = \Bbb Q'_{(T^*,{\bold I}^*,{\bold d}^*)}$ parallelly. \nl
4) For $p \in \Bbb Q$ (or $p \in \Bbb Q'$) we write $p = (T^p,\bold I^p,\bold d^p)$.
In this case for $\eta \in T^p$ we define $q = p^{[\eta]}$ by: \nl
$T^{[q]} = \{\nu \in T^p:\nu \trianglelefteq \eta$ or $\eta \trianglelefteq
\nu\},\bold I^q = \bold I^p \restriction T^q,\bold d^q = \bold d^p 
\restriction T^p$.  Clearly $p \in \Bbb Q \Rightarrow p \le q \in \Bbb Q$ and $p \in
\Bbb Q' \Rightarrow p \le q \in \Bbb Q'$.
\enddefinition
\bigskip

\proclaim{\stag{0.C} Claim}  Let $(T,\bold I^*,\bold d^*)$ and $\Bbb
Q,\Bbb Q'$ be as in
\scite{0.B}.  A sufficient condition for ``$\aleph_1$ not collapsed" is:
\mr
\item "{$(a)$}"  for $\Bbb Q$: $(**)$ below
\sn
\item "{$(b)$}"  for $\Bbb Q':(*) + (**)$ below where
\sn
{\roster
\itemitem{ $(*)$ }  $\Bbb I$ has $\aleph_1$-bigness: 
$$
\align
(\forall {\frak c} \in \Bbb I)(\forall x \in D^{\frak c})
\biggl[ \text{nor}^{\frak c}(x) > 0 \rightarrow &(\forall h \in
{}^{(\lambda^{\frak c})}\omega_1)(\exists y) \\
  &[x \le^{\frak c} y \and \text{ nor}^{\frak c}(y) \ge \text{ nor}^{\frak c}
(x)-1 \and \\
  &(h \restriction \text{ val}^{\frak c}(y) \text{ is constant}] \biggr]
\endalign
$$

\itemitem{ $(**)$ }  $\Bbb I$ is $(\aleph_1,\aleph_1)$-indecomposable where
$\Bbb I$ is $(\mu,\kappa)$-indecomposable means:
\mn
\block
$\boxtimes_{\Bbb I,\mu,\kappa} \qquad$ if ${\frak c} \in \Bbb I$ and $x \in 
D^{\frak c}$ satisfies nor$^{\frak c}(x) > 2$ and $A_\alpha \subseteq
\lambda^{\frak c}$ for $\alpha < \mu$ are such that $(\forall y)(x \le y \in 
D^{\frak c} \wedge \text{ val}^{\frak c}(y) \subseteq A_\alpha \rightarrow 
\text{ nor}^{\frak c}(y)+2 \le \text{ nor}^{\frak c}(x))$, \ub{then} we can
find $u \subseteq \lambda^{\frak c}$ of cardinality $< \mu$ such that for
every large enough $\alpha < \mu$ we have $u \nsubseteq A_\alpha$.
\endblock
\endroster}
\endroster
\endproclaim
\bigskip

\demo{Proof for $\Bbb Q$}  Lets use given $p = (T,\bold I,\bold d) \in
\Bbb Q$ and 
$\Bbb Q$-name $\underset\tilde {}\to \tau$ such that $\Vdash 
{\underset\tilde {}\to \tau}:\omega \rightarrow \omega_1$.  Now we 
choose by induction on $n,p^n,A_n$ such that:
\medskip
\roster
\item "{$(a)$}"  $p^n \le p^{n+1}$
\sn
\item "{$(b)$}"  $A_0,\dotsc,A_n$ are fronts of $p^n$ which means
$(\forall \eta \in \lim^{p^n})(\exists!n)(\eta \restriction n \in A_\ell)$
\sn
\item "{$(c)$}"  $A_\ell$ below $A_{\ell + 1}$ which means \nl 
$(\forall \eta \in A_{\ell +1})(\exists \nu \triangleleft \eta)
\nu \in A_\ell$ \newline
(so $A_n \subseteq T^{p^{n+1}},(\forall \eta \in T^{p^n} \backslash
T^{p^{n+1}})(\exists \nu \triangleleft \eta)(\nu \in A_n))$
\sn
\item "{$(d)$}"  $(\forall \nu \in A_n)(\forall \eta \in \text{ Suc}_{T^{p_n}}
(\nu))(p^{[\eta]}_n$ forces a value to $\underset\tilde {}\to \tau(n))$
\sn
\item "{$(e)$}"  $\eta \in A_n \Rightarrow \text{ nor}^{{\bold I}^*_\eta}
(\bold d^{p_n}_\eta) \ge n \and \bold d^{p_n}_\eta = \bold d^{p_{n+1}}_\eta$,
it follows that $\ell < n \and \eta \in A_\ell \Rightarrow \bold d^{p_\ell}
_\eta = \bold d^{p_n}_\eta \and \text{ Suc}_{T_{p_n}}(\eta) =
\text{ Suc}_{T^{p_\ell}}(\eta)$.
\ermn
So $p^*$ is defined by $T^{p^*} = \dbca_{n < \omega} T^{p_n},\bold I^{p^*} =
\bold I \restriction T^{p^*},\bold d^{p^*} = \dbcu_n\{\bold d^{p_n}
\restriction \{\eta \restriction \ell:\eta \in A_n$ and $\ell \le \ell g
(\eta)\}:n < \omega\}$ belong to $\Bbb Q$ and is an upper bound of $\{p_n:n <
\omega\}$. 

We define $h:T^{p^*} \rightarrow \omega_1$ as follows: if 
$\eta \in T^{p^*},\nu \triangleleft \eta \trianglelefteq \nu',\nu \in 
A_{\ell -1},\nu' \in A_\ell$ \newline
(if $\ell = 0$ omit $\nu'$, so just $\eta \trianglelefteq \nu'$), then 
$p^{*^{[\eta]}}$ forces value to $\underset\tilde {}\to \tau \restriction 
\ell$ call it $(\tau \restriction \ell)^{p^{*^{[\eta]}}}$ and let

$$
h(\eta) = \text{ Sup Rang}(\underset\tilde {}\to \tau \restriction 
\ell)^{p^{*^{[\eta]}}}.
$$
\medskip

\noindent
For notational simplicity $A_n = T^{p^*} \cap {}^n\text{Ord}$. \newline
We now define a game $\Game = \Game^\alpha_{T^{p^*}}$ for each $\alpha < 
\omega_1$:
\sn

A play of the game last $\omega$ moves, in the $(n-1)$-th move a member
$\eta_n$ of $A_n$ is chosen such that $m < n \Rightarrow \eta_m
\triangleleft \eta_n$, and fixing some $\eta_{-1} \in A_n$.
\sn
In the $n$-the move:
\mr
\item "{$(a)$}"  the anti-decidability player chooses a set $A_n \subseteq
\text{ Suc}_{T^{p^*}}(\eta_{n-1})$ such that
$$
B_n = \emptyset \vee (\text{nor}^{{\bold I}_\eta}(A_n) \le n-2,n \ge 2)
$$
{\roster
\itemitem{ $\boxtimes_1$ }  $B_n \ne \emptyset \text{ \ub{or} } n \ge 3$
and for no $d$, satisfying $\bold d^{p^*}(\eta) \le d \wedge
\text{ nor}^{{\bold I}^*_2}(d) \ge n-2$ do we have $\eta \char 94 \langle
\alpha \rangle \in A_n \Rightarrow \alpha \in \text{ val}^{{\bold I}^*_\eta}
(d)$
\endroster}
\item "{$(b)$}"  the decidability player chooses $\eta_n \in A_n$ such that
$\neg(\exists \nu \in B_n)(\nu \trianglelefteq \eta_n)$ and $n \ge 1
\Rightarrow (\eta_n \restriction (\ell g(\eta_{n-1})-1)) \le \alpha$.
\ermn
Without loss of generality $A_0 = \{<>\}$.

If for some $\alpha$ and decidability player has a winning strategy, we can
produce a condition as required.

If not, for every $\alpha < \omega_1$ the antidecidability player has a
winning strategy {\bf St}$_\alpha$.  For each $\eta \in T^{p^*}$ and $\alpha
< \omega_1$, we consider the play of the game in which the antidecidability
player has winning strategy {\bf St}$_\alpha$ and in some move $n$ the
decidability player chooses $\eta_n = \eta$.  Reflecting there is no freedom
left so there is at most one such play and $n$ and let the antidecidability
player choose set $B_{\eta,\alpha}$ as there (if no such game let
$B_{\eta,\alpha} = \emptyset$).
\bn
\ub{Case 1}:  For some $n$ and $\eta \in A_n$, we have: there is no
countable $u \subseteq \text{ Suc}_{T^{p^*}}(\eta)$ such that for every
large enough $\alpha < \omega_1,u \nsubseteq B_{\eta,\alpha}$ so by the
assumption $(**)$ we get a contradiction.
\bn
\ub{Case 2}:  Not Case 1.

We can choose by induction on $n$, a countable subset $u_n \subseteq A_n$
such that: $u_0 = \{<>\}$

$$
\align
\text{if } \eta \in A_n &\text{ then for some } \alpha_\eta < \omega_1
\text{ for every } \alpha \in [\alpha_\eta,\omega_1), \\
  &\text{ if in the play in which the antidecidability player uses }
{\bold {St}}_\alpha \\
  &\text{ and they arrive to } \eta, \text{ there is } \eta',\eta
\triangleleft \eta' \in u_n +1 \\
  &\text{ which is a legal response of the decidability player}.
\endalign
$$
mn
Now let

$$
\align
\alpha^* = &\sup\{h(\eta)+1:\text{ for some } \nu \in \dbcu_n u_n,
\eta \triangleleft \nu,\eta \in \text{ Dom}(h)\} + \\
  &\sup\{\alpha_\eta:\eta \in \dbcu_n u_n\}
\endalign
$$
\mn
and we can find a play of $\Game^{\alpha^*}$ as above where the
decidability player chooses $\eta$'s from $\dbcu_{n < \omega} u_n$.  We get
a contradiction.
\enddemo
\bigskip

\demo{Proof for $\Bbb Q'$}?  We should make changes: in $p^{n+1}$ we shrink
$p^{[\eta]}_n$ for each $\eta \in T^{p^n} \cap {}^n\text{Ord}$, to $q_\eta,
p^{[\eta]}_n \le_{pr} [\eta]$ and for each $\ell \le n$, if possible, 
$q_\eta$ forces a bound to $\underset\tilde {}\to \tau(\ell)$
and, of course, $p^{[\eta]}_{n+1} = q_\eta$ for each such $\eta$ and
$T^{p^{n+1}} \cap {}^{\eta \ge}\text{Ord } = 
T^{p^n} \cap {}^{n \ge}\text{Ord}$ and 
$\bold d^{p^{n+1}} \restriction {}^{n \ge}\text{Ord } = \bold d^{p^n}
\restriction {}^{n \ge}\text{Ord}$.  
So let $p^* = \dbca_n p_n$ be naturally defined, and we use 2-bigness to
prove enough times $q_\eta$ forces a bound.
\enddemo
\bn
Now we give details.  
\demo{Proof for $\Bbb Q$}  Given $p = (T,\bold I,\bold d)$, for notational
simplicity tr$(p) = <>$ and nor$^{{\bold d}^p_\eta}(\text{Suc}_{T^p}(\eta))
> 2$ and $\Bbb P$-name $\underset\tilde {}\to \tau$ such that 
$\Vdash {\underset\tilde {}\to \tau}:\omega \rightarrow \omega_1$ we 
choose by induction on $n,p^n$ such that:
\medskip
\roster
\item "{$(a)$}"  $p^n \le p^{n+1}$, and $p^n$ has trunk $n$
\sn
\item "{$(b)$}"  $A_0,\dotsc,A_n$ are fronts of $p^n$
\sn
\item "{$(c)$}"  $A_\ell$ below $A_{\ell + 1}$ which means \nl 
$(\forall \eta \in A_{\ell +1})(\exists \nu \triangleleft \eta)
\nu \in A_\ell$ \newline
(so $A_n \subseteq T^{p^{n+1}},(\forall \eta \in T^{p^n} \backslash
T^{p^{n+1}})(\exists \nu \triangleleft \eta)(\nu \in A_n))$
\sn
\item "{$(d)$}"  $A_0 = \{<>\}[\eta \in A_n \wedge \eta \trianglelefteq
\nu \in T^{p^n} \Rightarrow \text{ nor}^{{\bold I}_\nu}(\bold d^{p_n}_\nu)
> n+2]$
\sn
\item "{$(e)$}"  when $\eta \in A_n$ let $\ell = \ell_\eta \le n$ be
maximal such that there are $\alpha_m < \omega_1,m < \ell$ and $q$ satisfying
$\text{tr}q = \eta,p^{[\eta]}_n \le q \in P,q \Vdash ``\dsize \bigwedge_{m <
\ell} \underset\tilde {}\to \tau(m) < \alpha'_n$ and $\eta \trianglelefteq \nu
\in T^q \Rightarrow \text{ nor}^{{\bold d}^*_\nu}(\bold d^q_\nu) \ge n$ and
we demand: $p^{[\eta]}_{n+1}$ satisfies the demand on $q$ for some $\langle
\alpha_m:m < \ell_\eta \rangle$, note possible $\ell_\eta = \nu$ then we are
left with demand on norm.
\ermn
So $p^*,T^{p^*} = \dbca_{n < \omega} T^{p_n}$ is an upper bound of 
$\{p_n:n < \omega\}$.

Clearly $p^* \in \Bbb P$ and $n < \omega \Rightarrow p_n \le p^*$.  Let for
$\eta \in T^{p^*}$, let $n(\eta) = \text{ Max}\{n:\text{there is } \nu
\triangleleft \eta,\nu \in A_n\}$ and $\nu_\eta \triangleleft \eta$ be in
$A_n$ and $\beta_\eta < \omega_1$ be minimal such that $p^{[\nu_\eta]}_{n+1}$
forces $\tau(0),\dotsc,\tau(\ell_{\nu_\eta}-1) < \beta_\eta$.  Using games as
in the proof for $\Bbb Q'$ there is $p^+$ such that:
\mr
\item "{$(a)$}"  $p^* \le p^* \in \Bbb P$
\sn
\item "{$(b)$}"  $\rho \in T^{p^+} \Rightarrow \text{ nor}^{{\bold I}^*_\rho}
(\bold d^{p^+}_\rho) \ge \text{ nor}^{{\bold I}^*_\rho}(\bold d^{p^*}_\rho)-1$
\sn
\item "{$(c)$}"  $\beta^* = \sup\{\beta_\eta:\eta \in T^{p^+}\} < \omega_1$.
\ermn
We continue as in \cite[Ch.XIV,\S5]{Sh:f}.
\enddemo
\bn
\centerline {$* \qquad * \qquad *$}
\bn
\ub{Discussion}:  We can continue to do iteration.  \newline
But more urgent: can $\Bbb Q,\Bbb Q'$ like this do anything not already covered by
composition?

A natural thought is splitting or reaping numbers.  We can think of the
tree splitting in $T^*$ as a list of the reals.  BUT, what is the norm?
\bigskip

\centerline {$* \qquad * \qquad *$}
\bigskip

Not finished...check the better's theorem proof?
\bigskip

\noindent
\underbar{Assignment}: \cite[XIII,XVI]{Sh:f} and here put together, so does
the reflection $Pr_a(\lambda,f)$ works for ???
\bn
\ub{97/2/2 - Discussion}:

Saying a creature is $\mu$-complete meana that for pure extensions, increasing
chains of length $< \mu$ have pure upper bounds?  Probably pure means not
changing the norm; maybe the $\aleph_1$-indecomposable is enough. \nl
So the $\Bbb I$-th condition has a new meaning. 
\mn
\ub{Question}:  Does the theorem here hold?
\mn
\ub{Question}:  Does this new context have real applications?

The first result to be discussed is moving from 
${\underset\tilde {}\to {\Bbb I}}$ to one in the ground model. \nl
The second are \scite{5.2}, \sciteu{6.2} preservation of $N$ being suitable.
\newpage

\head{\S11 Variants of the iteration} \endhead  \resetall \sectno=11
\bigskip

As mentioned in \S1, we can consider $\kappa$-RS iteration and
variants of the Sp$_e$ iteration.
\bigskip

\definition{\stag{11.13} Definition/Claim}
Let $\kappa$ be a successor cardinal or
an infinite ordinal not a cardinal but an ordinal of power $|\kappa|,\kappa$ 
fix \footnote{For $\kappa$ inaccessible, see \sciteu{1.22}.}.  
We define and prove the following by induction on $\alpha$ (here 
$e = \{3,4,5,6\}$).  If $\kappa = \aleph_1$, we may omit it and this is 
the main case.  

We repeat \scite{1.13} with the following changes in the proof and
definition:
\mr
\item "{$(B)$}"  We say $\underset\tilde {}\to \zeta$ is a simple
$\bar{\Bbb Q}$-named$_e$ \, $[j,\beta)$-ordinal if
{\roster
\itemitem{ $(*)_1$ }  $\underset\tilde {}\to \zeta$ is a simple
$\bar{\Bbb Q}$-named$^1$ \, $[j,\beta)$-ordinal and may restrict ourselves
to $\kappa = \aleph_1 \Rightarrow e \in \{3,4\}$
\sn
\itemitem{ $(*)_2$ }  if $e \in \{5,6\}$ and $\kappa = \aleph_1$, \ub{then}
$\underset\tilde {}\to \zeta$ is a simple $\bar{\Bbb Q}$-named$^2$ \, 
$[j,\beta)$-ordinal.
\endroster}
\item "{$(F)(a)$}"  in (v) replace the remark in the end by: \nl
``if $e \in \{5,6\},\alpha \in w$ then this demand follows by \scite{1.6A}
\nl
and add:
{\roster
\itemitem{ $(vii)$ }  if $e=3,5$ then for some $n < \omega$ and simple
$\bar{\Bbb Q}$-named $[0,\ell g(\bar{\Bbb Q}))$-ordinals
${\underset\tilde {}\to \xi_1},\dotsc,{\underset\tilde {}\to \xi_n}$ 
we have, for every $\xi < \ell g(\bar{\Bbb Q}) 
\Vdash_{{\Bbb P}_\xi}$ ``if for $\ell = 1,\dotsc,n$ we have
${\underset\tilde {}\to \xi_\ell}
[{\underset\tilde {}\to G_{{\Bbb P}_\xi}}] \ne \xi$
(for example ${\underset\tilde {}\to \xi_\ell}[G_{{\Bbb P}_\xi}]$ 
not well defined) then 
$\emptyset_{\underset\tilde {}\to {\Bbb Q}_\xi} \le_{\text{pr}}
p \restriction \{\xi\}$ in ${\underset\tilde {}\to {\hat{\Bbb Q}}_\xi}"$ 
\endroster}
\item "{$(F)(e)(iii)$}"   inside change $p_2 \restriction \xi
\Vdash_{{\Bbb P}_\xi} ``\ldots"$ by \nl
$p_2 \restriction \xi \Vdash_{{\Bbb P}_\xi}$ ``if
$\xi \ne {\underset\tilde {}\to \xi_\ell}
[{\underset\tilde {}\to G_{{\Bbb P}_\xi}}]$ for \nl

$\qquad \qquad \ell = 1,\dotsc,n$ and \nl

$\qquad \qquad [e=4 \vee e=6 \Rightarrow 
\neg(\emptyset_{\underset\tilde {}\to {\Bbb Q}_\xi} \le_{vp} p^1 
\restriction \{\xi\})]$  \nl

$\qquad \qquad$ then: 
${\underset\tilde {}\to {\hat{\Bbb Q}}_\xi} 
\models p^1 \restriction
\{\xi\} \le_{pr} p^2 \restriction \{\xi\}"$.
\endroster
\enddefinition
\bigskip

\proclaim{\stag{11.14} Claim}  1) As in \scite{1.14} adding $\kappa \ne
\aleph_1 \vee e \in \{3,4\}$. \nl
2)  If $\bar{\Bbb Q}$ is an $\kappa-\text{Sp}_e(W)$-iteration, 
and for each $i$ the quasi-order $\le^{{\Bbb Q}_i}_{\text{pr}}$ is 
equality hence $\le^{{\Bbb Q}_i}_{\text{vpr}}$ is equality,
\ub{then} $\bar{\Bbb Q}$ is essentially a finite support iteration. 
\nl
[Saharon:  maybe restrict yourself above the constantly function
$\zeta \mapsto \emptyset_{{\Bbb Q}_\zeta}$, so we have to use $\kappa >
\ell g(\bar{\Bbb Q})$.]
\endproclaim
\bigskip

\proclaim{\stag{11.15} Claim}  1)  Add
\mr
\item "{$(d)$}"  if $e=3,5$ then $r$ is pure outside
$\{\xi_1,\dotsc,\xi_n\}$. 
\ermn
2) In the proof on ``$\xi^*_{\ell \times 1}$" ?? we say, i.e., simple$^1$
if $e \in \{3,4\}$ and simple$^2$ if $e \in \{5,6\}$.
\endproclaim
\bigskip

\citewarning{\noindent \llap{---$\!\!>$} MARTIN WARNS: Label 1.16 on next line is also used somewhere else (Perhaps should have used scite instead of stag?)}
\proclaim{\stag{1.16} Claim}
5) If $e \in \{3,4\}$ and \footnote{for the parallel for $e \in \{5,6\},
\kappa = \aleph_1$ we need pure decidability and restrict ourselves to
``above $p$" 
for purely dense sets of $p-s$} for each $\beta < \ell g(\bar{\Bbb Q}),
{\underset\tilde {}\to {\bold t}_\beta}$ is a $\Bbb P_\beta$-name of a truth
value, \ub{then} there is a simple $(\bar{\Bbb Q},W)$-named 
$[0,\alpha)$-ordinal
$\underset\tilde {}\to \zeta$ such that $\underset\tilde {}\to \zeta[G_\beta]
= \beta$ iff ${\underset\tilde {}\to {\bold t}_\beta}[G_\beta] =$ truth and
$\gamma < \beta \Rightarrow {\underset\tilde {}\to {\bold t}_\gamma}
[G_\beta] =$ false for any subset $G_\beta$ of $\Bbb P_\beta$ generic
over $\bold V$.
\endproclaim
\bigskip
 
We can deal with the parallel of hereditarily countable names.  This is not
used in later sections.
\definition{\stag{11.21} Definition}  We define for an 
$\kappa-\text{Sp}_e(W)$-iteration $\bar{\Bbb Q}$, and 
cardinal $\mu$ ($\mu$ regular), when is a $(\bar{\Bbb Q},W)$-name 
hereditarily $< \mu$, and in particular when a
$(\bar{\Bbb Q},W)$-named $[j,\alpha)$-ordinal is
hereditarily $< \mu$ and a $(\bar{\Bbb Q},W)$-named 
$[j,\alpha)$-atomic condition hereditarily $< \mu$, 
and which conditions of $\text{Sp}_e(W)$-Lim$_\kappa(\bar{\Bbb Q})$ 
are hereditarily $< \mu$.  For simplicity we are assuming that the set of 
members of $\Bbb Q_i$ is in $\bold V$.  This is done by induction on $\alpha = 
\ell g(\bar{\Bbb Q})$.
\mn
\underbar{First Case}:  $\alpha = 0$.

Trivial.
\mn
\underbar{Second Case}:  $\alpha > 0$.
\mr
\item "{$(A)$}"  A $\bar{\Bbb Q}$-named $[j,\alpha)$-ordinal 
$\underset\tilde {}\to \xi$ 
hereditarily $< \mu$ is a $(\bar{\Bbb Q},W)$-named
$[j,\alpha)$-ordinal which can be represented as follows: there is
$\langle (p_i,\xi_i):i < i^* \rangle,i^* < \mu$, each $\xi_i$ an ordinal in
$[j,\alpha),p_i \in \Bbb P_{\xi_i}$ is a member of $\Bbb P_{\xi_i}$ 
hereditarily $< \mu$ and for any $G \in \text{ Gen}^r(\bar{\Bbb Q}),
\underset\tilde {}\to \zeta[G]$ is $\zeta$ iff for some $i$ we have:
{\roster
\itemitem{ $(a)$ }  $p_i \in G,\zeta_i = \zeta$
\sn
\itemitem{ $(b)$ }  if $p_j \in G$ then $\zeta_i < \zeta_j \vee (\zeta_i =
\zeta_j \and i < j)$
\endroster}
\item "{$(B)$}"  A $(\bar{\Bbb Q},W)$-named $[j,\alpha)$-atomic condition 
$\underset\tilde {}\to q$ hereditarily $< \mu$, is a $(\bar{\Bbb Q},W)$-named
$[j,\alpha)$-atomic condition which can be represented as follows: there is
$\langle (p_i,\zeta_i,q_i):i < i^* \rangle,i^* < \mu,\zeta_i \in [j,\alpha),
p_i \in \Bbb P_{\zeta_i},q_i \in \bold V$, and for any $G \in \text{ Gen}^r(\bar{\Bbb Q}),
\underset\tilde {}\to q[G]$ is $q$ iff for some $i$ we have:
\sn
{\roster
\itemitem{ $(a)$ }  $p_i \in G,q = q_i$, and $p_i \Vdash_{{\Bbb P}_{\zeta_i}}
``q \in {\underset\tilde {}\to {\Bbb Q}_{\zeta_i}}"$
\sn
\itemitem{ $(b)$ }  if $p_j \in G$ then $\zeta_i < \zeta_j \vee (\zeta_i =
\zeta_j \and i < j)$
\endroster} 
\sn
\item "{$(C)$}"  A member $p$ of $\Bbb P_\alpha =
Sp_e(W)$-Lim$_\kappa(\bar{\Bbb Q})$ is 
hereditarily $< \mu$ if each member of 
$r$ is a $(\bar{\Bbb Q},W)$-named atomic
condition hereditarily $< \mu$
\sn
\item "{$(D)$}"  A $(\bar{\Bbb Q},W)$-name of a member of 
$\bold V$ hereditarily $< \kappa$
is defined as in clause $(B)$, similarly for member $x \in \bold
V^{{\Bbb P}_\alpha}$ such
that $y \in$ transitive closure of $x \Rightarrow |y| < \mu$.
\endroster
\enddefinition
\bigskip

\remark{\stag{11.22} Concluding Remarks}  1) We 
have not really dealt with the 
case $\kappa$ is inaccessible.  The point is that in this case, we do not 
know a priori the length of the list of the members of a condition 
(which are atomic conditions).  It is natural to work on it together 
with ``decidability on bound
on $\underset\tilde {}\to \alpha < \kappa$ by pure extensions",
see \scite{1.25} below. \newline
2) We can think of putting together \cite[Ch.XIV]{Sh:f} and \cite{Sh:587}.
\newline
3) We can ask:  Does ``Souslin forcing notions" help?
\endremark
\bigskip

\citewarning{\noindent \llap{---$\!\!>$} MARTIN WARNS: Label 11.27 on next line is also used somewhere else (Perhaps should have used scite instead of stag?)}
\proclaim{\stag{11.27} Claim}  $e \in \{4,5\}$ is O.K.
\endproclaim
\bigskip

\citewarning{\noindent \llap{---$\!\!>$} MARTIN WARNS: Label 11.27 on next line is also used somewhere else (Perhaps should have used scite instead of stag?)}
\proclaim{\stag{11.27} Claim}  In the proof of \scite{1.27}, in case 3
add:
\nl
(the point is that $e \in \{4,6\}$).  Instead $e \in \{4,6\}$ it is enough to
assume:
\mr
\item "{$\boxtimes_{\underset\tilde {}\to {\Bbb Q}_\beta}$}"  for every
$q',q'' \in {\underset\tilde {}\to {\Bbb Q}_{\beta_0}}$ we have \nl
$\emptyset_{{\Bbb Q}_{\beta_0}} \le_{\text{vpr}} q' \le q'' \Rightarrow q' 
\le_{\text{pr}} q''$.
\endroster
\endproclaim
\bigskip

\remark{\stag{11.27B} Remark}  1) Add: \nl
but for $e=4$ we could use appropriate $p_1 = p \cup \{
{\underset\tilde {}\to r_1}\},\underset\tilde {}\to r$ an atomic
$(\bar{\Bbb Q},W)$-named condition,
${\underset\tilde {}\to \zeta_{\underset\tilde {}\to r}} =
\underset\tilde {}\to \zeta$, see \scite{1.6}(5).  \nl
2) Holds for $e \in \{4,6\}$.
\endremark
\bigskip

\definition{\stag{n.2} Definition}  0) Let $\bar{\Bbb Q}$ be a
$\kappa_1$-Sp$_e(W)$-iteration of length $\alpha$.  Let 
$\underset\tilde {}\to \zeta$ denote a simple $\bar{\Bbb Q}$-named
$[0,\alpha)$-ordinal \ub{or} a simple $\bar{\Bbb Q}$-named$^2 
[0,\alpha)$-ordinals and $\Xi$ a countable set of such objects. \nl
1) For an atomic simple $\bar{\Bbb Q}$-named condition $r,r \restriction
\underset\tilde {}\to \zeta$ is defined by $r \restriction
\underset\tilde {}\to \zeta[G] = r^* \in \Bbb P_\zeta$ if
$\underset\tilde {}\to \zeta[G] \ge {\underset\tilde {}\to \zeta_r}[G],
r[G] = r^*$ and $\emptyset_{p_\zeta}$ otherwise. \nl
2) For $q \in \Bbb P_\alpha,q \restriction \underset\tilde {}\to \zeta =
\{r \restriction \zeta:r \in q\}$ and $q \restriction \Xi =
\dbcu_{{\underset\tilde {}\to \zeta} \in \Xi} q \restriction \zeta$. \nl
3) $\Bbb P_{\underset\tilde {}\to \zeta} 
= \{p \in \Bbb P_\alpha:p \restriction
\underset\tilde {}\to \zeta = p$, i.e., for every $G \subseteq \Bbb P_\alpha$
generic over $\bold V,p \restriction \underset\tilde {}\to \zeta[G] = p[G]\}$

$$
P_\Xi = \{p \in \Bbb P_\alpha:p \restriction \Xi = p\}
$$
\mn
both with the order inherited from $\Bbb P_\alpha$.
\enddefinition
\bigskip

\proclaim{\stag{n.3} Claim}  Let $\bar{\Bbb Q}$ be an
$\aleph_1-{\text{\rm Sp\/}}_e(W)$-iteration. \nl
1) If ${\underset\tilde {}\to \zeta_1}$ is a simple $\bar{\Bbb Q}$-named
$[\beta_1,\beta_2)$-ordinal, ${\underset\tilde {}\to \zeta_2}$ is a simple
$\bar{\Bbb Q}$-named$^1\,[\beta_1,\beta_2)$-ordinal, 
\ub{then} there is a simple
$\bar{\Bbb Q}$-named 
$[\beta_1,\beta_2)$-ordinal $\underset\tilde {}\to \zeta$
such that for $G \subseteq \Bbb P_\alpha$ is generic over $\bold V$:
\mr
\item "{$(a)$}"  if ${\underset\tilde {}\to \zeta_1}[G] =
{\underset\tilde {}\to \zeta_1}[G \cap \Bbb P_\xi] = \xi$ 
and ${\text{\rm Min\/}}\{\varepsilon$:
some $p \in G \cap \Bbb P_\varepsilon$ decided to be $\varepsilon$ or be
undefined$\} > \varepsilon$ \ub{then} $\underset\tilde {}\to \zeta[G] =
\underset\tilde {}\to \zeta[G \cap \Bbb P_\xi] = \xi$
\sn
\item "{$(b)$}"  otherwise undefined.
\ermn
2) Let 
$\underset\tilde {}\to \zeta$ be a simple $\bar{\Bbb Q}$-named$^1$ ordinal.
For $r$ an atomic $\bar{\Bbb Q}$-named condition $r \restriction
\underset\tilde {}\to \zeta$ is an atomic $\bar{\Bbb Q}$-named condition. \nl
3) For $q \in \Bbb P_\alpha$ we have 
$q \restriction \underset\tilde {}\to \zeta \in \Bbb P_\alpha$. \nl
4) For $q_1,q_2 \in \Bbb P_\alpha,q_1 \le q_2 \Rightarrow q_1 \restriction
\underset\tilde {}\to \zeta \le q_2 \restriction 
\underset\tilde {}\to \zeta$.
\nl
5) If $q 
\in \Bbb P_{\underset\tilde {}\to \zeta},p \in \Bbb P_\alpha,p \restriction
\zeta \le q$ \ub{then} $p \cup q \in \Bbb P_\alpha$ is a lub of $p$
and $q$. 
\nl
6) $\Bbb P_{\underset\tilde {}\to \zeta} \lessdot \Bbb P_\alpha$. \nl
7) If $G \subseteq \Bbb P_{\underset\tilde {}\to \alpha}$ is generic
over $\bold V,\xi = \underset\tilde {}\to \zeta[G]$ \ub{then} $G \cap
P_{\underset\tilde {}\to \zeta},G \cap \Bbb P_\xi$ are 
essentially the same. \nl
8) The parallel statements with $\Xi$ instead of $\underset\tilde {}\to 
\zeta$.
\endproclaim
\bigskip

\remark{Remark}  In fact 
by part (1), part (6) follows from the parts (2)-(5).
\endremark
\bigskip

\proclaim{\stag{n.4} Claim}  Let $e = 4(2)$.
Assume ${\underset\tilde {}\to \zeta_n}$ is a
simple $\bar{\Bbb Q}$-named for $n < \omega,
{\underset\tilde {}\to \zeta_n} < 
{\underset\tilde {}\to \zeta_{n+1}}$ and 
for every $G \subseteq \Bbb P_\alpha$ generic over $\bold V$, 
for some $n,\varphi(n,G \cap 
\Bbb P_{{\underset\tilde {}\to \zeta_n}[G]})$.  \ub{Then} for some simple
$\bar{\Bbb Q}$-name ordinal $\underset\tilde {}\to \xi$, we have

$$
\Vdash_{{\Bbb P}_\alpha} ``\text{for some } n,\underset\tilde {}\to \xi
[{\underset\tilde {}\to G_{{\Bbb P}_\alpha}}] = 
{\underset\tilde {}\to \zeta_n}
[{\underset\tilde {}\to G_{{\Bbb P}_\alpha}}]
\text{ and } \varphi(n,{\underset\tilde {}\to G_{{\Bbb P}_\alpha}} \cap
\Bbb P_{{\underset\tilde {}\to \zeta_n}[G_{{\Bbb P}_\alpha}]})".
$$
\endproclaim
\bigskip

\newpage

     \shlhetal 

\newpage
    
REFERENCES.  
\bibliographystyle{lit-plain}
\bibliography{lista,listb,listx,listf,liste}

\enddocument